\documentclass[10pt]{amsart}

\usepackage{amsthm,amsmath,amssymb}
\usepackage[right=1.13in,left=1.13in,bottom=1in,top=1in]{geometry}
\usepackage[numbers,sort & compress]{natbib}
\usepackage{datetime}
\usepackage{multirow}
\usepackage{multicol}
\usepackage{graphicx}
\usepackage{booktabs}
\usepackage{caption}
\usepackage{listings}
\usepackage{makecell}
\usepackage{float}
\usepackage{array}
\usepackage{longtable}
\usepackage{lmodern}

\newcolumntype{L}[1]{>{\raggedright\let\newline\\\arraybackslash\hspace{0pt}}m{#1}}
\newcolumntype{C}[1]{>{\centering\let\newline\\\arraybackslash\hspace{0pt}}m{#1}}
\newcolumntype{R}[1]{>{\raggedleft\let\newline\\\arraybackslash\hspace{0pt}}m{#1}}

\apptocmd{\thebibliography}{\raggedright}{}{}

\allowdisplaybreaks

\newcommand\Tstrut{\makebox[0pt][c]{\rule{0pt}{2.6ex}}} 
\newcommand\Bstrut{\makebox[0pt][c]{\rule[-1.2ex]{0pt}{0pt}}} 
\def\mystrut(#1,#2){\vrule height #1pt depth #2pt width 0pt}

\makeatletter
\patchcmd{\@settitle}{\uppercasenonmath\@title}{}{}{}
\patchcmd{\@setauthors}{\MakeUppercase}{}{}{}
\patchcmd{\section}{\scshape}{}{}{}
\makeatother

\usepackage[table]{xcolor}
\definecolor{solv}{rgb}{0.65,1,0.65}
\definecolor{impo}{rgb}{1,0.65,1}
\definecolor{tent}{rgb}{1,1,0.65}
\definecolor{real}{rgb}{0.65,1,1}
\definecolor{unkn}{gray}{0.85}
\definecolor{newt}{rgb}{0.5,0.5,1}
\definecolor{na}{gray}{0.1}
\definecolor{bug}{rgb}{1,0.5,0.5}

\title{Some spherical codes in S2 and their algebraic numbers}
\author{Randall L. Rathbun}
\author{Wesley J.M. Ridgway}
\email{randallrathbun@protonmail.com}
\email{wjr704@mail.usask.ca}
\keywords{Thomson problem, global minima for points on S2, Smale Problem \#7}

\begin{document}

\begin{abstract}
The first 195 spherical codes for the global minima of 1-65 points on S2 have been obtained for 3 types of potentials: logarithmic $log(r)$, Coulomb $1/r$, called the \textit{Thomson problem}, and the inverse square law $1/{r^2}$, with 77, 38, and 38 digits precision respectively. Solutions for the logarithmic potential are solutions to Smale's Problem \#7.

It was discovered that certain point sets have embedded polygonal structures, constraining the points, enabling them to be parameterized and to successfully recover the algebraic polynomial.

So far 49 algebraic number sets have been recovered, but 109 more remain to be recovered from their 1,585 parameters of which 1,124 are known to at least 50,014 digit precision. The very high algebraic degree of these minimal polynomials eludes finding the algebraic numbers from the spherical codes and requires new mathematical tools to meet this challenge.
\end{abstract}

\maketitle

\setcounter{section}{0}
\raggedbottom

\section{Energy Minimization of Points on a Sphere}
\subsection{Introduction}
The distribution of points on the S2 sphere is considered so important a mathematical problem that Steve Smale made it unsolved problem \#7 on his list\cite{55,56} of mathematical problems for the 21st century.

While extensive literature\cite{1,2,3,4,5,6,7,8,9,10,11,12,13,14,15,16,17,18,19,20,21,22,23,24,25,26,27,28,29,30,31,32,33,34,35,36,37,38,39,40,41,42,43,44,45,46,47,48,49,50,51,52,53,54,55,56,57,58,59,60,61,62} exists on optimizing spherical codes in $S2$ under various constraints, such as energy minimization for a given potential (force-law field), very little has been published\cite{60} regarding the algebraic numbers of the global minima.

It was decided to take a 3-prong approach to obtain these numbers by using two computer algorithms, \textit{percolating annealing}, and \textit{gradient descent}, in order to obtain high precision $E^3$ real numbers for the coordinates, and then utilizing certain polygon structures embedded in the cluster, to parameterize the configuration, then use a Jacobian matrix of these parameters to successfully recover the spherical code if possible, from the very high precision values (50,014 digits). As a final check, a Hessian matrix was run on the configuration to verify that a global minima had been reached.

\subsection{Point Set on Unit Sphere}

Let $\vec{x}$ be a unit vector in $R^3$ with real coordinates $\vec{x}_{1..3}\in{R}$ such that $|x|=1$, thus the point $x$ is on the unit sphere S2. We define a point set of size $n$ as
\begin{equation}
P_n = \{x_1,x_2,...x_n\}
\end{equation}

\subsection{Potentials}
We seek to minimize the energy under a potential. We need to consider the distance between a point pair first.

Let $d = |x_i - x_j|$ be the Euclidean distance $d$ between the two points $x_i$ and $x_j$($i\ne j$). We define the potential as a force depending upon a function of the distance $d$. The potentials considered
are

\noindent
\textit{Logarithmic Potential}
\begin{equation}
\mathcal{P}_{log} = -log(d) \quad (d\ne 0)
\end{equation}

\noindent
\textit{Coulomb Potential}
\begin{equation}
\mathcal{P}_{\frac{1}{r}} = \frac{1}{d} \quad (d\ne 0)
\end{equation}

\noindent
\textit{Inverse Square Law}
\begin{equation}
\mathcal{P}_{\frac{1}{r^2}} = \frac{1}{d^2} \quad (d\ne 0)
\end{equation}
\subsection{Riesz-$s$ kernel}
It is important to point out that the 3 potentials considered are 3 cases of the Reisz-$s$ kernel defined as
\begin{equation}
K_s(x,y)=\frac{1}{|x-y|^s} \text{ for } s>0 \text{ and } K_{log}(x,y)= log\frac{1}{|x-y|}
\end{equation}
for $s=1,2$ and the $log$ case. Again $|x-y| \ne 0$.
\subsection{Energy Minimization}
The goal is to find the minimal energy for a point set $P_n$ under these potentials.
\noindent
For the logarithmic potential we have:
\begin{equation}
\mathcal{E}_{log}(n) = -\sum_{i\ne j}^{n}{log(|\vec{x}_i - \vec{x}_j|)}
\label{eq:log}
\end{equation}
\noindent
For the Coulomb potential we have:
\begin{equation}
\mathcal{E}_{\frac{1}{r}}(n) = \sum_{i\ne j}^{n}{\frac{1}{|\vec{x}_i - \vec{x}_j|}}
\label{eq:coul}
\end{equation}
This energy minimization is the one considered in the \textit{Thomson problem}.

\noindent
For the inverse square law we have:
\begin{equation}
\mathcal{E}_{\frac{1}{r^2}}(n) = \sum_{i\ne j}^{n}{\frac{1}{|\vec{x}_i - \vec{x}_j|^2}}
\label{eq:invsq}
\end{equation}

\section{Computational Algorithms to find minimal energy}
It was decided that the way to find the minimal solutions was to use a computational program to actually search out an optimal $P_n$ configuration starting with a known configuration of points.

\subsection{Exhaustive search}
In the exhaustive search approach, each point position is carefully perturbed in a controlled manner, for all $n$ points, and after each perturbation of all points, the energy was determined. If a lower value of energy was found, the minimal energy positions were recorded, then used as the starting place for another search. This search continues until the desired precision is reached.

\noindent
\textbf{Exhaustive Search Algorithm:}
\begin{lstlisting}
points_type pts; /* predetermined point set */
mpf_class scale = 1.0e-1
mpf_class final_precision = 1.0e-40;
bool done = false;
bool found = false;
points_type qpts = pts;
while (! done); do {
  found = false;
  for(int i=\-10; i<=10; i++) {
    qpts[i] = scale*i + pts[0];
    for(int j=-10; j<=10; j++) {
      qpts[j] = scale*j + pts[1];
      for(int k=-10; k<=10; k++) {
        qpts[k] = scale*k + pts[2];
           [etc, sub nesting loops for all n points]
        mpf_class diff = energy(law,qpts) - energy(law,pts);
        if (diff < 0) {
          pts = qpts;
          found = true;
        }
      }
    }
  }
  if ( found ) {
    scale /= 10; /* advance to the next decade */
  }
  if ( diff < final_precision ) { done = true; }
}
store(pts);
\end{lstlisting}
It can be seen that this exhaustive algorithm will indeed
find the global minima, but the order of the searching is
$\mathcal{O}^n$ which grows exponentially in size as $n$
increases. The algorithm becomes too cumbersome and
inefficient, even for moderate sizes of $n$ points.

\noindent
This approach was abandoned and two more efficient approaches were tried.

\subsection{Percolating Annealing}
The main idea of this algorithm is to add random changes to the coordinates of the points, test for a lower minimum, but keep trying until no more lower energy configurations were found. Then the scaling factor "\texttt{scaling}" would be decreased by a ratio, and the algorithm repeated again, ad infinitum, until the desired accuracy was found. The original value of \texttt{scaling} is found manually to start the algorith off, depending upon the initial accuracy of the starting point set.

\noindent
\textbf{Percolating Annealing Algorithm:}
\begin{lstlisting}
percolating_anneal(points_type P) {
  const int total_passes = 100000;
  mpf_class scaling = 1.0e-1;
  mpf_class final_precision = 1.0e-44;
  points_type Q = P;
  bool found = false;
  bool done = false
  while ( ! done ); do {
    for (int passes = 0; passes < total_passes; passes+=1;) {
      Q = jiggle(scaling,P);
      diff = energy(law,Q) - energy(law,P)
      if (diff < 0) {
        P = Q;
        found = true;
      }
    }
    if (found) {
      scaling *= 0.8; /* decrease the scaling factor */
    }
    if (diff < final_precision) {
      done = true;
    }
  }
  return(P);
}
\end{lstlisting}
\texttt{Q=Jiggle(scaling,P)} is a function call that takes the
current scaling value and points and slightly randomly jiggles
each point to obtain Q. Its net effect is to slightly jostle
(boil) the points, in hopes of finding a lower energy value set.
The amount of variation in the jiggle is $v = scaling*random()$
where $-1 \le random() \le +1$. Each $x,y,z$ coordinate in a
point is randomly varied, for all $n$ points.

Please notice that if no lower energy configuration is
found during the first 100,000 passes, it repeats again in hopes of finding a lower energy configuration. This is similar to percolating coffee on the stove, hence the term \textit{percolating annealing}. This algorithm can lock up if no lower energy is found, but that only indicates that a local minimum has already been reached.

\subsection{A comment on isometric equivalency}
An isometry is a mapping transform which preserves distances in the object. Since the Gram matrix is uniquely determined by every possible dot-product pair in the object, then two objects are considered geometrically equivalent if their Gram matrix is equivalent.
\begin{equation}
\text{\textit{ If }} Gram(P) = Gram(Q) \text{ \textit{then }} isometry(P) = Q
\end{equation}
i.e., an isometric mapping can be found to map the
coordinates of one object to the other and vice-versa.

To test that this \textit{percolating annealing} algorithm was finding at least a local minimum, 5 different clusters of 19 points were created, and using $s=1$ or Coulomb repulsion, the same isometries
of points were obtained to 57 digits (energy to 54 digits) indicating that the same (global) minimum had been reached in each case.

This checking was repeated for other $n$-points sets, and in every case, the same energy minimum was found for a set, and double-checking the Gram matrix of the point set, an isometric configuration had been found each time.

\subsection{A comment on the scaling ratio decrease of 0.8}
The ability of this algorithm to detect the global minimum depends critically on the scaling ratio (0.8) which was found by numerical simulations. Too small a value would indicate violent boiling, thus possibly missing a true global minimum, but too large a value would simmer the points, taking an excessive amount of time to find the global minimum.

Some point clusters took much time to converge, 53 points and 36 points, in particular, took several weeks to settle down to their minimum energy value to 38 digits accuracy.

\subsection{Gradient Descent}

Let the force $f_i(s)$ on a point $P_i$ be defined as
\begin{equation}
F_i(s) = \sum_{i\ne j}^{n} \nabla_i \frac{1}{|P_j-P_i|^s} \quad s>0 \text{ or } \
F_i(s) = \sum_{i\ne j}^{n} \nabla_i log(\frac{1}{|P_j-P_i|}) \quad s=0
\end{equation}
where $s$ is the Riesz-$s$ kernel parameter or force law. $\nabla_i$ is the gradient operator at point $i$.

The force $f_i(s)$ on each point has 2 components, one is
normal to the S2 unit sphere at the point, the other one
is tangential and it is the one we are interested in.
\begin{equation}
\vec{F}(i) = \vec{F}_N(i) + \vec{F}_T(i)
\end{equation}
The gradient operator upon a point $P_i$ is defined as
the tangential component of the derivative of the force upon the point.
\begin{equation}
\nabla_i = {\frac{d(F_i)}{ds}}_T
\end{equation}
\noindent
\textbf{Descent Algorithm:}
\begin{lstlisting}
descent(points_type P) {
  points_type Q;
  unsigned long size = P.size();
  double alpha = 0.5 / (double) size;
  for(unsigned long i = 0; i < size; i++) {
    coordinate_type F = force(power-1,P,i); /* derivative */
    coordinate_type T = tangential_component(P[i],F);
    coordinate_type R = P[i] + alpha*T;
    coordinate_type S = unit_vector(R);
    Q[i] = S
  }
  return(Q)
}
\end{lstlisting}
We notice that derivative of force is the typical one
taught in college, for example for $s=1$ we have
\begin{equation}
\frac{d(\frac{1}{s})}{ds} = -\frac{1}{s^2}
\end{equation}
similarly for all $s>0$, while for the log condition we have
\begin{equation}
\frac{d(log(s)}{ds} = \frac{1}{s}
\end{equation}
hence the appearance of \texttt{force(power-1,P,i)} in the descent procedure listed above, to find the gradient $\nabla_i$ for a point under the $s$ potential condition.

The $\alpha$ factor of $\mathbf{0.5}$ was found by experimentation
to be the typical factor to slightly diminish the tangential force $f_T$ contribution when finding the newer point $Q_i = ||P_i+\alpha \* T_i||$. Once or twice, it was found necessary to adjust the $\alpha$ value to around $\mathbf{0.1}$ or so, in order to proceed on the descent, as lock up would occur with the $\mathbf{0.5}$ value.

\subsection{Why 2 algorithms?}
The reason for 2 algorithms is simple, each is a check upon the other. After obtaining a point set for $n$, the
alternate algorithm was used to double-check the results by comparing the minimal energy values.

In every case, both converged to the same isometric configuration of points and returned the same energy value to within the precision of the points (often 38 - 115 digits).

\section{The Gram Matrix of the point set $P_n$}
The Gram matrix is defined as
\begin{equation}
G = V^TV
\end{equation}
where $\vec{V}$ is a $n \times n$ matrix whose columns are the coordinates of the points. This means that $G$
\begin{equation}
G_{ij} = V_i \bullet V_j
\end{equation}
is the matrix of dot products of the vectors of the point set.

\textit{It is important to recognize that the Gram matrix G is the unique signature of the point set, and can be used to identity identical point sets, though an isometry mapping, to within limits of precision.}

Furthermore, the Gram matrix $G$ will be algebraic numbers for points on the unit sphere S2, provided that at least 1 point of the set has been embedded into algebraic coordinates, either rational $\mathcal{Q}$ or algebraic $\mathcal{A}$.

\subsection{Coordinates for the spherical code instead of the Gram Matrix}
While it could have been proper to only give the Gram
matrix $\boldsymbol{G}$ for a point set as the unique signature of a global minimal configuration, it was decided to deal with the coordinates of a configuration instead, in order to obtain a spherical code, because an {\bf SVD} matrix decomposition to very high precision $2n$ digits is needed to successfully recover the coordinates into algebraic numbers of precision $n$ digits.

\section{Embedded geometric structures}
It was soon discovered during the early part of the search phase, that optimal placing of the points for the minimal energy often led to a configuration of some of the points being on certain $2d$ planes as polygons.

Often these polygons would be oriented along a normal
of a plane and would give a strong hint on how to rotate
and orient the minimal set into the "best" orientation.

While any 3 points of a point set are most certainly located upon a $2d$ triangle, not every triad was considered when looking for embedded geometric structures.

Consider the case of 5 points given in the figure \ref{fig:5pts_x} below.
\begin{figure}[h]
	\begin{center}
		\includegraphics[type=pdf,ext=pdf,read=pdf,height=1in,width=1in,angle=0]{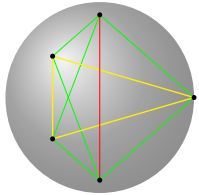}
		\caption{5 points.}
		\label{fig:5pts_x}
	\end{center}
\end{figure}
It is easy to see that 2 of the points are poles, while the other 3 points are located on the equator. They are the vertices of an equilateral triangle shown in yellow.

Similarly for other point sets, alignments of polygons
did occur, and rotations were employed in order to successfully recover the spherical code.

\section{Using the Jacobian matrix of the parameters}

After discovering the embedded polygons for a particular point set, it was possible to derive an expression for the points utilizing parameters.

\subsection{Parameterizing a configuration}
Consider the following, suppose that the configuration contained an embedded equilateral triangle. The following equation gives the algebraic coordinates of an equilateral triangle at different heights $n$ (algebraic) along the z-axis
for the point set.
\begin{equation}
a=\sqrt{1-n^2} \; ; \quad
\triangle(n)=[a,0,n],\left[-\frac{a}{2},a\frac{\sqrt{3}}{2},n \right],\left[-\frac{a}{2},-a\frac{\sqrt{3}}{2},n \right]
\end{equation}
A point set would be carefully examined, and then the positions of the coordinates would be expressed as a function of the needed parameters to express the embedded polygons found.

For example, the search for 27 points was optimized into two pole points and 5 pentagons at different heights along the pole axis. Thus it initially took 5 parameters to express this polyhedron

Actually, due to balance required for minimal energy, it was discovered that the 5 parameters dropped to just 2 with the 5 pentagons at levels $[a,b,0,-b,-a]$ along the pole axis.

This parametric equation for the 27-point set under the $1/r^2$ law is given:
\begin{align*}
P_{1} &= \left[\; 0,\; 0,\; 1 \; \right] \\
P_{2} &= \left[\sqrt{1-a^2},\; 0,\; a \; \right] \\
P_{3} &= \left[\frac{\sqrt{5}-1}{4}\sqrt{1-a^2},\; \frac{\sqrt{10 + 2\sqrt{5}}}{4}\sqrt{1-a^2},\; a \; \right] \\
P_{4} &= \left[-\frac{\sqrt{5}+1}{4}\sqrt{1-a^2},\; \frac{\sqrt{10 - 2\sqrt{5}}}{4}\sqrt{1-a^2},\; a \; \right] \\
P_{5} &= \left[-\frac{\sqrt{5}+1}{4}\sqrt{1-a^2},\; -\frac{\sqrt{10 - 2\sqrt{5}}}{4}\sqrt{1-a^2},\; a \; \right] \\
P_{6} &= \left[\frac{\sqrt{5}-1}{4}\sqrt{1-a^2},\; -\frac{\sqrt{10 + 2\sqrt{5}}}{4}\sqrt{1-a^2},\; a \; \right] \\
P_{7} &= \left[-\sqrt{1-b^2},\; 0,\; b \; \right] \\
P_{8} &= \left[-\frac{\sqrt{5}-1}{4}\sqrt{1-b^2},\; -\frac{\sqrt{10 + 2\sqrt{5}}}{4}\sqrt{1-b^2},\; b \; \right] \\
P_{9} &= \left[\frac{\sqrt{5}+1}{4}\sqrt{1-b^2},\; -\frac{\sqrt{10 - 2\sqrt{5}}}{4}\sqrt{1-b^2},\; b \; \right] \\
P_{10} &= \left[\frac{\sqrt{5}+1}{4}\sqrt{1-b^2},\; \frac{\sqrt{10 - 2\sqrt{5}}}{4}\sqrt{1-b^2},\; b \; \right] \\
P_{11} &= \left[-\frac{\sqrt{5}-1}{4}\sqrt{1-b^2},\; \frac{\sqrt{10 + 2\sqrt{5}}}{4}\sqrt{1-b^2},\; b \; \right] \\
P_{12} &= \left[\; 1,\; 0,\; 0 \; \right] \\
P_{13} &= \left[\frac{\sqrt{5}-1}{4},\; \frac{\sqrt{10+2\sqrt{5}}}{4},\; 0 \; \right] \\
P_{14} &= \left[-\frac{\sqrt{5}+1}{4},\; \frac{\sqrt{10-2\sqrt{5}}}{4},\; 0 \; \right] \\
P_{15} &= \left[-\frac{\sqrt{5}+1}{4},\; -\frac{\sqrt{10-2\sqrt{5}}}{4},\; 0 \; \right] \\
P_{16} &= \left[\frac{\sqrt{5}-1}{4},\; -\frac{\sqrt{10+2\sqrt{5}}}{4},\; 0 \; \right] \\
P_{17} &= \left[-\sqrt{1-b^2},\; 0,\; -b \; \right] \\
P_{18} &= \left[-\frac{\sqrt{5}-1}{4}\sqrt{1-b^2},\; -\frac{\sqrt{10 + 2\sqrt{5}}}{4}\sqrt{1-b^2},\; -b \; \right] \\
P_{19} &= \left[\frac{\sqrt{5}+1}{4}\sqrt{1-b^2},\; -\frac{\sqrt{10 - 2\sqrt{5}}}{4}\sqrt{1-b^2},\; -b \; \right] \\
P_{20} &= \left[\frac{\sqrt{5}+1}{4}\sqrt{1-b^2},\; \frac{\sqrt{10 - 2\sqrt{5}}}{4}\sqrt{1-b^2},\; -b \; \right] \\
P_{21} &= \left[-\frac{\sqrt{5}-1}{4}\sqrt{1-b^2},\; \frac{\sqrt{10 + 2\sqrt{5}}}{4}\sqrt{1-b^2},\; -b \; \right] \\
P_{22} &= \left[\sqrt{1-a^2},\; 0,\; -a \; \right] \\
P_{23} &= \left[\frac{\sqrt{5}-1}{4}\sqrt{1-a^2},\; \frac{\sqrt{10 + 2\sqrt{5}}}{4}\sqrt{1-a^2},\; -a \; \right] \\
P_{24} &= \left[-\frac{\sqrt{5}+1}{4}\sqrt{1-a^2},\; \frac{\sqrt{10 - 2\sqrt{5}}}{4}\sqrt{1-a^2},\; -a \; \right] \\
P_{25} &= \left[-\frac{\sqrt{5}+1}{4}\sqrt{1-a^2},\; -\frac{\sqrt{10 - 2\sqrt{5}}}{4}\sqrt{1-a^2},\; -a \; \right] \\
P_{26} &= \left[\frac{\sqrt{5}-1}{4}\sqrt{1-a^2},\; -\frac{\sqrt{10 + 2\sqrt{5}}}{4}\sqrt{1-a^2},\; -a \; \right] \\
P_{27} &= \left[\; 0,\; 0,\; -1 \; \right]
\end{align*}
It can be quickly seen that if $a,b \in \mathcal{A}$ are algebraic, then we have the spherical code for 27 points under the $\frac{1}{r^2}$ potential.

Letting the search programs run for awhile and then checking the configuration of points led to the discovery of the embedded polygons and enabled a parameterization to be created for many of the point sets, but not all of them.

\subsection{Constraints upon points}
We come to a critical part of the global minimal search. After obtaining a parameterization, such as the one given for 27 points $\frac{1}{r^2}$ above, how can the energy be minimized?

This was done by setting up a Jacobian matrix for the parameters under the force law, and using the Newton approximation to derive ultra-high precision (50,014 digits) solution so that the spherical code can be successfully recovered.

Using the parameters $a,b,c,\dots$ then let the first derivatives of the energy function $\mathcal{E}(P_n,r)$ for $n$ points be collected as a vector function
\begin{equation}
V(a,b,c,\dots) = \left[
\frac{d\mathcal{E}}{da}, \quad
\frac{d\mathcal{E}}{db}, \quad
\frac{d\mathcal{E}}{dc}, \quad \dots
\right]
\end{equation}

for parameters $a,b,c,\dots$ then the
Jacobian is defined as:

\begin{equation}
J(p)_{a,b,c,\dots} =
\begin{bmatrix}
\frac{\partial\left(\frac{d\mathcal{E}}{da}\right)}{\partial a} &
\frac{\partial\left(\frac{d\mathcal{E}}{da}\right)}{\partial b} &
\frac{\partial\left(\frac{d\mathcal{E}}{da}\right)}{\partial c} & \dotsc \, \\[0.5em]
\frac{\partial\left(\frac{d\mathcal{E}}{db}\right)}{\partial a} &
\frac{\partial\left(\frac{d\mathcal{E}}{db}\right)}{\partial b} &
\frac{\partial\left(\frac{d\mathcal{E}}{db}\right)}{\partial c} & \dotsc \, \\[0.5em]
\frac{\partial\left(\frac{d\mathcal{E}}{dc}\right)}{\partial a} &
\frac{\partial\left(\frac{d\mathcal{E}}{dc}\right)}{\partial b} &
\frac{\partial\left(\frac{d\mathcal{E}}{dc}\right)}{\partial c} & \dotsc \, \\
\vdots & \vdots & \vdots & \ddots
\end{bmatrix}
\end{equation}
\subsection{Newton's method}
Newton's method is used to converge quadratically to the desired precision to constrain the points.

\noindent
Repeatedly solve the matrix equation \textbf{MX=B} for the incremental $h$, if given an initial $V_0$
\begin{equation}
\begin{split}
J_p \, h & = -V_i^T \\
V_{i+1} & = V_{i} + h
\end{split}
\end{equation}
until the desired precision (50,014 digits) is reached.

Once the desired precision was obtained, an adaption of the arbitrary precision GP-Pari number theory calculator command $algdep(n,d)$ was used in an attempt to find the correct algebraic polynomial of degree $d$, the spherical code.

\section{Hessian matrix check on the global minimum}
Having found at least the local minimum, it is still necessary to check the Hessian matrix to ensure that the global minimum has actually been found.

The Hessian matrix is defined as
\begin{equation}
H(a,b,c,\dots) =
\begin{bmatrix}
\frac{\partial^2 \mathcal{E}}{\partial a^2} &
\frac{\partial^2 \mathcal{E}}{\partial a\partial b} &
\frac{\partial^2 \mathcal{E}}{\partial a\partial c} & \dotsc \, \\[0.5em]
\frac{\partial^2 \mathcal{E}}{\partial a\partial b} &
\frac{\partial^2 \mathcal{E}}{\partial b^2} &
\frac{\partial^2 \mathcal{E}}{\partial b\partial c} & \dotsc \, \\[0.5em]
\frac{\partial^2 \mathcal{E}}{\partial a\partial c} &
\frac{\partial^2 \mathcal{E}}{\partial b\partial c} &
\frac{\partial^2 \mathcal{E}}{\partial c^2} & \dotsc \, \\[0.5em]
\vdots & \vdots & \vdots & \ddots
\end{bmatrix}
\end{equation}
For a global minimum to occur, all the eigenvalues of the Hessian matrix must be positive or zero. It was discovered by checking the actual matrices, using the GP-Pari \textit{mateigen()} command to find both the eigenvalues and eigenvectors, that the first 3 values of the eigenvalues were 0 and the next were positive, thus confirming that the Hessian matrix is positive-definite. The first 3 values of 0 are probably due to the fact that an infinite number of 3d isometry maps exist which are rotations or reflections of a point set.

The Hessian check was used both for points found by a spherical code and points found by direct searching, using either \textit{percolating annealing} or \textit{descent} algorithms, in order to double check the global minimum point configuration.

We used Maple to symbolically compute the Hessian\cite{47}, and the eigenvalues thereof, up to the desired accuracy(digits) of the point set.

\section{Alignments and Gram matrix}
Carefully checking the polygon alignments along a common normal of a plane enabled the preferred orientation axis to be discovered and exploited for parameterizing the polyhedron.

When the Gram matrix was obtained, the components did show evidence of groups. However, in pursuing the spherical code, it was decided to work with the coordinates of the points, not the Gram matrix itself, due to the difficulty of recovering the coordinates accurately from the \textbf{SVD} matrix decomposition in an algebraic form which requires $2n$ precision due to the square root operation inside the algorithm.

In the solutions for points 1-65, the aligned planes, the Gram matrix groups, and the polygons are listed. A pictorial figure is also provided.

\section{Spherical codes and Algebraic Numbers}
The question is asked, how do we know that the spherical codes found represent an algebraic number?

First of all, the spherical code is a solution of the energy equations (\ref{eq:log}), (\ref{eq:coul}), or (\ref{eq:invsq}) which are rational polynomials for the $1/r$ and $1/r^2$ potentials, and we can take the antilog of the $log$ potential to recover the rational polynomial for the $log$ potential, since a summation of logs is equivalent to a product, then the summation is also an algebraic number, since an algebraic field is closed for the 4 types of common arithmetic operators (+,-,*,/) by definition.

Secondly, when using the {\it algdep(n,deg)} command of GP-Pari, it was noticed that when the precision of solution was increased, that at a certain point the integer coefficients for the algebraic number suddenly "snapped into place" or dropped in size and also the powers of the polynomial usually became even, strongly hinting that the correct characteristic polynomial had been located.

Thirdly, once a characteristic polynomial was found, a few runs of increasing precision found that the spherical code stayed strictly on the S2 unit sphere, which the code must, if a true solution actually was found. Similarly the algebraic number for the energy was checked by increasing the precision of the spherical code to see if they matched.

In this paper, where the algebraic number is given, that result has been checked. Where the spherical codes exist, to 50,014 digits of accuracy, it is asserted that they represent an algebraic number of degree $>360$ (in some cases, $>480$).

In all cases where the algebraic number was found, the assertion that the spherical code was algebraic has been shown to be true.

It is strongly conjectured that all spherical codes in the {\it Current Status of Unknown Parameters} section are algebraic, but current mathematical tools are unable to uncover the high degree of the polynomial $(>360)$ in a reasonable time. It is hoped that further progress will resolve these unknown parameters for the global minima solutions and successfully resolve the algebraic numbers.

\section{Globally optimal results}
A configuration of the global minimal solution for points 1 - 65 on the unit sphere $S2$ is provided for the 3 types of potentials, \textit{logarithmic, Coulomb $1/r$, and inverse square $1/r^2$}.

If an algebraic solution exists, the configuration has been rotated (often twice) into a preferred orientation.

For points found under the \textit{log} potential, the accuracy of their coordinates is 77 digits (often 80). For points found under the $1/r$ or $1/r^2$ potential, the accuracy is 38 digits (often 40).

For algebraic solutions, or its spherical code, the accuracy is unlimited.

The minimal energy found by these solutions has been compared to those discovered by Ridgway and Cheviakov\cite{47} and the results match for all 195 configurations to within 4.94707e-08.

\noindent
\textbf{Important Note on the \textit{Coulomb} or \textit{J. J. Thomson Problem}\cite{57} solutions - }
\begin{itemize}
	\setlength\itemsep{0.7em}
	\item The solutions for the \textit{Coulomb} or $1/r$ solutions are improvements of Neil J.A. Sloane's \textit{"Spherical codes for minimal energy"} database\cite{53}. They have improved his 12 digits accuracy to at least 38 digits accuracy. His putative solutions are the global minimum, as found out in Ridgway \& Cheviakov\cite{47} and here, for points 1 to 65 inclusive.
\end{itemize}

\subsection{Comment on Points 1 to 8}
Points of size 1 to 8 have had their globally optimal configurations determined and known, in some cases for hundred of years, such as 3, 4 or 6 points. The solutions are the equilateral triangle, the tetrahedron, and octahedron, and the last 2 are Platonic solids.

\noindent
The global solution for 5 points has only been recently proven\cite{51,52}.

\noindent
The solutions for 7 and 8 points are relatively new.

It must be noted that the solution for 8 points is NOT the cube in figure \ref{fig:cube}, despite the common assumption. It is the square antiprism shown in figure \ref{fig:square_antiprism} which has the lowest energy value for 8 points on S2 for all 3 potentials. One face is rotated $45^{\circ}$ with respect to the parallel face as shown in yellow outlines in both figures.
\begin{figure}[ht]
	\centering
	\begin{minipage}{0.45\textwidth}
		\centering
		\includegraphics[type=pdf,ext=pdf,read=pdf,height=1in,width=1in,angle=0]{square.}
		\caption{Cube.}
		\label{fig:cube}
	\end{minipage}\hfill
	\begin{minipage}{0.45\textwidth}
		\centering
		\includegraphics[type=pdf,ext=pdf,read=pdf,height=1in,width=1in,angle=0]{square-antiprism.}
		\caption{Square antiprism.}
		\label{fig:square_antiprism}
	\end{minipage}
\end{figure}

\subsection{1 point}
The optimal solution for 1 point is trivial. Since only 1 point is involved, there is no minimal energy.
\begin{figure}[ht]
	\begin{center}
		\includegraphics[type=pdf,ext=pdf,read=pdf,height=1in,width=1in,angle=0]{r-1.1pt.}
		\caption{1 point.}
		\label{fig:1pt}
	\end{center}
\end{figure}
\begin{center}
	\begin{tabular}{c|ccc}
		\multicolumn{4}{c}{Spherical code} \\
		\hline\Tstrut
		pt & x & y & z \\
		\hline\Tstrut
		1 & 1 & 0 & 0
	\end{tabular}
\end{center}
\begin{center}
\begin{tabular}{l|l||l|l}
\multicolumn{2}{c}{Symmetries - 1 point} & \multicolumn{2}{c}{Minimal Energy} \\
\hline\Tstrut
planes & [] & log & 0 \\
\hline\Tstrut
Gram groups & [[1, 1]] & $1/r$ & 0 \\
\hline\Tstrut
Polygons & [] & $1/r^2$ & 0
\end{tabular}
\end{center}

\subsection{2 points}
The optimal solution for 2 points is on the diameter.
\begin{figure}[h]
	\begin{center}
		\includegraphics[type=pdf,ext=pdf,read=pdf,height=1in,width=1in,angle=0]{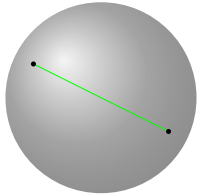}
		\caption{2 points.}
		\label{fig:2pts}
	\end{center}
\end{figure}
\begin{center}
	\begin{tabular}{c|ccc}
		\multicolumn{4}{c}{Spherical code} \\
		\hline\Tstrut
		pt & x & y & z \\
		\hline\Tstrut
		1 & 1 & 0 & 0 \\
		2 & -1 & 0 & 0
	\end{tabular}
\end{center}
\begin{center}
	\begin{tabular}{l|l||l|l}
		\multicolumn{2}{c}{Symmetries - 2 points} & \multicolumn{2}{c}{Minimal Energy} \\
		\hline\Tstrut
		planes & [] & log & $-log(2)$ \\
		\hline\Tstrut
		Gram groups & [[2, 2]] & $1/r$ & $1/2$ \\
		\hline\Tstrut
		Polygons & [] & $1/r^2$ & $1/4$
	\end{tabular}
\end{center}

\subsection{3 points}
The optimal solution for 3 points is an equilateral triangle, with the triangle along the equator.
\begin{figure}[ht]
	\begin{center}
		\includegraphics[type=pdf,ext=pdf,read=pdf,height=1in,width=1in,angle=0]{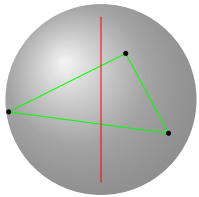}
		\caption{3 points.}
		\label{fig:3pts}
	\end{center}
\end{figure}
\begin{center}
	\begin{tabular}{c|ccc}
		\multicolumn{4}{c}{Spherical code} \\
		\hline\Tstrut
		pt & x & y & z \\
		\hline\Tstrut
		1 & 1 & 0 & 0 \\
		2 & $-1/2$ & $\sqrt{3}/2$ & 0 \\
		3 & $-1/2$ & $-\sqrt{3}/2$ & 0
	\end{tabular}
\end{center}
\begin{center}
	\begin{tabular}{l|l||l|l}
		\multicolumn{2}{c}{Symmetries - 3 points} & \multicolumn{2}{c}{Minimal Energy} \\
		\hline\Tstrut
		planes & [] & log & $-log(27)/2$ \\
		\hline\Tstrut
		Gram groups & [[3, 1], [6, 1]] & $1/r$ & $\sqrt{3}$ \\
		\hline\Tstrut
		Polygons & [] & $1/r^2$ & $1$
	\end{tabular}
\end{center}

\subsection{4 points}
The optimal solution for 4 points is the tetrahedron with vertices on the S2 sphere.
\begin{figure}[h]
	\begin{center}
		\includegraphics[type=pdf,ext=pdf,read=pdf,height=1in,width=1in,angle=0]{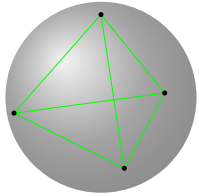}
		\caption{4 points.}
		\label{fig:4pts}
	\end{center}
\end{figure}
\begin{center}
	\begin{tabular}{c|ccc}
		\multicolumn{4}{c}{Spherical code} \\
		\hline\Tstrut
		pt & x & y & z \\
		\hline\Tstrut
		1 & 0 & 0 & 1 \\
		2 & $\frac{2\sqrt{2}}{3}$ & 0 & $-\frac{1}{3}$ \\
		3 & $-\frac{\sqrt{2}}{3}$ & $\sqrt{\frac{2}{3}}$ & $-\frac{1}{3}$ \\
		4 & $-\frac{\sqrt{2}}{3}$ & $-\sqrt{\frac{2}{3}}$ & $-\frac{1}{3}$
	\end{tabular}
\end{center}
\begin{center}
	\begin{tabular}{l|l||l|l}
		\multicolumn{2}{c}{Symmetries - 4 points} & \multicolumn{2}{c}{Minimal Energy} \\
		\hline\Tstrut
		planes & [] & log & $log(27/512)$\Bstrut \\
		\hline\Tstrut
		\Tstrut Gram groups & [[4, 1], [12, 1]] & $1/r$ & $\frac{3\sqrt{6}}{2}$ \mystrut(12.5,5) \\
		\hline\Tstrut
		Polygons & [] & $1/r^2$ & $9/4$
	\end{tabular}
\end{center}

\subsection{5 points}
The optimal solution for 5 points is an equilateral triangle on the equator and the two poles, north and south, the triangular bipyramid. Ron Schwartz proved\cite{51} that this configuration was optimal for the Coulomb potential and inverse square law potential. Searching for the logarithmic potential also showed this to be true.
\begin{figure}[ht]
	\begin{center}
		\includegraphics[type=pdf,ext=pdf,read=pdf,height=1in,width=1in,angle=0]{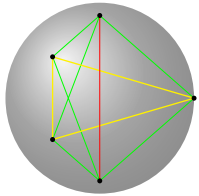}
		\caption{5 points.}
		\label{fig:5pts}
	\end{center}
\end{figure}
\begin{center}
	\begin{tabular}{c|ccc}
		\multicolumn{4}{c}{Spherical code} \\
		\hline\Tstrut
		pt & x & y & z \\
		\hline\Tstrut
		1 & 0 & 0 & 1 \\
		2 & 1 & 0 & 0 \\
		3 & $-\frac{1}{2}$ & $\frac{\sqrt{3}}{2}$ & 0 \\[0.7ex]
		4 & $-\frac{1}{2}$ & $-\frac{\sqrt{3}}{2}$ & 0 \\[0.7ex]
		5 & 0 & 0 & -1
	\end{tabular}
\end{center}
\begin{center}
	\begin{tabular}{l|l||l|l}
		\multicolumn{2}{c}{Symmetries - 5 points} & \multicolumn{2}{c}{Minimal Energy} \\
		\hline
		planes \rule{0pt}{3ex} & [1,1] & log & $-\frac{log(6912)}{2}$ \Bstrut \\
		\hline
		Gram groups \rule{0pt}{4ex} & [[2, 1], [5, 1], [6, 1], [12, 1]] &
		$1/r$ & \makecell[l]{6.474691494688162439\dots a root of \\
		$16x^4 - 32x^3 - 648x^2 + 664x + 3433$} \\
		\hline\Tstrut
		Polygons & [3,1] & $1/r^2$ & $17/4$
	\end{tabular}
\end{center}

\subsection{6 points}
The optimal solution for 6 points is an octahedron. A good alignment for this figure is along the coordinate axis.
\begin{figure}[h]
	\begin{center}
		\includegraphics[type=pdf,ext=pdf,read=pdf,height=1in,width=1in,angle=0]{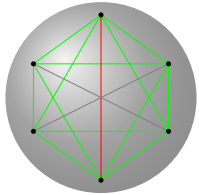}
		\caption{6 points.}
		\label{fig:6pts}
	\end{center}
\end{figure}
\begin{center}
	\begin{tabular}{c|ccc}
		\multicolumn{4}{c}{Spherical code} \\
		\hline\Tstrut
		pt & x & y & z \\
		\hline\Tstrut
		1 & 0 & 0 & 1 \\
		2 & 1 & 0 & 0 \\
		3 & -1 & 0 & 0 \\
		4 & 0 & 1 & 0 \\
		5 & 0 & -1 & 0 \\
		6 & 0 & 0 & -1
	\end{tabular}
\end{center}
\begin{center}
	\begin{tabular}{l|l||l|l}
		\multicolumn{2}{c}{Symmetries - 6 points} & \multicolumn{2}{c}{Minimal Energy} \\
		\hline\Tstrut
		planes & [[2, 4], [4, 3]] & log & $-log(512)$\Bstrut \\
		\hline\Tstrut
		\Tstrut Gram groups & [[6, 2], [24, 1]] &
		$1/r$ & $\frac{3}{2} + 6\sqrt{2} $ \Bstrut \\
		\hline\Tstrut
		Polygons & [[3, 8], [4, 3]] & $1/r^2$ & $27/4$
	\end{tabular}
\end{center}

\subsection{7 points}
The optimal solution for 7 points is a pentagon at the equator and two poles, or the pentagonal bipyramid. The minimal energy codes are given.
\begin{figure}[ht]
	\begin{center}
		\includegraphics[type=pdf,ext=pdf,read=pdf,height=1in,width=1in,angle=0]{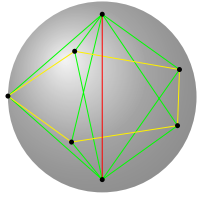}
		\caption{7 points.}
		\label{fig:7pts}
	\end{center}
\end{figure}
\begin{center}
	\begin{tabular}{c|ccc}
		\multicolumn{4}{c}{Spherical code -- 7 points} \\
		\hline\Tstrut
		pt & x & y & z \\
		\hline\Tstrut
		1 & 0 & 0 & 1 \\
		2 & 0 & 1 & 0 \\
		3 & $-\frac{\sqrt{10+2\sqrt{5}}}{4}$ & $\frac{\sqrt{5}-1)}{4}$ & 0\\[0.7ex]
		4 & $-\frac{\sqrt{10-2\sqrt{5}}}{4}$ & $-\frac{\sqrt{5}+1)}{4}$ & 0 \\[0.7ex]
		5 & $\frac{\sqrt{10-2\sqrt{5}}}{4}$ & $-\frac{\sqrt{5}+1)}{4}$ & 0 \\[0.7ex]
		6 & $\frac{\sqrt{10+2\sqrt{5}}}{4}$ & $\frac{\sqrt{5}-1)}{4}$ & 0 \\[0.7ex]
		7 & 0 & 0 & -1
	\end{tabular}
\end{center}
\begin{center}
	\begin{tabular}{l|l||l|l}
		\multicolumn{2}{c}{Symmetries - 7 points} & \multicolumn{2}{c}{Minimal Energy} \\
		\hline
		\rule{0pt}{3ex} planes & [[10, 1]] & log & $-\frac{log(12800000)}{2}$ \Bstrut \\
		\hline\Tstrut
		Gram groups & [[2, 1], [7, 1], [10, 2], [20, 1]] &
		$1/r$ & \rule{0pt}{6.7ex} \makecell[l]{14.45297741422134293\dots a root of \\ $256x^8 - 1024x^7 - 75008x^6 + 228608x^5$ \\ $+ 5537120x^4 - 11456448x^3 - 103335888x^2$ \\ $+ 109102384x - 23637199$ } \\
		\hline\Tstrut
		Polygons & [[5, 1]] & $1/r^2$ & $41/4$
	\end{tabular}
\end{center}

\subsection{8 points}
The optimal solution for 8 points is the square antiprism, or a cube with one face twisted 45 degrees with respect of the opposite face. It is not a cube, contrary to popular belief.

\begin{figure}[!ht]
	\begin{center}
		\includegraphics[type=pdf,ext=pdf,read=pdf,height=1in,width=1in,angle=0]{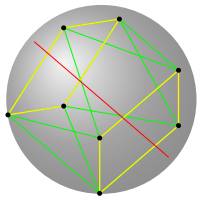}
		\caption{8 points.}
		\label{fig:8pts}
	\end{center}
\end{figure}

The 8 coordinates of the points for the global minimum are all different under the 3 potentials, but the spherical code can be easily expressed as a function of one algebraic parameter $a$ and two derived parameters $b = \sqrt{\frac{1-a^2}{2}}$ and $c = \sqrt{1-a^2}$.

\noindent
The parametric expression for 8 points is listed:
\begin{center}
	\begin{tabular}{c|ccc}
		\multicolumn{4}{c}{Spherical code} \\
		\hline\Tstrut
		pt & x & y & z \\
		\hline\Tstrut
		1 & $a$ & $b$ & $b$ \\
		2 & $a$ & $b$ & $-b$ \\
		3 & $a$ & $-b$ & $b$ \\
		4 & $a$ & $-b$ & $-b$ \\
		5 & $-a$ & 0 & $c$ \\
		6 & $-a$ & $c$ & 0 \\
		7 & $-a$ & 0 & $-c$ \\
		8 & $-a$ & $-c$ & 0
	\end{tabular}
\end{center}
We express the $a$ parameter and $energy$ as algebraic numbers as follows:

\noindent
\textit{\textbf{logarithmic potential -}}
\begin{align*}
a & = 0.5646169639331753669\dots \text{a root of } 7x^4 + 26x^2 - 9 \\
energy & = log\left(\frac{3114459466\sqrt{58}+23719027063}{1603087953297408}\right)
\end{align*}

\textit{\textbf{Coulomb $\mathbf{1/r}$ potential -}}
\begin{align*}
a = {} & 0.5604367652904311982\dots \text{ a root of }
 1227233x^{48} + 10470984x^{46} + 209556948x^{44} \\ &
 + 1689426104x^{42} + 15775473090x^{40} + 115047892248x^{38} + 620026175876x^{36} \\ & + 2363468454888x^{34} + 6203218661391x^{32} + 11155824048720x^{30} \\ & + 15270599090856x^{28} + 21228525244848x^{26} + 33209646716572x^{24} \\ & + 31438858405296x^{22} + 12676840005288x^{20} + 27288086224464x^{18} \\ & + 32464672800015x^{16} - 71697151292952x^{14} + 119763272586116x^{12} \\ & - 143004483949800x^{10} + 112981931600322x^8 - 72347775538504x^6 \\ & + 32198152383444x^4 - 7663346651064x^2 + 712364630753 \\
energy = {} & 19.67528786123276226\dots \text{ a root of }
x^{48} - 1008x^{46} + 403128x^{44} - 87692448x^{42} \\ & + 12008348040x^{40} - 1120658189376x^{38} + 74545764373344x^{36} \\ & - 3606531729475968x^{34} + 129269831646915312x^{32} - 4164949576246842880x^{30} \\ & + 141942516961279104768x^{28} - 6079710008497542226944x^{26} \\ & + 171029843453720580773632x^{24} - 10734215529571167885035520x^{22} \\ & + 119485611963873574006938624x^{20} - 7158961249326434415489593344x^{18} \\ & + 223410229304763643888886542080x^{16} + 1665579124457906626567818989568x^{14} \\ & + 360289207228162192431513880221696x^{12} \\ & + 4138484064890644161343737472278528x^{10} \\ & + 152548694458573411883645969649469440x^8 \\ & + 1581691619429822432587864824950079488x^6 \\ & + 6289231802234486775557087083681996800x^4 \\ & + 12513023795837806250674990759664713728x^2 \\ & + 16204114268229804228795137561699028992
\end{align*}

\textit{\textbf{Inverse square law $\mathbf{1/r^2}$ potential -}}

\begin{align*}
a & = 0.5563309621802899475\dots \text{ a root of } 11x^8 - 60x^6 - 158x^4 - 188x^2 + 75 \\
energy & = 14.33679108359450163\dots \text{ a root of } 64x^4 - 1408x^3 + 8144x^2 - 16432x + 6897
\end{align*}
The group symmetries and Gram matrix groups for the 8 points are identical for all 3 potentials.

\noindent
\textit{\textbf{Symmetries -}}
\begin{center}
	\begin{tabular}{l|l}
		\multicolumn{2}{c}{Symmetries - 8 points} \\
		\hline
		\rule{0pt}{3ex} planes & [[4, 8], [8, 1]] \\
		\hline
		\rule{0pt}{3ex} Gram groups & [[8, 2], [16, 3]] \\
		\hline
		\rule{0pt}{3ex} Polygons & [[4, 10]]
	\end{tabular}
\end{center}

\subsection{9 points}
The optimal solution for 9 points has 3 embedded equilateral triangles aligned along a common axis, the middle triangle is rotated 180$^\circ$ with respect to the other 2. This is shown in figure \ref{fig:9pts} below.
\begin{figure}[ht]
	\begin{center}
		\includegraphics[type=pdf,ext=pdf,read=pdf,height=1in,width=1in,angle=0]{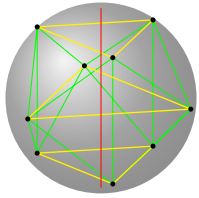}
		\caption{9 points.}
		\label{fig:9pts}
	\end{center}
\end{figure}
Like the optimal 8-point solution, the global solutions for 9 points, under all 3 potentials, are determinable from one algebraic parameter $a$.

\begin{longtable}[c]{c|ccc}
	\caption{Spherical code for 9 points} \\
	pt & $x$ & $y$ & $z$ \\
	\hline\vspace*{-2.2ex}
	\endfirsthead
	\multicolumn{4}{c}%
	{\tablename\ \thetable\ -- 9 points code -- \textit{continued}} \\
	pt & $x$ & $y$ & $z$ \\
	\hline\vspace*{-2.2ex}
	\endhead
	1 & 0 & $\sqrt{1-a^2}$ & $a$ \\[0.5ex]
	2 & $\frac{\sqrt{3}}{2}\sqrt{1-a^2}$ & $-\frac{\sqrt{1-a^2}}{2}$ & $a$ \\[1ex]
	3 & $-\frac{\sqrt{3}}{2}\sqrt{1-a^2}$ & $-\frac{\sqrt{1-a^2}}{2}$ & $a$ \\[1ex]
	4 & 0 & $-1$ & 0 \\[1ex]
	5 & $\frac{\sqrt{3}}{2}$ & $\frac{1}{2}$ & 0 \\[1ex]
	6 & $-\frac{\sqrt{3}}{2}$ & $\frac{1}{2}$ & 0 \\[1ex]
	7 & 0 & $\sqrt{1-a^2}$ & $-a$ \\[1ex]
	8 & $\frac{\sqrt{3}}{2}\sqrt{1-a^2}$ & $-\frac{\sqrt{1-a^2}}{2}$ & $-a$ \\[1ex]
	9 & $-\frac{\sqrt{3}}{2}\sqrt{1-a^2}$ & $-\frac{\sqrt{1-a^2}}{2}$ & $-a$
\end{longtable}

Again we express the $a$ parameter and $energy$ as algebraic numbers as follows:

\noindent
\textit{\textbf{logarithmic potential -}}
\begin{align*}
a & = 0.7031106068430678248\dots \text{a root of } 64x^8 + 105x^6 - 87x^4 - 45x^2 + 27 \\
\begin{split}
energy & = -12.88775272575927896\dots log( \text{a root of }
38311009236316727915402195098082798279121/ \\
& \qquad 16073449381704719530872738340965868632806598377472x^8 - 943440073904807815958/ \\
& \qquad 67136200068986103420023732084394928240591343540991466405888x^6 - 187182785938/ \\
& \qquad 03854967001728213902036260549626384901267496901275549696x^4 + 297604557206794/ \\
& \qquad 7154685831689083848055683186077146195085747x^2 + 5444517870735015415413993718/ \\
& \qquad 908291383296)
\end{split}
\end{align*}

\textit{\textbf{Coulomb $\mathbf{1/r}$ potential -}}
\textbf{Special Note:} The algebraic number $a$ for the $1/r$ Coulomb potential is degree 208. The Maxima program, when creating the Jacobians, indeed shows this exponential increase due to the pairing of vertices in the distance function, which increases as $2^n$. This makes finding the spherical codes for $n>50$ very difficult.
\begin{longtable}{l}
$a = 0.7036483958041317758\dots \text{ a root of }$ \\
\hphantom{000} $5227573613485916806405226496x^{208}$ \\
\hphantom{000} $ - \; \; 202492098767053611550667833344x^{206}$ \\
\hphantom{000} $ + \; 1635712891968416454830997356544x^{204}$ \\
\hphantom{000} $ + \; 16396822731607253984747139784704x^{202}$ \\
\hphantom{000} $ + \; 193015835879039113070530950538752x^{200}$ \\
\hphantom{000} $ - \; 12740560120959318844079308513602048x^{198}$ \\
\hphantom{000} $ + \; 66065492671088058986802507829937088x^{196}$ \\
\hphantom{000} $ + \; 1282498725954877734448708630566514848x^{194}$ \\
\hphantom{000} $ - \; 2666472074456316879582985542457095455x^{192}$ \\
\hphantom{000} $ - \; 163740764791297575959260674167115526560x^{190}$ \\
\hphantom{000} $ + \; 1841697492573886723673740348130760851760x^{188}$ \\
\hphantom{000} $ - \; 8753150354792477208666776190919886368032x^{186}$ \\
\hphantom{000} $ + \; 35831730265780789420875186338849301572040x^{184}$ \\
\hphantom{000} $ - \; 418511150503375126753022400479683759736736x^{182}$ \\
\hphantom{000} $ + \; 4723504967176370835104259160065678571230672x^{180}$ \\
\hphantom{000} $ - \; 32896125445825134349104369238966235831958432x^{178}$ \\
\hphantom{000} $ + \; 160250395374898385556443296103966472162186996x^{176}$ \\
\hphantom{000} $ - \; 777210588907367878918618525113530436484946208x^{174}$ \\
\hphantom{000} $ + \; 5176086281012420705340086488834699674179439408x^{172}$ \\
\hphantom{000} $ - \; 35552759164704704283191811880888377104290509984x^{170}$ \\
\hphantom{000} $ + \; 189526330881341940897249634575046057994136907128x^{168}$ \\
\hphantom{000} $ - \; 722623783766351451872044447499678949493274213664x^{166}$ \\
\hphantom{000} $ + \; 1850803721338323365546008024433094520683517724880x^{164}$ \\
\hphantom{000} $ - \; 2577173145419445019237889952625679507466949498464x^{162}$ \\
\hphantom{000} $ - \; 425041937830163455053965077556698888986767539422x^{160}$ \\
\hphantom{000} $ + \; 5370230173271408315640214033457807427976515499680x^{158}$ \\
\hphantom{000} $ + \; 21760215587746247543245341346365971920452921930192x^{156}$ \\
\hphantom{000} $ - \; 156825796042113370198023834121320626652051537915360x^{154}$ \\
\hphantom{000} $ + \; 370449477545798979633880321793729090073446999492376x^{152}$ \\
\hphantom{000} $ - \; 219416915478692787229076000486747272946597890447200x^{150}$ \\
\hphantom{000} $ - \; 833092552148915429448195890142602845769411197593552x^{148}$ \\
\hphantom{000} $ + \; 1265124935209099217546644317610758941413802465679264x^{146}$ \\
\hphantom{000} $ + \; 4239300915803771925710919681385862952162488372288484x^{144}$ \\
\hphantom{000} $ - \; 17543658013800437066549973649703060716901195687369952x^{142}$ \\
\hphantom{000} $ + \; 20167156394458672311134716750820304725731030381481040x^{140}$ \\
\hphantom{000} $ + \; 18251304228562827610817108694285513284780411444521632x^{138}$ \\
\hphantom{000} $ - \; 64486472782511145039963013337957007436040935686512472x^{136}$ \\
\hphantom{000} $ - \; 34221823328110813472282323726112977108481747606157024x^{134}$ \\
\hphantom{000} $ + \; 355271861877389569697471816974230362310237820234334832x^{132}$ \\
\hphantom{000} $ - \; 531883414015416640515093499120235800871433337922647296x^{130}$ \\
\hphantom{000} $ + \; 7618794661746601378495855728613132091024993696014383x^{128}$ \\
\hphantom{000} $ + \; 686847173756211027571126511361241055869984457265977280x^{126}$ \\
\hphantom{000} $ + \; 385275947456624085919698538812732099789374318046766560x^{124}$ \\
\hphantom{000} $ - \; 3538958792111643280196572711964415980913531236673579840x^{122}$ \\
\hphantom{000} $ + \; 4822133300609051851223198162901744259078518692895752016x^{120}$ \\
\hphantom{000} $ - \; 1819722467118346025132284468382013863814836964520653888x^{118}$ \\
\hphantom{000} $ + \; 2841842743985437493828416350489217420015259834924145696x^{116}$ \\
\hphantom{000} $ - \; 17703044162803129430928364050130001119004634487500480576x^{114}$ \\
\hphantom{000} $ + \; 31985658681155073097593351770669127329487040991197264872x^{112}$ \\
\hphantom{000} $ - \; 16597778835609587800637859688604858745601682258649464640x^{110}$ \\
\hphantom{000} $ - \; 6068940895157028997083965036209830269280825681306815520x^{108}$ \\
\hphantom{000} $ - \; 55918758584975433685426356887540492594585261693393965632x^{106}$ \\
\hphantom{000} $ + \; 228687951599110057364295036151588223530487086210413527856x^{104}$ \\
\hphantom{000} $ - \; 338277782357120087083167260352343443139757224028752717632x^{102}$ \\
\hphantom{000} $ + \; 184305228984461880679819790638385922484592597690395797024x^{100}$ \\
\hphantom{000} $ + \; 66739548908477144631366510678059578891857737702280159552x^{98}$ \\
\hphantom{000} $ + \; 124480365335767082131200311512379063357138542150867168860x^{96}$ \\
\hphantom{000} $ - \; 1013076507855007754531811501435122490370229930891625021120x^{94}$ \\
\hphantom{000} $ + \; 2012235163511936218327204670103911138526362729213761720224x^{92}$ \\
\hphantom{000} $ - \; 2190911192270389121819696508413705831193066778686328588224x^{90}$ \\
\hphantom{000} $ + \; 1459820418662214876314521025318535189821091567344251510704x^{88}$ \\
\hphantom{000} $ - \; 939619881171602265931059851392357296043319397249640868544x^{86}$ \\
\hphantom{000} $ + \; 1864327746590088208838083582979487273215357691638565451872x^{84}$ \\
\hphantom{000} $ - \; 4244529302038621191829193851648105579864919551414919559360x^{82}$ \\
\hphantom{000} $ + \; 6835473113085243964189148714277305428142602650922297577832x^{80}$ \\
\hphantom{000} $ - \; 8324565069765753066822132585523019815086477177179356608960x^{78}$ \\
\hphantom{000} $ + \; 8326493984070990445522216822953509110505650985699955057312x^{76}$ \\
\hphantom{000} $ - \; 7311965981067155607159912730202729914820872285506746989248x^{74}$ \\
\hphantom{000} $ + \; 5927571174969833388720376314540740448829553204930434513872x^{72}$ \\
\hphantom{000} $ - \; 4553994310183118854339980896662437176818602051368763688896x^{70}$ \\
\hphantom{000} $ + \; 3325341324561030971898949608384203806911064049152084584096x^{68}$ \\
\hphantom{000} $ - \; 2294648307177577403462538674106109621459203318233928490848x^{66}$ \\
\hphantom{000} $ + \; 1496520887549301632413059460021541708796960049747780536399x^{64}$ \\
\hphantom{000} $ - \; 929825137724837322362910689792431285720014699744646024224x^{62}$ \\
\hphantom{000} $ + \; 554859595505928725401929812437565280626006496483922989552x^{60}$ \\
\hphantom{000} $ - \; 318519891905174659227900924338757331942199310241038005664x^{58}$ \\
\hphantom{000} $ + \; 175404209026249043653246616132110843946909011595668442472x^{56}$ \\
\hphantom{000} $ - \; 92576142500067287170026103983098408338555844334337550880x^{54}$ \\
\hphantom{000} $ + \; 46964989531868527353444339376719480199872687998910807312x^{52}$ \\
\hphantom{000} $ - \; 22966026250475205780793863125483405647954136242788123168x^{50}$ \\
\hphantom{000} $ + \; 10819237422448745973316816227342475473966937437632779044x^{48}$ \\
\hphantom{000} $ - \; 4900655137945533650990114110775531980538337802035248544x^{46}$ \\
\hphantom{000} $ + \; 2134626644986201436591167418773300351739262588843288048x^{44}$ \\
\hphantom{000} $ - \; 895715034377167915797208433182517146030784943633667872x^{42}$ \\
\hphantom{000} $ + \; 362082992975930132665633494632415656328666325990751448x^{40}$ \\
\hphantom{000} $ - \; 140686161111196995693429547017865511161994453691879840x^{38}$ \\
\hphantom{000} $ + \; 52467164469723304606637106819226508827129657143182352x^{36}$ \\
\hphantom{000} $ - \; 18795329736364554003190711487462442933293199468946656x^{34}$ \\
\hphantom{000} $ + \; 6466978544030595942068892686027940180178737224698530x^{32}$ \\
\hphantom{000} $ - \; 2130601040898348166948434663171521723223101658698208x^{30}$ \\
\hphantom{000} $ + \; 670001508774632296475359960495536298619995433354640x^{28}$ \\
\hphantom{000} $ - \; 201046445760791244182987316490920919204103664565344x^{26}$ \\
\hphantom{000} $ + \; 57537623918123597895479298052392721942825071305400x^{24}$ \\
\hphantom{000} $ - \; 15626203658099356726334009357516964024057920165856x^{22}$ \\
\hphantom{000} $ + \; 4001518292751511821571280885373479037410443182192x^{20}$ \\
\hphantom{000} $ - \; 964644360566438025076459824632638299883394695392x^{18}$ \\
\hphantom{000} $ + \; 218721025628246802900104514738861883223145440436x^{16}$ \\
\hphantom{000} $ - \; 46139836976724455428599917340630156184877716320x^{14}$ \\
\hphantom{000} $ + \; 8900758654816724806973784311616075816075766800x^{12}$ \\
\hphantom{000} $ - \; 1559355579948554407769028730729076204079084000x^{10}$ \\
\hphantom{000} $ + \; 247661035207059081823969484559025591399125000x^8$ \\
\hphantom{000} $ - \; 34270581496560378496423284570998324895900000x^6$ \\
\hphantom{000} $ + \; 3706547514018742501148033218396620906750000x^4$ \\
\hphantom{000} $ - \; 262831328024290368138135058539014145000000x^2$ \\
\hphantom{000} $ + \; 8800156072241865004625057763583062890625$
\end{longtable}
$energy = {} 25.75998653126983156\ldots = \text{algebraic number, degree} > 208\text{ ?}$

\textit{\textbf{Inverse square law $\mathbf{1/r^2}$ potential -}}
\begin{align*}
\begin{split}
a & = 0.7046074370271068597\dots \text{ a root of } 324x^{18} + 9973x^{16} - 20004x^{14} - 24768x^{12} - 11484x^{10} \\ & \qquad + 187650x^8 - 159084x^6 + 44712x^4 - 17496x^2 + 6561 \\
energy & = 19.25286878398789640\dots \text{ a root of } 65536x^9 - 6488064x^8 + 270663680x^7 - 6200782848x^6 \\ & \qquad + 85555689984x^5 - 738608792832x^4 + 4038778392192x^3 - 13674172605264x^2 \\ & \qquad + 25867905614397x - 19342011471804
\end{split}
\end{align*}

\noindent
\textit{\textbf{Symmetries -}}
\begin{center}
	\begin{tabular}{l|l}
		\multicolumn{2}{c}{Symmetries - 9 points} \\
		\hline\Tstrut
		planes & [[3, 1], [4, 9]] \\[0.2ex]
		\hline\Tstrut
		Gram groups & [[6, 2], [9, 1], [12, 3], [24, 1]] \\
		\hline\Tstrut
		Polygons & [[3, 3], [4, 9]]
	\end{tabular}
\end{center}

\subsection{10 points}
The optimal solution for 10 points is a polyhedron which contains 2 embedded parallel squares and two end points or poles, which are aligned along the normal axis of the two squares. The squares are rotated $90^{\circ}$ with respect to each other's vertices, making an antiprism.

The optimal alignment axis was found by examining the sets of parallel-plane squares.
\begin{figure}[ht]
	\begin{center}
		\includegraphics[type=pdf,ext=pdf,read=pdf,height=1in,width=1in,angle=0]{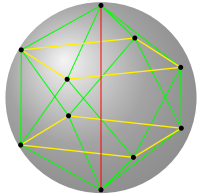}
		\caption{10 points.}
		\label{fig:10pts}
	\end{center}
\end{figure}
As in the case for 8 and 9 points, the optimal 10-point solutions, under all 3 potentials, are determinable from one algebraic parameter $a$.

The spherical code for 10 points, assuming an algebraic $a\in\mathcal{A}$, is given below.

\begin{longtable}[c]{c|ccc}
	\caption{Algebraic code for 10 points} \\
	pt & $x$ & $y$ & $z$ \\
	\hline\vspace*{-2.2ex}
	\endfirsthead
	\multicolumn{4}{c}%
	{\tablename\ \thetable\ -- Algebraic code 10 points -- \textit{continued}} \\
	pt & $x$ & $y$ & $z$ \\
	\hline\vspace*{-2.2ex}
	\endhead
	1 & 0 & 0 & 1 \\
	2 & 0 & $\sqrt{1-a^2}$ & $a$ \\[0.5ex]
	3 & 0 & $-\sqrt{1-a^2}$ & $a$ \\[0.5ex]
	4 & $\sqrt{1-a^2}$ & 0 & $a$ \\[0.5ex]
	5 & $-\sqrt{1-a^2}$ & 0 & $a$ \\[0.5ex]
	6 & $\sqrt{\frac{1-a^2}{2}}$ & $\sqrt{\frac{1-a^2}{2}}$ & $-a$ \\[0.7ex]
	7 & $\sqrt{\frac{1-a^2}{2}}$ & $-\sqrt{\frac{1-a^2}{2}}$ & $-a$ \\[0.7ex]
	8 & $-\sqrt{\frac{1-a^2}{2}}$ & $\sqrt{\frac{1-a^2}{2}}$ & $-a$ \\[0.7ex]
	9 & $-\sqrt{\frac{1-a^2}{2}}$ & $-\sqrt{\frac{1-a^2}{2}}$ & $-a$ \\[0.7ex]
	10 & 0 & 0 & -1
\end{longtable}

The algebraic numbers for the $a$ parameter and $energy$ are given next:

\noindent
\textit{\textbf{logarithmic potential -}}
\begin{align*}
a & = 0.4204838855379730022\dots \text{a root of } 9x^4 + 38x^2 - 7 \\
\begin{split}
energy & = -15.56312338902193911\dots log( \text{a root of }
50706024009129176059868128215040000000000x^2 \\ & \quad - 8832445142182596378429651924549632x + 150094635296999121)
\end{split}
\end{align*}

\textit{\textbf{Coulomb $\mathbf{1/r}$ potential -}}
The algebraic number $a$ for the $1/r$ Coulomb potential is degree 224. The degree of the algebraic number for the energy is unknown.
\begin{longtable}{l}
	$a = 0.4226874240439860662\dots \text{ a root of }$ \\
	\hphantom{000} $ 4294967296x^{224}$ \\
	\hphantom{000} $ + \; 2458868776960x^{222}$ \\
	\hphantom{000} $ + \; 1228624183164928x^{220}$ \\
	\hphantom{000} $ + \; 393776570801586176x^{218}$ \\
	\hphantom{000} $ + \; 97567522372779900928x^{216}$ \\
	\hphantom{000} $ + \; 10416129061243948859392x^{214}$ \\
	\hphantom{000} $ + \; 672902660086655438019584x^{212}$ \\
	\hphantom{000} $ + \; 30111292352687221496374400x^{210}$ \\
	\hphantom{000} $ + \; 1068457611997512842663975617x^{208}$ \\
	\hphantom{000} $ + \; 31676697179069531247453849952x^{206}$ \\
	\hphantom{000} $ + \; 846094263704413332601486841296x^{204}$ \\
	\hphantom{000} $ + \; 21002899296096951269176738959968x^{202}$ \\
	\hphantom{000} $ + \; 481140804389145695369741359799400x^{200}$ \\
	\hphantom{000} $ + \; 9943523475991007656584012813224864x^{198}$ \\
	\hphantom{000} $ + \; 182294560088379576330705880425562480x^{196}$ \\
	\hphantom{000} $ + \; 2942546681779552607510655549727930080x^{194}$ \\
	\hphantom{000} $ + \; 41775448400487384505354245671797251860x^{192}$ \\
	\hphantom{000} $ + \; 522691883378839588430975033199887622304x^{190}$ \\
	\hphantom{000} $ + \; 5782227663387798904865923846487902190160x^{188}$ \\
	\hphantom{000} $ + \; 56751855616590100909028140821772000760800x^{186}$ \\
	\hphantom{000} $ + \; 495851939646000872744516744964533213876120x^{184}$ \\
	\hphantom{000} $ + \; 3868398756758797125597333941912969768214048x^{182}$ \\
	\hphantom{000} $ + \; 27018057375015865054899175027727076029994864x^{180}$ \\
	\hphantom{000} $ + \; 169265797647881768079165707334155069510701920x^{178}$ \\
	\hphantom{000} $ + \; 951893815198086820225796622184185279723382018x^{176}$ \\
	\hphantom{000} $ + \; 4796738551580872064182344002463481878510303328x^{174}$ \\
	\hphantom{000} $ + \; 21512647204482169427708388955431594087607988528x^{172}$ \\
	\hphantom{000} $ + \; 84409180380589476711576891091370221650078967968x^{170}$ \\
	\hphantom{000} $ + \; 278178477826317517616598344024655600814433021816x^{168}$ \\
	\hphantom{000} $ + \; 687660031690059811886749816622978755339355986272x^{166}$ \\
	\hphantom{000} $ + \; 687211663675273626333324170743800306840925771152x^{164}$ \\
	\hphantom{000} $ - \; 4590648804119487190034633037004472546126223230944x^{162}$ \\
	\hphantom{000} $ - \; 35033566058817210905131026928635038139477395944508x^{160}$ \\
	\hphantom{000} $ - \; 140294888863422320715147328463554282237107325257376x^{158}$ \\
	\hphantom{000} $ - \; 352430014234513986404022520172282351966538368561360x^{156}$ \\
	\hphantom{000} $ - \; 306477899230089265204616954705696149017425001790432x^{154}$ \\
	\hphantom{000} $ + \; 2019517627337898274559234203278267267966273277417096x^{152}$ \\
	\hphantom{000} $ + \; 12450571073991588564435898598623150920882784228714976x^{150}$ \\
	\hphantom{000} $ + \; 40024530453463030626839489236569653734857796726271504x^{148}$ \\
	\hphantom{000} $ + \; 82937651675831053559936135615603933793670592226669856x^{146}$ \\
	\hphantom{000} $ + \; 92293661613905299543060384563270794177783627492836911x^{144}$ \\
	\hphantom{000} $ - \; 53810200002206943411201101986259954762297718337768704x^{142}$ \\
	\hphantom{000} $ - \; 457186029843595686973226500422842071231577644334686432x^{140}$ \\
	\hphantom{000} $ - \; 919722879740359331619486552347987101265569439199381568x^{138}$ \\
	\hphantom{000} $ - \; 786685992076974047916907957113632317699877970620534384x^{136}$ \\
	\hphantom{000} $ + \; 398069072912919404130496617795497527525257905335985216x^{134}$ \\
	\hphantom{000} $ + \; 1635992559313482562953250722348820536543468538828349536x^{132}$ \\
	\hphantom{000} $ + \; 1082201488171269375598004998864829312919938889786905792x^{130}$ \\
	\hphantom{000} $ + \; 46534424483997810750785921878271725584482576873879848x^{128}$ \\
	\hphantom{000} $ + \; 4060999172285328139690270679634291506666126858118249024x^{126}$ \\
	\hphantom{000} $ + \; 12514944804257162029459072218968128554908255378809118496x^{124}$ \\
	\hphantom{000} $ + \; 9946276477842177111594332751804629055365839262743743168x^{122}$ \\
	\hphantom{000} $ - \; 13424082347695448789217825475890077109593049160049628048x^{120}$ \\
	\hphantom{000} $ - \; 36288282801695549983355930360583417973205344236718315200x^{118}$ \\
	\hphantom{000} $ - \; 41147334428448495670257448086834058097031385976859517344x^{116}$ \\
	\hphantom{000} $ - \; 48991162685311338099819120926865316972528591294224250432x^{114}$ \\
	\hphantom{000} $ - \; 28124016532815903322786368239095790960430785203864361956x^{112}$ \\
	\hphantom{000} $ + \; 116078127584179642701535728548632504960048794375598869184x^{110}$ \\
	\hphantom{000} $ + \; 253595295998036067047089597081044449688841056826525479520x^{108}$ \\
	\hphantom{000} $ + \; 92597825762910220427728363565922673388681046683422948672x^{106}$ \\
	\hphantom{000} $ - \; 137361163526500395779452006471767279326226476035404014480x^{104}$ \\
	\hphantom{000} $ + \; 644899702101790082817343240432460560596356419098969792x^{102}$ \\
	\hphantom{000} $ - \; 9330922704029552738205487163805417075010586027973679328x^{100}$ \\
	\hphantom{000} $ - \; 575757471746248826408062171398756294630320914095133291456x^{98}$ \\
	\hphantom{000} $ - \; 696592896261855298134588128636607064932319425832102180568x^{96}$ \\
	\hphantom{000} $ - \; 15874336691325881235781357479724507841230033840723447104x^{94}$ \\
	\hphantom{000} $ + \; 906537547679743507972036231261750143810023975136036001888x^{92}$ \\
	\hphantom{000} $ + \; 1684398895125160168531270816721441458150424113874344312896x^{90}$ \\
	\hphantom{000} $ + \; 479973121410231599069457466506033248373523201326802780560x^{88}$ \\
	\hphantom{000} $ - \; 1307756259338011300798931695518884582412823667625507415104x^{86}$ \\
	\hphantom{000} $ - \; 565113266111116835016564448182314677333008705800032265696x^{84}$ \\
	\hphantom{000} $ - \; 1627385428030381626788372867771144139957172351791073509696x^{82}$ \\
	\hphantom{000} $ - \; 688446724905494068394615343618906547290594605162339433937x^{80}$ \\
	\hphantom{000} $ + \; 2733950882606103239928731707597901164994741104207316913504x^{78}$ \\
	\hphantom{000} $ - \; 1632402289087606611335203055390685104151753502675618716912x^{76}$ \\
	\hphantom{000} $ + \; 2894190988896132768982689299303050431326559617076040822752x^{74}$ \\
	\hphantom{000} $ + \; 1982610080883240293441144156426932051204710099304364678792x^{72}$ \\
	\hphantom{000} $ - \; 10539662547455100693577419890771076371996721589230309243360x^{70}$ \\
	\hphantom{000} $ + \; 18118544960814181129928639610705841404599988861865542574896x^{68}$ \\
	\hphantom{000} $ - \; 29206604067745604715995115213371883246366666487232027519392x^{66}$ \\
	\hphantom{000} $ + \; 38519849314899652729701590249900875045418761733805389374404x^{64}$ \\
	\hphantom{000} $ - \; 42834768091091310096268230277103406849805858704138109648096x^{62}$ \\
	\hphantom{000} $ + \; 44648015999758750386311763988728987779765371875397212414352x^{60}$ \\
	\hphantom{000} $ - \; 42386750797107057003830169862590804041490523607155394081952x^{58}$ \\
	\hphantom{000} $ + \; 36453222423641619782737742887442316053279861345211899524984x^{56}$ \\
	\hphantom{000} $ - \; 29417666627083637623467401586445185290280328990840821993824x^{54}$ \\
	\hphantom{000} $ + \; 22092115031686153975799887009050699133792082895567728593712x^{52}$ \\
	\hphantom{000} $ - \; 15300523238889514546211602293886036861602743446733449120544x^{50}$ \\
	\hphantom{000} $ + \; 9917085582769303612510060944737969459717315661498587067650x^{48}$ \\
	\hphantom{000} $ - \; 6003259073427820950219536946205265231495240713846152100128x^{46}$ \\
	\hphantom{000} $ + \; 3357031939680247464479895859787191675165596621077637070192x^{44}$ \\
	\hphantom{000} $ - \; 1740589188214655897138171596860282199704712736477922037728x^{42}$ \\
	\hphantom{000} $ + \; 836822020936025813494453975967351783808669933148480628632x^{40}$ \\
	\hphantom{000} $ - \; 367970482712309166305388753700392239947769079581914813984x^{38}$ \\
	\hphantom{000} $ + \; 146656157391374477619360138979290311037452981938710314576x^{36}$ \\
	\hphantom{000} $ - \; 52814763361922332001294142442829579865480741641174361696x^{34}$ \\
	\hphantom{000} $ + \; 16883286677068925213661898054627689453103721361205929748x^{32}$ \\
	\hphantom{000} $ - \; 4614101566794821969796668544492339160449171715787278880x^{30}$ \\
	\hphantom{000} $ + \; 1034600431715262197835110627360222379900830013858930032x^{28}$ \\
	\hphantom{000} $ - \; 184326743519123553080740385630059705658629387570654816x^{26}$ \\
	\hphantom{000} $ + \; 25516906250852194407965339672129141594482307113972840x^{24}$ \\
	\hphantom{000} $ - \; 2698277947807530325223307345346774484494006895885216x^{22}$ \\
	\hphantom{000} $ + \; 214553027716969793135534686425158221720204521134800x^{20}$ \\
	\hphantom{000} $ - \; 12620001246779136961826002476430533582706257275488x^{18}$ \\
	\hphantom{000} $ + \; 540051885378925822270452451876948851438441621185x^{16}$ \\
	\hphantom{000} $ - \; 16580231821234125940385941529627306493735261120x^{14}$ \\
	\hphantom{000} $ + \; 362222806767905100981163591225811065476663040x^{12}$ \\
	\hphantom{000} $ - \; 5593939457139973612037719110654963251580928x^{10}$ \\
	\hphantom{000} $ + \; 60348011980880725501571311652806659432448x^{8}$ \\
	\hphantom{000} $ - \; 443647540128519185709379165378893316096x^{6}$ \\
	\hphantom{000} $ + \; 2114266757533773534832009081416318976x^{4}$ \\
	\hphantom{000} $ - \; 5881401167590744690846765982679040x^{2}$ \\
	\hphantom{000} $ + \; 7243565730450179340973882998784$ \\
	\\
	$energy = 32.71694946014758813\dots\text{high degree algebraic number}$
\end{longtable}

\textit{\textbf{Inverse square law $\mathbf{1/r^2}$ potential -}}
\begin{align*}
\begin{split}
a & = 0.4242756082881730876\dots \text{ a root of } 3x^8 - 156x^6 - 462x^4 - 284x^2 + 67 \\
energy & = 25.04135972210499009\dots \text{ a root of } 256x^4 - 9984x^3 + 102432x^2 - 327568x + 82881
\end{split}
\end{align*}

\noindent
\textit{\textbf{Symmetries -}}

\begin{center}
	\begin{longtable}{l|l}
		\caption{Symmetries -- 10 points} \\
		\hline\Tstrut
		planes & [[4, 12], [8, 1]] \\
		\hline\Tstrut
		Gram groups & [[2, 1], [8, 1], [10, 1], [16, 5]] \\
		\hline\Tstrut
		Polygons & [[4, 14]]
	\end{longtable}
\end{center}

\subsection{11 points}
The optimal solution for 11 points is the first one for which no algebraic number has been yet found for the coordinates, the algebraic degree is too high for current mathematical tools to discover, although 50,014 digits of the parameters for the polyhedron have been computed (making the spherical codes 50,014 digits in size).

\begin{figure}[ht]
	\begin{center}
		\includegraphics[type=pdf,ext=pdf,read=pdf,height=1in,width=1in,angle=0]{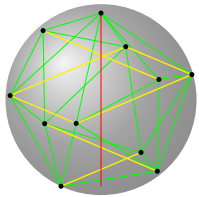}
		\caption{11 points.}
		\label{fig:11pts}
	\end{center}
\end{figure}

The suggested optimal axis was found by examining the parallel planes, aligning the squares and the dipoles in the figure \ref{fig:11pts}.

The figure contains one isolated pole, one dipole, one rectangle, and two dipoles as one proceeds away from the isolated pole along the preferred axis.
The top and middle dipole are aligned parallel with each other, the bottom dipole is rotated $90^{\circ}$ with respect to the other two.

In parameterizing this configuration, the distance between the dipoles and rectangle were used, as well as the height/width ratio of the rectangle, resulting in 5 parameters:$a$, $b$, $c$, $d$, and $e$.

The resulting parameterization, in terms of $a$, $b$, $c$, $d$, and $e$, is
\begin{center}
	\begin{tabular}{c|ccc}
		\multicolumn{4}{c}{Spherical code for 11 points} \\
		\hline\Tstrut
		pt & x & y & z \\
		\hline\Tstrut
		1 & 0 & 0 & 1 \\[0.5ex]
		2 & $\sqrt{1-a^2}$ & 0 & a \\[0.5ex]
		3 & -$\sqrt{1-a^2}$ & 0 & a \\[0.5ex]
		4 & c & $\sqrt{1-b^2-c^2}$ & b \\[0.5ex]
		5 & c & -$\sqrt{1-b^2-c^2}$ & b \\[0.5ex]
		6 & -c & $\sqrt{1-b^2-c^2}$ & b \\[0.5ex]
		7 & -c & -$\sqrt{1-b^2-c^2}$ & b \\[0.5ex]
		8 & $\sqrt{1-d^2}$ & 0 & d \\[0.5ex]
		9 & -$\sqrt{1-d^2}$ & 0 & d \\[0.5ex]
		10 & 0 & $\sqrt{1-e^2}$ & e \\[0.5ex]
		11 & 0 & -$\sqrt{1-e^2}$ & e
	\end{tabular}
\end{center}
The decimal values for the parameters $a$, $b$, $c$, $d$, and $e$ were optimized by using the Jacobian matrix and Newton's method to 50,014 digits and attempts were made to locate the algebraic numbers, but they was unsuccessful to a degree of 360.
{\small
\begin{center}
	\begin{tabular}{l|ccccc||cl}
		Pot+ & \multicolumn{5}{c}{Parameters - 11 points} & \multicolumn{2}{c}{Minimal Energy} \\
		\hline\Tstrut
		Law & a & b & c & d & e & energy \\[0.1ex]
		\hline\Tstrut
		-log & 0.5176071807 & 0.1729166133 & 0.4856631610 & -0.5583040204 & -0.8051363870 & -18.42047972 \\
		\hline\Tstrut
		$1/r$ & 0.5153071412 & 0.1683597004 & 0.4889224036 & -0.5526814629 & -0.8059548965 & 40.59645050 \\
		\hline\Tstrut
		$1/r^2$ & 0.5131231311 & 0.1619382013 & 0.4919211169 & -0.5483209203 & -0.8074504105 & 31.83472163
	\end{tabular}
\end{center}
}

\noindent
\textit{\textbf{Symmetries -}}
\begin{center}
	\begin{tabular}{l|l}
		\multicolumn{2}{c}{Symmetries - 11 points} \\
		\hline\Tstrut
		planes & [[4, 7], [10, 1]] \\[0.2ex]
		\hline\Tstrut
		Gram groups & [[2, 3], [4, 8], [8, 9], [11, 1]] \\
		\hline\Tstrut
		Polygons & [[4, 7], [5, 1]]
	\end{tabular}
\end{center}

\subsection{12 points}
The optimal solution for 12 points is the regular icosahedron, one of the Platonic solids. The multiple symmetries embedded in this polyhedron are noticeable when examining the group symmetries and Gram matrix, for example, there are 40 triangles, 15 squares and 12 pentagons embedded in the configuration.

\begin{figure}[ht]
	\begin{center}
		\includegraphics[type=pdf,ext=pdf,read=pdf,height=1in,width=1in,angle=0]{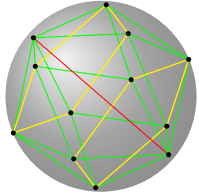}
		\caption{12 points - icosahedron.}
		\label{fig:12pts}
	\end{center}
\end{figure}

The spherical code for all 3 potentials is the same:

\noindent
\textit{\textbf{Spherical Code -}}
\begin{longtable}[c]{c|ccc}
	\caption{Algebraic code for 12 points} \\
	pt & $x$ & $y$ & $z$ \\
	\hline\vspace*{-2.2ex}
	\endfirsthead
	\multicolumn{4}{c}%
	{\tablename\ \thetable\ -- Algebraic code 12 points -- \textit{continued}} \\
	pt & $x$ & $y$ & $z$ \\
	\hline\vspace*{-2.2ex}
	\endhead
 	1 & 1 & 0 & 0 \\[0.5ex]
 	2 & $\sqrt{\frac{1}{5}}$ & 0 & $\frac{2}{\sqrt{5}}$ \\[0.7ex]
	3 & $\sqrt{\frac{1}{5}}$ & -$\sqrt{\frac{5+\sqrt{5}}{10}}$ & $\frac{5-\sqrt{5}}{10}$ \\[0.7ex]
	4 & $\sqrt{\frac{1}{5}}$ & -$\sqrt{\frac{5-\sqrt{5}}{10}}$ & $\frac{-5-\sqrt{5}}{10}$ \\[0.7ex]
	5 & $\sqrt{\frac{1}{5}}$ & $\sqrt{\frac{5-\sqrt{5}}{10}}$ & $\frac{-5-\sqrt{5}}{10}$ \\[0.7ex]
	6 & $\sqrt{\frac{1}{5}}$ & $\sqrt{\frac{5+\sqrt{5}}{10}}$ & $\frac{5-\sqrt{5}}{10}$ \\[0.7ex]
	7 & -$\sqrt{\frac{1}{5}}$ & $\sqrt{\frac{5-\sqrt{5}}{10}}$ & $\frac{5+\sqrt{5}}{10}$ \\[0.7ex]
	8 & -$\sqrt{\frac{1}{5}}$ & -$\sqrt{\frac{5-\sqrt{5}}{10}}$ & $\frac{5+\sqrt{5}}{10}$ \\[0.7ex]
	9 & -$\sqrt{\frac{1}{5}}$ & -$\sqrt{\frac{5+\sqrt{5}}{10}}$ & $\frac{-5+\sqrt{5}}{10}$ \\[0.7ex]
	10 & -$\sqrt{\frac{1}{5}}$ & 0 & -$\frac{2}{\sqrt{5}}$ \\[0.7ex]
	11 & -$\sqrt{\frac{1}{5}}$ & $\sqrt{\frac{5+\sqrt{5}}{10}}$ & $\frac{-5+\sqrt{5}}{10}$ \\[0.7ex]
	12 & -1 & 0 & 0
\end{longtable}

\begin{center}
	\begin{tabular}{l|c}
		Law & Minimal energy \\[0.5ex]
		\hline\Tstrut
		log & $log(30517578125/73786976294838206464)$ \\
		\hline
		$1/r$ & $3+15\sqrt{2\sqrt{5}+5}$ \mystrut(13,0) \\
		\hline\Tstrut
		$1/r^2$ & 39
	\end{tabular}
\end{center}

\noindent
\textit{\textbf{Symmetries -}}
\begin{center}
	\begin{tabular}{l|l}
		\multicolumn{2}{c}{Symmetries - 12 points} \\
		\hline\Tstrut
		planes & [[4, 25], [20, 6]] \\[0.2ex]
		\hline\Tstrut
		Gram groups & [[12, 2], [60, 2]] \\
		\hline\Tstrut
		Polygons & [[3, 40], [4, 15], [5, 12]]
	\end{tabular}
\end{center}

\subsection{13 points}
The optimal configuration for 13 points is interesting, there is one isolated pole, two dipoles, then a rectangle, then two more dipoles, moving away from the isolated pole. The dipoles, if next to each other, are rotated $90^{\circ}$ with respect to each other. The arrangement is [1:2:2:4:2:2].

\begin{figure}[!ht]
	\begin{center}
		\includegraphics[type=pdf,ext=pdf,read=pdf,height=1in,width=1in,angle=0]{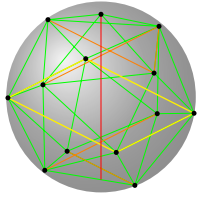}
		\caption{13 points.}
		\label{fig:13pts}
	\end{center}
\end{figure}

Six parameters, $a$, $b$, $c$, $d$, $e$, $f$, were used in an attempt to find the algebraic numbers.
\begin{longtable}[c]{c|ccc}
	\caption{Spherical code for 13 points} \\
	pt & $x$ & $y$ & $z$ \\
	\hline\vspace*{-2.2ex}
	\endfirsthead
	\multicolumn{4}{c}%
	{\tablename\ \thetable\ -- Spherical code 13 points -- \textit{continued}} \\
	pt & $x$ & $y$ & $z$ \\
	\hline\vspace*{-2.2ex}
	\endhead
	1 & 0 & 0 & 1 \\[0.5ex]
	2 & $\sqrt{1-a^2}$ & 0 & a \\[0.5ex]
	3 & -$\sqrt{1-a^2}$ & 0 & a \\[0.5ex]
	4 & 0 & $\sqrt{1-b^2}$ & b \\[0.5ex]
	5 & 0 & $-\sqrt{1-b^2}$ & b \\[0.5ex]
	6 & d & $\sqrt{1-c^2-d^2}$ & c \\[0.5ex]
	7 & d & -$\sqrt{1-c^2-d^2}$ & c \\[0.5ex]
	8 & -d & $\sqrt{1-c^2-d^2}$ & c \\[0.5ex]
	9 & -d & -$\sqrt{1-c^2-d^2}$ & c \\[0.5ex]
	10 & 0 & $\sqrt{1-e^2}$ & e \\[0.5ex]
	11 & 0 & -$\sqrt{1-e^2}$ & e \\[0.5ex]
	12 & $\sqrt{1-f^2}$ & 0 & e \\[0.5ex]
	13 & -$\sqrt{1-f^2}$ & 0 & e
\end{longtable}

The values to 19 digits of the parameters optimized for the minimal solution are:

\begin{center}
	\begin{tabular}{c|c|c|c}
		Parameter & log & 1/r & $1/r^2$ \\
		\hline\Tstrut
		a & 0.6104352212015100601 & 0.6112932937024949408 & 0.6107598627096058981 \\
		b & 0.5000426887660601445 & 0.5006051736217342854 & 0.5006390420498682103 \\
		c & -0.1021508091983655130 & -0.1043503753542226428 & -0.1073338258910784390 \\
		d & 0.8062500285826657195 & 0.8078028912837569615 & 0.8086348537465551129 \\
		e & -0.5441180397131400020 & -0.5446310313223069086 & -0.5458793128928694460 \\
		f & -0.8620582518576991767 & -0.8629768688323714373 & -0.8639426921905183161 \\
		\hline\Tstrut
		energy & -24.86672187550609084 & 58.85323061170244816 & 47.77330898085704154
	\end{tabular}
\end{center}
They are all known to 50,014 digits precision, for all 3 potentials, but again, attempts to find the algebraic numbers for the parameters failed, the algebraic degree $> 360$.

\noindent
\textit{\textbf{Symmetries -}}
\begin{center}
	\begin{tabular}{l|l}
		\multicolumn{2}{c}{Symmetries - 13 points} \\
		\hline\Tstrut
		planes & [[4, 9], [10, 2]] \\[0.2ex]
		\hline\Tstrut
		Gram groups & [[2, 4], [4, 11], [8, 13], [13, 1]] \\
		\hline\Tstrut
		Polygons & [[4, 9], [5, 2]]
	\end{tabular}
\end{center}

\subsection{14 points}
The optimal solution for 14 points is interesting, it is two poles and two parallel-plane hexagons whose normal is oriented along the pole axis. The 2 hexagons are rotated $30^{\circ}$ with respect to each other. The arrangement is [1:6:6:1].

\begin{figure}[h]
	\begin{center}
		\includegraphics[type=pdf,ext=pdf,read=pdf,height=1in,width=1in,angle=0]{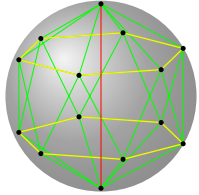}
		\caption{14 points.}
		\label{fig:14pts}
	\end{center}
\end{figure}
Only 3 parameters, $a$, $b$, and $c$ were required to parameterize the polyhedron in order to locate the algebraic numbers.
\begin{longtable}[c]{c|ccc}
	\caption{Algebraic code for 14 points} \\
	pt & $x$ & $y$ & $z$ \\
	\hline\vspace*{-2.2ex}
	\endfirsthead
	\multicolumn{4}{c}%
	{\tablename\ \thetable\ -- Algebraic code 14 points -- \textit{continued}} \\
	pt & $x$ & $y$ & $z$ \\
	\hline\vspace*{-2.2ex}
	\endhead
	1 & 0 & 0 & 1 \\[0.5ex]
	2 & $b$ & $-\sqrt{1-a^2-b^2}$ & $a$ \\[0.5ex]
	3 & $-b$ & $\sqrt{1-a^2-b^2}$ & $a$ \\[0.5ex]
	4 & $b$ & $\sqrt{1-a^2-b^2}$ & $a$ \\[0.5ex]
	5 & $-b$ & $-\sqrt{1-a^2-b^2}$ & $a$ \\[0.5ex]
	6 & $\sqrt{1-a^2}$ & 0 & $a$ \\[0.5ex]
	7 & $-\sqrt{1-a^2}$ & 0 & $a$ \\[0.5ex]
	8 & 0 & $-\sqrt{1-a^2}$ & $-a$ \\[0.5ex]
	9 & 0 & $\sqrt{1-a^2}$ & $-a$ \\[0.5ex]
	10 & $c$ & $-\sqrt{1-a^2-c^2}$ & $-a$ \\[0.5ex]
	11 & $-c$ & $\sqrt{1-a^2-c^2}$ & $-a$ \\[0.5ex]
	12 & $c$ & $\sqrt{1-a^2-c^2}$ & $-a$ \\[0.5ex]
	13 & $-c$ & $-\sqrt{1-a^2-c^2}$ & $-a$ \\[0.5ex]
	14 & 0 & 0 & -1
\end{longtable}
While the algebraic numbers were successfully recovered for the \textit{logarithmic} and \textit{inverse square law} potentials, the attempts to locate the algebraic polynomials for the $1/r$ potential was unsuccessful, the degree is $> 360$. A more efficient LLL matrix factoring program is needed.

\noindent
\textit{\textbf{logarithmic potential -}}
\begin{align*}
a & = 0.4591508204907729375\dots \text{a root of } 13x^6 + 159x^4 + 75x^2 - 23 \\
b & = 0.4441791654396483527\dots \text{a root of } 26x^6 - 99x^4 + 54x^2 - 7 \\
c & = 0.7693408822050128806\dots \text{a root of } 26x^6 - 297x^4 + 486x^2 - 189
\end{align*}
\begin{align*}
\begin{split}
energy & = -28.40781300924246057\dots log( \text{a root of }
4169296310946615347081896890345435474217/ \\ & \qquad 625419757666259830598370852784315823192280014758230283939998256898644349920/ \\ & \qquad 9490432x^3 - 19173819271508477744387320310627874199431928762462502133020630604/ \\ & \qquad 428505798411098847020422132399316243664863232x^2 + 4671594209038318175565468/ \\ & \qquad 7111676478249362798441801853034094278841304536804968718712384258048x \\ & \qquad - 27783742160348572763840067510872319734178277)
\end{split}
\end{align*}

\noindent
\textit{\textbf{Coulomb $\mathbf{1/r}$ potential -}}

No algebraic numbers were recovered for the 3 parameters under the Coulomb potential, to algebraic degree 360. Neither was the energy polynomial recovered.
\begin{center}
	\begin{tabular}{l|ccc}
		Potential & \multicolumn{3}{c}{Parameters - 14 points} \\
		\hline\Tstrut
		Law & $a$ & $b$ & $c$ \\[0.1ex]
		\hline\Tstrut
		$1/r$ & 0.4553677951624630035 & 0.4451517076033959054 & 0.7710253746451266125
	\end{tabular}
\end{center}
\indent
$energy = {} 69.30636329662642126\ldots = \text{algebraic number, degree} > 360\text{ ?}$
The decimal values for $a$, $b$, $c$ and $energy$ are known to 50,014 digits.

\noindent
\textit{\textbf{Inverse square law $\mathbf{1/r^2}$ potential -}}
\begin{align*}
\begin{split}
a & = 0.4518625916952697588\dots \text{ a root of } 49x^{12} - 1266x^{10} - 12273x^8 - 31132x^6 \\ & \qquad - 15729x^4 - 1266x^2 + 1201 \\
b & = 0.4460437753815296572\dots \text{ a root of } 196x^{12} + 972x^{10} - 4467x^8 + 5744x^6 - 3042x^4 \\ & \qquad + 708x^2 - 59 \\
c & = 0.7725704813606493514\dots \text{ a root of } 196x^{12} + 2916x^{10} - 40203x^8 + 155088x^6 \\ & \qquad - 246402x^4 + 172044x^2 - 43011 \\
energy & = 57.12120859837665927\dots \text{ a root of } 4096x^6 - 516096x^5 + 22827648x^4 - 446078592x^3 \\ & \qquad + 3791794905x^2 - 12207846852x + 5692229883
\end{split}
\end{align*}

\noindent
\textit{\textbf{Symmetries -}}
\begin{center}
	\begin{tabular}{l|l}
		\multicolumn{2}{c}{Symmetries - 14 points} \\
		\hline\Tstrut
		planes & [[4, 42], [40, 1]] \\[0.2ex]
		\hline\Tstrut
		Gram groups & [[2, 1], [12, 1], [14, 1], [24, 7]] \\
		\hline\Tstrut
		Polygons & [[4, 42], [6, 2]]
	\end{tabular}
\end{center}

\subsection{15 points}
The optimal configuration for 15 points is quite interesting, as it consists of 5 sets of parallel-plane equilateral triangles, alternately rotated $180^{\circ}$ with respect to their neighbor along the axis normal to them all. The arrangement is [3:3:3:3:3].

\begin{figure}[ht]
	\begin{center}
		\includegraphics[type=pdf,ext=pdf,read=pdf,height=1in,width=1in,angle=0]{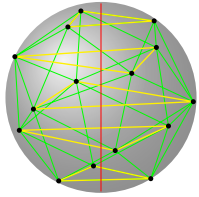}
		\caption{15 points.}
		\label{fig:15pts}
	\end{center}
\end{figure}

The configuration of 15 points is parameterized by 4 parameters, $a$, $b$, $c$, and $d$ as given below:

\noindent
\textit{\textbf{Parameterized Structure -}}
\begin{longtable}[c]{c|ccc}
	\caption{Spherical code for 15 points} \\
	pt & $x$ & $y$ & $z$ \\
	\hline\vspace*{-2.2ex}
	\endfirsthead
	\multicolumn{4}{c}%
	{\tablename\ \thetable\ -- Spherical code 15 points -- \textit{continued}} \\
	pt & $x$ & $y$ & $z$ \\
	\hline\vspace*{-2.2ex}
	\endhead
	1 & $a$ & $b$ & $\sqrt{1-a^2-b^2}$ \\[0.5ex]
	2 & $\frac{-a+b\sqrt{3}}{2}$ & $\frac{-b-a\sqrt{3}}{2}$ & $\sqrt{1-a^2-b^2}$ \\[0.5ex]
	3 & $\frac{-a-b\sqrt{3}}{2}$ & $\frac{-b+a\sqrt{3}}{2}$ & $\sqrt{1-a^2-b^2}$ \\[0.5ex]
	4 & $-c$ & $d$ & $\sqrt{1-c^2-d^2}$ \\[0.5ex]
	5 & $\frac{c+d\sqrt{3}}{2}$ & $\frac{-d+c\sqrt{3}}{2}$ & $\sqrt{1-c^2-d^2}$ \\[0.5ex]
	6 & $\frac{c-d\sqrt{3}}{2}$ & $\frac{-d-c\sqrt{3}}{2}$ & $\sqrt{1-c^2-d^2}$ \\[0.5ex]
	7 & $\frac{\sqrt{3}}{2}$ & $-\frac{1}{2}$ & 0 \\[0.5ex]
	8 & 0 & 1 & 0 \\[0.5ex]
	9 & $-\frac{\sqrt{3}}{2}$ & $-\frac{1}{2}$ & 0 \\[0.5ex]
	10 & $c$ & $d$ & $-\sqrt{1-c^2-d^2}$ \\[0.5ex]
	11 & $\frac{-c+d\sqrt{3}}{2}$ & $\frac{-d-c\sqrt{3}}{2}$ & $-\sqrt{1-c^2-d^2}$ \\[0.5ex]
	12 & $\frac{-c-d\sqrt{3}}{2}$ & $\frac{-d+c\sqrt{3}}{2}$ & $-\sqrt{1-c^2-d^2}$ \\[0.5ex]
	13 & $-a$ & $b$ & $-\sqrt{1-a^2-b^2}$ \\[0.5ex]
	14 & $\frac{a+b\sqrt{3}}{2}$ & $\frac{-b+a\sqrt{3}}{2}$ & $-\sqrt{1-a^2-b^2}$ \\[0.5ex]
	15 & $\frac{a-b\sqrt{3}}{2}$ & $\frac{-b-a\sqrt{3}}{2}$ & $-\sqrt{1-a^2-b^2}$
\end{longtable}

The values of the 4 parameters optimized for the minimal solution are:
\begin{center}
	\begin{tabular}{c|c|c|c}
		Parameter & log & 1/r & $1/r^2$ \\
		\hline\Tstrut
		$a$ & 0.07491147758374605908 & 0.06812067962095187165 & 0.06471950175770332021 \\
		$b$ & 0.5594393279679571222 & 0.5577319898936203239 & 0.5555619674382903803 \\
		$c$ & 0.7403487282949123332 & 0.7425947069353261638 & 0.7425546518580456061 \\
		$d$ & 0.5376117024907701347 & 0.5320796700922574934 & 0.5300812467220562797 \\
		\hline\Tstrut
		energy & -32.14787628384166228 & 80.67024411429390062 & 67.48618473529605946
	\end{tabular}
\end{center}
The decimal values for the 4 parameters are known to 50,014 digits, as usual, and the energy values likewise. Unfortunately the algebraic degree for all is $> 360$ as searching has shown.

\noindent
\textit{\textbf{Symmetries -}}

\begin{center}
	\begin{tabular}{l|l}
		\multicolumn{2}{c}{Symmetries - 15 points} \\
		\hline\Tstrut
		planes & [[5, 1]] \\[0.2ex]
		\hline\Tstrut
		Gram groups & [[6, 7], [12, 14], [15, 1]] \\
		\hline\Tstrut
		Polygons & [[3, 5]]
	\end{tabular}
\end{center}

\subsection{16 points}
The optimal solution for 16 points is a configuration with 5 embedded parallel-plane triangles along an axis, with a single pole at one end. Each triangle is rotated $180^{\circ}$ with respect to the neighboring triangle(s) and the pole is aligned with a line through the triangle centroids. The arrangement is [1:3:3:3:3:3].

\begin{figure}[ht]
	\begin{center}
		\includegraphics[type=pdf,ext=pdf,read=pdf,height=1in,width=1in,angle=0]{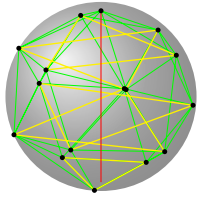}
		\caption{16 points.}
		\label{fig:16pts}
	\end{center}
\end{figure}

\noindent
\textit{\textbf{Parameterized Structure -}}
\begin{longtable}[c]{c|ccc}
	\caption{Spherical code for 16 points} \\
	pt & $x$ & $y$ & $z$ \\
	\hline\vspace*{-2.2ex}
	\endfirsthead
	\multicolumn{4}{c}%
	{\tablename\ \thetable\ -- Spherical code 16 points -- \textit{continued}} \\
	pt & $x$ & $y$ & $z$ \\
	\hline\vspace*{-2.2ex}
	\endhead
	1 & 0 & 0 & 1 \\[0.5ex]
	2 & $\sqrt{1-a^2-b^2}$ & $b$ & $a$ \\[0.7ex]
	3 & $\frac{b\sqrt{3}-\sqrt{1-a^2-b^2}}{2}$ & $\frac{\left(-\sqrt{3}\sqrt{1-a^2-b^2}\right)-b}{2}$ & $a$ \\[0.7ex]
	4 & $\frac{-b\sqrt{3}-\sqrt{1-a^2-b^2}}{2}$ & $\frac{\sqrt{3}\sqrt{1-a^2-b^2}-b}{2}$ & $a$ \\[0.7ex]
	5 & $\sqrt{1-c^2-d^2}$ & $d$ & $c$ \\[0.7ex]
	6 & $\frac{d\sqrt{3}-\sqrt{1-c^2-d^2}}{2}$ & $\frac{\left(-\sqrt{3}\sqrt{1-c^2-d^2}\right)-d}{2}$ & $c$ \\[0.7ex]
	7 & $\frac{-d\sqrt{3}-\sqrt{1-c^2-d^2}}{2}$ & $\frac{\sqrt{3}\sqrt{1-c^2-d^2}-d}{2}$ & $c$ \\[0.7ex]
	8 & $\sqrt{1-e^2-f^2}$ & $f$ & $e$ \\[0.7ex]
	9 & $\frac{f\sqrt{3}-\sqrt{1-e^2-f^2}}{2}$ & $\frac{\left(-\sqrt{3}\sqrt{1-e^2-f^2}\right)-f}{2}$ & $e$ \\[0.7ex]
	10 & $\frac{-f\sqrt{3}-\sqrt{1-e^2-f^2}}{2}$ & $\frac{\sqrt{3}\sqrt{1-e^2-f^2}-f}{2}$ & $e$ \\[0.7ex]
	11 & $\frac{\sqrt{2}}{3}$ & $-\sqrt{\frac{2}{3}}$ & $-\frac{1}{3}$ \\[0.7ex]
	12 & $-\frac{2\sqrt{2}}{3}$ & 0 & $-\frac{1}{3}$ \\[0.7ex]
	13 & $\frac{\sqrt{2}}{3}$ & $\sqrt{\frac{2}{3}}$ & $-\frac{1}{3}$ \\[0.7ex]
	14 & $\sqrt{1-g^2-h^2}$ & $h$ & $g$ \\[0.7ex]
	15 & $\frac{h\sqrt{3}-\sqrt{1-g^2-h^2}}{2}$ & $\frac{\left(-\sqrt{3}\sqrt{1-g^2-h^2}\right)-h}{2}$ & $g$ \\[0.7ex]
	16 & $\frac{-h\sqrt{3}-\sqrt{1-g^2-h^2}}{2}$ & $\frac{\sqrt{3}\sqrt{1-g^2-h^2}-h}{2}$ & $g$
\end{longtable}
\noindent

\textit{\textbf{logarithmic potential -}}
\begin{align*}
\text{let } A(x) & = 89915392x^{36} - 539492352x^{34} + 838139904x^{32} + 3993407488x^{30} - 21748887552x^{28} \\ & \quad + 24754297344x^{26} + 60398596528x^{24} - 215912904384x^{22} + 121681883520x^{20} \\ & \quad + 335366512064x^{18} - 140298350304x^{16} - 255855385152x^{14} + 57862428200x^{12} \\ & \quad + 79451789520x^{10} - 12304178880x^{8} - 8598909328x^{6} + 2428094280x^{4} \\ & \quad - 177973632x^{2} + 1022283 \\
\text{and} \\
\text{let } B(x) & = 1701x^{12} - 3402x^{10} - 1134x^9 - 8019x^8 + 1512x^7 - 10980x^6 - 5796x^5 + 26199x^4 \\ & \quad + 840x^3 - 11298x^2 + 1506x + 679
\end{align*}
be the two defining algebraic polynomials, then the solutions for the \textit{logarithmic} potential are:

\begin{tabular}{l}
$a = \;\ 0.5899523558141528103 = \text{a root of B(x)}$ \\
$b = \,\;\ 0.4001422226961329020 = \text{a root of A(x)}$ \\
$c = \,\;\ 0.4606669288440935079 = \text{a root of B(x)}$ \\
$d = -0.4793130545097588813 = \text{a root of A(x)}$ \\
$e = -0.1927625845883329931 = \text{a root of B(x)}$ \\
$f = \;\ 0.07917083181362597929 = \text{a root of A(x)}$ \\
$g = -0.8578567000699133251 = \text{a root of B(x)}$ \\
$h = -0.07917083181362597929 = \text{a root of A(x)}$ \\
\end{tabular}
\begin{align*}
energy & = -36.1061521619506068725\ldots = log(\text{a root of } 76957043352332967211482500195592995713/ \\ & \qquad 046365762627825523336510555167425334955489475418488779072100860950445293568x^3 \\ & \qquad - 16052612687571957819904801038866776948912323897668308741713914864782031146/ \\ & \qquad 426236031726716444475392x^2 + 173088986728844883554508658463496666984905048975/ \\ & \qquad 6484549909556494336x - 4259245936403381258086839702909)
\end{align*}

\noindent
\textit{\textbf{Coulomb $\mathbf{1/r}$ potential -}}

For the Coulomb potential, the algebraic degree for the polynomials $A(x)$ and $B(x)$ is $>480$ and the attempts to find them have failed.

The first 19 digits of the real spherical code and minimal energy is given below for the Coulomb $1/r$ potential:

\begin{tabular}{l}
	$a = \;\ 0.5784553451430919801\ldots$ \\
	$b = \,\;\ 0.4066356292761919381\ldots$ \\
	$c = \,\;\ 0.4725437786058905049\ldots$ \\
	$d = -0.4714929532449352116\ldots$ \\
	$e = -0.1914894233689195990\ldots$ \\
	$f = \,\;\ 0.06485732396874327353\ldots$ \\
	$g = -0.8595097003800628861\ldots$ \\
	$h = -0.06485732396874327353\ldots$ \\
	\\
	$energy = 92.9116553025449436907\ldots$
\end{tabular}

\noindent
\textit{\textbf{Inverse square law $\mathbf{1/r^2}$ potential -}}

\begin{longtable}[l]{l}
$\text{let } B(x) = 1990280943581607x^{60} - 19902809435816070x^{58} $ \\
\hphantom{000} $ - \; 27421648556013252x^{57} - 43941234574323423x^{56} $ \\
\hphantom{000} $ + \; 255935386522790352x^{55} + 1250459032205919348x^{54} $ \\
\hphantom{000} $ + \; 332161988700866328x^{53} - 4502047973778290349x^{52} $ \\
\hphantom{000} $ - \; 9254712809334718512x^{51} + 1848754561419189630x^{50} $ \\
\hphantom{000} $ + \; 40043464834629679956x^{49} + 20754711971072321661x^{48} $ \\
\hphantom{000} $ - \; 127224497718745869600x^{47} - 129032333365902031800x^{46} $ \\
\hphantom{000} $ + \; 366503839233648588144x^{45} + 655441487568573180795x^{44} $ \\
\hphantom{000} $ - \; 669111487290513977760x^{43} - 1899095746076333915958x^{42} $ \\
\hphantom{000} $ + \; 707866494468313444380x^{41} + 3503698955684621738685x^{40} $ \\
\hphantom{000} $ - \; 585877852174542474384x^{39} - 5139653744209018233780x^{38} $ \\
\hphantom{000} $ + \; 317649251715694661736x^{37} + 6437655525601892490471x^{36} $ \\
\hphantom{000} $ + \; 1419741771734647269360x^{35} - 3770384156108979284610x^{34} $ \\
\hphantom{000} $ - \; 721906282040742833580x^{33} + 2079936844030153363833x^{32} $ \\
\hphantom{000} $ + \; 4586289138687879233856x^{31} + 7330669816248030061296x^{30} $ \\
\hphantom{000} $ + \; 1632127979399261997600x^{29} - 1689853496295756171195x^{28} $ \\
\hphantom{000} $ + \; 13509139969273281053760x^{27} + 23634992650336254976158x^{26} $ \\
\hphantom{000} $ + \; 3980029437337343519700x^{25} - 5876462430313445547405x^{24} $ \\
\hphantom{000} $ + \; 13372329567591304865712x^{23} + 13363907663844619776972x^{22} $ \\
\hphantom{000} $ - \; 3008289018150579923352x^{21} - 2399390462232922407831x^{20} $ \\
\hphantom{000} $ + \; 3811959959294391158448x^{19} + 2201597932109314943466x^{18} $ \\
\hphantom{000} $ - \; 1611682502421498836580x^{17} - 794300844233841993369x^{16} $ \\
\hphantom{000} $ + \; 1267106177356657640928x^{15} + 285708020409215700552x^{14} $ \\
\hphantom{000} $ - \; 488638078721611154064x^{13} + 3723688174313436633x^{12} $ \\
\hphantom{000} $ + \; 64019467833519357792x^{11} - 83645669146145463522x^{10} $ \\
\hphantom{000} $ - \; 60003466723495376172x^{9} + 38993648108050287x^{8} $ \\
\hphantom{000} $ + \; 26543401733915042064x^{7} + 8319347141860947444x^{6} $ \\
\hphantom{000} $ - \; 4410051391500266216x^{5} - 1134983421447927267x^{4}$ \\
\hphantom{000} $ + \; 311613104829005712x^{3} + 9048953409216282x^{2}$ \\
\hphantom{000} $ + \; 1091221162080444x + 2242635315117795$
\end{longtable}
{\footnotesize
\begin{longtable}[l]{l}
$\text{let } A(x) = 2224225399066765295616x^{180} - 66726761972002958868480x^{178} $ \\
\hphantom{000} $ + \; 808032965090852688887808x^{176} - 4654845173544489647603712x^{174} $ \\
\hphantom{000} $ - \; 4419320395191815937982464x^{172} + 404990473123255093309734912x^{170} $ \\
\hphantom{000} $ - \; 4687291122972092315385987072x^{168} + 29615087152686382580005601280x^{166} $ \\
\hphantom{000} $ - \; 61370785356292588245765586944x^{164} - 927319976276608044803886678016x^{162} $ \\
\hphantom{000} $ + \; 13640375637643154991589388976128x^{160} - 112342344680441260173498058276864x^{158} $ \\
\hphantom{000} $ + \; 695261621972565687452794831765504x^{156} - 3452220355326481081987587859546112x^{154} $ \\
\hphantom{000} $ + \; 13846016116460829738689498772930560x^{152} - 43494742725799721407279302954713088x^{150} $ \\
\hphantom{000} $ + \; 96632258756361759744332622299398144x^{148} - 87715458048721096084831245365673984x^{146} $ \\
\hphantom{000} $ - \; 391887640892939402993602581853372416x^{144} + 2389708442838665095909906267525611520x^{142} $ \\
\hphantom{000} $ - \; 7277484680738157005869242025095725056x^{140} + 14032533223633918680693595327205474304x^{138} $ \\
\hphantom{000} $ - \; 12649369803012158857725355017421062144x^{136} - 22104955528279903875489542700320456704x^{134} $ \\
\hphantom{000} $ + \; 119948218840100264697623835629514726400x^{132} - 265259355575627302186858919256618797056x^{130} $ \\
\hphantom{000} $ + \; 320449646844816073631014310806800430080x^{128} - 16714832260731330260912421855176255168x^{126} $ \\
\hphantom{000} $ - \; 874130832127977138118422842465301982464x^{124} + 2327376117261221684578048963785044036320x^{122} $ \\
\hphantom{000} $ - \; 4433713725153147105072468551667909183495x^{120} + 8875277675960918149827347910957913730060x^{118} $ \\
\hphantom{000} $ - \; 20298796240860198470466563590205710380492x^{116} + 44000187813236574306008823203915156448320x^{114} $ \\
\hphantom{000} $ - \; 78322494671113161408191942221196431588752x^{112} + 110950401783144589788868131580230578861664x^{110} $ \\
\hphantom{000} $ - \; 139417656543358097210112892880939673992680x^{108} + 222320772412817409533305242285753009403008x^{106} $ \\
\hphantom{000} $ - \; 514986098995667877201161801070296382403768x^{104} + 1196056835811287060147635348292175592275600x^{102} $ \\
\hphantom{000} $ - \; 2245070332192896261272362433577301752784224x^{100} + 3217906404958508850886819594145597799495568x^{98} $ \\
\hphantom{000} $ - \; 3330083041064200873965959609788770455844784x^{96} + 2027008933097900091853206993852172116289856x^{94} $ \\
\hphantom{000} $ + \; 256065553555719775292531588418658385770880x^{92} - 1856086362727517132425898513800135128546480x^{90} $ \\
\hphantom{000} $ + \; 678095060308511112369756613047305458330320x^{88} + 4309721555117202397046698884867341068762496x^{86} $ \\
\hphantom{000} $ - \; 11993202813895890039643215430532572127173696x^{84} + 19427182193009887599857075901316636644199296x^{82} $ \\
\hphantom{000} $ - \; 23466816306690069757245762218306738739357824x^{80} + 22573622221660248732658163508617872867676416x^{78} $ \\
\hphantom{000} $ - \; 17556191070081111444182790750049657241355456x^{76} + 10863184497747124945087395366323868618766144x^{74} $ \\
\hphantom{000} $ - \; 5064240365235337480607452464015486595252096x^{72} + 1580691017159711034842565710290326979140096x^{70} $ \\
\hphantom{000} $ - \; 366390472221472176838182697574745923281280x^{68} + 469697799095354002027496601893676734093056x^{66} $ \\
\hphantom{000} $ - \; 876121332260276891126942598208770537248000x^{64} + 1037077879669904537587853257182182274686464x^{62} $ \\
\hphantom{000} $ - \; 884414647349094376394129678572833314886848x^{60} + 577590787311914206969257934265879502570368x^{58} $ \\
\hphantom{000} $ - \; 292420238420369318905025484827152244445952x^{56} + 117407554793610301839928183155521862047744x^{54} $ \\
\hphantom{000} $ - \; 41059736004841117709565818211344171608576x^{52} + 15380612534157064709782526267052470477824x^{50} $ \\
\hphantom{000} $ - \; 9076258057827374205521188380760928027136x^{48} + 7768428535099042491739439960861805510656x^{46} $ \\
\hphantom{000} $ - \; 5315010300779892593249393495752551416832x^{44} + 2397294446334587683469708207701702647808x^{42} $ \\
\hphantom{000} $ - \; 875932110351698986977659308562374891520x^{40} + 361671770185444073383475648938402447360x^{38} $ \\
\hphantom{000} $ - \; 113792111590768835700323376990156344320x^{36} + 9794291383785225940117182938899073024x^{34} $ \\
\hphantom{000} $ + \; 587642877485643607168837384742584320x^{32} + 1636844792111403858152998077200187392x^{30} $ \\
\hphantom{000} $ - \; 73048126031844571941679458503655424x^{28} - 174401845158783490349073877898051584x^{26} $ \\
\hphantom{000} $ - \; 20863370618668547224552034971025408x^{24} + 10732008152387774649408151150919680x^{22} $ \\
\hphantom{000} $ + \; 2129681450563201799814129672650752x^{20} - 441801998942565500176793936527360x^{18} $ \\
\hphantom{000} $ - \; 48713092709549829100631320625152x^{16} + 7098447961372115198142368972800x^{14} $ \\
\hphantom{000} $ + \; 254830098732817001004489768960x^{12} - 25700118754142025187453829120x^{10} $ \\
\hphantom{000} $ + \; 1117483756902132381044441088x^{8} - 82695495956318028130418688x^{6} $ \\
\hphantom{000} $ + \; 13377928794696579343712256x^{4} - 185520202801158295650304x^{2} $ \\
\hphantom{000} $ + \; 420353570822297747456$
\end{longtable}
}
Then the spherical code for the 16 point Inverse Square law $1/r^2$ potential is found:

\begin{tabular}{l}
	$a = \;\ 0.5689949134399233803 = \text{a root of B(x)}$ \\
	$b = \,\;\ 0.4116447726156913951 = \text{a root of A(x)}$ \\
	$c = \,\;\ 0.4820317283041705916 = \text{a root of B(x)}$ \\
	$d = -0.4648986301131366439 = \text{a root of A(x)}$ \\
	$e = -0.1901813704907959089 = \text{a root of B(x)}$ \\
	$f = \,\;\ 0.05325385749744524880 = \text{a root of A(x)}$ \\
	$g = -0.8608452712532980630 = \text{a root of B(x)}$ \\
	$h = -0.05325385749744524880 = \text{a root of A(x)}$ \\
\end{tabular}

\noindent
likewise the minimal energy is given:
\begin{align*}
energy & = 78.7726490309895138199\ldots = \text{ a root of }
1048576x^{15} - 334757888x^{14} + 47184216064x^{13} \\ & \qquad - 3921290297344x^{12} + 215897097752576x^{11} - 8369979026080768x^{10} \\ & \qquad + 236600677270422528x^9 - 4973176072973820672x^8 + 78418910577341147904x^7 \\ & \qquad - 927517684105981411776x^6 + 8154887946779466170424x^5 \\ & \qquad - 52268128066305944923824x^4 + 236014485410944499420934x^3 \\ & \qquad - 707824749993648711558966x^2 + 1260534253142866044717126x \\ & \qquad - 1007561675892684901708755
\end{align*}

\noindent
\textit{\textbf{Symmetries -}}

\begin{center}
	\begin{tabular}{l|l}
		\multicolumn{2}{c}{Symmetries - 16 points} \\
		\hline\Tstrut
		planes & [[5, 4]] \\[0.2ex]
		\hline\Tstrut
		Gram groups & [[12, 4], [16, 1], [24, 8]] \\
		\hline\Tstrut
		Polygons & [[3, 20]]
	\end{tabular}
\end{center}

\subsection{17 points}
\begin{figure}[ht]
	\centering
	\begin{minipage}{0.9\textwidth}
		\centering
		\includegraphics[type=pdf,ext=pdf,read=pdf,height=1in,width=1in,angle=0]{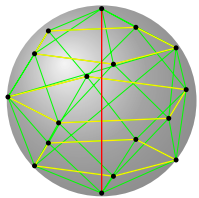}
		\caption{17 points.}
		\label{fig:17pts.log}
	\end{minipage}\hfill
\end{figure}

The configuration for 17 points has the points arrange themselves into a dipole, and 3 parallel pentagons, with the normal to the parallel-plane pentagons as the dipole axis. See figure \ref{fig:17pts.log}. This arrangement is [1:5:5:5:1], where 1 denotes a single point and 5 denotes a pentagon.

\noindent
\textit{\textbf{Parametric Solution --}}

\begin{center}
	\begin{longtable}{c|ccc}
		\caption{Spherical code for 17 points} \\
		pt & $x$ & $y$ & $z$ \\
		\hline\vspace*{-2.2ex}
		\endfirsthead
		\multicolumn{4}{c}%
		{\tablename\ \thetable\ -- 17 points spherical code -- \textit{cont\ldots}} \\
		pt & $x$ & $y$ & $z$ \\
		\hline\vspace*{-2.2ex}
		\endhead
		1 & 0 & 0 & 1 \\[0.5ex]
		2 & $-\sqrt{1-a^2}$ & 0 & $a$ \\[0.7ex]
		3 & $-\frac{\left(\sqrt{5}-1\right)\sqrt{1-a^2}}{4}$ & $\frac{\sqrt{2\sqrt{5}+10}\sqrt{1-a^2}}{4}$ & $a$ \\[0.7ex]
		4 & $\frac{\left(\sqrt{5}+1\right)\sqrt{1-a^2}}{4}$ & $\frac{\sqrt{10-2\sqrt{5}}\sqrt{1-a^2}}{4}$ & $a$ \\[0.7ex]
		5 & $\frac{\left(\sqrt{5}+1\right)\sqrt{1-a^2}}{4}$ & $-\frac{\sqrt{10-2\sqrt{5}}\sqrt{1-a^2}}{4}$ & $a$ \\[0.7ex]
		6 & $-\frac{\left(\sqrt{5}-1\right)\sqrt{1-a^2}}{4}$ & $-\frac{\sqrt{2\sqrt{5}+10}\sqrt{1-a^2}}{4}$ & $a$ \\[0.7ex]
		7 & 1 & 0 & 0 \\[0.7ex]
		8 & $\frac{\sqrt{5}-1}{4}$ & $\frac{\sqrt{2\sqrt{5}+10}}{4}$ & 0 \\[0.7ex]
		9 & $-\frac{\sqrt{5}+1}{4}$ & $\frac{\sqrt{10-2\sqrt{5}}}{4}$ & 0 \\[0.7ex]
		10 & $-\frac{\sqrt{5}+1}{4}$ & $-\frac{\sqrt{10-2\sqrt{5}}}{4}$ & 0 \\[0.7ex]
		11 & $\frac{\sqrt{5}-1}{4}$ & $-\frac{\sqrt{2\sqrt{5}+10}}{4}$ & 0 \\[0.7ex]
		12 & $-\sqrt{1-a^2}$ & 0 & $-a$ \\[0.7ex]
		13 & $-\frac{\left(\sqrt{5}-1\right)\sqrt{1-a^2}}{4}$ & $\frac{\sqrt{2\sqrt{5}+10}\sqrt{1-a^2}}{4}$ & $-a$ \\[0.7ex]
		14 & $\frac{\left(\sqrt{5}+1\right)\sqrt{1-a^2}}{4}$ & $\frac{\sqrt{10-2\sqrt{5}}\sqrt{1-a^2}}{4}$ & $-a$ \\[0.7ex]
		15 & $\frac{\left(\sqrt{5}+1\right)\sqrt{1-a^2}}{4}$ & $-\frac{\sqrt{10-2\sqrt{5}}\sqrt{1-a^2}}{4}$ & $-a$ \\[0.7ex]
		16 & $-\frac{\left(\sqrt{5}-1\right)\sqrt{1-a^2}}{4}$ & $-\frac{\sqrt{2\sqrt{5}+10}\sqrt{1-a^2}}{4}$ & $-a$ \\[0.7ex]
		17 & 0 & 0 & -1
	\end{longtable}
\end{center}
\noindent
Please note that the parameterization uses only 1 parameter, $a$.

\noindent
\textit{\textbf{logarithmic potential -}}

The global minimum for the \textit{logarithmic} potential has been determined.

Let $a = 0.6076810889242587549\ldots = \text{ a root of } 256x^{12} + 3545x^{10} + 11335xx^8 - 4470x^6 - 2050x^4 + 525x^2 + 75$.

The spherical code for the $energy$ is also known:
\begin{longtable}{l}
$energy = -40.2730669611811632688\ldots log(\text{ a root of }$ 16012754614967250797153798078/ \\
09264529669290183544775350331701073809543109113900584706040601860365657301144/ \\
95697054184260158283544631485901043893918796476892459555975251860397378921329/ \\
71004650239322719387721326595979010744368498782391378404208279897929676570447/ \\
41817257803953686812869344133570887680000000000000000000000000000000000000000/ \\
00000000000000000000000000000000000000000000000000000000000000000000000000000/ \\
0000000000000000000000000000000000000000000000000000000000000000000000000$x^{12}$ \\
$ + \; $377334792306556323678291825408953346835679595542161491732944851627138337355/ \\
36920118345413442481166282186917443618678834130306531841601354913282415746027/ \\
42984817259503666508681133162384521236708766464018845099896556636384809713286/ \\
39317785421109660416808307250098202311857787110626173283446284712039685082316/ \\
80000000000000000000000000000000000000000000000000000000000000000000000000000/ \\
00000000000000000000000000000000000000000000000000000000000000000000000000$x^{10}$ \\
$ - \; $144692374548576282746422991383115216477646062545643846512248403699542458201/ \\
58483392912221917153018803369845861061464070702822650153011276189565292852314/ \\
18267334657295033112788384238822540650787962034790328112735488875338076436951/ \\
79714406615429448430831521046528276370129510422219042377435822581383686109931/ \\
13661440000000000000000000000000000000000000000000000000000000000000000000000/ \\
0000000000000000000000000000000000000000$x^8 \;+\; $96468947935755003581836553712501/ \\
23077866708662553804882588896702637156405899244361486524549751584650689817437/ \\
35104081177131698109367301315515499165946856009873141343944455609919091792415/ \\
04776756226136204905288352875223687307472233624954900626819664727856576612634/ \\
36558940835753634660720283284763102630559353340847718400000000000000000000000/ \\
00000000000000000000000000000000000000000000000$x^6 \;+\; $1416545651654981248994952/ \\
53351663652489405460495773587222627003959706683912402319587637007842352012499/ \\
77153623020433987638394461778517400543413561974931363911869296115633277552873/ \\
37689273185817730372301229959921086554583166296073883742434306077385855257934/ \\
71331051347323439106054315785579825241588446003200000000000000000000000000000/ \\
0000000000000000000$x^4 \;+\; $41041765620529242133611691056883041736616521262187380/ \\
01693921773618916140585455112841800710702201488638383670501319602887938491450/ \\
78877437731470009943587939240379895685794943269641182065282516728632198605949/ \\
4286877893220716049098699633151314459420780627320110797882080078125$x^2 \;+\; $10130/ \\
65324433836171511818326096474890383898005918563696288002277756507034036354527/ \\
92961597874685151227739206216096210673398319118052045295602706905129735441578/ \\
6421338721071661056)
\end{longtable}

\noindent
\textit{\textbf{Coulomb $1/r$ potential -}}

The global minimum energy and coordinates for the \textit{Coulomb $\frac{1}{r}$} potential has been determined to 50,014 digits, however the degree of the algebraic polynomial is $>600$ for parameter $a$.

\begin{longtable}[c]{cr|l}
	\caption{17 points - Coulomb $1/r$ potential} \\
	\toprule
	Parameter & $a$ & $0.6095575990554807772\ldots$ \\
	\midrule
	Minimal & $energy$ & $106.0504048286187048\ldots$
\end{longtable}

\noindent
\textit{\textbf{Inverse Square $1/r^2$ potential -}}

The algebraic code has been obtained for both the coordinates and the minimal energy.

\noindent
\begin{align*}
a \; & = 0.6117975792003008025\ldots = \text{ a root of } 100x^{26} + 7981x^{24} + 23840x^{22} - 143330x^{20} - 1346500x^{18} \\ & \quad + 2043775x^{16} + 6878400x^{14} + 15610900x^{12} - 1916100x^{10} - 4286725x^8 - 580000x^6 \\ & \quad + 456750x^4 + 22500x^2 + 5625
\end{align*}
\begin{align*}
energy \; & = 91.04731931420740319\ldots = \text{ a root of } 1024x^{13} - 612608x^{12} + 162461952x^{11} \\ & \quad - 25242603136x^{10} + 2563215005760x^9 - 179852969538112x^8 + 8988670083974544x^7 \\ & \quad - 324447567015364856x^6 + 8463548460720194888x^5 - 157615754889421646115x^4 \\ & \quad + 2037337935965765371114x^3 - 17321248088620678256598x^2 \\ & \quad + 86909807362528534375692x - 194568980729327870281101
\end{align*}

\noindent
\textit{\textbf{Symmetries -}}

17 points has the same symmetry groups for all 3 potentials.

\begin{center}
	\begin{tabular}{l|l}
		\multicolumn{2}{c}{Symmetries - 17 points} \\
		\hline\Tstrut
		planes & [[4, 50], [8, 5], [10, 5], [30, 1]] \\[0.2ex]
		\hline\Tstrut
		Gram groups & [[2, 1], [10, 3], [17, 1], [20, 8], [40, 2]] \\
		\hline\Tstrut
		Polygons & [[4, 60], [5, 8]] \\[0.2ex]
	\end{tabular}
\end{center}

\subsection{18 points}

The optimal configuration for 18 points, in all 3 potentials, consists of 4 parallel-plane squares, each rotated $90^{\circ}$ with respect to its neighbor, and 2 poles, located along the line passing through all centroids of the squares. The arrangement is [1:4:4:4:4:1].

\begin{figure}[H]
	\begin{center}
		\includegraphics[type=pdf,ext=pdf,read=pdf,height=1in,width=1in,angle=0]{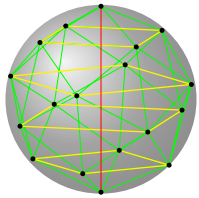}
		\caption{18 points.}
		\label{fig:18pts}
	\end{center}
\end{figure}

\noindent
\textit{\textbf{Parameterized Structure -}}
\begin{center}
	\begin{longtable}[c]{l|c c c}
		\caption{Parameterized structure for 18 points} \\
		pt & $x$ & $y$ & $z$ \\
		\hline\vspace*{-2.2ex}
		\endfirsthead
		\multicolumn{4}{c}%
		{\tablename\ \thetable\ -- 18 points parameters -- \textit{cont.}} \\
		pt & $x$ & $y$ & $z$ \\
		\hline\vspace*{-2.2ex}
		\endhead
		1 & 0 & 0 & 1 \\[0.5ex]
		2 & $\sqrt{1-a^2}$ & 0 & $a$ \\[0.5ex]
		3 & $-\sqrt{1-a^2}$ & 0 & $a$ \\[0.5ex]
		4 & 0 & $\sqrt{1-a^2}$ & $a$ \\[0.5ex]
		5 & 0 & $-\sqrt{1-a^2}$ & $a$ \\[0.5ex]
		6 & $\sqrt{\frac{1-b^2}{2}}$ & $\sqrt{\frac{1-b^2}{2}}$ & $b$ \\[0.5ex]
		7 & $\sqrt{\frac{1-b^2}{2}}$ & $-\sqrt{\frac{1-b^2}{2}}$ & $b$ \\[0.5ex]
		8 & $-\sqrt{\frac{1-b^2}{2}}$ & $\sqrt{\frac{1-b^2}{2}}$ & $b$ \\[0.5ex]
		9 & $-\sqrt{\frac{1-b^2}{2}}$ & $-\sqrt{\frac{1-b^2}{2}}$ & $b$ \\[0.5ex]
		10 & 0 & $\sqrt{1-b^2}$ & $-b$ \\[0.5ex]
		11 & 0 & $-\sqrt{1-b^2}$ & $-b$ \\[0.5ex]
		12 & $\sqrt{1-b^2}$ & 0 & $-b$ \\[0.5ex]
		13 & $-\sqrt{1-b^2}$ & 0 & $-b$ \\[0.5ex]
		14 & $\sqrt{\frac{1-a^2}{2}}$ & $\sqrt{\frac{1-a^2}{2}}$ & $-a$ \\[0.5ex]
		15 & $\sqrt{\frac{1-a^2}{2}}$ & $-\sqrt{\frac{1-a^2}{2}}$ & $-a$ \\[0.5ex]
		16 & $-\sqrt{\frac{1-a^2}{2}}$ & $\sqrt{\frac{1-a^2}{2}}$ & $-a$ \\[0.5ex]
		17 & $-\sqrt{\frac{1-a^2}{2}}$ & $-\sqrt{\frac{1-a^2}{2}}$ & $-a$ \\[0.5ex]
		18 & 0 & 0 & -1
	\end{longtable}
\end{center}
\noindent
\textit{\textbf{logarithmic potential -}}
\begin{align*}
\text{let } f(x) & = 6975757441x^{40} + 208693421574x^{38} + 2727780575657x^{36} + 20307644474000x^{34} \\ & \qquad + 93075033736436x^{32} + 258946513145640x^{30} + 373983661039812x^{28} \\ & \qquad + 55293016177264x^{26} - 546051138638354x^{24} - 383999611370844x^{22} \\ & \qquad + 342483271515086x^{20} + 240930872783728x^{18} - 107690500640380x^{16} \\ & \qquad - 82973538391256x^{14} + 39303141877044x^{12} + 479477056528x^{10} \\ & \qquad - 1194833123559x^{8} - 72796608506x^{6} + 980311185x^{4} + 115941504x^{2} + 2560000
\end{align*}
For the spherical code, we let $a$ and $b$ be roots of $f(x)$:
\begin{align*}
a & = 0.6754406562091057220\ldots = \text{ a root of f(x)} \hspace{2.5in} \\
b & = 0.2063761761970050338\ldots = \text{ a root of f(x)}
\end{align*}
The minimal energy polynomial has been found, but it is excessively long. A precision of 12,002 digits had to be employed to successfully recover this polynomial.
\begin{longtable}[l]{l}
	$energy = -44.6502872592307253555\ldots = log(\text{a root of } 349517029005437161283359578148251/ $ \\
	$ \; 3476901109690634419058795694716815043148971274171932580984121459006491231497041542/ $ \\
	$ \quad 313380560502302063652326536276730384167219482425213883791784359267965624864476143/ $ \\
	$ \quad 48188935756915368751266109911460576609233308993820766582549613502551845177680401/ $ \\
	$ \quad 33399072906144382981170642512621029317210689563707184529409262882208418202138290/ $ \\
	$ \quad 23049915771457356905713901068702774827736025047937311524392733763487449667843374/ $ \\
	$ \quad 49787403959956951462529143257666741808373719075235663641578953833996132107193050/ $ \\
	$ \quad 34630039617468429506356895753741337878060982409701341313152243864260254012971853/ $ \\
	$ \quad 318587192815332236473292711526400000000000000000000000000000000000000000000000000x^{10}$ \\
	$ \quad - 52551132042685849026065501179765397186197506119258920304392355633440530153395782/ $ \\
	$ \quad 45140911657432298197888180253272066534620410066019948852205915732529159595545247/ $ \\
	$ \quad 04699047577407761852301343921408888491106707676760325977305060254767296767854402/ $ \\
	$ \quad 75383578067239402597389082037257859934443457625169087128394601538167237636391020/ $ \\
	$ \quad 08101786020267973189620110005270889427465404641676106116626767243128514272579306/ $ \\
	$ \quad 07505915768419748368630097794376062624152380571701567886704377943031006348655214/ $ \\
	$ \quad 22632164354284461946438633090646494600464973960701546937636677147555683451503358/ $ \\
	$ \quad 10074284283554173656383383569695739333839812760043058626560000000000000000000000/ $ \\
	$ \quad 000000000000000000x^9 + 2675742080694921408517522202613315311040307885149900415453/ $ \\
	$ \quad 21817506117008012082287024592094096611691607907148747251032006934537698323232417/ $ \\
	$ \quad 38733884975945952338022881759151505137670921401986626454912191682569285160484726/ $ \\
	$ \quad 81082525982250387006920502897065765026048802135085926066553419152138557952445705/ $ \\
	$ \quad 20718751584517547430212281096443677928960837103713753961156691874278793660327915/ $ \\
	$ \quad 52067435484828853969835923238817677520628054474352025209917606118352844697537140/ $ \\
	$ \quad 18130059483429324075964121688947795490754799879091608413776504194751712012349838/ $ \\
	$ \quad 55254067462222253979982650319953606760380898972546904409253517834846208000000000/ $ \\
	$ \quad 000000000000000000000x^8 - 4540411473798394290854644920459525725968730563069224991/ $ \\
	$ \quad 03778881133741940276589586468710810431043970427392621472014097378401081211901124/ $ \\
	$ \quad 23645455497507751153305393939038718549486513545437429161627322773663604670486543/ $ \\
	$ \quad 63161724836957956339496720219517054422361007349260790720840096878604448790558647/ $ \\
	$ \quad 16939448346345061360813847396887545349807385153026090174005871302534541453136521/ $ \\
	$ \quad 47088226454121935775626044968993856838134978435545374328202917926740267648595722/ $ \\
	$ \quad 64045283507844243215197185578700478268008733032060741227301096124419889176325309/ $ \\
	$ \quad 04972871277103049704924709557793428281714001496853898988533841920000000000000000/ $ \\
	$ \quad 0000x^7 + 378135800462764072947621624377493985460375577930933181986970433306104710/ $ \\
	$ \quad 51618022963132766158259343101725899131179432057605421875235374798021734215332834/ $ \\
	$ \quad 68318448873655653921763788472572444972375865885611381008974159810830322835853112/ $ \\
	$ \quad 13947915543531808833272173225599614126193911741779595787972295924708388262134306/ $ \\
	$ \quad 37761445054972394126499593489783681941816407430851001444167529832296251077892261/ $ \\
	$ \quad 20914939019900668090010831008952196831386340055980768854374652827665944932564469/ $ \\
	$ \quad 08194299838206920299246692881193074299771049275078970400174330833019796468278756/ $ \\
	$ \quad 45429109168530631287786962944000000000000x^6 - 21446725848125520216440937147047308/ $ \\
	$ \quad 52489294718815277291554977218650229512350097217110212549930470771761034527078706/ $ \\
	$ \quad 34647766578318797662966981317392349046018408664543166551236932060045813395057249/ $ \\
	$ \quad 38846892131212958612857329209642735245671997829403284716719500204555380999903572/ $ \\
	$ \quad 03225892831588086277932323235088587802980700801336702306052130859981571733234114/ $ \\
	$ \quad 84415821489296158694001616729309584643135134731440110813403562712821080912160814/ $ \\
	$ \quad 00950589425557300072604312447709314263971893487864227851561401878516511159324653/ $ \\
	$ \quad 0603751973642572539433429420615902390494531944448x^5 + 9158078299321659636551986/ $ \\
	$ \quad 43241555688072970094824596789791938994909150684094442702621562512207278923758581/ $ \\
	$ \quad 80019668821741597238414923796303117699341275392045457010341797810384896019036773/ $ \\
	$ \quad 56975317243743308771955968774528858805265500572299168446940912926055007752708653/ $ \\
	$ \quad 45807372701538867413535202471998142108960253070085481141401978781319625337601660/ $ \\
	$ \quad 76320742643883727922519520631607910859742664336416529066267333956352285102764895/ $ \\
	$ \quad 85107854964628351179792967789288063028162503802038574566810696418649714404595051/ $ \\
	$ \quad 5646563397740265472x^4 - 204101629117025960998019700094879039496278036551923062343/ $ \\
	$ \quad 02317808077907179387023756903547804544133578811190216716313963136974236221106971/ $ \\
	$ \quad 26607087400960409315441055410826639134527298332566960099136729559592055434748313/ $ \\
	$ \quad 49610537053069778179142769789682414166133351411671035565661904087124150324254856/ $ \\
	$ \quad 58053655910580747595247748086483105118269464552843912747133799766536863255658053/ $ \\
	$ \quad 83338368843434364325322382985695681503233778775503697738224034598056426035746224/ $ \\
	$ \quad 58704914858822935760601088x^3 + 14296451571354724028108133137020864406647059956678/ $ \\
	$ \quad 66495163170979231468505015383935788852819383935836729414000709482076539890944591/ $ \\
	$ \quad 70701775730143912941769291551426622078320550841966485715380534201998437393354419/ $ \\
	$ \quad 93063081322652570355877017431053908805882925195011379884751380968594548835651612/ $ \\
	$ \quad 20194073591401483664221376889165338926108644596635581304989531003297062197999608/ $ \\
	$ \quad 09543351053911716063975926443992660186862404274531941179427623163396096x^2 - 20103/ $ \\
	$ \quad 15677609239296966224595263682760068958175286548363376432846552165814201211023396/ $ \\
	$ \quad 90514905917185821646465510220351536083883466404520623626726869265388688248774859/ $ \\
	$ \quad 29393638173799198411263760686319013949082028813006925483717571236490229507563565/ $ \\
	$ \quad 24566321719914951115285032222223348361395117542822642964058138081839739632527380/ $ \\
	$ \quad 144735873405737955147423956235277717509212828702688985756144238592x + 480957404123/ $ \\
	$ \quad 03883405904058368793114819338547249811295617352588070406967510296802602942856634/ $ \\
	$ \quad 74952287554532805259845369359443286525460197855205982877052379844503781444429263/ $ \\
	$ \quad 63638928549735332604850561599810492196356796818556639656990812164802450218345851/ $ \\
	$ \quad 79435651099671173276655522316999170350787534643330580951332145555259884723651938561) $
\end{longtable}

\noindent\textbf{NOTE:
The excessively long polynomial for this minimal energy for 18 points provides a strong hint that obtaining the algebraic code for minimal energy for points $>18$ is going to quickly become intractable.} GP-Pari found the 8th root of the polroots($poly(x_{10})$) command to be the correct root of this polynomial with huge coefficients only after extending precision to over 12,000 digits.

\noindent
\textit{\textbf{Coulomb 1/r potential -}}

For the spherical code, we let $a$ and $b$ be roots of $A(x)$ where $A(x)$ is an unknown algebraic polynomial of degree $>480$:
\begin{align*}
a & = 0.6751471684502996248\ldots = \text{ a root of A(x)} \hspace{2.5in} \\
b & = 0.2034104243431649960\ldots = \text{ a root of A(x)}
\end{align*}
As stated before, more capable mathematical tools are needed to discover the algebraic number, as the decimal values of $a$, and $b$ are known to 50,014 digits, similarly for the minimum energy value which is

$\quad\ energy = 120.08446744749231864308798638022609878\ldots$.

\noindent
\textit{\textbf{Inverse square law $\mathbf{1/r^2}$ potential -}}
\begin{longtable}[l]{l}
	$\text{let} f(x) = 31403151600910336x^{196} + 8407642195544244224x^{194}$ \\
	$ + \; 1055673537097791176704x^{192} + 83007599346596586455040x^{190}$ \\
	$ + \; 4575809607849871764275456x^{188} + 184773687073221458433242788x^{186}$ \\
	$ + \; 5453323019930964616155245808x^{184} + 112282175519762189381115322289x^{182}$ \\
	$ + \; 1404317605011834971424585433582x^{180} + 4612270393131975161590664549920x^{178}$ \\
	$ - \; 152041554742835474938694521144415x^{176} - 1705724149568660789414932318156291x^{174}$ \\
	$ + \; 31411130116374507601322712986646320x^{172} + 1151783771973171002091177068360467930x^{170}$ \\
	$ + \; 16901988691407695224508671314405219101x^{168}$ \\
	$ + \; 146927867806382634761500824029519682618x^{166}$ \\
	$ + \; 723564620168144139994720282028588896164x^{164}$ \\
	$ + \; 494621875394290750415059292082509853736x^{162}$ \\
	$ - \; 21087835370504638426036070550811722340902x^{160}$ \\
	$ - \; 146018239458231958209935508555225147164742x^{158}$ \\
	$ - \; 129435685596309111204026664034450592844392x^{156}$ \\
	$ + \; 4394476714345846528918794750472415532345520x^{154}$ \\
	$ + \; 34295600723992661893110749908639250152030314x^{152}$ \\
	$ + \; 115115757789617056925271733324745722112198777x^{150}$ \\
	$ - \; 20826140611135746163117950963461859219441962x^{148}$ \\
	$ - \; 2402088105828107290716948901325146443501248312x^{146}$ \\
	$ - \; 14174310758244536458122934035158295318193888343x^{144}$ \\
	$ - \; 45243003599993284917166048595126386870639078755x^{142}$ \\
	$ - \; 60330137534019447428474166229686096615583939960x^{140}$ \\
	$ + \; 284772359001968581704583748303277453591459731130x^{138}$ \\
	$ + \; 2809762647471768742513911640722869150699519918109x^{136}$ \\
	$ + \; 12846576865618020863689870206973461712862894111256x^{134}$ \\
	$ + \; 28314837260485826501380584778302040636544739184304x^{132}$ \\
	$ - \; 31531768959658682365071948955819686531823359642848x^{130}$ \\
	$ - \; 455044641189813586346582978670085628086406852865640x^{128}$ \\
	$ - \; 1517700025477634021901415121486639734870962919414440x^{126}$ \\
	$ - \; 1945299294414137544252113427684591406293338320213024x^{124}$ \\
	$ + \; 3429821206324237430347472850309835354605848483538904x^{122}$ \\
	$ + \; 22889520301360565571093859974234546326725790843126840x^{120}$ \\
	$ + \; 63182729955768082432027148831287197786852838412457186x^{118}$ \\
	$ + \; 100939707327275642192630316582423425244822607078586556x^{116}$ \\
	$ - \; 7623448590860374037781849843254904698515131096048416x^{114}$ \\
	$ - \; 436934444797563522383376424786920673244094811746969662x^{112}$ \\
	$ - \; 1058609212978076401907057246886409838431531348244057638x^{110}$ \\
	$ - \; 1846153652534479444518815546985819741433540477330630464x^{108}$ \\
	$ - \; 2397518204892649068499671503031942058153853899701327820x^{106}$ \\
	$ - \; 748330347705709468962379792834194696923745922122927974x^{104}$ \\
	$ + \; 1805146686805877620811771097938480395244808962619975036x^{102}$ \\
	$ + \; 4148988744949951304737953263702433126468947530543173464x^{100}$ \\
	$ + \; 8801007357603004190621940240212634462885107958626683440x^{98}$ \\
	$ + \; 8813312244443671436483147008717720122337576967897027388x^{96}$ \\
	$ + \; 8531777971693161870526811464047429153943597326587730876x^{94}$ \\
	$ + \; 7613866437438510237986395927291631935573790541646633552x^{92}$ \\
	$ + \; 5572112553043805388789724561386116235421178061979180080x^{90}$ \\
	$ + \; 2877928053107318689162832268311997493479846394015755164x^{88}$ \\
	$ + \; 417905659286082904047192874560774465363794879808450874x^{86}$ \\
	$ - \; 2324083806652453756619288828311046346372552858438129732x^{84}$ \\
	$ - \; 3945192449222137281348580368397621991797257266844244496x^{82}$ \\
	$ - \; 3490642653752956730177550370448442743648003418580109286x^{80}$ \\
	$ - \; 926966753064093933784947931443682359060427678014379902x^{78}$ \\
	$ + \; 1403503058343378576147577172904375968089026441565845552x^{76}$ \\
	$ + \; 1579037784594663168943242060732394993777179843216931780x^{74}$ \\
	$ + \; 432031423945254437660964083489074303813047867966030338x^{72}$ \\
	$ - \; 354583804260482701683816886466720396046889414113862440x^{70}$ \\
	$ - \; 321604193875919768208827743173613506911764258629202256x^{68}$ \\
	$ - \; 20017456253538234599991402289606018044615677558023264x^{66}$ \\
	$ + \; 67758385409108994167983652366590627431842121481613400x^{64}$ \\
	$ + \; 19170367561415241352926833792268084778708067087646360x^{62}$ \\
	$ - \; 8028178840295601196205816315781710793186228516584864x^{60}$ \\
	$ - \; 4707753025103569948375293068939389281465561974745052x^{58}$ \\
	$ + \; 1578753999220962191096691180804903472253493676527304x^{56}$ \\
	$ + \; 762903773327851618260770699245222884519252399665805x^{54}$ \\
	$ - \; 537337071602136181543273846496183641425230738326698x^{52}$ \\
	$ + \; 70591370889160371438198220242108166216521526260160x^{50}$ \\
	$ + \; 93442967125071686844813881033770961089346191381437x^{48}$ \\
	$ - \; 56680913017080902891230978126279827147692836511991x^{46}$ \\
	$ + \; 1940245742734836507760230175801457053437394966864x^{44}$ \\
	$ + \; 8348295262152541572595060509667583635011275513266x^{42}$ \\
	$ - \; 3040670849353779058972452270592708183915913993175x^{40}$ \\
	$ + \; 75605731159025143084715987638607460194812826970x^{38}$ \\
	$ + \; 253347375477301158064652433029885532035583385508x^{36}$ \\
	$ - \; 90445524107163445513346531337783920980910132504x^{34}$ \\
	$ + \; 12422725628364640073285569731186459157000582906x^{32}$ \\
	$ + \; 1551436239041858782152684282890857894995966618x^{30}$ \\
	$ - \; 1185450973792042073653842883393420461369593640x^{28}$ \\
	$ + \; 285363828544049124572623679533967615248908416x^{26}$ \\
	$ - \; 35689087546960003826485053301819428528986614x^{24}$ \\
	$ + \; 1691707659682436108262791425958734494534957x^{22}$ \\
	$ + \; 83122040668556750206757648171547561201518x^{20}$ \\
	$ - \; 2534610747340695703497535590133538220920x^{18}$ \\
	$ - \; 1274019922501404641912766197935315809507x^{16}$ \\
	$ - \; 6348609513647865081388460945493362687x^{14} + 794317057079693927169728572567890952x^{12}$ \\
	$ + \; 33213294747380755318486312961516994x^{10} + 247210670109825215212920653560833x^{8}$ \\
	$ + \; 22789416663355912783375401624576x^{6} + 128531996519460884825736871936x^{4}$ \\
	$ + \; 2725494326569170291690307584x^{2} + 37706458790564181808513024$
\end{longtable}
For the spherical code, we let $a$ and $b$ be roots of $f(x)$:
\begin{align*}
a & = 0.6743335122024262360\ldots = \text{ a root of f(x)} \hspace{2.5in} \\
b & = 0.2007314823505518450\ldots = \text{ a root of f(x)}
\end{align*}
The minimal energy is only known as a real decimal value, the degree of its algebraic polynomial $>196$.

$energy = 104.31468528482482860938250508403941892\ldots$.

\noindent
\textit{\textbf{Symmetries -}}

\begin{center}
	\begin{tabular}{l|l}
		\multicolumn{2}{c}{Symmetries - 18 points} \\
		\hline\Tstrut
		planes & [[4, 48], [16, 1], [20, 4]]] \\[0.2ex]
		\hline\Tstrut
		Gram groups & [[2, 1], [8, 2], [16, 12], [18, 1], [32, 3]] \\
		\hline\Tstrut
		Polygons & [[4, 52], [6, 4]]
	\end{tabular}
\end{center}

\subsection{19 points}
The configuration for 19 points is significant, as a change occurs for the \textit{Inverse square law $1/r^2$} configuration, when all 19 points are moved in such a way, that no planer figure can be found among them, except trivial triangles for any 3 points.

For the \textit{logarithmic} and \textit{Coulomb $1/r$} configurations, which are essentially the same structure, embedded polygons are found, enabling the polyhedra to be parameterized.

\begin{figure}[ht]
	\centering
	\begin{minipage}{0.45\textwidth}
		\centering
		\includegraphics[type=pdf,ext=pdf,read=pdf,height=1in,width=1in,angle=0]{r-1.19pts.}
		\caption{19 points - Log or Coulomb potential.}
		\label{fig:19pts}
	\end{minipage}\hfill
	\begin{minipage}{0.45\textwidth}
		\centering
		\includegraphics[type=pdf,ext=pdf,read=pdf,height=1in,width=1in,angle=0]{r-2.19pts.}
		\caption{19 points - Inverse square law potential.}
		\label{fig:19ptsr2}
	\end{minipage}
\end{figure}

For those two potentials, we discover that the odd number of points produces one lone pole, the other 18 points are balanced with respect to the equatorial plane, either as part of a dipole pair, or one of the 4 points of a rectangle.

Specifically, proceeding from the one lone pole, we find a rectangle, then a dipole, then another rectangle, the largest, then two dipoles, oriented $90^{\circ}$ with respect to each other, the nearest dipole to the rectangle also oriented $90^{\circ}$ to the first-most dipole which is next to the pole. After the two dipoles, we come to another rectangle finishing the configuration. The 3 centroids of the rectangles are aligned exactly with the axis of the lone pole. The arrangement is [1:4:2:4:2:2:4].

19 point configurations were used as a check upon the \textit{percolating anneal} algorithm, as 5 different random groups of 19 points were used as a starting configuration, and in all 5 cases, the point set and energy converged to 57 decimal points accuracy. A check of the Gram matrix also showed the quintuples were identical after convergence.

It was a surprise to discover, however, that the \textit{Inverse square law $1/r^2$} potential resulted in no embedded polygons in the 19 points and cannot be parameterized.

\noindent
\textit{\textbf{Parametric structure -}}

\begin{longtable}[c]{l|ccc}
	\caption{Parameterized structure for 19 points} \\
	pt & $x$ & $y$ & $z$ \\
	\hline\vspace*{-2.2ex}
	\endfirsthead
	\multicolumn{4}{c}%
	{\tablename\ \thetable\ -- 19 points parameters -- \textit{cont.}} \\
	pt & $x$ & $y$ & $z$ \\
	\hline\vspace*{-2.2ex}
	\endhead
	1 & $0$ & $0$ & $1$ \\[0.5ex]
	2 & $b$ & $\sqrt{1-a^2-b^2}$ & $a$ \\[0.5ex]
	3 & $-b$ & $-\sqrt{1-a^2-b^2}$ & $a$ \\[0.5ex]
	4 & $b$ & $-\sqrt{1-a^2-b^2}$ & $a$ \\[0.5ex]
	5 & $-b$ & $\sqrt{1-a^2-b^2}$ & $a$ \\[0.5ex]
	6 & $\sqrt{1-c^2}$ & $0$ & $c$ \\[0.5ex]
	7 & $-\sqrt{1-c^2}$ & $0$ & $c$ \\[0.5ex]
	8 & $e$ & $\sqrt{1-d^2-e^2}$ & $d$ \\[0.5ex]
	9 & $-e$ & $-\sqrt{1-d^2-e^2}$ & $d$ \\[0.5ex]
	10 & $e$ & $-\sqrt{1-d^2-e^2}$ & $d$ \\[0.5ex]
	11 & $-e$ & $\sqrt{1-d^2-e^2}$ & $d$ \\[0.5ex]
	12 & $0$ & $-\sqrt{1-f^2}$ & $f$ \\[0.5ex]
	13 & $0$ & $\sqrt{1-f^2}$ & $f$ \\[0.5ex]
	14 & $\sqrt{1-g^2}$ & $0$ & $g$ \\[0.5ex]
	15 & $-\sqrt{1-g^2}$ & $0$ & $g$ \\[0.5ex]
	16 & $i$ & $\sqrt{1-h^2-i^2}$ & $h$ \\[0.5ex]
	17 & $-i$ & $-\sqrt{1-h^2-i^2}$ & $h$ \\[0.5ex]
	18 & $i$ & $-\sqrt{1-h^2-i^2}$ & $h$ \\[0.5ex]
	19 & $-i$ & $\sqrt{1-h^2-i^2}$ & $h$
\end{longtable}

The parametric structure uses 9 parameters, $a$, $b$, $c$, $d$, $e$, $f$, $g$, $h$, and $i$.

\noindent
\textit{\textbf{logarithmic potential -}}

The values for the parameters have been found to very high precision (50,014 digits), however the algebraic degree $>360$ for this potential.

\begin{longtable}[l]{l}
$ \qquad a =  0.6381528951337401855\ldots$ \\
$ \qquad b =  0.3831842970606347904\ldots$ \\
$ \qquad c =  0.5467365671168718879\ldots$ \\
$ \qquad d = -0.07231369021952538767\ldots$ \\
$ \qquad e = -0.7150404304460191670\ldots$ \\
$ \qquad f = -0.2005761421600021806\ldots$ \\
$ \qquad g = -0.3263507617273051298\ldots$ \\
$ \qquad h = -0.8257440365289970866\ldots$ \\
$ \qquad i =  0.3790800442392165742\ldots$ \\
\\
$\text{likewise for the energy:}$ \\
$\qquad energy = -49.19989156578651401\ldots$
\end{longtable}

\noindent
\textit{\textbf{Coulomb $\mathbf{1/r}$ potential -}}

Similarly the parameters have been located for the Coulomb potential:

\begin{longtable}[l]{l}
	$ \qquad a =  0.6384202837804736562\ldots$ \\
	$ \qquad b =  0.3839817226765534235\ldots$ \\
	$ \qquad c =  0.5402937261933615502\ldots$ \\
	$ \qquad d = -0.07074005957424242974\ldots$ \\
	$ \qquad e =  0.7175494647203290165\ldots$ \\
	$ \qquad f = -0.1966975337708868157\ldots$ \\
	$ \qquad g = -0.3327331842299849652\ldots$ \\
	$ \qquad h = -0.8230779376210749545\ldots$ \\
	$ \qquad i =  0.3819553076876148924\ldots$ \\
	\\
	$\text{and for the energy:}$ \\
	\\
	$ \qquad energy = 135.0894675566793420\ldots$
\end{longtable}
\vspace{-1.0em}
and the algebraic degree for the spherical code is $>480$.

\noindent
\textit{\textbf{Inverse square law $\mathbf{1/r^2}$ potential -}}

The \textit{Inverse square law $1/r^2$} potential moves the points slightly in such a way that no parameterization became possible. A direct search using the points themselves has to be done.

{
\textit{This is the first configuration of points where a parameterization apparently is impossible for the minimal energy under the Inverse Square $1/r^2$ potential.}
}

Using the parameterized structure, as in the \textit{logarithmic} or \textit{Coulomb} potentials, the energy turned out to be slightly higher than the actual minimum:

\begin{align*}
energy(parameterized)\ \mathcal{E}_p & = 118.8255644829909889\ldots \\
energy(minimum)\ \mathcal{E}_m & = 118.8255636487880665\ldots \\
\triangle (energy) = \mathcal{E}_p - \mathcal{E}_m & = 8.342029224022331538\ldots\text{e-7}
 \end{align*}

The decimal values for the coordinates of the 19 points under the \textit{Inverse square law $1/r^2$} potential had to be obtained by direct search Jacobian matrix -- Newton's method instead. It simply would be too tedious to obtain these values by the \textit{descent.3d} or \textit{percolating anneal} algorithms.

This search takes 38 parameters, and the point configuration is given below:

\noindent
\textit{\textbf{Direct Search -- 19 points --}}

\begin{longtable}[c]{l|ccc}
	\caption{Direct search for 19 points} \\
	pt & $x$ & $y$ & $z$ \\
	\hline\vspace*{-2.2ex}
	\endfirsthead
	\multicolumn{4}{c}%
	{\tablename\ \thetable\ -- 19 points search -- \textit{cont.}} \\
	pt & $x$ & $y$ & $z$ \\
	\hline\vspace*{-2.2ex}
	\endhead
	1 & $b$ & $\sqrt{1-a^2-b^2}$ & $a$ \\[0.3ex]
	2 & $d$ & $-\sqrt{1-c^2-d^2}$ & $c$ \\[0.3ex]
	3 & $f$ & $\sqrt{1-e^2-f^2}$ & $e$ \\[0.3ex]
	4 & $h$ & $-\sqrt{1-g^2-h^2}$ & $g$ \\[0.3ex]
	5 & $j$ & $-\sqrt{1-i^2-j^2}$ & $i$ \\[0.3ex]
	6 & $l$ & $\sqrt{1-k^2-l^2}$ & $k$ \\[0.3ex]
	7 & $n$ & $\sqrt{1-m^2-n^2}$ & $m$ \\[0.3ex]
	8 & $p$ & $-\sqrt{1-o^2-p^2}$ & $o$ \\[0.3ex]
	9 & $r$ & $-\sqrt{1-q^2-r^2}$ & $q$ \\[0.3ex]
	10 & $t$ & $-\sqrt{1-s^2-t^2}$ & $s$ \\[0.3ex]
	11 & $v$ & $\sqrt{1-u^2-v^2}$ & $u$ \\[0.3ex]
	12 & $x$ & $-\sqrt{1-w^2-x^2}$ & $w$ \\[0.3ex]
	13 & $z$ & $\sqrt{1-y^2-z^2}$ & $y$ \\[0.3ex]
	14 & $B$ & $\sqrt{1-A^2-B^2}$ & $A$ \\[0.3ex]
	15 & $D$ & $-\sqrt{1-C^2-D^2}$ & $C$ \\[0.3ex]
	16 & $F$ & $\sqrt{1-E^2-F^2}$ & $E$ \\[0.3ex]
	17 & $H$ & $\sqrt{1-G^2-H^2}$ & $G$ \\[0.3ex]
	18 & $J$ & $-\sqrt{1-I^2-J^2}$ & $I$ \\[0.3ex]
	19 & $L$ & $\sqrt{1-K^2-L^2}$ & $K$
\end{longtable}

The values of coordinates have been obtained to 60,013 digits precision, available in the r2.3.19.parm file.

\noindent
\textit{\textbf{Symmetries -}}

The symmetry groups for the \textit{logarithmic} and \textit{Coulomb $1/r$} potentials is known:

\begin{center}
	\begin{tabular}{l|l}
		\multicolumn{2}{c}{Symmetries - 19 points} \\
		\hline\Tstrut
		planes & [[4, 42], [10, 1], [12, 1]] \\[0.2ex]
		\hline\Tstrut
		Gram groups & [[2, 3], [4, 14], [8, 35], [19, 1]] \\
		\hline\Tstrut
		Polygons & [[4, 45], [5, 1]]
	\end{tabular}
\end{center}

The symmetry groups for the \textit{Inverse square law $1/r^2$} potential changes to:

\begin{center}
	\begin{tabular}{l|l}
		\multicolumn{2}{c}{Symmetries - 19 points} \\
		\hline\Tstrut
		planes & [] \\[0.2ex]
		\hline\Tstrut
		Gram groups & [[2, 9], [4, 81], [19, 1]] \\
		\hline\Tstrut
		Polygons & []
	\end{tabular}
\end{center}

Please note that the 81 groups of tetrads in the Gram matrix hints at the regular order of this configuration, but no planar structures are embedded among the points.

\subsection{20 points}
The optimal arrangement for 20 points contains 4 triangles, one hexagon and two poles. The hexagon is on the equator, 2 sets of triangles, rotated $180^{\circ}$ with respect to each other, on one side, then a pole, similarly for the opposite side. The arrangement is [1:3:3:6:3:3:1].

\begin{figure}[H]
	\begin{center}
		\includegraphics[type=pdf,ext=pdf,read=pdf,height=1in,width=1in,angle=0]{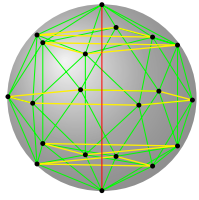}
		\caption{20 points.}
		\label{fig:20pts}
	\end{center}
\end{figure}

\noindent
\textit{\textbf{Parametric structure -}}

\begin{longtable}[c]{r|ccc}
	\caption{Parameterization for 20 points} \\
	pt & $x$ & $y$ & $z$ \\
	\hline\vspace*{-2.2ex}
	\endfirsthead
	\multicolumn{4}{c}%
	{\tablename\ \thetable\ -- 20 points parameters -- \textit{continued}} \\
	pt & $x$ & $y$ & $z$ \\
	\hline\vspace*{-2.2ex}
	\endhead
	1 & $0$ & $0$ & $1$ \\[0.5ex]
	2 & $\sqrt{1-a^2}$ & $0$ & $a$ \\[0.5ex]
	3 & $-\frac{\sqrt{1-a^2}}{2}$ & $\frac{\sqrt{3}\sqrt{1-a^2}}{2}$ & $a$ \\[0.5ex]
	4 & $-\frac{\sqrt{1-a^2}}{2}$ & $-\frac{\sqrt{3}\sqrt{1-a^2}}{2}$ & $a$ \\[0.5ex]
	5 & $-\sqrt{1-b^2}$ & $0$ & $b$ \\[0.5ex]
	6 & $\frac{\sqrt{1-b^2}}{2}$ & $\frac{\sqrt{3}\sqrt{1-b^2}}{2}$ & $b$ \\[0.5ex]
	7 & $\frac{\sqrt{1-b^2}}{2}$ & $-\frac{\sqrt{3}\sqrt{1-b^2}}{2}$ & $b$ \\[0.5ex]
	8 & $c$ & $\sqrt{1-c^2}$ & $0$ \\[0.5ex]
	9 & $d$ & $\sqrt{1-d^2}$ & $0$ \\[0.5ex]
	10 & $e$ & $\sqrt{1-e^2}$ & $0$ \\[0.5ex]
	11 & $e$ & $-\sqrt{1-e^2}$ & $0$ \\[0.5ex]
	12 & $d$ & $-\sqrt{1-d^2}$ & $0$ \\[0.5ex]
	13 & $c$ & $-\sqrt{1-c^2}$ & $0$ \\[0.5ex]
	14 & $\frac{\sqrt{1-b^2}}{2}$ & $-\frac{\sqrt{3}\sqrt{1-b^2}}{2}$ & $-b$ \\[0.5ex]
	15 & $\frac{\sqrt{1-b^2}}{2}$ & $\frac{\sqrt{3}\sqrt{1-b^2}}{2}$ & $-b$ \\[0.5ex]
	16 & $-\sqrt{1-b^2}$ & $0$ & $-b$ \\[0.5ex]
	17 & $-\frac{\sqrt{1-a^2}}{2}$ & $-\frac{\sqrt{3}\sqrt{1-a^2}}{2}$ & $-a$ \\[0.5ex]
	18 & $-\frac{\sqrt{1-a^2}}{2}$ & $\frac{\sqrt{3}\sqrt{1-a^2}}{2}$ & $-a$ \\[0.5ex]
	19 & $\sqrt{1-a^2}$ & $0$ & $-a$ \\[0.5ex]
	20 & $0$ & $0$ & $-1$
\end{longtable}

The parametric structure uses 5 parameters, $a$, $b$, $c$, $d$, and $e$.

\noindent
\textit{\textbf{logarithmic potential -}}

The decimal values for the parameters have been determined to 50,014 digits. However the algebraic degree is $>360$ for this potential.

\begin{longtable}[l]{l}
	$ \qquad a = 0.6980302156823009886\ldots$ \\
	$ \qquad b = 0.5371537469422876921\ldots$ \\
	$ \qquad c = 0.9069309469841543857\ldots$ \\
	$ \qquad d = -0.08862688856468946150\ldots$ \\
	$ \qquad e = -0.8183040584194649242\ldots$ \\
	\\
	$\qquad energy = -54.01112997458416356\ldots$
\end{longtable}

\noindent
\textit{\textbf{Coulomb $\mathbf{1/r}$ potential -}}

Similarly the parameters have been located for the Coulomb potential:

\begin{tabular}{l}
	$ \qquad a = 0.6934858852542112728\ldots$ \\
	$ \qquad b = 0.5422222757303874383\ldots$ \\
	$ \qquad c = 0.9062454955978355516\ldots$ \\
	$ \qquad d = -0.08700893177789533511\ldots$ \\
	$ \qquad e = -0.8192365638199402165\ldots$ \\
	$\text{and for the energy:}$ \\
	$ \qquad energy = 150.8815683337565034\ldots$
\end{tabular}

and the algebraic degree for the spherical code is $>360$.

\noindent
\textit{\textbf{Inverse square law $\mathbf{1/r^2}$ potential -}}

The parameters also have been determined for the Inverse square law potential:

\begin{tabular}{l}
	$ \qquad a = 0.6894142938281278435\ldots$ \\
	$ \qquad b = 0.5468551803398145473\ldots$ \\
	$ \qquad c = 0.9057332382851492562\ldots$ \\
	$ \qquad d = -0.08580330607639148262\ldots$ \\
	$ \qquad e = -0.8199299322087577735\ldots$ \\
	$\text{and for the energy:}$ \\
	$ \qquad energy = 133.9369785684329726\ldots$
\end{tabular}

\noindent
While the decimal precision of the parameter values is accurately known to 50,014 digits, the algebraic polynomials were not recovered. The polynomials must have degree $>360$, if they are found.

\noindent
\textit{\textbf{Symmetries -}}

The symmetry groups for all 3 potentials is identical:
\begin{center}
	\begin{tabular}{l|l}
		\multicolumn{2}{c}{Symmetries - 20 points} \\
		\hline\Tstrut
		planes & [[4, 54], [8, 3], [20, 3], [24, 1]] \\[0.2ex]
		\hline\Tstrut
		Gram groups & [[2, 1], [6, 5], [12, 11], [20, 1], [24, 9]] \\
		\hline\Tstrut
		Polygons & [[3, 4], [4, 60], [6, 4]]
	\end{tabular}
\end{center}

\subsection{21 points}
The polyhedra for 21 points has partial symmetry, but it is rather unique and required 18 parameters to fully constrain the figure.
While the symmetry polygon groups showed squares and 2 heptagons, trying to set up the figure aligned on a square showed too many parameters would have to be used. However aligning upon a heptagon revealed that the heptagon is NOT regular, and hence 7 parameters were required for that alone.

\begin{figure}[H]
	\begin{center}
		\includegraphics[type=pdf,ext=pdf,read=pdf,height=1in,width=1in,angle=0]{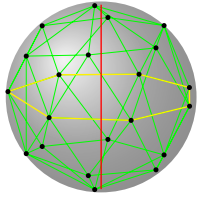}
		\caption{21 points.}
		\label{fig:21pts}
	\end{center}
\end{figure}

The heptagon is highlighted in yellow.

\noindent
\textit{\textbf{Algebraic parameterization -}}

The algebraic parameterized structure is given below:

\begin{longtable}[c]{r|ccc}
	\caption{Parameterization for 21 points} \\
	pt & $x$ & $y$ & $z$ \\
	\hline\vspace*{-2.2ex}
	\endfirsthead
	\multicolumn{4}{c}%
	{\tablename\ \thetable\ -- 21 points parameters -- \textit{continued}} \\
	pt & $x$ & $y$ & $z$ \\
	\hline\vspace*{-2.2ex}
	\endhead
	1 & $b$ & $-\sqrt{1-a^2-b^2}$ & $a$ \\[0.5ex]
	2 & $d$ & $\sqrt{1-c^2-d^2}$ & $c$ \\[0.5ex]
	3 & $\sqrt{1-e^2}$ & $0$ & $e$ \\[0.5ex]
	4 & $f$ & $\sqrt{1-e^2-f^2}$ & $e$ \\[0.5ex]
	5 & $h$ & $-\sqrt{1-g^2-h^2}$ & $g$ \\[0.5ex]
	6 & $i$ & $-\sqrt{1-g^2-i^2}$ & $g$ \\[0.5ex]
	7 & $k$ & $-\sqrt{1-j^2-k^2}$ & $j$ \\[0.5ex]
	8 & $l$ & $\sqrt{1-l^2}$ & $0$ \\[0.5ex]
	9 & $m$ & $-\sqrt{1-m^2}$ & $0$ \\[0.5ex]
	10 & $n$ & $\sqrt{1-n^2}$ & $0$ \\[0.5ex]
	11 & $o$ & $\sqrt{1-o^2}$ & $0$ \\[0.5ex]
	12 & $p$ & $-\sqrt{1-p^2}$ & $0$ \\[0.5ex]
	13 & $q$ & $\sqrt{1-q^2}$ & $0$ \\[0.5ex]
	14 & $r$ & $-\sqrt{1-r^2}$ & $0$ \\[0.5ex]
	15 & $k$ & $-\sqrt{1-j^2-k^2}$ & $-j$ \\[0.5ex]
	16 & $h$ & $-\sqrt{1-g^2-h^2}$ & $-g$ \\[0.5ex]
	17 & $i$ & $-\sqrt{1-g^2-i^2}$ & $-g$ \\[0.5ex]
	18 & $f$ & $\sqrt{1-e^2-f^2}$ & $-e$ \\[0.5ex]
	19 & $\sqrt{1-e^2}$ & $0$ & $-e$ \\[0.5ex]
	20 & $d$ & $\sqrt{1-c^2-d^2}$ & $-c$ \\[0.5ex]
	21 & $b$ & $-\sqrt{1-a^2-b^2}$ & $-a$
\end{longtable}

\noindent
\textbf{NOTE:} The parameterization above suggests that ALL point configurations could be set up this way and only require 2 parameters per points. Indeed this is the case, however that reflects in spherical coordinates $\left(\phi,\theta,\rho\right)$, since $\rho=1$, then only $\theta$ and $\phi$ change. We see this behavior here in the parameter tables since all points are on the S2 sphere.

\noindent
\textit{\textbf{Parameterization values --}}

\begin{longtable}[c]{c|c|c|c}
	\caption{Parameter values for 21 points} \\
	Parameter & log & 1/r & $1/r^2$ \\
	\hline\vspace*{-2.2ex}
	\endfirsthead
	\multicolumn{4}{c}%
	{\tablename\ \thetable\ -- 21 points parameters -- \textit{continued}} \\
	Parameter & log & 1/r & $1/r^2$ \\
	\hline\vspace*{-2.2ex}
	\endhead
	$a$ & 0.9972153605278213274 & 0.9973393074145580182 & 0.9974615439622348757 \\
	$b$ & -0.03489481469930750138 & -0.03390188354653576863 & -0.03293150836912674948 \\
	$c$ & 0.7018770304414993892 & 0.7012708203071860937 & 0.7008902972449793164 \\
	$d$ & 0.3332926741502225372 & 0.3315324960089965675 & 0.3298686089650181290 \\
	$e$ & 0.6609394617323170835 & 0.6611862673594078435 & 0.6616005559006836318 \\
	$f$ & -0.4218347726064498261 & -0.4257172448949750270 & -0.4290941670736423810 \\
	$g$ & 0.6525255598288283041 & 0.6544907612808201166 & 0.6553227532392476260 \\
	$h$ & 0.3769293850991284010 & 0.3751544353166537904 & 0.3745038347929542487 \\
	$i$ & -0.7555609085128865075 & -0.7533910088183089325 & -0.7522598853817577594 \\
	$j$ & 0.5343775883874167618 & 0.5288815289849839779 & 0.5259574583890506802 \\
	$k$ & -0.3955004858544503296 & -0.3946867403884559194 & -0.3933394069973158405 \\
	$l$ & 0.9527159863609153685 & 0.9512398555464450847 & 0.9501902243866868108 \\
	$m$ & 0.8479301445451463908 & 0.8510104453928346264 & 0.8533520217043249261 \\
	$n$ & 0.4679117597536944551 & 0.4650509689418503278 & 0.4624742817702815718 \\
	$o$ & -0.2842264525816205634 & -0.2858051644714194939 & -0.2881343801698224016 \\
	$p$ & 0.1078533280398298287 & 0.1205546286123695020 & 0.1296658538440677401 \\
	$q$ & -0.9150659309084105854 & -0.9153200328835932232 & -0.9158595061128591815 \\
	$r$ & -0.8828594370392068651 & -0.8858082887690374015 & -0.8873655851079998991 \\
	\hline\Tstrut
	$energy$ & -59.00091213514238479 & 167.6416223992704410 & 150.3251227433946007
\end{longtable}

\noindent
\textit{\textbf{Symmetries -}}

The symmetry groups for all 3 potentials are identical:
\begin{center}
	\begin{tabular}{l|l}
		\multicolumn{2}{c}{Symmetries - 20 points} \\
		\hline\Tstrut
		planes & [[4, 32], [8, 1], [35, 2]] \\[0.2ex]
		\hline\Tstrut
		Gram groups & [[2, 6], [4, 24], [8, 39], [21, 1]] \\
		\hline\Tstrut
		Polygons & [[4, 34], [7, 2]]
	\end{tabular}
\end{center}

\subsection{22 points}
The optimal configuration for 22 points consists of one pole, then two acute triangles, then a stretched hexagon near the equator, then three more triangles, the triangles and hexagon are parallel-plane and the line from the pole passes through the centroids of all the polygons. The arrangement is [1:3:3:6:3:3:3].

\begin{figure}[ht]
	\begin{center}
		\includegraphics[type=pdf,ext=pdf,read=pdf,height=1in,width=1in,angle=0]{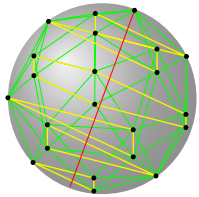}
		\caption{22 points.}
		\label{fig:22pts}
	\end{center}
\end{figure}

\noindent
\textit{\textbf{Parametric structure -}}
\begin{center}
	\begin{longtable}[c]{l|ccc}
		\caption{Parameterized structure for 22 points} \\
		pt & $x$ & $y$ & $z$ \\
		\hline\vspace*{-2.2ex}
		\endfirsthead
		\multicolumn{4}{c}%
		{\tablename\ \thetable\ -- 22 points parameters -- \textit{cont.}} \\
		pt & $x$ & $y$ & $z$ \\
		\hline\vspace*{-2.2ex}
		\endhead
		1 & $0$ & $0$ & $1$ \\[0.5ex]
		2 & $\sqrt{1-a^2}$ & $0$ & $a$ \\[0.5ex]
		3 & $-\frac{\sqrt{1-a^2}}{2}$ & $\frac{\sqrt{3-3a^2}}{2}$ & $a$ \\[0.5ex]
		4 & $-\frac{\sqrt{1-a^2}}{2}$ & $-\frac{\sqrt{3-3a^2}}{2}$ & $a$ \\[0.5ex]
		5 & $-\sqrt{\frac{2}{3}}$ & $0$ & $\sqrt{\frac{1}{3}}$ \\[0.5ex]
		6 & $\sqrt{\frac{1}{6}}$ & $\sqrt{\frac{1}{2}}$ & $\sqrt{\frac{1}{3}}$ \\[0.5ex]
		7 & $\sqrt{\frac{1}{6}}$ & $-\sqrt{\frac{1}{2}}$ & $\sqrt{\frac{1}{3}}$ \\[0.5ex]
		8 & $c$ & $\sqrt{1-b^2-c^2}$ & $b$ \\[0.5ex]
		9 & $c$ & $-\sqrt{1-b^2-c^2}$ & $b$ \\[0.5ex]
		10 & $d$ & $\sqrt{1-b^2-d^2}$ & $b$ \\[0.5ex]
		11 & $d$ & $-\sqrt{1-b^2-d^2}$ & $b$ \\[0.5ex]
		12 & $e$ & $\sqrt{1-b^2-e^2}$ & $b$ \\[0.5ex]
		13 & $e$ & $-\sqrt{1-b^2-e^2}$ & $b$ \\[0.5ex]
		14 & $-\frac{2}{3}\sqrt{2}$ & $0$ & $-\frac{1}{3}$ \\[0.5ex]
		15 & $\frac{\sqrt{2}}{3}$ & $\sqrt{\frac{2}{3}}$ & $-\frac{1}{3}$ \\[0.5ex]
		16 & $\frac{\sqrt{2}}{3}$ & $-\sqrt{\frac{2}{3}}$ & $-\frac{1}{3}$ \\[0.5ex]
		17 & $\sqrt{\frac{2}{3}}$ & $0$ & $-\sqrt{\frac{1}{3}}$ \\[0.5ex]
		18 & $-\sqrt{\frac{1}{6}}$ & $\sqrt{\frac{1}{2}}$ & $-\sqrt{\frac{1}{3}}$ \\[0.5ex]
		19 & $-\sqrt{\frac{1}{6}}$ & $-\sqrt{\frac{1}{2}}$ & $-\sqrt{\frac{1}{3}}$ \\[0.5ex]
		20 & $-\sqrt{1-f^2}$ & $0$ & $f$ \\[0.5ex]
		21 & $\frac{\sqrt{1-f^2}}{2}$ & $\frac{\sqrt{3-3f^2}}{2}$ & $f$ \\[0.5ex]
		22 & $\frac{\sqrt{1-f^2}}{2}$ & $-\frac{\sqrt{3-3f^2}}{2}$& $f$
	\end{longtable}
\end{center}

The parametric structure uses 6 parameters, $a$, $b$, $c$, $d$, $e$, and $f$ which were discovered later on to be derived from 4 algebraic polynomials for all 3 potentials.

\noindent
\textit{\textbf{logarithmic potential -}}

The algebraic spherical code has been found for this potential. Let the following polynomials be defined:
\begin{align*}
a(x) & = 1701x^7 + 108x^6 - 1431x^5 + 36x^4 + 63x^3 - 12x^2 + 51x - 4 \\
b(x) & = 2893401x^{14} - 9759852x^{12} + 9017244x^{10} - 44064x^8 - 2244240x^6 + 210240x^4 \\ & \quad + 160320x^2 - 2048 \\
c(x) & = 964467x^{14} - 122472x^{13} - 2866428x^{12} + 1003104x^{11} + 5935761x^{10} - 2245320x^9 - 9147978x^8 \\ & \quad + 2623104x^7 + 9881973x^6 - 1596024x^5 - 6175224x^4 + 613152x^3 + 1588647x^2 \\ & \quad - 177240x - 17378 \\
d(x) & = 133948309548816x^{28} - 615760615236312x^{26} + 1389731144428089x^{24} - 2506767990117372x^{22} \\ & \quad + 4520740051544868x^{20} - 6995578973321016x^{18} + 8332668829477872x^{16} \\ & \quad - 8304982392648096x^{14} + 7666397260136208x^{12} - 5961362623108128x^{10} \\ & \quad + 3199163617746432x^8 - 987900759644544x^6 + 133179792383232x^4 \\ & \quad - 1396375163904x^2 + 85525504 \\
\end{align*}
then let\par
\begin{tabular}{l}
	$ \qquad a = 0.7293924592514163468\ldots \text{ a root of } c(x)$\\
	$ \qquad b = 0.07935488721519961829\ldots \text{ a root of } a(x)$ \\
	$ \qquad c = 0.9139184782999168457\ldots \text{ a root of } d(x)$ \\
	$ \qquad d = -0.1122247577403226282\ldots \text{ a root of } b(x)$ \\
	$ \qquad e = -0.8016937205595942175\ldots \text{ a root of } d(x)$ \\
	$ \qquad f = -0.8881022336818155834\ldots \text{ a root of } c(x)$
\end{tabular}

\noindent
The energy polynomial has been found:
\begin{align*}
energy = & -64.20600776166054730\ldots = log(\text{a root of } 2588053261692056235625455974689531391243/ \\ & \ 1729557072883075695847564932007266772650321354613741961057191648260644369903826/ \\ & \ 006401246678598678738568569941919302870264850253367573009297895302498414553140/ \\ & \ 715440546094805040283591990253290657818454660532363918392827227506117379543067/ \\ & \ 166594839312530692887313254776929075324443352473531187552110145853486609344736/ \\ & \ 361946730348409764237396151682952596933211179572367284208094842802928765546806/ \\ & \ 25858280803615156808318543420811026366464x^7 - 258028599427320297565699822417878/ \\ & \ 447414691835527081431325168993128750160036338518659856934382190603642270863346/ \\ & \ 108471826446558913066148748906215193922277172723835258631643699007100118228277/ \\ & \ 307107176977154627362562534083776103716345014198813808230155103921214056324271/ \\ & \ 679761585269071526574448029725938494727306709384779694353828365374588764979336/ \\ & \ 364110771171344656502540597236541165437073505810744330741120340378895083612242/ \\ & \ 528193550395377763483648x^6 + 16833282612848153737290918125434464491942266388350/ \\ & \ 063506058411602844457031962546204470767807612832847818490632918706019685899720/ \\ & \ 478717587816267814606221508725783202555102779848946153953414282359984674012671/ \\ & \ 830346899901759915525272529886184809508829060974026570487244125706262152709017/ \\ & \ 548498401973383944954018105304534936430336391043520261057499541831165904310451/ \\ & \ 1023786382888434231475229042172214113076565877068453117952x^5 - 1758359160359315/ \\ & \ 881447270522902988344766529958360141706021975333634653367844975788187213732493/ \\ & \ 039342891147403469763264262820694566409942862895403973570618792818840536105398/ \\ & \ 268601426712915228853975899026832652727068525835689511275443987117173607327548/ \\ & \ 644130650836480438126489421171311131821429621046034129052984180270902758227169/ \\ & \ 3566141062580155723029502264571381508822277004111501237568929792x^4 + 3172291523/ \\ & \ 382983292164676075503792085008254704263465632846105905501620096565469972884311/ \\ & \ 378618877459115742400928628284322509755841948885019557561266995149028203640910/ \\ & \ 057572827164810635181651172723119434416807729960949538982964274143068968902469/ \\ & \ 116438128903797968256606353818987765604548162094470864517689071224034805745239/ \\ & \ 30686856373379678151892981284627021824x^3 - 163150664589794832134640591295747341/ \\ & \ 139032499953049357543403623959927532919596842015578557816050381431625915414479/ \\ & \ 204103240689019729655550684702506709990692154740800116036066948086410247520686/ \\ & \ 551825672302693414195174450065608876539509710874427490665566300517907919486998/ \\ & \ 9636628149094345949957845048383082533393369808370160107520x^2 - 745186319275479/ \\ & \ 304669131773526962850703619650044711896100641648895935753066488375470439014266/ \\ & \ 856206874061568588015127457699728364662297717750070591134596183486051122673652/ \\ & \ 008324964740436041232187298205484732348512298701699190277606460293450517423757/ \\ & \ 831900640511006856968077791413554626166784x - 21117289084166899734739778811894/ \\ & \ 824981395151344306134044029420973675659755198350503440549428467901301953548795/ \\ & \ 984157044867930729643440269086562234868314375084777836591270518332457756122607/ \\ & \ 9019503576722231301101408152581092064329226639241590648203)
\end{align*}

\noindent
\textit{\textbf{Coulomb $\mathbf{1/r}$ potential -}}

Similarly the parameters have been located for the Coulomb potential to 50,014 digits:

\begin{tabular}{l}
	$ \qquad a = 0.7277487398168151126\ldots$ \\
	$ \qquad b = 0.08072697305800621906\ldots$ \\
	$ \qquad c = 0.9146005785716013866\ldots$ \\
	$ \qquad d = -0.1141651801479598416\ldots$ \\
	$ \qquad e = -0.8004353984236415450\ldots$ \\
	$ \qquad f = -0.8892026859328275507\ldots$ \\
	$\text{and for the energy:}$ \\
	$ \qquad energy = 185.2875361493075820\ldots$
\end{tabular}

and the algebraic degree for the spherical code is $>480$.

\noindent
\textit{\textbf{Inverse square law $\mathbf{1/r^2}$ potential -}}

The algebraic spherical code has also been determined for the Inverse square law $1/r^2$ potential.

\noindent
Let the following polynomials be defined:
\begin{align*}
a(x) & = 472392x^{16} + 505197x^{15} + 1259712x^{14} + 5024997x^{13} - 1166400x^{12} - 7260111x^{11} + 31104x^{10} \\ & \quad + 3979449x^9 + 255312x^8 - 1115289x^7 - 53568x^6 + 185391x^5 - 16128x^4 - 6597x^3 + 3840x^2 \\ & \quad - 2317x + 168 \\
b(x) & = 24794911296x^{32} - 53033560272x^{31} - 311643911133x^{30} + 2743587545904x^{29} \\ & \quad + 2721936639564x^{28} - 38400761741328x^{27} + 83815977635109x^{26} + 141551317071216x^{25} \\ & \quad - 532899996509466x^{24} - 207662118306192x^{23} + 1403642340204867x^{22} + 52257770730864x^{21} \\ & \quad - 2128972490243136x^{20} + 166165430721840x^{19} + 2135902577651325x^{18} - 23798429240400x^{17} \\ & \quad - 1705479562133982x^{16} - 393602125086576x^{15} + 1430294663939265x^{14} + 579607833079440x^{13} \\ & \quad - 1171805322609252x^{12} - 442997359888176x^{11} + 656576540961831x^{10} + 284119432143696x^9 \\ & \quad - 204171773263830x^8 - 188047841542704x^7 + 67923527451993x^6 + 86965369183440x^5 \\ & \quad - 52148498654184x^4 - 15314781956592x^3 + 17324074418359x^2 \\ & \quad - 1472704852080x - 100684279602 \\
c(x) & = 111577100832x^{32} + 934931733399x^{30} - 9184697106948x^{28} - 83297167393356x^{26} \\ & \quad + 564986627305056x^{24} - 1456891190591760x^{22} + 2177887691709504x^{20} \\ & \quad - 2162107552423104x^{18} + 1503531434382336x^{16} - 724165828326144x^{14} + 218558329242624x^{12} \\ & \quad - 27591723795456x^{10} - 5709220651008x^8 + 2759555469312x^6 - 174030077952x^4 \\ & \quad - 66818031616x^2 + 924844032 \\
d(x) & = 12748236216396078174437376x^{64} - 492893683534492079224135680x^{62} \\ & \quad + 10429642314758465150355744912x^{60} - 144031497264720456920422650648x^{58} \\ & \quad + 1321175051794900688676895926513x^{56} - 8475789436765902736771956104484x^{54} \\ & \quad + 128969331786335907773221327035588x^{52} - 1025931364596101003898190164957792x^{50} \\ & \quad + 4496214128620565028367348543722528x^{48} - 12778042681906716967988820283777968x^{46} \\ & \quad + 25992180246053610277742651860107360x^{44} - 40292904661774663631462671553774784x^{42} \\ & \quad + 49687263157134060295692125439007872x^{40} - 50379261715587487079230525462747392x^{38} \\ & \quad + 43424825233410185308185561018219072x^{36} - 33205796278489660249408815854964864x^{34} \\ & \quad + 23682804205076694615078533889372672x^{32} - 16353428969235057759406169225339904x^{30} \\ & \quad + 10975670896492433808386988393784320x^{28} - 6990470597711596867787824319784960x^{26} \\ & \quad + 4066461559876684057471735940874240x^{24} - 2016358740422299334549403763826688x^{22} \\ & \quad + 736015327835052723058113304707072x^{20} - 147019341960177733715688227930112x^{18} \\ & \quad + 30052963295632558124257958166528x^{16} - 65159815573468810216690855575552x^{14} \\ & \quad + 70107867952828981554785010057216x^{12} - 41513362888712907881184617299968x^{10} \\ & \quad + 16231925015609001203188676689920x^8 - 4129883670351075080509077848064x^6 \\ & \quad + 524271703275498208600803770368x^4 - 7431861139235886307309780992x^2 \\ & \quad + 3509541404990926303002624
\end{align*}
then the parameters are roots of the polynomials:\par
\begin{tabular}{l}
	$ \qquad a = 0.7261505868988778041\ldots \text{a root of } b(x)$ \\
	$ \qquad b = 0.08205723688026544899\ldots \text{ a root of } a(x)$ \\
	$ \qquad c = 0.9152570778468481725\ldots \text{ a root of } d(x)$ \\
	$ \qquad d = -0.7992106205599150556\ldots \text{ a root of } d(x)$ \\
	$ \qquad e = -0.1160464572869331169\ldots \text{ a root of } c(x)$ \\
	$ \qquad f = -0.8902650606594087021\ldots \text{ a root of } b(x)$ \\
	\\
	$ \text{and for the energy:}$ \\
	$ \qquad energy = 167.6657856364984302\ldots \text{ a root of } \mathcal{E}_{22}$ \\
\end{tabular}

with the polynomial $\mathcal{E}_{22}$ for the minimum energy given by
\begin{align*}
\mathcal{E}_{22} & = 536870912x^{38} - 849061347328x^{37} + 617262264352768x^{36} - 273893662826430464x^{35} \\ & \quad + 83052981864535949312x^{34} - 18257258188385905803264x^{33} \\ & \quad + 3011522863648204458164224x^{32} - 380534292974393628982444032x^{31} \\ & \quad + 37268900572129290388293115904x^{30} - 2843172608546179940154814992384x^{29} \\ & \quad + 168800893079672767742190711363072x^{28} - 7746173601437329275248608481420544x^{27} \\ & \quad + 270958586078582827447953343855936736x^{26} - 7056005260481588606830526674498057488x^{25} \\ & \quad + 131522373719849859329220333364950365742x^{24} \\ & \quad - 1637183948660058102301401192337429178709x^{23} \\ & \quad + 11772210466900020878144287476608456973929x^{22} \\ & \quad - 30018118814568920564060273990791883857916x^{21} \\ & \quad - 81100981726919064693993121083775856258463x^{20} \\ & \quad - 35821444597459154524448623101842701906496x^{19} \\ & \quad + 224662365901753751865491463740755617957886x^{18} \\ & \quad + 28013628473795083245016389714934049266505x^{17} \\ & \quad - 9024876107074338167148106855556229836109x^{16} \\ & \quad + 22979673943661815022852220965556253219862x^{15} \\ & \quad + 110222888609848042598535263326630674717362x^{14} \\ & \quad - 31315146307987784058942022287800793532353x^{13} \\ & \quad - 110723897466037210791418425047730181551966x^{12} \\ & \quad - 113040734960226389531886021838560912902286x^{11} \\ & \quad - 271462393896255999103805115003516413714103x^{10} \\ & \quad - 11155703477528266731656505030457778824008x^{9} \\ & \quad + 46014173721142348782036310108424710409637x^{8} \\ & \quad + 55641130745524838538237349526328622733349x^{7} \\ & \quad - 23530383295416886415683809682696370548770x^{6} \\ & \quad + 123772707320625250642615457982390862333671x^{5} \\ & \quad - 57664040277072034991139555118230893984481x^{4} \\ & \quad + 217836134163059146817210517684826659158946x^{3} \\ & \quad - 78680956102867977550797844174034524314987x^{2} \\ & \quad + 18660531450590089490994811367532675538869x \\ & \quad - 91437892323110947361247296887555519526883
\end{align*}

\noindent
\textit{\textbf{Symmetries -}}

The symmetry groups for all 3 potentials is identical:
\begin{center}
	\begin{tabular}{l|l}
		\multicolumn{2}{c}{Symmetries - 22 points} \\
		\hline\Tstrut
		planes & [[4, 120], [12, 3], [20, 6], [25, 4]] \\[0.2ex]
		\hline\Tstrut
		Gram groups & [[6, 1], [12, 2], [22, 1], [24, 10], [48, 4]] \\
		\hline\Tstrut
		Polygons & [[3, 20], [4, 129], [6, 10]]
	\end{tabular}
\end{center}

\subsection{23 points}
For 23 points, when looking at the only symmetry, 7 parallel triangles, more carefully, it was discovered that this figure has a balanced [1:3:3:3:3:3:3:3:1] arrangement, with an equilateral triangle on the equator, and 3 more on either side, but significantly rotated with respect to each other.

\begin{figure}[ht]
	\begin{center}
		\includegraphics[type=pdf,ext=pdf,read=pdf,height=1in,width=1in,angle=0]{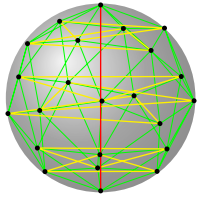}
		\caption{23 points.}
		\label{fig:23pts}
	\end{center}
\end{figure}

\noindent
\textit{\textbf{Parameterized Structure --}}

After noticing this key fact of 7 parallel triangles, a parameterization was successfully obtained.

\begin{longtable}[c]{r|ccc}
	\caption{Parameterization for 23 points --} \\
	pt & $x$ & $y$ & $z$ \\
	\hline\vspace*{-2.2ex}
	\endfirsthead
	\multicolumn{4}{c}%
	{\tablename\ \thetable\ -- 23 points -- \textit{continued}} \\
	pt & $x$ & $y$ & $z$ \\
	\hline\vspace*{-2.2ex}
	\endhead
	1 & $0$ & $0$ & $1$ \\[0.65ex]
	2 & $b$ & $\sqrt{1-a^2-b^2}$ & $a$ \\[0.65ex]
	3 & $\frac{-b-\sqrt{3}}{2}\sqrt{1-a^2-b^2}$ & $\frac{b\sqrt{3}-\sqrt{1-a^2-b^2}}{2}$ & $a$ \\[0.65ex]
	4 & $\frac{-b+\sqrt{3}}{2}\sqrt{1-a^2-b^2}$ & $\frac{-b\sqrt{3}-\sqrt{1-a^2-b^2}}{2}$ & $a$ \\[0.65ex]
	5 & $d$ & $-\sqrt{1-c^2-d^2}$ & $c$ \\[0.65ex]
	6 & $\frac{-d+\sqrt{3}}{2}\sqrt{1-c^2-d^2}$ & $\frac{d\sqrt{3}+\sqrt{1-c^2-d^2}}{2}$ & $c$ \\[0.65ex]
	7 & $\frac{-d-\sqrt{3}}{2}\sqrt{1-c^2-d^2}$ & $\frac{-d\sqrt{3}+\sqrt{1-c^2-d^2}}{2}$ & $c$ \\[0.65ex]
	8 & $f$ & $\sqrt{1-e^2-f^2}$ & $e$ \\[0.65ex]
	9 & $\frac{-f-\sqrt{3}}{2}\sqrt{1-e^2-f^2}$ & $\frac{f\sqrt{3}-\sqrt{1-e^2-f^2}}{2}$ & $e$ \\[0.65ex]
	10 & $\frac{-f+\sqrt{3}}{2}\sqrt{1-e^2-f^2}$ & $\frac{-f\sqrt{3}-\sqrt{1-e^2-f^2}}{2}$ & $e$ \\[0.65ex]
	11 & $1$ & $0$ & $0$ \\[0.65ex]
	12 & $-\frac{1}{2}$ & $\frac{\sqrt{3}}{2}$ & $0$ \\[0.65ex]
	13 & $-\frac{1}{2}$ & $-\frac{\sqrt{3}}{2}$ & $0$ \\[0.65ex]
	14 & $f$ & $-\sqrt{1-e^2-f^2}$ & $-e$ \\[0.65ex]
	15 & $\frac{-f+\sqrt{3}}{2}\sqrt{1-e^2-f^2}$ & $\frac{f\sqrt{3}+\sqrt{1-e^2-f^2}}{2}$ & $-e$ \\[0.65ex]
	16 & $\frac{-f-\sqrt{3}}{2}\sqrt{1-e^2-f^2}$ & $\frac{-f\sqrt{3}+\sqrt{1-e^2-f^2}}{2}$ & $-e$ \\[0.65ex]
	17 & $d$ & $\sqrt{1-c^2-d^2}$ & $-c$ \\[0.65ex]
	18 & $\frac{-d-\sqrt{3}}{2}\sqrt{1-c^2-d^2}$ & $\frac{d\sqrt{3}-\sqrt{1-c^2-d^2}}{2}$ & $-c$ \\[0.65ex]
	19 & $\frac{-d+\sqrt{3}}{2}\sqrt{1-c^2-d^2}$ & $\frac{-d\sqrt{3}-\sqrt{1-c^2-d^2}}{2}$ & $-c$ \\[0.65ex]
	20 & $b$ & $-\sqrt{1-a^2-b^2}$ & $-a$ \\[0.65ex]
	21 & $\frac{-b+\sqrt{3}}{2}\sqrt{1-a^2-b^2}$ & $\frac{b\sqrt{3}+\sqrt{1-a^2-b^2}}{2}$ & $-a$ \\[0.65ex]
	22 & $\frac{-b-\sqrt{3}}{2}\sqrt{1-a^2-b^2}$ & $\frac{-b\sqrt{3}+\sqrt{1-a^2-b^2}}{2}$ & $-a$ \\[0.65ex]
	23 & $0$ & $0$ & $-1$
\end{longtable}

The values to 19 digits of the parameters optimized for the minimal solutions of 23 points are:

\pagebreak
\begin{longtable}{c|c|c|c}
	Parameter & log & 1/r & $1/r^2$ \\
	\hline\Tstrut
	$a$ & 0.7340917690492719012 & 0.7310318038872114760 & 0.7281967383910171980 \\
	$b$ & 0.6394869826704419718 & 0.6438361580150366677 & 0.6480289594830250432 \\
	$c$ & 0.6136317746144812205 & 0.6135985939094442390 & 0.6132564278157956058 \\
	$d$ & 0.6293719233011343503 & 0.6249423636891649152 & 0.6206910170628525391 \\
	$e$ & 0.1418558536006295830 & 0.1406492270986577138 & 0.1397764134860746589 \\
	$f$ & 0.7497476107917851120 & 0.7491741097328619614 & 0.7487568007366313359 \\
	\hline\Tstrut
	$energy$ & -69.57838259252510490 & 203.9301906628785362 & 186.4037128696907775
\end{longtable}

They are all known to 50,014 digits precision, for all 3 potentials, but again, attempts to find the algebraic numbers for the parameters have failed, the algebraic degree $> 360$.

\noindent
\textit{\textbf{Symmetries -}}

The symmetry groups for all 3 potentials are identical for 23 points.

\begin{center}
	\begin{tabular}{l|l}
		\multicolumn{2}{c}{Symmetries - 23 points} \\
		\hline\Tstrut
		planes & [[7, 1]] \\[0.2ex]
		\hline\Tstrut
		Gram groups & [[2, 1], [6, 10], [12, 37], [23, 1]] \\
		\hline\Tstrut
		Polygons & [[3, 7]]
	\end{tabular}
\end{center}

The 7 parallel triangles are the only polygons embedded in this configuration, interestingly enough.

\subsection{24 points}
The optimal solution for 24 points contains embedded squares, it was discovered that a workable alignment contains 6 squares with the $z$-axis aligned through the centroids of the squares. The arrangement is [4:4:4:4:4:4], where 4 denotes a square.

\begin{figure}[ht]
	\begin{center}
		\includegraphics[type=pdf,ext=pdf,read=pdf,height=1in,width=1in,angle=0]{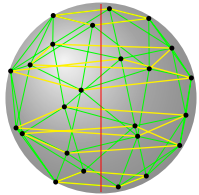}
		\caption{24 points.}
		\label{fig:24pts}
	\end{center}
\end{figure}

\noindent
\textit{\textbf{Parametric structure -}}

A parametric structure was found for 24 points, using 9 parameters, $a$, $b$, $c$, $d$, $e$, $f$, $g$, $h$, and $i$:
\begin{longtable}[c]{r|ccc}
	\caption{Parameterization for 24 points} \\
	pt & $x$ & $y$ & $z$ \\
	\hline\vspace*{-2.2ex}
	\endfirsthead
	\multicolumn{4}{c}%
	{\tablename\ \thetable\ -- 24 points -- \textit{continued\ldots}} \\
	pt & $x$ & $y$ & $z$ \\
	\hline\vspace*{-2.2ex}
	\endhead
	1 & $b$ & $\sqrt{1-a^2-b^2}$ & $a$ \\[0.5ex]
	2 & $\sqrt{1-a^2-b^2}$ & $-b$ & $a$ \\[0.5ex]
	3 & $-b$ & $-\sqrt{1-a^2-b^2}$ & $a$ \\[0.5ex]
	4 & $-\sqrt{1-a^2-b^2}$ & $b$ & $a$ \\[0.5ex]
	5 & $d$ & $\sqrt{1-c^2-d^2}$ & $c$ \\[0.5ex]
	6 & $\sqrt{1-c^2-d^2}$ & $-d$ & $c$ \\[0.5ex]
	7 & $-d$ & $-\sqrt{1-c^2-d^2}$ & $c$ \\[0.5ex]
	8 & $-\sqrt{1-c^2-d^2}$ & $d$ & $c$ \\[0.5ex]
	9 & $f$ & $\sqrt{1-e^2-f^2}$ & $e$ \\[0.5ex]
	10 & $\sqrt{1-e^2-f^2}$ & $-f$ & $e$ \\[0.5ex]
	11 & $-f$ & $-\sqrt{1-e^2-f^2}$ & $e$ \\[0.5ex]
	12 & $-\sqrt{1-e^2-f^2}$ & $f$ & $e$ \\[0.5ex]
	13 & $g$ & $\sqrt{1-e^2-g^2}$ & $-e$ \\[0.5ex]
	14 & $\sqrt{1-e^2-g^2}$ & $-g$ & $-e$ \\[0.5ex]
	15 & $-g$ & $-\sqrt{1-e^2-g^2}$ & $-e$ \\[0.5ex]
	16 & $-\sqrt{1-e^2-g^2}$ & $g$ & $-e$ \\[0.5ex]
	17 & $h$ & $\sqrt{1-c^2-h^2}$ & $-c$ \\[0.5ex]
	18 & $\sqrt{1-c^2-h^2}$ & $-h$ & $-c$ \\[0.5ex]
	19 & $-h$ & $-\sqrt{1-c^2-h^2}$ & $-c$ \\[0.5ex]
	20 & $-\sqrt{1-c^2-h^2}$ & $h$ & $-c$ \\[0.5ex]
	21 & $i$ & $\sqrt{1-a^2-i^2}$ & $-a$ \\[0.5ex]
	22 & $\sqrt{1-a^2-i^2}$ & $-i$ & $-a$ \\[0.5ex]
	23 & $-i$ & $-\sqrt{1-a^2-i^2}$ & $-a$ \\[0.5ex]
	24 & $-\sqrt{1-a^2-i^2}$ & $i$ & $-a$
\end{longtable}

Unfortunately for all 3 potentials, no algebraic polynomial of degree $< 360$ was discovered.

\noindent
\textit{\textbf{Parameterization values --}}

\begin{longtable}[c]{c|c|c|c}
	\caption{Parameter values for 24 points} \\
	Parameter & log & 1/r & $1/r^2$ \\
	\hline\vspace*{-2.2ex}
	\endfirsthead
	\multicolumn{4}{c}%
	{\tablename\ \thetable\ -- 24 points parameters -- \textit{continued}} \\
	Parameter & log & 1/r & $1/r^2$ \\
	\hline\vspace*{-2.2ex}
	\endhead
	$a$ & 0.8624744511248889895 & 0.8616157394691719940 & 0.8607395756746204637 \\
	$b$ & 0.3562946772141632259 & 0.3588999286005723978 & 0.3599495679032510744 \\
	$c$ & 0.4395840425428552288 & 0.4415679900467134637 & 0.4434758265238369016 \\
	$d$ & 0.01050597855609840377 & 0.01159706995212176130 & 0.009278990123764351111 \\
	$e$ & 0.2508060818610719374 & 0.2502719074628582630 & 0.2499131332212735685 \\
	$f$ & 0.7131286545805561954 & 0.7127751741154940637 & 0.7098035240999151707 \\
	$g$ & 0.1923792391932357068 & 0.1959747357332816245 & 0.1957137667180801688 \\
	$h$ & 0.7637175417435373665 & 0.7637562252327480142 & 0.7618476430257817818 \\
	$i$ & 0.1238325294140447517 & 0.1235412117021585841 & 0.1214706723830686456 \\
	\hline\Tstrut
	$energy$ & -75.21398478862851694 & 223.3470740518051052 & 205.6584379977072105
\end{longtable}

\noindent
It is hoped that the algebraic polynomials might be recovered for these configurations of 24 points, at least under 1 potential. Some of the algebraic numbers were recovered, but not all.

\noindent
\textit{\textbf{Symmetries -}}

The symmetry groups are identical for 24 points under all 3 potentials:

\begin{center}
	\begin{tabular}{l|l}
		\multicolumn{2}{c}{Symmetries - 24 points} \\
		\hline\Tstrut
		planes & [[8, 4], [24, 3]] \\[0.2ex]
		\hline\Tstrut
		Gram groups & [[24, 10], [48, 7]] \\
		\hline\Tstrut
		Polygons & [[3, 32], [4, 18]]
	\end{tabular}
\end{center}

\subsection{25 points}
After some consideration of the symmetry groups, it was decided to parameterize the point set, but the embedded pentagon was irregular, however it was on an equator. It turns out that 24 parameters are necessary to constrain the structure of the point configuration.

\begin{figure}[ht]
	\begin{center}
		\includegraphics[type=pdf,ext=pdf,read=pdf,height=1in,width=1in,angle=0]{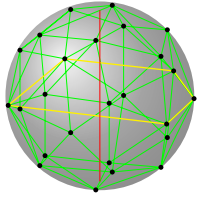}
		\caption{25 points.}
		\label{fig:25pts}
	\end{center}
\end{figure}
The embedded pentagon is highlighted in yellow in figure \ref{fig:25pts}.

\noindent
\textit{\textbf{Parametric structure -}}

A parametric structure was found for 25 points, using 24 parameters, $a - x$, but no polygon was used as an embedded constraint. The situation is actually a direct point constrained type of search.
\begin{longtable}[c]{r|ccc}
	\caption{Parameterization for 25 points} \\
	pt & $x$ & $y$ & $z$ \\
	\hline\vspace*{-2.2ex}
	\endfirsthead
	\multicolumn{4}{c}%
	{\tablename\ \thetable\ -- 25 points -- \textit{continued\ldots}} \\
	pt & $x$ & $y$ & $z$ \\
	\hline\vspace*{-2.2ex}
	\endhead
	1 & $b$ & $\sqrt{1-a^2-b^2}$ & $a$ \\[0.5ex]
	2 & $d$ & $-\sqrt{1-c^2-d^2}$ & $c$ \\[0.5ex]
	3 & $f$ & $-\sqrt{1-e^2-f^2}$ & $e$ \\[0.5ex]
	4 & $h$ & $\sqrt{1-g^2-h^2}$ & $g$ \\[0.5ex]
	5 & $j$ & $\sqrt{1-i^2-j^2}$ & $i$ \\[0.5ex]
	6 & $l$ & $\sqrt{1-k^2-l^2}$ & $k$ \\[0.5ex]
	7 & $n$ & $\sqrt{1-m^2-n^2}$ & $m$ \\[0.5ex]
	8 & $p$ & $-\sqrt{1-o^2-p^2}$ & $o$ \\[0.5ex]
	9 & $r$ & $-\sqrt{1-q^2-r^2}$ & $q$ \\[0.5ex]
	10 & $t$ & $-\sqrt{1-s^2-t^2}$ & $s$ \\[0.5ex]
	11 & $1$ & $0$ & $0$ \\[0.5ex]
	12 & $u$ & $\sqrt{1-u^2}$ & $0$ \\[0.5ex]
	13 & $v$ & $\sqrt{1-v^2}$ & $0$ \\[0.5ex]
	14 & $w$ & $\sqrt{1-w^2}$ & $0$ \\[0.5ex]
	15 & $x$ & $-\sqrt{1-x^2}$ & $0$ \\[0.5ex]
	16 & $t$ & $-\sqrt{1-s^2-t^2}$ & $-s$ \\[0.5ex]
	17 & $r$ & $-\sqrt{1-q^2-r^2}$ & $-q$ \\[0.5ex]
	18 & $p$ & $-\sqrt{1-o^2-p^2}$ & $-o$ \\[0.5ex]
	19 & $n$ & $\sqrt{1-m^2-n^2}$ & $-m$ \\[0.5ex]
	20 & $l$ & $\sqrt{1-k^2-l^2}$ & $-k$ \\[0.5ex]
	21 & $j$ & $\sqrt{1-i^2-j^2}$ & $-i$ \\[0.5ex]
	22 & $h$ & $\sqrt{1-g^2-h^2}$ & $-g$ \\[0.5ex]
	23 & $f$ & $-\sqrt{1-e^2-f^2}$ & $-e$ \\[0.5ex]
	24 & $d$ & $-\sqrt{1-c^2-d^2}$ & $-c$ \\[0.5ex]
	25 & $b$ & $\sqrt{1-a^2-b^2}$ & $-a$
\end{longtable}

\noindent
\textit{\textbf{Parameterization values --}}

\begin{longtable}[c]{c|c|c|c}
	\caption{Parameter values for 25 points} \\
	Parameter & log & 1/r & $1/r^2$ \\
	\hline\vspace*{-2.2ex}
	\endfirsthead
	\multicolumn{4}{c}%
	{\tablename\ \thetable\ -- 25 points parameters -- \textit{continued}} \\
	Parameter & log & 1/r & $1/r^2$ \\
	\hline\vspace*{-2.2ex}
	\endhead
	$a$ & 0.9748400482776597105 & 0.9751821992500684939 & 0.9755092092721821243 \\
	$b$ & -0.03422946851960271674 & -0.03123602076239851574 & -0.02863828272626720468 \\
	$c$ & 0.8738307687010033467 & 0.8720231895972047146 & 0.8703290253723059420 \\
	$d$ & 0.3119389150694342939 & 0.3114825717454717616 & 0.3114200305543287559 \\
	$e$ & 0.7632084283804168806 & 0.7640457322810302437 & 0.7651370478385824907 \\
	$f$ & -0.3705972251105773610 & -0.3727771061620920591 & -0.3744006206848329942 \\
	$g$ & 0.7199369033886241052 & 0.7195890586039713739 & 0.7188736896484927827 \\
	$h$ & -0.6580466682056721329 & -0.6588336948752226961 & -0.6598248588308915649 \\
	$i$ & 0.6485304089108194685 & 0.6502690064498257011 & 0.6515986500287267937 \\
	$j$ & 0.7410345585219819283 & 0.7399075515191026582 & 0.7390887130742449337 \\
	$k$ & 0.6296976591376905007 & 0.6273440654023089354 & 	0.6253307435492714428 \\
	$l$ & 0.2901253289916902046 & 0.2902387849243536460 & 0.2905173872831476453 \\
	$m$ & 0.4053729526955878345 & 0.4060219631153653671 & 0.4058193304142052376 \\
	$n$ & -0.4435906951787374045 & -0.4423711901494016314 & -0.4414501091940891950 \\
	$o$ & 0.3588865023298887721 & 0.3597520916705570228 & 0.3608411398452656937 \\
	$p$ & 0.1061337532401854996 & 0.1069139035885959770 & 0.1073401811458837796 \\
	$q$ & 0.3483221777082064666 & 0.3482226170044710232 & 		0.3481786113138171417 \\
	$r$ & 0.7289803113831196293 & 0.7312643708675318560 & 0.7333559841824886316 \\
	$s$ & 0.3436290173807913696 & 0.3419594895646984177 & 0.3407312707655539995 \\
	$t$ & -0.8865497960927564253 & -0.8873815641901340489 & -0.8878527312909185546 \\
	$u$ & 0.7548398673944465477 & 0.7544231966150258870 & 0.7540570111554236954 \\
	$v$ & 0.1113197258584147103 & 0.1107650830530035312 & 0.1108931283047425362 \\
	$w$ & -0.9192023080455600329 & -0.9199013678639954006 & -0.9208416685184383615 \\
	$x$ & -0.5173553134054322557 & -0.5187105602738805705 & -0.5198136192911186486 \\
	\hline\Tstrut
	$energy$ & -80.99750999019678051 & 243.8127602987657369 & 226.5450772595475468
\end{longtable}

For all 3 laws, the parameters were obtained to 50,014 digits but no algebraic numbers were recovered, their degree is $>420$.

\noindent
\textit{\textbf{Minimal Energy values -}}

The coordinates for 25 points are known to 77 digits for the \textit{log} potential and 38 digits for the other two. The minimal energies have been determined for all 3 potentials as well.
\begin{center}
	\begin{tabular}{l|l}
		\multicolumn{2}{c}{Minimal Energy - 25 points} \\
		\hline\Tstrut
		logarithmic & -80.99750999019678051\ldots \\[0.2ex]
		\hline\Tstrut
		Coulomb & 243.8127602987657369\ldots \\[0.2ex]
		\hline\Tstrut
		Inverse square law & 226.5450772595475468\ldots
	\end{tabular}
\end{center}

\noindent
\textit{\textbf{Symmetries -}}

The symmetry groups are identical for 25 points under all 3 potentials:

\begin{center}
	\begin{tabular}{l|l}
		\multicolumn{2}{c}{Symmetries - 25 points} \\
		\hline\Tstrut
		planes & [[4, 45], [10, 1]] \\[0.2ex]
		\hline\Tstrut
		Gram groups & [[2, 20], [4, 140], [25, 1]] \\
		\hline\Tstrut
		Polygons & [[4, 45], [5, 1]]
	\end{tabular}
\end{center}

\subsection{26 points}
The configuration for 26 points does not admit a parameterization for all 3 potentials. This is the first time that this occurs for all 3 potentials, although 19 points for the \textit{Inverse square law $1/r^2$} potential does not admit a parameterization.

\textit{It is quite an interesting discovery that minimal energy configurations can contain no embedded polygons (except the trivial triangles).}

\begin{figure}[ht]
	\begin{center}
		\includegraphics[type=pdf,ext=pdf,read=pdf,height=1in,width=1in,angle=0]{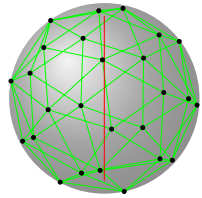}
		\caption{26 points.}
		\label{fig:26pts}
	\end{center}
\end{figure}
\noindent
\textit{\textbf{Minimal Energy values -}}

The coordinates for 26 points are known to 77 digits for the \textit{log} potential and 38 digits for the other two. The minimal energies have been determined for all 3 potentials as well.
\begin{center}
	\begin{tabular}{l|l}
		\multicolumn{2}{c}{Minimal Energy - 26 points} \\
		\hline\Tstrut
		logarithmic & -87.00942305704718702\ldots \\[0.2ex]
		\hline\Tstrut
		Coulomb & 265.1333263173565392\ldots \\[0.2ex]
		\hline\Tstrut
		Inverse square law & 248.2671389221036336\ldots
	\end{tabular}
\end{center}

\noindent
\textit{\textbf{Symmetries -}}

The symmetry groups are identical for 26 points under all 3 potentials:

\begin{center}
	\begin{tabular}{l|l}
		\multicolumn{2}{c}{Symmetries - 26 points} \\
		\hline\Tstrut
		planes & [] \\[0.2ex]
		\hline\Tstrut
		Gram groups & [[2, 13], [4, 156], [26, 1]] \\
		\hline\Tstrut
		Polygons & []
	\end{tabular}
\end{center}

\subsection{27 points}
Initially, the optimal configuration for 27 points was ignored for parameterization, due to the odd number of points, but after examining the symmetry groups, it was decided to check all 3 potentials anyways.

\begin{figure}[ht]
	\begin{center}
		\includegraphics[type=pdf,ext=pdf,read=pdf,height=1in,width=1in,angle=0]{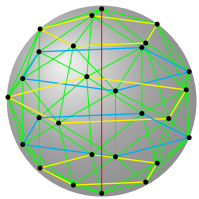}
		\caption{27 points.}
		\label{fig:27pts}
	\end{center}
\end{figure}

A preferable alignment was obtained, which put 3 pentagons in a parallel-plane alignment, with 2 poles, the pole axis passing through the centers of the pentagons. Conveniently, one pentagon is located on the equator. This arrangement is [1:5:5:5:5:5:1]. The pentagons in figure \ref{fig:27pts} are alternately colored yellow and cyan.

\noindent
\textit{\textbf{Constrained parameterization -}}

\noindent
A parameterization was thus obtained:

\noindent
Let $a$ and $b$ be 2 algebraic parameters, and $m$ and $n$ from $a$ and $b$ respectively, then let $c_1$, $c_2$, $s_1$ and $s_2$ be constants for creating a regular pentagon; all defined as:

\begin{align*}
c_1 = & \frac{\sqrt{5}-1}{4} \qquad\qquad c_2 = \frac{\sqrt{5}+1}{4} \\
s_1 = & \frac{\sqrt{10+2\sqrt{5}}}{4} \qquad s_2 = \frac{\sqrt{10-2\sqrt{5}}}{4} \\
m = & \sqrt{1-a^2} \qquad \qquad n = \sqrt{1-b^2}
\end{align*}

then the parameterization of 27 points is successfully accomplished:

\begin{longtable}[c]{r|ccc}
	\caption{Parameters for 27 points} \\
	pt & $x$ & $y$ & $z$ \\
	\hline\vspace*{-2.2ex}
	\endfirsthead
	\multicolumn{4}{c}%
	{\tablename\ \thetable\ -- 27 points -- \textit{cont.}} \\
	pt & $x$ & $y$ & $z$ \\
	\hline\vspace*{-2.2ex}
	\endhead
	1 & $0$ & $0$ & $1$ \\[0.5ex]
	2 & $m$ & $0$ & $a$ \\[0.5ex]
	3 & $c_1m$ & $s_1m$ & $a$ \\[0.5ex]
	4 & $-c_2m$ & $s_2m$ & $a$ \\[0.5ex]
	5 & $-c_2m$ & $-s_2m$ & $a$ \\[0.5ex]
	6 & $c_1m$ & $-s_1m$ & $a$ \\[0.5ex]
	7 & $-n$ & $0$ & $b$ \\[0.5ex]
	8 & $-c_1n$ & $-s_1n$ & $b$ \\[0.5ex]
	9 & $c_2n$ & $-s_2n$ & $b$ \\[0.5ex]
	10 & $c_2n$ & $s_2n$ & $b$ \\[0.5ex]
	11 & $-c_1n$ & $s_1n$ & $b$ \\[0.5ex]
	12 & $1$ & $0$ & $0$ \\[0.5ex]
	13 & $c_1$ & $s_1$ & $0$ \\[0.5ex]
	14 & $-c_2$ & $s_2$ & $0$ \\[0.5ex]
	15 & $-c_2$ & $-s_2$ & $0$ \\[0.5ex]
	16 & $c_1$ & $-s_1$ & $0$ \\[0.5ex]
	17 & $-n$ & $0$ & $-b$ \\[0.5ex]
	18 & $-c_1n$ & $-s_1n$ & $-b$ \\[0.5ex]
	19 & $c_2n$ & $-s_2n$ & $-b$ \\[0.5ex]
	20 & $c_2n$ & $s_2n$ & $-b$ \\[0.5ex]
	21 & $-c_1n$ & $s_1n$ & $-b$ \\[0.5ex]
	22 & $m$ & $0$ & $-a$ \\[0.5ex]
	23 & $c_1m$ & $s_1m$ & $-a$ \\[0.5ex]
	24 & $-c_2m$ & $s_2m$ & $-a$ \\[0.5ex]
	25 & $-c_2m$ & $-s_2m$ & $-a$ \\[0.5ex]
	26 & $c_1m$ & $-s_1m$ & $-a$ \\[0.5ex]
	27 & $0$ & $0$ & $-1$ \\[0.5ex]
\end{longtable}

\noindent
\textit{\textbf{Spherical codes -}}

The algebraic number was successfully recovered for the {\it logarithmic} potential. Unfortunately a search for the algebraic polynomials for $a$ and $b$ for the {\it Coulomb} or {\it Inverse square} potentials has failed to uncover any with degree $<360$ even though the spherical codes are known to 50,014 digits.

\noindent
\textit{\textbf{Logarithmic potential -}}

\noindent
Let the following polynomial be defined:
\begin{align*}
f(x) = & \quad 1216317681592371164083264000000x^{140} + 183491886158513882495311143040000x^{138} \\
 & + 13266184192578719000088327853017600x^{136} + 611762340031238753408375379716968320x^{134} \\
 & + 20201573276377960197232174927951327696x^{132} \\
 & + 508353855838310930890336987690699932740x^{130} \\
 & + 10127332828682321341107824026483487761120x^{128} \\
 & + 163786222858509079620015377666032552374745x^{126} \\
 & + 2187432611814481471102174015520709971340725x^{124} \\
 & + 24407581558045048436881425026354620909720525x^{122} \\
 & + 229272412585134595674886547591938043326494025x^{120} \\
 & + 1820938484893431876744710991658567249750376625x^{118} \\
 & + 12244835621182215840770032562854184289078310125x^{116} \\
 & + 69602799529249419835432371939603727135740379125x^{114} \\
 & + 332753300752855480356375966235683293774425060625x^{112} \\
 & + 1325100965006541694499567599022061306975075460625x^{110} \\
 & + 4322096782936211238772813545214881816648977498125x^{108} \\
 & + 11202058862824085405799672389489230227577511053125x^{106} \\
 & + 21674706885351371743421252214906730124373489790625x^{104} \\
 & + 26211143785941643622578669540447994973164956015625x^{102} \\
 & + 1838885029290492457203494075160532587521168953125x^{100} \\
 & - 66686794117260215673929263194680026847805110171875x^{98} \\
 & - 132531054679314755872135450541220685808481887984375x^{96} \\
 & - 75433370596586409657703414134832648760026569921875x^{94} \\
 & + 146399199122142718186911729756393179707018534765625x^{92} \\
 & + 300308440042103407071017086560910129969471386640625x^{90} \\
 & + 78656966610742652309399575542921847102732973828125x^{88} \\
 & - 327283816765513265019942969882270930401144873046875x^{86} \\
 & - 302475277152452143843445538855397932006814780859375x^{84} \\
 & + 164422817079319302304675011749617840358935634765625x^{82} \\
 & + 319296354219545048743140058053743513996761939453125x^{80} \\
 & - 23632805922699457178896077047015540322188591796875x^{78} \\
 & - 203122992361078029592371918490161017461587099609375x^{76} \\
 & - 7263564453786815622851688344978213997208115234375x^{74} \\
 & + 91460763616108429681346043043185102836924814453125x^{72} \\
 & - 8796233145525578575542719936660712504981982421875x^{70} \\
 & - 29949858575466304731146311719218672962180615234375x^{68} \\
 & + 14272825484594435996976907220005401710124072265625x^{66} \\
 & + 6284403663119039717575910679099902379655517578125x^{64} \\
 & - 8033869182023186100769828097704587979626220703125x^{62} \\
 & - 202810736109082381235495363009932498261962890625x^{60} \\
 & + 2474562728919976271999614644088630834818115234375x^{58} \\
 & - 453996948315498890152072177159472477325439453125x^{56} \\
 & - 420587072016337581007972361541790783468017578125x^{54} \\
 & + 217986131173780918244923793256115400091552734375x^{52} \\
 & + 21194259600009861439220138770294016253662109375x^{50} \\
 & - 62379087403209195452479582153446949053955078125x^{48} \\
 & + 7048575661875721592884782523980382781982421875x^{46} \\
 & + 12752340697779295911795027057829865325927734375x^{44} \\
 & - 1803929464584725776386519939595644195556640625x^{42} \\
 & - 1910234029290474519011255422721430206298828125x^{40} \\
 & + 189548963796025381200024614427562713623046875x^{38} \\
 & + 205934103851505136470360192012946624755859375x^{36} \\
 & - 7469727589898383872940354909903717041015625x^{34} \\
 & - 15368380455734066708067584954566192626953125x^{32} \\
 & - 398604557454446026883205555271148681640625x^{30} \\
 & + 748229752710951140728646615764617919921875x^{28} \\
 & + 61306377944844628929407601757049560546875x^{26} \\
 & - 21566591372675028417421342563629150390625x^{24} \\
 & - 2974856922287427599845699405670166015625x^{22} \\
 & + 295211366268221293687586498260498046875x^{20} \\
 & + 68500283682788640215350627899169921875x^{18} - 72850670006427596361637115478515625x^{16} \\
 & - 700617088321181672937870025634765625x^{14} - 37352991647147328751087188720703125x^{12} \\
 & + 1774553986930528628826141357421875x^{10} + 234814160952770841121673583984375x^{8} \\
 & + 8729774004736840724945068359375x^{6} + 143172778267085552215576171875x^{4} \\
 & + 957711039483547210693359375x^{2} + 1319848001003265380859375
\end{align*}
then let:
\begin{longtable}[l]{l}
	$ \qquad a = 0.7538089984441335383\ldots$ \text{a root of }f(x) \\
	$ \qquad b = 0.3604942753234939635\ldots$ \text{a root of }f(x)\\
	$\text{and for the energy:}$ \\
	$ \qquad energy = -93.25198640000452027\ldots$
\end{longtable}

\noindent
\textit{\textbf{Coulomb $1/r$ potential -}}

\begin{longtable}[l]{l}
	$ \qquad a = 0.7538564449703482744\ldots$ \\
	$ \qquad b = 0.3589242703564896574\ldots$ \\
	$\text{and for the energy:}$ \\
	$ \qquad energy = 287.3026150330391631\ldots$
\end{longtable}
Algebraic polynomial $>360$ degree.

\noindent
\textit{\textbf{Inverse square $1/r^2$ potential -}}

\begin{longtable}[l]{l}
	$ \qquad a = 0.7539171374221273508\ldots$ \\
	$ \qquad b = 0.3574199262261141346\ldots$ \\
	$\text{and for the energy:}$ \\
	$ \qquad energy = 270.7984042081812851\ldots$
\end{longtable}
Algebraic polynomial $>360$ degree.

\noindent
\textit{\textbf{Symmetries -}}

The symmetry groups are identical for 27 points under all 3 potentials:

\begin{center}
	\begin{tabular}{l|l}
		\multicolumn{2}{c}{Symmetries - 27 points} \\
		\hline\Tstrut
		planes & [[4, 200], [16, 5], [35, 5], [50, 1]] \\[0.2ex]
		\hline\Tstrut
		Gram groups & [[2, 1], [10, 4], [20, 17], [27, 1], [40, 8]] \\
		\hline\Tstrut
		Polygons & [[4, 220], [5, 5], [7, 5]]
	\end{tabular}
\end{center}

\subsection{28 points}
The constrained configuration for 28 points is interesting, it consists of 1 pole and 9 triangles in a 1:3:3:3:3:3:3:3:3:3 arrangement. It requires 16 parameters $a$ - $o$ to constrain the arrangement.

\begin{figure}[ht]
	\begin{center}
		\includegraphics[type=pdf,ext=pdf,read=pdf,height=1in,width=1in,angle=0]{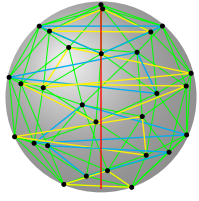}
		\caption{28 points.}
		\label{fig:28pts}
	\end{center}
\end{figure}
In this figure \ref{fig:28pts}, the embedded 9 co-planar triangles are shown in alternating yellow and cyan colors.

\noindent
\textit{\textbf{Constrained parameterization -}}

\begin{longtable}[c]{r|ccc}
	\caption{Parameterization for 28 points} \\
	pt & $x$ & $y$ & $z$ \\
	\hline\vspace*{-2.2ex}
	\endfirsthead
	\multicolumn{4}{c}%
	{\tablename\ \thetable\ -- 28 points parameterization -- \textit{continued\ldots}} \\
	pt & $x$ & $y$ & $z$ \\
	\hline\vspace*{-2.2ex}
	\endhead
	1 & $0$ & $0$ & $1$ \\[0.7ex]
	2 & $b$ & $-\sqrt{1-a^2-b^2}$ & $a$ \\[0.7ex]
	3 & $\frac{-\sqrt{3}\sqrt{1-a^2-b^2}-b}{2}$ & $-\frac{\sqrt{3}b-\sqrt{1-a^2-b^2}}{2}$ & $a$ \\[0.7ex]
	4 & $\frac{\sqrt{3}\sqrt{1-a^2-b^2}-b}{2}$ & $-\frac{-\sqrt{1-a^2-b^2}-\sqrt{3}b}{2}$ & $a$ \\[0.7ex]
	5 & $d$ & $\sqrt{1-c^2-d^2}$ & $c$ \\[0.7ex]
	6 & $\frac{-\sqrt{3}\sqrt{1-c^2-d^2}-d}{2}$ & $\frac{\sqrt{3}d-\sqrt{1-c^2-d^2}}{2}$ & $c$ \\[0.7ex]
	7 & $\frac{\sqrt{3}\sqrt{1-c^2-d^2}-d}{2}$ & $\frac{-\sqrt{1-c^2-d^2}-\sqrt{3}d}{2}$ & $c$ \\[0.7ex]
	8 & $f$ & $\sqrt{1-e^2-f^2}$ & $e$ \\[0.7ex]
	9 & $\frac{-\sqrt{3}\sqrt{1-e^2-f^2}-f}{2}$ & $\frac{\sqrt{3}f-\sqrt{1-e^2-f^2}}{2}$ & $e$ \\[0.7ex]
	10 & $\frac{\sqrt{3}\sqrt{1-e^2-f^2}-f}{2}$ & $\frac{-\sqrt{1-e^2-f^2}-\sqrt{3}f}{2}$ & $e$ \\[0.7ex]
	11 & $h$ & $-\sqrt{1-g^2-h^2}$ & $g$ \\[0.7ex]
	12 & $\frac{-\sqrt{3}\sqrt{1-g^2-h^2}-h}{2}$ & $-\frac{\sqrt{3}h-\sqrt{1-g^2-h^2}}{2}$ & $g$ \\[0.7ex]
	13 & $\frac{\sqrt{3}\sqrt{1-g^2-h^2}-h}{2}$ & $-\frac{-\sqrt{1-g^2-h^2}-\sqrt{3}h}{2}$ & $g$ \\[0.7ex]
	14 & $j$ & $\sqrt{1-i^2-j^2}$ & $i$ \\[0.7ex]
	15 & $\frac{-\sqrt{3}\sqrt{1-i^2-j^2}-j}{2}$ & $\frac{\sqrt{3}j-\sqrt{1-i^2-j^2}}{2}$ & $i$ \\[0.7ex]
	16 & $\frac{\sqrt{3}\sqrt{1-i^2-j^2}-j}{2}$ & $\frac{-\sqrt{1-i^2-j^2}-\sqrt{3}j}{2}$ & $i$ \\[0.7ex]
	17 & $\frac{2}{3}\sqrt{2}$ & $0$ & $-\frac{1}{3}$ \\[0.7ex]
	18 & $-\frac{\sqrt{2}}{3}$ & $\sqrt{\frac{2}{3}}$ & $-\frac{1}{3}$ \\[0.7ex]
	19 & $-\frac{\sqrt{2}}{3}$ & $-\sqrt{\frac{2}{3}}$ & $-\frac{1}{3}$ \\[0.7ex]
	20 & $l$ & $-\sqrt{1-k^2-l^2}$ & $k$ \\[0.7ex]
	21 & $\frac{-\sqrt{3}\sqrt{1-k^2-l^2}-l}{2}$ & $-\frac{\sqrt{3}l-\sqrt{1-k^2-l^2}}{2}$ & $k$ \\[0.7ex]
	22 & $\frac{\sqrt{3}\sqrt{1-k^2-l^2}-l}{2}$ & 	$-\frac{-\sqrt{1-k^2-l^2}-\sqrt{3}l}{2}$ & $k$ \\[0.7ex]
	23 & $n$ & $\sqrt{1-m^2-n^2}$ & $m$ \\[0.7ex]
	24 & $\frac{-\sqrt{3}\sqrt{1-m^2-n^2}-n}{2}$ & $\frac{\sqrt{3}n-\sqrt{1-m^2-n^2}}{2}$ & $m$ \\[0.7ex]
	25 & $\frac{\sqrt{3}\sqrt{1-m^2-n^2}-n}{2}$ & $\frac{-\sqrt{1-m^2-n^2}-\sqrt{3}n}{2}$ & $m$ \\[0.7ex]
	26 & $p$ & $-\sqrt{1-o^2-p^2}$ & $o$ \\[0.7ex]
	27 & $\frac{-\sqrt{3}\sqrt{1-o^2-p^2}-p}{2}$ & $-\frac{\sqrt{3}p-\sqrt{1-o^2-p^2}}{2}$ & $o$ \\[0.7ex]
	28 & $\frac{\sqrt{3}\sqrt{1-o^2-p^2}-p}{2}$ & $-\frac{-\sqrt{1-o^2-p^2}-\sqrt{3}p}{2}$ & $o$
\end{longtable}

\noindent
\textit{\textbf{Parameterization values --}}

\begin{longtable}[c]{c|c|c|c}
	\caption{Parameter values for 28 points} \\
	Parameter & log & 1/r & $1/r^2$ \\
	\hline\vspace*{-2.2ex}
	\endfirsthead
	\multicolumn{4}{c}%
	{\tablename\ \thetable\ -- 28 points parameters -- \textit{continued}} \\
	Parameter & log & 1/r & $1/r^2$ \\
	\hline\vspace*{-2.2ex}
	\endhead
	$a$ & 0.7907708157515443835 & 0.7899009401486432163 & 0.7890262437836545087 \\
	$b$ & 0.5947229186382124928 & 0.5957394888529532691 & 0.5967238615644449130 \\
	$c$ & 0.6679507326864148493 & 0.6688958370340451721 & 0.6695846149511056347 \\
	$d$ & 0.4959390043217711336 & 0.4962265397070834826 & 0.4966443708261904194 \\
	$e$ & 0.2971198730109991089 & 0.2983682631351601216 & 0.2995879040828373313 \\
	$f$ & 0.9437868477891841413 & 0.9433055779439641978 & 0.9428090305395777610 \\
	$g$ & 0.2449255331189668743 & 0.2448815892974706984 & 0.2450459316153989162 \\
	$h$ & 0.7950629915487061674 & 0.7960498896013280029 & 0.7968385526888769798 \\
	$i$ & -0.003384550570525765522 & -0.004976051626784320694 & -0.006389697496857813468 \\
	$j$ & 0.7072721195417864797 & 0.7077118734869209250 & 0.7079426335002857609 \\
	$k$ & -0.4256627251957432557 & -0.4253884685862338868 & -0.4249805760090653975 \\
	$l$ & 0.6882446095313745062 & 0.6874189314791589315 & 0.6866353876860705716 \\
	$m$ & -0.6622279635668002367 & -0.6628807346975694511 & -0.6636335718574264425 \\
	$n$ & 0.6046061674087812837 & 0.6034527355558161951 & 0.6022588118286400375 \\
	$o$ & -0.9094917152348559581 & -0.9088013747047315498 & -0.9082408490696467374 \\
	$p$ & 0.3869148592338547215 & 0.3881613660086515981 & 0.3890901010512777794 \\
	\hline\Tstrut
	$energy$ & -99.65860938412481523 & 310.4915423582018490 & 294.8784716100522203
\end{longtable}
After checking the 16 parameters to 50,014 digits, the algebraic degree for the minimal polynomials is $>360$.

\noindent
\textit{\textbf{Symmetries -}}

The symmetry groups are identical for 28 points under all 3 potentials:

\begin{center}
	\begin{tabular}{l|l}
		\multicolumn{2}{c}{Symmetries - 28 points} \\
		\hline\Tstrut
		planes & [[9, 4]] \\[0.2ex]
		\hline\Tstrut
		Gram groups & [[12, 7], [24, 28], [28, 1]] \\
		\hline\Tstrut
		Polygons & [[3, 36]]
	\end{tabular}
\end{center}

\subsection{29 points}
Although the \textit{logarithmic} potential for 29 points has no embedded polygons, nevertheless the other two potentials have a configuration which admits parameterization. These configurations for the \textit{Coulomb $1/r$} and \textit{Inverse Square $1/r^2$} consist of a pole and antipole and 9 co-planar triangles whose centroids are aligned with the poles. This is an 1:3:3:3:3:3:3:3:3:3:1 arrangement. The \textit{logarithmic} potential would require 58 parameters to constrain, using a direct method on the points.

\begin{figure}[ht]
	\begin{center}
		\includegraphics[type=pdf,ext=pdf,read=pdf,height=1in,width=1in,angle=0]{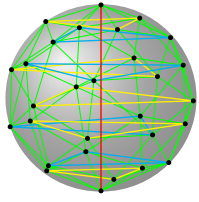}
		\caption{29 points.}
		\label{fig:29pts}
	\end{center}
\end{figure}
\noindent

\noindent
\textit{\textbf{Parameterization for the Coulomb and Inverse Square potentials}}

\begin{longtable}[c]{r|ccc}
	\caption{Parameterization for 29 points} \\
	pt & $x$ & $y$ & $z$ \\
	\hline\vspace*{-2.2ex}
	\endfirsthead
	\multicolumn{4}{c}%
	{\tablename\ \thetable\ -- 29 points parameterization -- \textit{continued\ldots}} \\
	pt & $x$ & $y$ & $z$ \\
	\hline\vspace*{-2.2ex}
	\endhead
	1 & $0$ & $0$ & $1$ \\[0.7ex]
	2 & $b$ & $\sqrt{1-a^2-b^2}$ & $a$ \\[0.7ex]
	3 & $\frac{-\sqrt{3}\sqrt{1-a^2-b^2}-b}{2}$ & $\frac{b\sqrt{3}-\sqrt{1-a^2-b^2}}{2}$ & $a$ \\[0.7ex]
	4 & $\frac{\sqrt{3}\sqrt{1-a^2-b^2}-b}{2}$ & $\frac{-\sqrt{1-a^2-b^2}-b\sqrt{3}}{2}$ & $a$ \\[0.7ex]
	5 & $d$ & $\sqrt{1-c^2-d^2}$ & $c$ \\[0.7ex]
	6 & $\frac{-\sqrt{3}\sqrt{1-c^2-d^2}-d}{2}$ & $\frac{d\sqrt{3}-\sqrt{1-c^2-d^2}}{2}$ & $c$ \\[0.7ex]
	7 & $\frac{\sqrt{3}\sqrt{1-c^2-d^2}-d}{2}$ & $\frac{-\sqrt{1-c^2-d^2}-d\sqrt{3}}{2}$ & $c$ \\[0.7ex]
	8 & $f$ & $-\sqrt{1-e^2-f^2}$ & $e$ \\[0.7ex]
	9 & $\frac{-\sqrt{3}\sqrt{1-e^2-f^2}-f}{2}$ & $-\frac{f\sqrt{3}-\sqrt{1-e^2-f^2}}{2}$ & $e$ \\[0.7ex]
	10 & $\frac{\sqrt{3}\sqrt{1-e^2-f^2}-f}{2}$ & $-\frac{-\sqrt{1-e^2-f^2}-f\sqrt{3}}{2}$ & $e$ \\[0.7ex]
	11 & $h$ & $\sqrt{1-g^2-h^2}$ & $g$ \\[0.7ex]
	12 & $\frac{-\sqrt{3}\sqrt{1-g^2-h^2}-h}{2}$ & $\frac{h\sqrt{3}-\sqrt{1-g^2-h^2}}{2}$ & $g$ \\[0.7ex]
	13 & $\frac{\sqrt{3}\sqrt{1-g^2-h^2}-h}{2}$ & $\frac{-\sqrt{1-g^2-h^2}-h\sqrt{3}}{2}$ & $g$ \\[0.7ex]
	14 & $1$ & $0$ & $0$ \\[0.7ex]
	15 & $-\frac{1}{2}$ & $\frac{\sqrt{3}}{2}$ & $0$ \\[0.7ex]
	16 & $-\frac{1}{2}$ & $-\frac{\sqrt{3}}{2}$ & $0$ \\[0.7ex]
	17 & $h$ & $-\sqrt{1-g^2-h^2}$ & $-g$ \\[0.7ex]
	18 & $\frac{-\sqrt{3}\sqrt{1-g^2-h^2}-h}{2}$ & $-\frac{h\sqrt{3}-\sqrt{1-g^2-h^2}}{2}$ & $-g$ \\[0.7ex]
	19 & $\frac{\sqrt{3}\sqrt{1-g^2-h^2}-h}{2}$ & $-\frac{-\sqrt{1-g^2-h^2}-h\sqrt{3}}{2}$ & $-g$ \\[0.7ex]
	20 & $f$ & $\sqrt{1-e^2-f^2}$ & $-e$ \\[0.7ex]
	21 & $\frac{-\sqrt{3}\sqrt{1-e^2-f^2}-f}{2}$ & $\frac{f\sqrt{3}-\sqrt{1-e^2-f^2}}{2}$ & $-e$ \\[0.7ex]
	22 & $\frac{\sqrt{3}\sqrt{1-e^2-f^2}-f}{2}$ & $\frac{-\sqrt{1-e^2-f^2}-f\sqrt{3}}{2}$ & $-e$ \\[0.7ex]
	23 & $d$ & $-\sqrt{1-c^2-d^2}$ & $-c$ \\[0.7ex]
	24 & $\frac{-\sqrt{3}\sqrt{1-c^2-d^2}-d}{2}$ & $-\frac{d\sqrt{3}-\sqrt{1-c^2-d^2}}{2}$ & $-c$ \\[0.7ex]
	25 & $\frac{\sqrt{3}\sqrt{1-c^2-d^2}-d}{2}$ & $-\frac{-\sqrt{1-c^2-d^2}-d\sqrt{3}}{2}$ & $-c$ \\[0.7ex]
	26 & $b$ & $-\sqrt{1-a^2-b^2}$ & $-a$ \\[0.7ex]
	27 & $\frac{-\sqrt{3}\sqrt{1-a^2-b^2}-b}{2}$ & $-\frac{b\sqrt{3}-\sqrt{1-a^2-b^2}}{2}$ & $-a$ \\[0.7ex]
	28 & $\frac{\sqrt{3}\sqrt{1-a^2-b^2}-b}{2}$ & $-\frac{-\sqrt{1-a^2-b^2}-b\sqrt{3}}{2}$ & $-a$ \\[0.7ex]
	29 & $0$ & $0$ & $-1$
\end{longtable}
\vspace{2.0em}

\noindent
\textit{\textbf{Parameterization values for the Coulomb and Inverse Square potentials--}}

It was found that only 8 parameters are needed to constrain the 29 points

\begin{longtable}[c]{c|c|c|c}
	\caption{Parameter values for 29 points} \\
	Parameter & log & 1/r & $1/r^2$ \\
	\hline\vspace*{-2.2ex}
	\endfirsthead
	\multicolumn{4}{c}%
	{\tablename\ \thetable\ -- 29 points parameters -- \textit{continued}} \\
	Parameter & log & 1/r & $1/r^2$ \\
	\hline\vspace*{-2.2ex}
	\endhead
	$a$ & - & 0.8049835831969647230 & 0.8039411075798132358 \\
	$b$ & - & 0.2794055079959605373 & 0.2808816029091532881 \\
	$c$ & - & 0.6677576661310917280 & 0.6703429832839722148 \\
	$d$ & - & 0.7435575225885066573 & 0.7413020298254826231 \\
	$e$ & - & 0.3208350404562844680 & 0.3222927507118441673 \\
	$f$ & - & 0.7622654414203364054 & -0.8675984615141489508 \\
	$g$ & - & 0.3015513170785488457 & 0.3029954853682648017 \\
	$h$ & - & 0.7247225871913261335 & 0.7246961356257399080 \\
	\hline\Tstrut
	$energy$ & -106.2545711708346378 & 334.6344399204156906 & 320.2160317560987546
\end{longtable}
After checking the 16 parameters to 50,014 digits, the algebraic degree for the minimal polynomials is $>360$.

The coordinates for 29 points are known to 77 digits for the \textit{logarithmic} potential.

\noindent
\textit{\textbf{Symmetries -}}

NOTE: The symmetries for the \textit{logarithmic} potential does not permit a parameterization.

\begin{center}
	\begin{tabular}{l|l}
		\multicolumn{2}{c}{Symmetries - 29 points - \textit{logarithmic}} \\
		\hline\Tstrut
		planes & [] \\[0.2ex]
		\hline\Tstrut
		Gram groups & [[2, 14], [4, 196], [29, 1]] \\
		\hline\Tstrut
		Polygons & []
	\end{tabular}
\end{center}

The symmetry groups for the \textit{Coulomb $1/r$} and \textit{Inverse square law $1/r^2$} potentials are identical.

\begin{center}
	\begin{tabular}{l|l}
		\multicolumn{2}{c}{Symmetries - 29 points - \textit{Coulomb} or \textit{Inverse sq}} \\
		\hline\Tstrut
		planes & [[9, 1]] \\[0.2ex]
		\hline\Tstrut
		Gram groups & [[2, 1], [6, 13], [12, 61], [29, 1]] \\
		\hline\Tstrut
		Polygons & [[3, 9]]
	\end{tabular}
\end{center}

\subsection{30 points}
The configuration of 30 points has a symmetric arrangement, and one quadrilateral appears at the equator, as shown by the symmetry groups, but it requires 20 parameters to constrain. The configuration is the same for all 3 potentials, although distances slightly change in each case.

\begin{figure}[ht]
	\begin{center}
		\includegraphics[type=pdf,ext=pdf,read=pdf,height=1in,width=1in,angle=0]{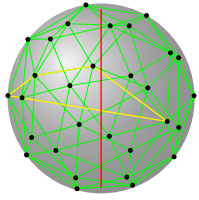}
		\caption{30 points.}
		\label{fig:30pts}
	\end{center}
\end{figure}
\noindent
The kite quadrilateral embedded in the configuration is outlined in yellow. This corresponds to the [4,14] polygon group shown in the symmetry groups. The other 26 points had to be individually constrained.

\noindent
\textit{\textbf{Parameterization for all 3 potentials}}

\begin{longtable}[c]{r|ccc}
	\caption{Parameterization for 30 points} \\
	pt & $x$ & $y$ & $z$ \\
	\hline\vspace*{-2.2ex}
	\endfirsthead
	\multicolumn{4}{c}%
	{\tablename\ \thetable\ -- 30 points parameterization -- \textit{continued\ldots}} \\
	pt & $x$ & $y$ & $z$ \\
	\hline\vspace*{-2.2ex}
	\endhead
	1 & $-b$ & $\sqrt{1-a^2-b^2}$ & $a$ \\[0.5em]
	2 & $b$ & $-\sqrt{1-b^2-c^2}$ & $c$ \\[0.5em]
	3 & $e$ & $\sqrt{1-d^2-e^2}$ & $d$ \\[0.5em]
	4 & $-e$ & $-\sqrt{1-e^2-f^2}$ & $f$ \\[0.5em]
	5 & $h$ & $\sqrt{1-g^2-h^2}$ & $g$ \\[0.5em]
	6 & $-h$ & $-\sqrt{1-h^2-i^2}$ & $i$ \\[0.5em]
	7 & $-k$ & $\sqrt{1-j^2-k^2}$ & $j$ \\[0.5em]
	8 & $k$ & $-\sqrt{1-k^2-l^2}$ & $l$ \\[0.5em]
	9 & $-n$ & $\sqrt{1-m^2-n^2}$ & $m$ \\[0.5em]
	10 & $n$ & $-\sqrt{1-n^2-o^2}$ & $o$ \\[0.5em]
	11 & $q$ & $\sqrt{1-p^2-q^2}$ & $p$ \\[0.5em]
	12 & $-q$ & $-\sqrt{1-q^2-r^2}$ & $r$ \\[0.5em]
	13 & $t$ & $\sqrt{1-s^2-t^2}$ & $s$ \\[0.5em]
	14 & $1$ & $0$ & $0$ \\[0.5em]
	15 & $-t$ & $\sqrt{1-t^2}$ & $0$ \\[0.5em]
	16 & $-t$ & $-\sqrt{1-t^2}$ & $0$ \\[0.5em]
	17 & $-1$ & $0$ & $0$ \\[0.5em]
	18 & $t$ & $-\sqrt{1-s^2-t^2}$ & $-s$ \\[0.5em]
	19 & $-q$ & $\sqrt{1-q^2-r^2}$ & $-r$ \\[0.5em]
	20 & $q$ & $-\sqrt{1-p^2-q^2}$ & $-p$ \\[0.5em]
	21 & $n$ & $\sqrt{1-n^2-o^2}$ & $-o$ \\[0.5em]
	22 & $-n$ & $-\sqrt{1-m^2-n^2}$ & $-m$ \\[0.5em]
	23 & $k$ & $\sqrt{1-k^2-l^2}$ & $-l$ \\[0.5em]
	24 & $-k$ & $-\sqrt{1-j^2-k^2}$ & $-j$ \\[0.5em]
	25 & $-h$ & $\sqrt{1-h^2-i^2}$ & $-i$ \\[0.5em]
	26 & $h$ & $-\sqrt{1-g^2-h^2}$ & $-g$ \\[0.5em]
	27 & $-e$ & $\sqrt{1-e^2-f^2}$ & $-f$ \\[0.5em]
	28 & $e$ & $-\sqrt{1-d^2-e^2}$ & $-d$ \\[0.5em]
	29 & $b$ & $\sqrt{1-b^2-c^2}$ & $-c$ \\[0.5em]
	30 & $-b$ & $-\sqrt{1-a^2-b^2}$ & $-a$
\end{longtable}

\noindent
\textit{\textbf{Parameterization values --}}

\begin{longtable}[c]{c|c|c|c}
	\caption{Parameter values for 30 points} \\
	Parameter & log & 1/r & $1/r^2$ \\
	\hline\vspace*{-2.2ex}
	\endfirsthead
	\multicolumn{4}{c}%
	{\tablename\ \thetable\ -- 30 points parameters -- \textit{continued}} \\
	Parameter & log & 1/r & $1/r^2$ \\
	\hline\vspace*{-2.2ex}
	\endhead
	$a$ & 0.9486527267413356788 & 0.9479246620637078755 & 0.9473564675102344637 \\
	$b$ & 0.3066734200644118059 & 0.3086684605898085157 & 0.3103675967559615939 \\
	$c$ & 0.9366966646890613759 & 0.9375329394737162094 & 0.9381957682673322194 \\
	$d$ & 0.8348763939425380052 & 0.8330520838115783220 & 0.8316900668842572719 \\
	$e$ & 0.1615635634965170551 & 0.1587264754655424326 & 0.1560736352850945125 \\
	$f$ & 0.7800307273926152722 & 0.7839095808737609976 & 0.7870515390893629380 \\
	$g$ & 0.6114046189205064800 & 0.6088759699493982358 & 0.6070039615681473432 \\
	$h$ & 0.7189758666937269015 & 0.7197644243912286318 & 0.7205727021183832515 \\
	$i$ & 0.5765436224713525535 & 0.5776317653443249921 & 0.5786664974548113327 \\
	$j$ & 0.5675096798101771037 & 0.5656310135090343176 & 0.5638093466300415529 \\
	$k$ & 0.7766621301043334901 & 0.7769588815656016999 & 0.7773464422959958165 \\
	$l$ & 0.5383884572786199876 & 0.5395046991159762670 & 0.5400351360335616450 \\
	$m$ & 0.5238322117415161400 & 0.5178503474943680813 & 0.5124590167803016716 \\
	$n$ & 0.3512803016237076140 & 0.3525736975276157293 & 0.3539157941793877699 \\
	$o$ & 0.4462680632513593951 & 0.4482405653180075719 & 0.4491261124958870240 \\
	$p$ & 0.2806873507348947443 & 0.2770686982856629910 & 0.2742751919187992182 \\
	$q$ & 0.2445616358795041037 & 0.2462481850409242538 & 0.2478166943029047182 \\
	$r$ & 0.1895437965497211218 & 0.1954064901966944200 & 0.2001027263629027045 \\
	$s$ & 0.06227250397460793596 & 0.05565485208568008268 & 	0.05040789193886369704 \\
	$t$ & 0.7655086170510579665 & 0.7673872585373448361 & 	0.7693445196233009931 \\
	\hline\Tstrut
	$energy$ & -113.0892554965139952 & 359.6039459037634590 & 	346.2636306407908126
\end{longtable}
After checking the 20 parameters to 50,014 digits, the algebraic degree for the minimal polynomials is $>420$.

\noindent
\textit{\textbf{Symmetries -}}

The symmetry groups for 30 points are identical for all 3 potentials.

\begin{center}
	\begin{tabular}{l|l}
		\multicolumn{2}{c}{Symmetries - 30 points} \\
		\hline\Tstrut
		planes & [[4, 14]] \\[0.2ex]
		\hline\Tstrut
		Gram groups & [[2, 1], [4, 21], [8, 98], [30, 1]] \\
		\hline\Tstrut
		Polygons & [[4, 14]]
	\end{tabular}
\end{center}

\subsection{31 points}
The configuration of 31 points is similar to 28 points, where one point is at a pole and the others are on co-planar triangles whose centroids are aligned with the pole. In the case of 31 points, there are 10 co-planar triangles. This is an 1:3:3:3:3:3:3:3:3:3:3 arrangement. This alignment was strongly suggested by the [46,1] group in the planes part of the symmetry groups.

\begin{figure}[ht]
	\begin{center}
		\includegraphics[type=pdf,ext=pdf,read=pdf,height=1in,width=1in,angle=0]{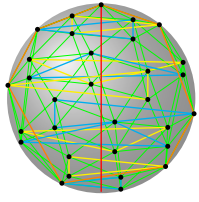}
		\caption{31 points.}
		\label{fig:31pts}
	\end{center}
\end{figure}
While there are 3 embedded heptagons in figure \ref{fig:31pts}, only one is shown in orange. The co-planar triangles are shown in cyan and yellow outlines.

\noindent
\textit{\textbf{Parameterization --}}

\begin{longtable}[c]{r|ccc}
	\caption{Parameterization for 31 points} \\
	pt & $x$ & $y$ & $z$ \\
	\hline\vspace*{-2.2ex}
	\endfirsthead
	\multicolumn{4}{c}%
	{\tablename\ \thetable\ -- 31 points parameterization -- \textit{continued\ldots}} \\
	pt & $x$ & $y$ & $z$ \\
	\hline\vspace*{-2.2ex}
	\endhead
	1 & $0$ & $0$ & $1$ \\[0.7em]
	2 & $b$ & $\sqrt{1-a^2-b^2}$ & $a$ \\[0.7em]
	3 & $\frac{-b-\sqrt{3}\sqrt{1-a^2-b^2}}{2}$ & $\frac{b\sqrt{3}-\sqrt{1-a^2-b^2}}{2}$ & $a$ \\[0.7em]
	4 & $\frac{-b+\sqrt{3}\sqrt{1-a^2-b^2}}{2}$ & $\frac{-b\sqrt{3}-\sqrt{1-a^2-b^2}}{2}$ & $a$ \\[0.7em]
	5 & $d$ & $\sqrt{1-c^2-d^2}$ & $c$ \\[0.7em]
	6 & $\frac{-d-\sqrt{3}\sqrt{1-c^2-d^2}}{2}$ & $\frac{d\sqrt{3}-\sqrt{1-c^2-d^2}}{2}$ & $c$ \\[0.7em]
	7 & $\frac{-d+\sqrt{3}\sqrt{1-c^2-d^2}}{2}$ & 	$\frac{-d\sqrt{3}-\sqrt{1-c^2-d^2}}{2}$ & $c$ \\[0.7em]
	8 & $f$ & $\sqrt{1-e^2-f^2}$ & $e$ \\[0.7em]
	9 & $\frac{-f-\sqrt{3}\sqrt{1-e^2-f^2}}{2}$ & $\frac{f\sqrt{3}-\sqrt{1-e^2-f^2}}{2}$ & $e$ \\[0.7em]
	10 & $\frac{-f+\sqrt{3}\sqrt{1-e^2-f^2}}{2}$ & $\frac{-f\sqrt{3}-\sqrt{1-e^2-f^2}}{2}$ & $e$ \\[0.7em]
	11 & $g$ & $\sqrt{1-e^2-g^2}$ & $e$ \\[0.7em]
	12 & $\frac{-g-\sqrt{3}\sqrt{1-e^2-g^2}}{2}$ & $\frac{g\sqrt{3}-\sqrt{1-e^2-g^2}}{2}$ & $e$ \\[0.7em]
	13 & $\frac{-g+\sqrt{3}\sqrt{1-e^2-g^2}}{2}$ & 	$\frac{-g\sqrt{3}-\sqrt{1-e^2-g^2}}{2}$ & $e$ \\[0.7em]
	14 & $i$ & $\sqrt{1-h^2-i^2}$ & $h$ \\[0.7em]
	15 & $\frac{-i-\sqrt{3}\sqrt{1-h^2-i^2}}{2}$ & $\frac{i\sqrt{3}-\sqrt{1-h^2-i^2}}{2}$ & $h$ \\[0.7em]
	16 & $\frac{-i+\sqrt{3}\sqrt{1-h^2-i^2}}{2}$ & $\frac{-i\sqrt{3}-\sqrt{1-h^2-i^2}}{2}$ & $h$ \\[0.7em]
	17 & $k$ & $\sqrt{1-j^2-k^2}$ & $j$ \\[0.7em]
	18 & $\frac{-k-\sqrt{3}\sqrt{1-j^2-k^2}}{2}$ & $\frac{k\sqrt{3}-\sqrt{1-j^2-k^2}}{2}$ & $j$ \\[0.7em]
	19 & $\frac{-k+\sqrt{3}\sqrt{1-j^2-k^2}}{2}$ & $\frac{-k\sqrt{3}-\sqrt{1-j^2-k^2}}{2}$ & $j$ \\[0.7em]
	20 & $m$ & $\sqrt{1-l^2-m^2}$ & $l$ \\[0.7em]
	21 & $\frac{-m-\sqrt{3}\sqrt{1-l^2-m^2}}{2}$ & $\frac{m\sqrt{3}-\sqrt{1-l^2-m^2}}{2}$ & $l$ \\[0.7em]
	22 & $\frac{-m+\sqrt{3}\sqrt{1-l^2-m^2}}{2}$ & $\frac{-m\sqrt{3}-\sqrt{1-l^2-m^2}}{2}$ & $l$ \\[0.7em]
	23 & $n$ & $\sqrt{1-l^2-n^2}$ & $l$ \\[0.7em]
	24 & $\frac{-n-\sqrt{3}\sqrt{1-l^2-n^2}}{2}$ & $\frac{n\sqrt{3}-\sqrt{1-l^2-n^2}}{2}$ & $l$ \\[0.7em]
	25 & $\frac{-n+\sqrt{3}\sqrt{1-l^2-n^2}}{2}$ & $\frac{-n\sqrt{3}-\sqrt{1-l^2-n^2}}{2}$ & $l$ \\[0.7em]
	26 & $p$ & $\sqrt{1-o^2-p^2}$ & $o$ \\[0.7em]
	27 & $\frac{-p-\sqrt{3}\sqrt{1-o^2-p^2}}{2}$ & $\frac{p\sqrt{3}-\sqrt{1-o^2-p^2}}{2}$ & $o$ \\[0.7em]
	28 & $\frac{-p+\sqrt{3}\sqrt{1-o^2-p^2}}{2}$ & $\frac{-p\sqrt{3}-\sqrt{1-o^2-p^2}}{2}$ & $o$ \\[0.7em]
	29 & $r$ & $\sqrt{1-q^2-r^2}$ & $q$ \\[0.7em]
	30 & $\frac{-r-\sqrt{3}\sqrt{1-q^2-r^2}}{2}$ & $\frac{r\sqrt{3}-\sqrt{1-q^2-r^2}}{2}$ & $q$ \\[0.7em]
	31 & $\frac{-r+\sqrt{3}\sqrt{1-q^2-r^2}}{2}$ & $\frac{-r\sqrt{3}-\sqrt{1-q^2-r^2}}{2}$ & $q$
\end{longtable}

\noindent
\textit{\textbf{Parameterization values --}}

It takes 18 parameters $a$ - $r$ to constrain the arrangement.

\begin{longtable}[c]{c|c|c|c}
	\caption{Parameter values for 31 points} \\
	Parameter & log & 1/r & $1/r^2$ \\
	\hline\vspace*{-2.2ex}
	\endfirsthead
	\multicolumn{4}{c}%
	{\tablename\ \thetable\ -- 30 points parameters -- \textit{continued}} \\
	Parameter & log & 1/r & $1/r^2$ \\
	\hline\vspace*{-2.2ex}
	\endhead
	$a$ & 0.7863298907020499427 & 0.7870999094095871799 & 0.7878563575150593603 \\
	$b$ & 0.6178068492567091888 & 0.6168255284984722928 & 0.6158590422516364810 \\
	$c$ & 0.7367473395695973885 & 0.7350765584881929921 & 0.7337481909540752649 \\
	$d$ & 0.3380840715107413775 & 0.3389920549073216064 & 0.3397107564795461042 \\
	$e$ & 0.3169221170283596425 & 0.3179017267194057892 & 0.3190569256801285921 \\
	$f$ & 0.8694417223274445707 & 0.8697325260339821684 & 0.8698546411206109922 \\
	$g$ & -0.1065081292185521326 & -0.1079444543524014153 & -0.1090950230169857251 \\
	$h$ & 0.1383690098767768550 & 0.1380042000429985443 & 0.1377248927162571981 \\
	$i$ & 0.4951903717525514523 & 0.4952158218318786815 & 0.4952352607413711422 \\
	$j$ & -0.1653379278962295080 & -0.1677988338361777688 & -0.1694832822408188735 \\
	$k$ & 0.9862369743621363602 & 0.9858212573094672332 & 0.9855330623783653306 \\
	$l$ & -0.4125042399687933459 & -0.4116653217317328665 & -0.4108163847952203757 \\
	$m$ & 0.7068688051657705802 & 0.7082299443358992576 & 0.7094355284719730024 \\
	$n$ & 0.7068688051657705802 & 0.7082299443358992576 & 0.7094355284719730024 \\
	$o$ & -0.7296161168785959254 & -0.7287802628567269414 & -0.7278462146559030662 \\
	$p$ & 0.6838569455602531569 & 0.6847476385285895863 & 0.6857403938890234842 \\
	$q$ & -0.9086612828260646794 & -0.9083394771095752304 & -0.9080016083542497285 \\
	$r$ & 0.2087670191223283802 & 0.2091168299804582061 & 0.2094833401646152920 \\
	\hline\Tstrut
	$energy$ & -120.1103466395223015 & 385.5308380632994252 & 373.5808689611434006
\end{longtable}
After checking the 18 parameters to 50,014 digits, the algebraic degree for the minimal polynomials is $>420$.

\noindent
\textit{\textbf{Symmetries -}}

The symmetry groups for 31 points are identical for all 3 potentials.

\begin{center}
	\begin{tabular}{l|l}
		\multicolumn{2}{c}{Symmetries - 31 points} \\
		\hline\Tstrut
		planes & [[4, 180], [35, 3], [46, 1]] \\[0.2ex]
		\hline\Tstrut
		Gram groups & [[6, 33], [12, 61], [31, 1]] \\
		\hline\Tstrut
		Polygons & [[3, 6], [4, 180], [6, 2], [7, 3]]
	\end{tabular}
\end{center}

\subsection{32 points}
The 32 vertices polyhedra has remarkable icosahedral symmetry over all 3 potentials. It was noticed during searching that the convergence was rapid, indicating a high degree of symmetry for the configuration. The arrangement is [2:2:4:2:2:8:2:2:4:2:2] which is balanced.

\begin{figure}[H]
	\begin{center}
		\includegraphics[type=pdf,ext=pdf,read=pdf,height=1in,width=1in,angle=0]{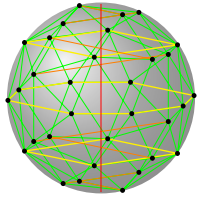}
		\caption{32 points.}
		\label{fig:32pts}
	\end{center}
\end{figure}

The simple algebraic spherical code obtained demonstrates the elegant symmetries found in this unique polyhedron.

\noindent
\textit{\textbf{Algebraic Spherical Code -}}

The algebraic spherical code for 32 points, for all 3 potentials, is:
\begin{longtable}[c]{r|ccc}
	\caption{Parameterization for 32 points} \\
	pt & $x$ & $y$ & $z$ \\
	\hline\vspace*{-2.2ex}
	\endfirsthead
	\multicolumn{4}{c}%
	{\tablename\ \thetable\ -- 32 points -- \textit{continued\ldots}} \\
	pt & $x$ & $y$ & $z$ \\
	\hline\vspace*{-2.2ex}
	\endhead
	1 & $\sqrt{\frac{3-\sqrt{5}}{6}}$ & $0$ & $\sqrt{\frac{3+\sqrt{5}}{6}}$ \\[0.7ex]
	2 & $-\sqrt{\frac{3-\sqrt{5}}{6}}$ & $0$ & $\sqrt{\frac{3+\sqrt{5}}{6}}$ \\[0.7ex]
	3 & $0$ & $\sqrt{\frac{5-\sqrt{5}}{10}}$ & $\sqrt{\frac{5+\sqrt{5}}{10}}$ \\[0.7ex]
	4 & $0$ & $-\sqrt{\frac{5-\sqrt{5}}{10}}$ & $\sqrt{\frac{5+\sqrt{5}}{10}}$ \\[0.7ex]
	5 & $\sqrt{\frac{1}{3}}$ & $\sqrt{\frac{1}{3}}$ & $\sqrt{\frac{1}{3}}$ \\[0.7ex]
	6 & $-\sqrt{\frac{1}{3}}$ & $\sqrt{\frac{1}{3}}$ & $\sqrt{\frac{1}{3}}$ \\[0.7ex]
	7 & $\sqrt{\frac{1}{3}}$ & $-\sqrt{\frac{1}{3}}$ & $\sqrt{\frac{1}{3}}$ \\[0.7ex]
	8 & $-\sqrt{\frac{1}{3}}$ & $-\sqrt{\frac{1}{3}}$ & $\sqrt{\frac{1}{3}}$ \\[0.7ex]
	9 & $\sqrt{\frac{5+\sqrt{5}}{10}}$ & $0$ & $\sqrt{\frac{5-\sqrt{5}}{10}}$ \\[0.7ex]
	10 & $-\sqrt{\frac{5+\sqrt{5}}{10}}$ & $0$ & $\sqrt{\frac{5-\sqrt{5}}{10}}$ \\[0.7ex]
	11 & $0$ & $\sqrt{\frac{3+\sqrt{5}}{6}}$ & $\sqrt{\frac{3-\sqrt{5}}{6}}$ \\[0.7ex]
	12 & $0$ & $-\sqrt{\frac{3+\sqrt{5}}{6}}$ & $\sqrt{\frac{3-\sqrt{5}}{6}}$ \\[0.7ex]
	13 & $-\sqrt{\frac{3+\sqrt{5}}{6}}$ & $\sqrt{\frac{3-\sqrt{5}}{6}}$ & $0$ \\[0.7ex]
	14 & $-\sqrt{\frac{3+\sqrt{5}}{6}}$ & $-\sqrt{\frac{3-\sqrt{5}}{6}}$ & $0$ \\[0.7ex]
	15 & $-\sqrt{\frac{5-\sqrt{5}}{10}}$ & $\sqrt{\frac{5+\sqrt{5}}{10}}$ & $0$ \\[0.7ex]
	16 & $-\sqrt{\frac{5-\sqrt{5}}{10}}$ & $-\sqrt{\frac{5+\sqrt{5}}{10}}$ & $0$ \\[0.7ex]
	17 & $\sqrt{\frac{5-\sqrt{5}}{10}}$ & $\sqrt{\frac{5+\sqrt{5}}{10}}$ & $0$ \\[0.7ex]
	18 & $\sqrt{\frac{5-\sqrt{5}}{10}}$ & $-\sqrt{\frac{5+\sqrt{5}}{10}}$ & $0$ \\[0.7ex]
	19 & $\sqrt{\frac{3+\sqrt{5}}{6}}$ & $\sqrt{\frac{3-\sqrt{5}}{6}}$ & $0$ \\[0.7ex]
	20 & $\sqrt{\frac{3+\sqrt{5}}{6}}$ & $-\sqrt{\frac{3-\sqrt{5}}{6}}$ & $0$ \\[0.7ex]
	21 & $0$ & $\sqrt{\frac{3+\sqrt{5}}{6}}$ & $-\sqrt{\frac{3-\sqrt{5}}{6}}$ \\[0.7ex]
	22 & $0$ & $-\sqrt{\frac{3+\sqrt{5}}{6}}$ & $-\sqrt{\frac{3-\sqrt{5}}{6}}$ \\[0.7ex]
	23 & $\sqrt{\frac{5+\sqrt{5}}{10}}$ & $0$ & $-\sqrt{\frac{5-\sqrt{5}}{10}}$ \\[0.7ex]
	24 & $-\sqrt{\frac{5+\sqrt{5}}{10}}$ & $0$ & $-\sqrt{\frac{5-\sqrt{5}}{10}}$ \\[0.7ex]
	25 & $\sqrt{\frac{1}{3}}$ & $\sqrt{\frac{1}{3}}$ & $-\sqrt{\frac{1}{3}}$ \\[0.7ex]
	26 & $-\sqrt{\frac{1}{3}}$ & $\sqrt{\frac{1}{3}}$ & $-\sqrt{\frac{1}{3}}$ \\[0.7ex]
	27 & $\sqrt{\frac{1}{3}}$ & $-\sqrt{\frac{1}{3}}$ & $-\sqrt{\frac{1}{3}}$ \\[0.7ex]
	28 & $-\sqrt{\frac{1}{3}}$ & $-\sqrt{\frac{1}{3}}$ & $-\sqrt{\frac{1}{3}}$ \\[0.7ex]
	29 & $0$ & $\sqrt{\frac{5-\sqrt{5}}{10}}$ & $-\sqrt{\frac{5+\sqrt{5}}{10}}$ \\[0.7ex]
	30 & $0$ & $-\sqrt{\frac{5-\sqrt{5}}{10}}$ & $-\sqrt{\frac{5+\sqrt{5}}{10}}$ \\[0.7ex]
	31 & $\sqrt{\frac{3-\sqrt{5}}{6}}$ & $0$ & $-\sqrt{\frac{3+\sqrt{5}}{6}}$ \\[0.7ex]
	32 & $-\sqrt{\frac{3-\sqrt{5}}{6}}$ & $0$ & $-\sqrt{\frac{3+\sqrt{5}}{6}}$
\end{longtable}

\noindent
\textit{\textbf{Minimal Energy values -}}

All 3 values of the minimal energy have been obtained as an algebraic number, from the 32 point configuration code listed above. The polynomial for the \textit{Coulomb $1/r$} potential has large coefficients, which was a bit surprising, when first discovered.
\begin{center}
	\begin{longtable}{l|l}
		\multicolumn{2}{c}{Minimal Energy - 32 points} \\
		\hline\mystrut(13,0)
		logarithmic & -127.3788676147802622\ldots = $log(
10515705214263760266382541387079545648/$ \\
		potential & $\qquad 545976643780561909556593187225764978620645706541836261749267578/$ \\
		$\; -log()$ & $\qquad 125/21967352512417951087942082557060458295262192960458577310062/$ \\
		& $\qquad 283069393738186872499366792190850116654575927348196452738745718/$ \\
		& $\qquad 7903615990932602002368905315681828864)$ \\[0.2ex]
		\hline\mystrut(13,0)
		Coulomb & 412.2612746505293101\ldots a root of $x^{32} - 256x^{31} - 321056x^{30}$ \\
		potential & $ + \; 56212480x^{29} + 45980440160x^{28} - 4662927392768x^{27} - 3725768316788352x^{26}$ \\
		$\quad\ (\frac{1}{r})$ & $ + \; 172160485617420288x^{25} + 187812559524684772800x^{24}$ \\
		& $ - \; 1188239661143788646400x^{23} - 6183823866230874346314240x^{22}$ \\
		& $ - \; 142087324208565134319452160x^{21} + 136198061614336090992556362240x^{20}$ \\
		& $ + \; 6203737940397149567528035123200x^{19}$ \\
		& $ - \; 2018494242800279884757142849177600x^{18}$ \\
		& $ - \; 130191651568644593241413290046423040x^{17}$ \\
		& $ + \; 19886979994385494880151097652299875840x^{16}$ \\
		& $ + \; 1651022156000651714852318455909034557440x^{15}$ \\
		& $ - \; 125191437238032671831800022541146854809600x^{14}$ \\
		& $ - \; 13308457372603210878758765690294294917939200x^{13}$ \\
		& $ + \; 453613065533944641684181743285017332487331840x^{12}$ \\
		& $ + \; 68053067238638400943741310984769323538878300160x^{11}$ \\
		& $ - \; 623287949142090854184722342779466792456777072640x^{10}$ \\
		& $ - \; 213183879390214546069562984506723843991442200985600x^{9}$ \\
		& $ - \; 1311815696392437117892800661653830993087821257523200x^{8}$ \\
		& $ + \; 384485209291462309415622819836536642510793775793569792x^{7}$ \\
		& $ + \; 5816470292685027671662290569991129400571245636014440448x^{6}$ \\
		& $ - \; 371131509602563559938383240089569329622889937518515453952x^{5}$ \\
		& $ - \; 7480021683546696333619583149415575181418775652922195312640x^{4}$ \\
		& $ + \; 164143334208143730520533868283744520771676778757659683717120x^{3}$ \\
		& $ + \; 3654063025430843865204601351924922535583012113242950851887104x^{2}$ \\
		& $ - \; 19859588203354051542890817248546892512751903736345856095289344x$ \\
		& $ - \; 375067918205041763744279215282103445336143562416034878481956864$ \\[0.2ex]
		\hline\mystrut(13,0)
		Inverse square & $401.5 = \frac{803}{2}$ \\
		potential $1/r^2$
	\end{longtable}
\end{center}

\noindent
\textit{\textbf{Symmetries -}}

The symmetry groups for 32 points are identical for all 3 potentials. This is quite remarkable as the 32-point set is not universally optimal. However this result that all 3 potentials have the same point set is also verified in \cite{47} to a precision of about 11 digits.

\begin{center}
	\begin{tabular}{l|l}
		\multicolumn{2}{c}{Symmetries - 32 points} \\
		\hline\Tstrut
		planes & [[2, 810], [4, 30], [8, 180], [46, 10], [60, 6], [64, 15]] \\[0.2ex]
		\hline\Tstrut
		Gram groups & [[32, 2], [60, 4], [120, 6]] \\
		\hline\Tstrut
		Polygons & [[3, 1680], [4, 420], [5, 36], [6, 20], [8, 15]]
	\end{tabular}
\end{center}

\subsection{33 points}
For 33 points, the embedded heptagon indicates an axis of symmetry, so the configuration was rotated into alignment on the equator and it was found that 13 points were above the equator and 13 points below, in a balanced arrangement, but no other embedded polygons occurred along the same alignment axis (quadrilaterals do exist, but not in alignment with the heptagon).

\begin{figure}[ht]
	\begin{center}
		\includegraphics[type=pdf,ext=pdf,read=pdf,height=1in,width=1in,angle=0]{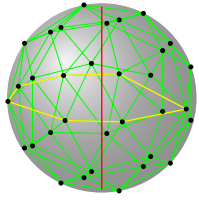}
		\caption{33 points.}
		\label{fig:33pts}
	\end{center}
\end{figure}
The heptagon is outlined in yellow in figure \ref{fig:33pts}.

\noindent
\textit{\textbf{Parameterization for all 3 potentials}}

\begin{longtable}[c]{r|ccc}
	\caption{Parameterization for 33 points} \\
	pt & $x$ & $y$ & $z$ \\
	\hline\vspace*{-2.2ex}
	\endfirsthead
	\multicolumn{4}{c}%
	{\tablename\ \thetable\ -- 33 points parameterization -- \textit{continued\ldots}} \\
	pt & $x$ & $y$ & $z$ \\
	\hline\vspace*{-2.2ex}
	\endhead
	1 & $-b$ & $\sqrt{1-a^2-b^2}$ & $a$ \\[0.5em]
	2 & $d$ & $-\sqrt{1-c^2-d^2}$ & $c$ \\[0.5em]
	3 & $-f$ & $-\sqrt{1-e^2-f^2}$ & $e$ \\[0.5em]
	4 & $h$ & $\sqrt{1-g^2-h^2}$ & $g$ \\[0.5em]
	5 & $j$ & $\sqrt{1-i^2-j^2}$ & $i$ \\[0.5em]
	6 & $-l$ & $\sqrt{1-k^2-l^2}$ & $k$ \\[0.5em]
	7 & $n$ & $-\sqrt{1-m^2-n^2}$ & $m$ \\[0.5em]
	8 & $-p$ & $-\sqrt{1-o^2-p^2}$ & $o$ \\[0.5em]
	9 & $-r$ & $\sqrt{1-q^2-r^2}$ & $q$ \\[0.5em]
	10 & $t$ & $-\sqrt{1-s^2-t^2}$ & $s$ \\[0.5em]
	11 & $-v$ & $-\sqrt{1-u^2-v^2}$ & $u$ \\[0.5em]
	12 & $x$ & $\sqrt{1-w^2-x^2}$ & $w$ \\[0.5em]
	13 & $z$ & $\sqrt{1-y^2-z^2}$ & $y$ \\[0.5em]
	14 & $1$ & $0$ & $0$ \\[0.5em]
	15 & $-A$ & $\sqrt{1-A^2}$ & $0$ \\[0.5em]
	16 & $-B$ & $\sqrt{1-B^2}$ & $0$ \\[0.5em]
	17 & $C$ & $-\sqrt{1-C^2}$ & $0$ \\[0.5em]
	18 & $-D$ & $\sqrt{1-D^2}$ & $0$ \\[0.5em]
	19 & $E$ & $-\sqrt{1-E^2}$ & $0$ \\[0.5em]
	20 & $-F$ & $-\sqrt{1-F^2}$ & $0$ \\[0.5em]
	21 & $z$ & $\sqrt{1-y^2-z^2}$ & $-y$ \\[0.5em]
	22 & $x$ & $\sqrt{1-w^2-x^2}$ & $-w$ \\[0.5em]
	23 & $-v$ & $-\sqrt{1-u^2-v^2}$ & $-u$ \\[0.5em]
	24 & $t$ & $-\sqrt{1-s^2-t^2}$ & $-s$ \\[0.5em]
	25 & $-r$ & $\sqrt{1-q^2-r^2}$ & $-q$ \\[0.5em]
	26 & $-p$ & $-\sqrt{1-o^2-p^2}$ & $-o$ \\[0.5em]
	27 & $n$ & $-\sqrt{1-m^2-n^2}$ & $-m$ \\[0.5em]
	28 & $-l$ & $\sqrt{1-k^2-l^2}$ & $-k$ \\[0.5em]
	29 & $j$ & $\sqrt{1-i^2-j^2}$ & $-i$ \\[0.5em]
	30 & $h$ & $\sqrt{1-g^2-h^2}$ & $-g$ \\[0.5em]
	31 & $-f$ & $-\sqrt{1-e^2-f^2}$ & $-e$ \\[0.5em]
	32 & $d$ & $-\sqrt{1-c^2-d^2}$ & $-c$ \\[0.5em]
	33 & $-b$ & $\sqrt{1-a^2-b^2}$ & $-a$
\end{longtable}

\noindent
\textit{\textbf{Parameterization values --}}

It takes 32 parameters, $a$ - $z$, $A$ - $F$, to constrain the polyhedron. The values given to 19 digits are listed below:

\begin{longtable}[c]{c|c|c|c}
	\caption{Parameter values for 33 points} \\
	Parameter & log & 1/r & $1/r^2$ \\
	\hline\vspace*{-2.2ex}
	\endfirsthead
	\multicolumn{4}{c}%
	{\tablename\ \thetable\ -- 33 points parameters -- \textit{continued}} \\
	Parameter & log & 1/r & $1/r^2$ \\
	\hline\vspace*{-2.2ex}
	\endhead
	$a$ & 0.9441354405857699281 & 0.9443379028688371297 & 0.9447610253876353098 \\
	$b$ & 0.2905640124100343812 & 0.3037040330136551239 & 0.3016836871749528318 \\
	$c$ & 0.9345998724028032232 & 0.9333938530416051179 & 0.9316825901300259677 \\
	$d$ & 0.3017306777032591648 & 0.3386313079897339111 & 0.3404703649925359162 \\
	$e$ & 0.8504519070490259610 & 0.8518529803217173233 & 0.8537453625341212703 \\
	$f$ & 0.2461935540949194669 & -0.2070212968704013093 & -0.2138902228328844522 \\
	$g$ & 0.8373194344583771142 & 0.8359803174718946552 & 0.8343922633990169106 \\
	$h$ & 0.2431523959305478594 & -0.2232902769227805514 & -0.2263045153295219886 \\
	$i$ & 0.6546496836754062898 & 0.6510565523039134609 & 0.6472905600791991335 \\
	$j$ & 0.7440042280355189756 & 0.3685982446524377337 & 0.3641511969357234569 \\
	$k$ & 0.6022814340652630615 & 0.6029568676883085885 & 0.6037529559528323030 \\
	$l$ & 0.3140334969622244110 & 0.7672333256310515118 & 0.7694927097706817583 \\
	$m$ & 0.5909892333166061138 & 0.5893326407922974706 & 0.5879108250331423993 \\
	$n$ & 0.2661015860213742934 & 0.7594740649765698504 & 0.7636135093616019175 \\
	$o$ & 0.5723203334842430753 & 0.5732358840938099116 & 0.5737728991063233710 \\
	$p$ & 0.7683290942541970021 & 0.2603038344456836052 & 0.2487583395856551723 \\
	$q$ & 0.5491395658902830166 & 0.5486023941900991534 & 0.5485488782362098275 \\
	$r$ & 0.7710177904574511397 & 0.7362560795962369714 & 	0.7305244960195997603 \\
	$s$ & 0.4442465871890066326 & 0.4420289651486112907 & 0.4361791084974002231 \\
	$t$ & 0.7922342695397467159 & 0.8260944317350690476 & 0.8248850189495450027 \\
	$u$ & 0.3706514795426231661 & 0.3724889251567467111 & 0.3762150755311320795 \\
	$v$ & 0.3907372189445071434 & -0.4041288502324899330 & -0.4107125501713670196 \\
	$w$ & 0.3324420333566057344 & 0.3313739092027784255 & 0.3300473247902436542 \\
	$x$ & 0.2581762659301997113 & -0.5451696716795546735 & -0.5602714990848728267 \\
	$y$ & 0.2916082231040449711 & 0.2922911885557252468 & 0.2935410680557640754 \\
	$z$ & 0.7573642041186359211 & 0.01965858954718621313 & 0.001443841219829245825 \\
	$A$ & 0.2823458487681062946 & -0.6267209832905903770 & -0.6101798489082580366 \\
	$B$ & 0.7745508124916516604 & -0.8171346843197990691 & -0.8186981292619572688 \\
	$C$ & 0.6364563713080943362 & -0.07019679826361042841 & -0.05740747625502538887 \\
	$D$ & 0.9997173400385455117 & 0.9234493530503589647 & 0.6376736480979082347 \\
	$E$ & 0.07395591266845277821 & -0.9789021592239545404 & -0.9752922714179999114 \\
	$F$ & 0.8175752029901418422 & 0.6486352294571746269 & 0.9298417669352073362 \\
	\hline\Tstrut
	$energy$ & -134.7478208243334052\ldots & 440.2040574476473575\ldots & 431.9318385881776597\ldots
\end{longtable}
After checking the 32 parameters to 50,014 digits, the algebraic degree for the minimal polynomials is $>420$.

\noindent
\textit{\textbf{Symmetries -}}

The symmetry groups for 33 points are identical for all 3 potentials.

\begin{center}
	\begin{tabular}{l|l}
		\multicolumn{2}{c}{Symmetries - 33 points} \\
		\hline\Tstrut
		planes & [[4, 78], [35, 1]] \\[0.2ex]
		\hline\Tstrut
		Gram groups & [[2, 34], [4, 247], [33, 1]] \\
		\hline\Tstrut
		Polygons & [[4, 78], [7, 1]]
	\end{tabular}
\end{center}

\subsection{34 points}
For 34 points, it was discovered that the arrangement for all 3 potentials is basically the same, 17 dipoles along one axis, with 8 dipoles above the xy-plane and 8 dipoles below, such that a pair of dipoles is the same opposing distance from the z-axis such that the arrangement is balanced. The red dipole on the equator is aligned along the \textit{x-axis}. This is an 2:2:2:2:2:2:2:2:2:2:2:2:2:2:2:2:2 arrangement.

\begin{figure}[ht]
	\begin{center}
		\includegraphics[type=pdf,ext=pdf,read=pdf,height=1in,width=1in,angle=0]{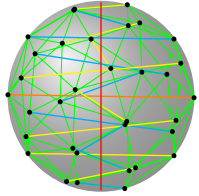}
		\caption{34 points.}
		\label{fig:34pts}
	\end{center}
\end{figure}
\noindent
\textit{\textbf{Algebraic Spherical Code -}}

The algebraic spherical code for 34 points, for all 3 potentials, created from 16 parameters, $a-p$ is:
\begin{longtable}[c]{r|ccc}
	\caption{Parameterization for 34 points} \\
	pt & $x$ & $y$ & $z$ \\
	\hline\vspace*{-2.2ex}
	\endfirsthead
	\multicolumn{4}{c}%
	{\tablename\ \thetable\ -- 34 points -- \textit{continued\ldots}} \\
	pt & $x$ & $y$ & $z$ \\
	\hline\vspace*{-2.2ex}
	\endhead
	1 & $b$ & $\sqrt{1-a^2-b^2}$ & $a$ \\[0.5ex]
	2 & $-b$ & $-\sqrt{1-a^2-b^2}$ & $a$ \\[0.5ex]
	3 & $d$ & $-\sqrt{1-c^2-d^2}$ & $c$ \\[0.5ex]
	4 & $-d$ & $\sqrt{1-c^2-d^2}$ & $c$ \\[0.5ex]
	5 & $f$ & $\sqrt{1-e^2-f^2}$ & $e$ \\[0.5ex]
	6 & $-f$ & $-\sqrt{1-e^2-f^2}$ & $e$ \\[0.5ex]
	7 & $h$ & $\sqrt{1-g^2-h^2}$ & $g$ \\[0.5ex]
	8 & $-h$ & $-\sqrt{1-g^2-h^2}$ & $g$ \\[0.5ex]
	9 & $j$ & $-\sqrt{1-i^2-j^2}$ & $i$ \\[0.5ex]
	10 & $-j$ & $\sqrt{1-i^2-j^2}$ & $i$ \\[0.5ex]
	11 & $l$ & $-\sqrt{1-k^2-l^2}$ & $k$ \\[0.5ex]
	12 & $-l$ & $\sqrt{1-k^2-l^2}$ & $k$ \\[0.5ex]
	13 & $n$ & $\sqrt{1-m^2-n^2}$ & $m$ \\[0.5ex]
	14 & $-n$ & $-\sqrt{1-m^2-n^2}$ & $m$ \\[0.5ex]
	15 & $p$ & $\sqrt{1-o^2-p^2}$ & $o$ \\[0.5ex]
	16 & $-p$ & $-\sqrt{1-o^2-p^2}$ & $o$ \\[0.5ex]
	17 & $1$ & $0$ & $0$ \\[0.5ex]
	18 & $-1$ & $0$ & $0$ \\[0.5ex]
	19 & $p$ & $-\sqrt{1-o^2-p^2}$ & $-o$ \\[0.5ex]
	20 & $-p$ & $\sqrt{1-o^2-p^2}$ & $-o$ \\[0.5ex]
	21 & $n$ & $-\sqrt{1-m^2-n^2}$ & $-m$ \\[0.5ex]
	22 & $-n$ & $\sqrt{1-m^2-n^2}$ & $-m$ \\[0.5ex]
	23 & $l$ & $\sqrt{1-k^2-l^2}$ & $-k$ \\[0.5ex]
	24 & $-l$ & $-\sqrt{1-k^2-l^2}$ & $-k$ \\[0.5ex]
	25 & $j$ & $\sqrt{1-i^2-j^2}$ & $-i$ \\[0.5ex]
	26 & $-j$ & $-\sqrt{1-i^2-j^2}$ & $-i$ \\[0.5ex]
	27 & $h$ & $-\sqrt{1-g^2-h^2}$ & $-g$ \\[0.5ex]
	28 & $-h$ & $\sqrt{1-g^2-h^2}$ & $-g$ \\[0.5ex]
	29 & $f$ & $-\sqrt{1-e^2-f^2}$ & $-e$ \\[0.5ex]
	30 & $-f$ & $\sqrt{1-e^2-f^2}$ & $-e$ \\[0.5ex]
	31 & $d$ & $\sqrt{1-c^2-d^2}$ & $-c$ \\[0.5ex]
	32 & $-d$ & $-\sqrt{1-c^2-d^2}$ & $-c$ \\[0.5ex]
	33 & $b$ & $-\sqrt{1-a^2-b^2}$ & $-a$ \\[0.5ex]
	34 & $-b$ & $\sqrt{1-a^2-b^2}$ & $-a$ \\[0.5ex]
\end{longtable}
The values for these 16 parameters to 19 digits is given below, along with the minimal energy for the points under the potential:

\begin{longtable}[c]{c|c|c|c}
	\caption{Parameter values for 34 points} \\
	Parameter & log & 1/r & $1/r^2$ \\
	\hline\vspace*{-2.2ex}
	\endfirsthead
	\multicolumn{4}{c}%
	{\tablename\ \thetable\ -- 34 points parameters -- \textit{continued}} \\
	Parameter & log & 1/r & $1/r^2$ \\
	\hline\vspace*{-2.2ex}
	\endhead
	$a$ & 0.9501746904318417810 & 0.9507345860669045433 & 0.9512670049190477768 \\
	$b$ & 0.2690296903895190684 & 0.2667993421427814118 & 0.2648543738542885809 \\
	$c$ & 0.8334906701638139447 & 0.8347171305054724743 & 0.8358386082813027937 \\
	$d$ & 0.3314547676755688974 & 0.3301374286789593185 & 0.3286706495008522924 \\
	$e$ & 0.6720105807965419493 & 0.6709162957697776888 & 0.6701241936721569551 \\
	$f$ & 0.3591611788440802614 & 0.3589314226631725173 & 0.3582682173168182203 \\
	$g$ & 0.6203239065468236863 & 0.6218438093075785121 & 0.6233799885625532601 \\
	$h$ & 0.7843451748727355187 & 0.7831345168065170310 & 0.7819041440050671677 \\
	$i$ & 0.4524090009853990815 & 0.4531721066862990811 & 0.4535566868316177422 \\
	$j$ & 0.1824031292301683503 & 0.1814743809559557118 & 0.1808113840729157640 \\
	$k$ & 0.3428005610247186752 & 0.3439704385614441753 & 0.3450194919610273067 \\
	$l$ & 0.7570163447483596234 & 0.7553751608582816613 & 0.7536679020995221525 \\
	$m$ & 0.2729534769375722201 & 0.2722721572403301061 & 0.2720230702769717561 \\
	$n$ & 0.8120626542505883933 & 0.8113438409105690329 & 0.8105442906535017426 \\
	$o$ & 0.09886379372735685423 & 0.1003489964767353087 & 0.1020597436732134232 \\
	$p$ & 0.3576913654812188770 & 0.3563114827352296087 & 0.3550466709305113318 \\
	\hline\Tstrut
	$energy$ & -142.3758522709015841 & 468.9048532813432615 & 462.7012364217078132
\end{longtable}
All 16 parameters for all 3 potentials have been found to 50,014 digits accuracy, but the degree of the algebraic polynomials $>360$.

\noindent
\textit{\textbf{Symmetries -}}

The symmetry groups for 34 points are identical for all 3 potentials.

\begin{center}
	\begin{tabular}{l|l}
		\multicolumn{2}{c}{Symmetries - 34 points} \\
		\hline\Tstrut
		planes & [[4, 16]] \\[0.2ex]
		\hline\Tstrut
		Gram groups & [[2, 1], [4, 24], [8, 128], [34, 1]] \\
		\hline\Tstrut
		Polygons & [[4, 16]]
	\end{tabular}
\end{center}

\subsection{35 points}
The optimal configuration for 35 points does not allow a parameterization to occur for all 3 potentials. This is the second set of points where this occurs.

\begin{figure}[ht]
	\begin{center}
		\includegraphics[type=pdf,ext=pdf,read=pdf,height=1in,width=1in,angle=0]{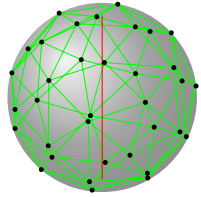}
		\caption{35 points.}
		\label{fig:35pts}
	\end{center}
\end{figure}
\noindent
\textit{\textbf{Minimal Energy values -}}

The coordinates for 35 points are known to 77 digits for the \textit{log} potential and 38 digits for the other two. The minimal energies have been determined for all 3 potentials as well.
\begin{center}
	\begin{tabular}{l|l}
		\multicolumn{2}{c}{Minimal Energy - 35 points} \\
		\hline\Tstrut
		logarithmic & -150.1920585107381031\ldots \\[0.2ex]
		\hline\Tstrut
		Coulomb & 498.5698724906454068\ldots \\[0.2ex]
		\hline\Tstrut
		Inverse square law & 494.8164319541329121\ldots
	\end{tabular}
\end{center}

\noindent
\textit{\textbf{Symmetries -}}

The symmetry groups are identical for 35 points under all 3 potentials.

\begin{center}
	\begin{tabular}{l|l}
		\multicolumn{2}{c}{Symmetries - 35 points} \\
		\hline\Tstrut
		planes & [] \\[0.2ex]
		\hline\Tstrut
		Gram groups & [[2, 17], [4, 289], [35, 1]] \\
		\hline\Tstrut
		Polygons & []
	\end{tabular}
\end{center}

\subsection{36 points}
The optimal configuration for 36 points does not allow a parameterization to be determined for all 3 potentials. This is the third set of points where this occurs (previous: 26 pts, 35 pts);

The algorithmic searches for the optimal configurations took an exceptionally long time, for the \textit{percolating anneal} program indicating the presence of many local minima. Some point sets converged quickly, under the potentials, but 36 points made the searches very arduous, the Coulomb $1/r$ search took over 6 weeks to converge, for billions of loops in the algorithm to converge to 38 digits accuracy.

\begin{figure}[ht]
	\begin{center}
		\includegraphics[type=pdf,ext=pdf,read=pdf,height=1in,width=1in,angle=0]{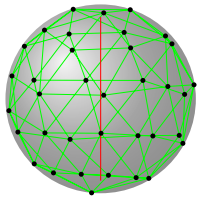}
		\caption{36 points.}
		\label{fig:36pts}
	\end{center}
\end{figure}
\noindent
\textit{\textbf{Minimal Energy values -}}

The coordinates for 36 points are known to 77 digits for the \textit{log} potential and 38 digits for the other two. The minimal energies have been determined for all 3 potentials as well.
\begin{center}
	\begin{tabular}{l|l}
		\multicolumn{2}{c}{Minimal Energy - 36 points} \\
		\hline\Tstrut
		logarithmic & -158.2240684255879753\ldots \\[0.2ex]
		\hline\Tstrut
		Coulomb & 529.1224083754138052\ldots \\[0.2ex]
		\hline\Tstrut
		Inverse square law & 527.9142565801658786\ldots
	\end{tabular}
\end{center}

\noindent
\textit{\textbf{Symmetries -}}

Interestingly, the Gram matrix differs for the \textit{Inverse Square $1/r^2$} potentials

\begin{center}
	\begin{tabular}{l|l}
		\multicolumn{2}{c}{Symmetries - 36 points - \textit{log} or \textit{Coulomb}} \\
		\hline\Tstrut
		planes & [] \\[0.2ex]
		\hline\Tstrut
		Gram groups & [[4, 27], [8, 144], [36, 1]] \\
		\hline\Tstrut
		Polygons & []
	\end{tabular}
\end{center}

\begin{center}
	\begin{tabular}{l|l}
		\multicolumn{2}{c}{Symmetries - 36 points - \textit{Inverse}$^2$} \\
		\hline\Tstrut
		planes & [] \\[0.2ex]
		\hline\Tstrut
		Gram groups & [[2, 18], [4, 306], [36, 1]] \\
		\hline\Tstrut
		Polygons & []
	\end{tabular}
\end{center}

\subsection{37 points}
The 9 pentagons embedded in the polyhedron suggest the preferred symmetry arrangement, and after 2 rotations, it was found that the figure is a pleasing arrangement of a pole axis and 7 co-planar pentagons arranged so that a line through their centroids is the pole axis. This enabled only 6 parameters necessary to constrain the arrangement. This is an 1:5:5:5:5:5:5:5:1 arrangement.

\begin{figure}[ht]
	\begin{center}
		\includegraphics[type=pdf,ext=pdf,read=pdf,height=1in,width=1in,angle=0]{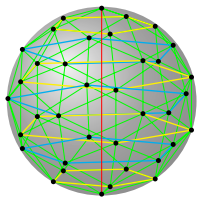}
		\caption{37 points.}
		\label{fig:37pts}
	\end{center}
\end{figure}
The embedded pentagons are outlined in yellow and cyan in figure \ref{fig:37pts}.

\noindent
\textit{\textbf{Parameterization for all 3 potentials}}

\begin{longtable}[c]{r|ccc}
	\caption{Parameterization for 37 points} \\
	pt & $x$ & $y$ & $z$ \\
	\hline\vspace*{-2.2ex}
	\endfirsthead
	\multicolumn{4}{c}%
	{\tablename\ \thetable\ -- 37 points parameterization -- \textit{continued\ldots}} \\
	pt & $x$ & $y$ & $z$ \\
	\hline\vspace*{-2.2ex}
	\endhead
	1 & $0$ & $0$ & $1$ \\[0.7em]
	2 & $b$ & $0$ & $a$ \\[0.7em]
	3 & $\frac{b\sqrt{5}-1}{4}$ & $\frac{b\sqrt{10+2\sqrt{5}}}{4}$ & $a$ \\[0.7em]
	4 & $-\frac{b\sqrt{5}+1}{4}$ & $\frac{b\sqrt{10-2\sqrt{5}}}{4}$ & $a$ \\[0.7em]
	5 & $-\frac{b\sqrt{5}+1}{4}$ & $-\frac{b\sqrt{10-2\sqrt{5}}}{4}$ & $a$ \\[0.7em]
	6 & $\frac{b\sqrt{5}-1}{4}$ & $-\frac{b\sqrt{10+2\sqrt{5}}}{4}$ & $a$ \\[0.7em]
	7 & $d$ & $0$ & $c$ \\[0.7em]
	8 & $\frac{d\sqrt{5}-1}{4}$ & $\frac{d\sqrt{10+2\sqrt{5}}}{4}$ & $c$ \\[0.7em]
	9 & $-\frac{d\sqrt{5}+1}{4}$ & $\frac{d\sqrt{10-2\sqrt{5}}}{4}$ & $c$ \\[0.7em]
	10 & $-\frac{d\sqrt{5}+1}{4}$ & $-\frac{d\sqrt{10-2\sqrt{5}}}{4}$ & $c$ \\[0.7em]
	11 & $\frac{d\sqrt{5}-1}{4}$ & $-\frac{d\sqrt{10+2\sqrt{5}}}{4}$ & $c$ \\[0.7em]
	12 & $f$ & $0$ & $e$ \\[0.7em]
	13 & $\frac{f\sqrt{5}-1}{4}$ & $\frac{f\sqrt{10+2\sqrt{5}}}{4}$ & $e$ \\[0.7em]
	14 & $-\frac{f\sqrt{5}+1}{4}$ & $\frac{f\sqrt{10-2\sqrt{5}}}{4}$ & $e$ \\[0.7em]
	15 & $-\frac{f\sqrt{5}+1}{4}$ & $-\frac{f\sqrt{10-2\sqrt{5}}}{4}$ & $e$ \\[0.7em]
	16 & $\frac{f\sqrt{5}-1}{4}$ & $-\frac{f\sqrt{10+2\sqrt{5}}}{4}$ & $e$ \\[0.7em]
	17 & $1$ & $0$ & $0$ \\[0.7em]
	18 & $\frac{\sqrt{5}-1}{4}$ & $\frac{\sqrt{10+2\sqrt{5}}}{4}$ & $0$ \\[0.7em]
	19 & $\frac{-\sqrt{5}-1}{4}$ & $\frac{\sqrt{10-2\sqrt{5}}}{4}$ & $0$ \\[0.7em]
	20 & $\frac{-\sqrt{5}-1}{4}$ & $-\frac{\sqrt{10-2\sqrt{5}}}{4}$ & $0$ \\[0.7em]
	21 & $\frac{\sqrt{5}-1}{4}$ & $-\frac{\sqrt{10+2\sqrt{5}}}{4}$ & $0$ \\[0.7em]
	22 & $f$ & $0$ & $-e$ \\[0.7em]
	23 & $\frac{f\sqrt{5}-1}{4}$ & $\frac{f\sqrt{10+2\sqrt{5}}}{4}$ & $-e$ \\[0.7em]
	24 & $-\frac{f\sqrt{5}+1}{4}$ & $\frac{f\sqrt{10-2\sqrt{5}}}{4}$ & $-e$ \\[0.7em]
	25 & $-\frac{f\sqrt{5}+1}{4}$ & $-\frac{f\sqrt{10-2\sqrt{5}}}{4}$ & $-e$ \\[0.7em]
	26 & $\frac{f\sqrt{5}-1}{4}$ & $-\frac{f\sqrt{10+2\sqrt{5}}}{4}$ & $-e$ \\[0.7em]
	27 & $d$ & $0$ & $-c$ \\[0.7em]
	28 & $\frac{d\sqrt{5}-1}{4}$ & $\frac{d\sqrt{10+2\sqrt{5}}}{4}$ & $-c$ \\[0.7em]
	29 & $-\frac{d\sqrt{5}+1}{4}$ & $\frac{d\sqrt{10-2\sqrt{5}}}{4}$ & $-c$ \\[0.7em]
	30 & $-\frac{d\sqrt{5}+1}{4}$ & $-\frac{d\sqrt{10-2\sqrt{5}}}{4}$ & $-c$ \\[0.7em]
	31 & $\frac{d\sqrt{5}-1}{4}$ & $-\frac{d\sqrt{10+2\sqrt{5}}}{4}$ & $-c$ \\[0.7em]
	32 & $b$ & $0$ & $-a$ \\[0.7em]
	33 & $\frac{b\sqrt{5}-1}{4}$ & $\frac{b\sqrt{10+2\sqrt{5}}}{4}$ & $-a$ \\[0.7em]
	34 & $-\frac{b\sqrt{5}+1}{4}$ & $\frac{b\sqrt{10-2\sqrt{5}}}{4}$ & $-a$ \\[0.7em]
	35 & $-\frac{b\sqrt{5}+1}{4}$ & $-\frac{b\sqrt{10-2\sqrt{5}}}{4}$ & $-a$ \\[0.7em]
	36 & $\frac{b\sqrt{5}-1}{4}$ & $-\frac{b\sqrt{10+2\sqrt{5}}}{4}$ & $-a$ \\[0.7em]
	37 & $0$ & $0$ & $-1$
\end{longtable}

\noindent
\textit{\textbf{Parameterization values --}}

It takes 6 parameters, $a$ - $f$, to constrain the polyhedron. The values given to 19 digits are listed below:

\begin{longtable}[c]{c|c|c|c}
	\caption{Parameter values for 37 points} \\
	Parameter & log & 1/r & $1/r^2$ \\
	\hline\vspace*{-2.2ex}
	\endfirsthead
	\multicolumn{4}{c}%
	{\tablename\ \thetable\ -- 37 points parameters -- \textit{continued}} \\
	Parameter & log & 1/r & $1/r^2$ \\
	\hline\vspace*{-2.2ex}
	\endhead
	$a$ & 0.8207614357951864012 & 0.8217755307731263974 & 0.8230114007595582573 \\
	$b$ & -0.5712710963381782287 & -0.5698113521355523059 & -0.5680248535229684736 \\
	$c$ & 0.5332210870197916288 & 0.5348307762724087591 & 0.5367177983270355014 \\
	$d$ & 0.8459759289468180863 & 0.8449591947259066592 & 0.8437618176707095150 \\
	$e$ & 0.2848889779325672258 & 0.2875615914220408508 & 0.2901640938863503629 \\
	$f$ & -0.9585605198695266650 & -0.9577621474765137333 & -0.9569769060009302293 \\
	\hline\Tstrut
	$energy$ & -166.4506975239976931 & 560.6188877310436776 & 562.2556382320535124
\end{longtable}
After checking the 6 parameters to 5,009 digits, the algebraic degree for the minimal polynomials is $>200$. There was a problem with the Jacobian matrix, making it impossible to find 50,014 digits of each parameter, instead a direct search to 5,009 digits was done for all 3 potentials.

\noindent
\textit{\textbf{Symmetries -}}

The symmetry groups are identical for 37 points under all 3 potentials.

This configuration offers promise to be parameterized.

\begin{center}
	\begin{tabular}{l|l}
		\multicolumn{2}{c}{Symmetries - 37 points} \\
		\hline\Tstrut
		planes & [[4, 450], [24, 5], [70, 1], [84, 5]] \\[0.2ex]
		\hline\Tstrut
		Gram groups & [[2, 1], [10, 5], [20, 28], [37, 1], [40, 18]] \\
		\hline\Tstrut
		Polygons & [[4, 480], [5, 7], [9, 5]]
	\end{tabular}
\end{center}

\subsection{38 points}
The configuration for 38 points has embedded polygons, a preferred orientation found an arrangement of two pole points, or dipole, and 6 parallel hexagons such that the pole axis passes through their centroids. The arrangement is [1:6:6:6:6:6:6:1].

\begin{figure}[ht]
	\begin{center}
		\includegraphics[type=pdf,ext=pdf,read=pdf,height=1in,width=1in,angle=0]{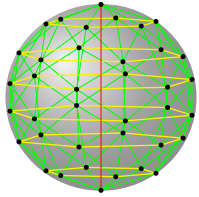}
		\caption{38 points.}
		\label{fig:38pts}
	\end{center}
\end{figure}

\noindent
\textit{\textbf{Parameterized Structure -}}

Upon careful analysis, only 3 algebraic parameters $a$, $b$, and $c$ are required to constrain the structure of the 38 points under minimal energy. This structure is given below:

\begin{longtable}[c]{r|ccc}
	\caption{Parameterization for 38 points} \\
	pt & $x$ & $y$ & $z$ \\
	\hline\vspace*{-2.2ex}
	\endfirsthead
	\multicolumn{4}{c}%
	{\tablename\ \thetable\ -- 38 points parameters -- \textit{continued}} \\
	pt & $x$ & $y$ & $z$ \\
	\hline\vspace*{-2.2ex}
	\endhead
	1 & $0$ & $0$ & $1$ \\[0.7ex]
	2 & $0$ & $\sqrt{1-a^2}$ & $a$ \\[0.7ex]
	3 & $\frac{\sqrt{3}}{2}\sqrt{1-a^2}$ & $\frac{\sqrt{1-a^2}}{2}$ & $a$ \\[0.7ex]
	4 & $\frac{\sqrt{3}}{2}\sqrt{1-a^2}$ & $-\frac{\sqrt{1-a^2}}{2}$ & $a$ \\[0.7ex]
	5 & $0$ & $-\sqrt{1-a^2}$ & $a$ \\[0.7ex]
	6 & $-\frac{\sqrt{3}}{2}\sqrt{1-a^2}$ & $-\frac{\sqrt{1-a^2}}{2}$ & $a$ \\[0.7ex]
	7 & $-\frac{\sqrt{3}}{2}\sqrt{1-a^2}$ & $\frac{\sqrt{1-a^2}}{2}$ & $a$ \\[0.7ex]
	8 & $\sqrt{1-b^2}$ & $0$ & $b$ \\[0.7ex]
	9 & $\frac{\sqrt{1-b^2}}{2}$ & $\frac{\sqrt{3}}{2}\sqrt{1-b^2}$ & $b$ \\[0.7ex]
	10 & $-\frac{\sqrt{1-b^2}}{2}$ & $\frac{\sqrt{3}}{2}\sqrt{1-b^2}$ & $b$ \\[0.7ex]
	11 & $-\sqrt{1-b^2}$ & $0$ & $b$ \\[0.7ex]
	12 & $-\frac{\sqrt{1-b^2}}{2}$ & $-\frac{\sqrt{3}}{2}\sqrt{1-b^2}$ & $b$ \\[0.7ex]
	13 & $\frac{\sqrt{1-b^2}}{2}$ & $-\frac{\sqrt{3}}{2}\sqrt{1-b^2}$ & $b$ \\[0.7ex]
	14 & $0$ & $\sqrt{1-c^2}$ & $c$ \\[0.7ex]
	15 & $\frac{\sqrt{3}}{2}\sqrt{1-c^2}$ & $\frac{\sqrt{1-c^2}}{2}$ & $c$ \\[0.7ex]
	16 & $\frac{\sqrt{3}}{2}\sqrt{1-c^2}$ & $-\frac{\sqrt{1-c^2}}{2}$ & $c$ \\[0.7ex]
	17 & $0$ & $-\sqrt{1-c^2}$ & $c$ \\[0.7ex]
	18 & $-\frac{\sqrt{3}}{2}\sqrt{1-c^2}$ & $-\frac{\sqrt{1-c^2}}{2}$ & $c$ \\[0.7ex]
	19 & $-\frac{\sqrt{3}}{2}\sqrt{1-c^2}$ & $\frac{\sqrt{1-c^2}}{2}$ & $c$ \\[0.7ex]
	20 & $\sqrt{1-c^2}$ & $0$ & $-c$ \\[0.7ex]
	21 & $\frac{\sqrt{1-c^2}}{2}$ & $\frac{\sqrt{3}}{2}\sqrt{1-c^2}$ & $-c$ \\[0.7ex]
	22 & $-\frac{\sqrt{1-c^2}}{2}$ & $\frac{\sqrt{3}}{2}\sqrt{1-c^2}$ & $-c$ \\[0.7ex]
	23 & $-\sqrt{1-c^2}$ & $0$ & $-c$ \\[0.7ex]
	24 & $-\frac{\sqrt{1-c^2}}{2}$ & $-\frac{\sqrt{3}}{2}\sqrt{1-c^2}$ & $-c$ \\[0.7ex]
	25 & $\frac{\sqrt{1-c^2}}{2}$ & $-\frac{\sqrt{3}}{2}\sqrt{1-c^2}$ & $-c$ \\[0.7ex]
	26 & $0$ & $\sqrt{1-b^2}$ & $-b$ \\[0.7ex]
	27 & $\frac{\sqrt{3}}{2}\sqrt{1-b^2}$ & $\frac{\sqrt{1-b^2}}{2}$ & $-b$ \\[0.7ex]
	28 & $\frac{\sqrt{3}}{2}\sqrt{1-b^2}$ & $-\frac{\sqrt{1-b^2}}{2}$ & $-b$ \\[0.7ex]
	29 & $0$ & $-\sqrt{1-b^2}$ & $-b$ \\[0.7ex]
	30 & $-\frac{\sqrt{3}}{2}\sqrt{1-b^2}$ & $-\frac{\sqrt{1-b^2}}{2}$ & $-b$ \\[0.7ex]
	31 & $-\frac{\sqrt{3}}{2}\sqrt{1-b^2}$ & $\frac{\sqrt{1-b^2}}{2}$ & $-b$ \\[0.7ex]
	32 & $\sqrt{1-a^2}$ & $0$ & $-a$ \\[0.7ex]
	33 & $\frac{\sqrt{1-a^2}}{2}$ & $\frac{\sqrt{3}}{2}\sqrt{1-a^2}$ & $-a$ \\[0.7ex]
	34 & $-\frac{\sqrt{1-a^2}}{2}$ & $\frac{\sqrt{3}}{2}\sqrt{1-a^2}$ & $-a$ \\[0.7ex]
	35 & $-\sqrt{1-a^2}$ & $0$ & $-a$ \\[0.7ex]
	36 & $-\frac{\sqrt{1-a^2}}{2}$ & $-\frac{\sqrt{3}}{2}\sqrt{1-a^2}$ & $-a$ \\[0.7ex]
	37 & $\frac{\sqrt{1-a^2}}{2}$ & $-\frac{\sqrt{3}}{2}\sqrt{1-a^2}$ & $-a$ \\[0.7ex]
	38 & $0$ & $0$ & $-1$
\end{longtable}

\noindent
\textit{\textbf{Algebraic parameters values -}}

The values to 19 digits of the 3 parameters $a-c$ which are optimal for the minimal solutions of 38 points are:

\begin{longtable}[c]{c|c|c|c}
	\caption{Parameter values for 38 points} \\
	Parameter & log & 1/r & $1/r^2$ \\
	\hline\vspace*{-2.2ex}
	\endfirsthead
	\multicolumn{4}{c}%
	{\tablename\ \thetable\ -- 38 points parameters -- \textit{continued}} \\
	Parameter & log & 1/r & $1/r^2$ \\
	\hline\vspace*{-2.2ex}
	\endhead
	$a$ & 0.8039422032494780264 & 0.8031706352420300965 & 0.8024795013067287797 \\
	$b$ & 0.4583733204758793321 & 0.4581934070811996284 & 0.4577217788720220286 \\
	$c$ & 0.1721603867747475720 & 0.1698720444743110028 & 0.1675495508816995084 \\
	\hline\Tstrut
	$energy$ & -174.8801971518150640 & 593.0385035664514014 & 597.7394530832607790
\end{longtable}

All 3 parameters have 50,014 digits precision, for all 3 potentials, but the degree of the algebraic polynomials are $>360$.

\noindent
\textit{\textbf{Symmetries -}}

The symmetry groups for 38 points are identical for all 3 potentials.

\begin{center}
	\begin{tabular}{l|l}
		\multicolumn{2}{c}{Symmetries - 38 points} \\
		\hline\Tstrut
		planes & [[4, 540], [56, 6], [120, 1]] \\[0.2ex]
		\hline\Tstrut
		Gram groups & [[2, 1], [12, 3], [24, 27], [38, 1], [48, 15]] \\
		\hline\Tstrut
		Polygons & [[4, 540], [6, 6], [8, 6]]
	\end{tabular}
\end{center}

\subsection{39 points}
The optimal solution for 39 points contains an interesting mix of parallel-plane polygons. A suitable arrangement consists of two triangles near a north pole, then a hexagon, then another triangle, then a hexagon at the equator, then another triangle, a hexagon, then two more triangles near the south pole. The arrangement of polygons is [3:3:6:3:6:3:6:3:3] which is balanced.

\begin{figure}[ht]
	\begin{center}
		\includegraphics[type=pdf,ext=pdf,read=pdf,height=1in,width=1in,angle=0]{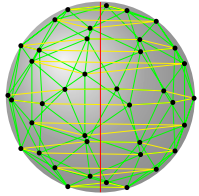}
		\caption{39 points.}
		\label{fig:39pts}
	\end{center}
\end{figure}

\noindent
\textit{\textbf{Parameterized Structure -}}

10 algebraic parameters are needed to constrain the 39 point minimal solution which is given below:

\begin{longtable}[c]{r|ccc}
	\caption{Parameterization for 39 points} \\
	pt & $x$ & $y$ & $z$ \\
	\hline\vspace*{-2.2ex}
	\endfirsthead
	\multicolumn{4}{c}%
	{\tablename\ \thetable\ -- 39 points parameters -- \textit{cont.}} \\
	pt & $x$ & $y$ & $z$ \\
	\hline\vspace*{-2.2ex}
	\endhead
	1 & $\sqrt{1-a^2}$ & $0$ & $a$ \\[0.7ex]
	2 & $-\frac{\sqrt{1-a^2}}{2}$ & $\frac{\sqrt{3-3a^2}}{2}$ & $a$ \\[0.7ex]
	3 & $-\frac{\sqrt{1-a^2}}{2}$ & $-\frac{\sqrt{3-3a^2}}{2}$ & $a$ \\[0.7ex]
	4 & $-\sqrt{1-b^2}$ & $0$ & $b$ \\[0.7ex]
	5 & $\frac{\sqrt{1-b^2}}{2}$ & $\frac{\sqrt{3-3b^2}}{2}$ & $b$ \\[0.7ex]
	6 & $\frac{\sqrt{1-b^2}}{2}$ & $-\frac{\sqrt{3-3b^2}}{2}$ & $b$ \\[0.7ex]
	7 & $\sqrt{1-c^2-d^2}$ & $d$ & $c$ \\[0.7ex]
	8 & $-\sqrt{1-c^2-e^2}$ & $e$ & $c$ \\[0.7ex]
	9 & $-\sqrt{1-c^2-f^2}$ & $f$ & $c$ \\[0.7ex]
	10 & $-\sqrt{1-c^2-f^2}$ & $-f$ & $c$ \\[0.7ex]
	11 & $-\sqrt{1-c^2-e^2}$ & $-e$ & $c$ \\[0.7ex]
	12 & $\sqrt{1-c^2-d^2}$ & $-d$ & $c$ \\[0.7ex]
	13 & $-\sqrt{1-g^2}$ & $0$ & $g$ \\[0.7ex]
	14 & $\frac{\sqrt{1-g^2}}{2}$ & $\frac{\sqrt{3-3g^2}}{2}$ & $g$ \\[0.7ex]
	15 & $\frac{\sqrt{1-g^2}}{2}$ & $-\frac{\sqrt{3-3g^2}}{2}$ & $g$ \\[0.7ex]
	16 & $1$ & $0$ & $0$ \\[0.7ex]
	17 & $h$ & $\sqrt{1-h^2}$ & $0$ \\[0.7ex]
	18 & $i$ & $\sqrt{1-i^2}$ & $0$ \\[0.7ex]
	19 & $-\frac{1}{2}$ & $\frac{\sqrt{3}}{2}$ & $0$ \\[0.7ex]
	20 & $-j$ & $\sqrt{1-j^2}$ & $0$ \\[0.7ex]
	21 & $-j$ & $-\sqrt{1-j^2}$ & $0$ \\[0.7ex]
	22 & $-\frac{1}{2}$ & $-\frac{\sqrt{3}}{2}$ & $0$ \\[0.7ex]
	23 & $i$ & $-\sqrt{1-i^2}$ & $0$ \\[0.7ex]
	24 & $h$ & $-\sqrt{1-h^2}$ & $0$ \\[0.7ex]
	25 & $-\sqrt{1-g^2}$ & $0$ & $-g$ \\[0.7ex]
	26 & $\frac{\sqrt{1-g^2}}{2}$ & $\frac{\sqrt{3-3g^2}}{2}$ & $-g$ \\[0.7ex]
	27 & $\frac{\sqrt{1-g^2}}{2}$ & $-\frac{\sqrt{3-3g^2}}{2}$ & $-g$ \\[0.7ex]
	28 & $\sqrt{1-c^2-d^2}$ & $d$ & $-c$ \\[0.7ex]
	29 & $-\sqrt{1-c^2-e^2}$ & $e$ & $-c$ \\[0.7ex]
	30 & $-\sqrt{1-c^2-f^2}$ & $f$ & $-c$ \\[0.7ex]
	31 & $-\sqrt{1-c^2-f^2}$ & $-f$ & $-c$ \\[0.7ex]
	32 & $-\sqrt{1-c^2-e^2}$ & $-e$ & $-c$ \\[0.7ex]
	33 & $\sqrt{1-c^2-d^2}$ & $-d$ & $-c$ \\[0.7ex]
	34 & $-\sqrt{1-b^2}$ & $0$ & $-b$ \\[0.7ex]
	35 & $\frac{\sqrt{1-b^2}}{2}$ & $\frac{\sqrt{3-3b^2}}{2}$ & $-b$ \\[0.7ex]
	36 & $\frac{\sqrt{1-b^2}}{2}$ & $-\frac{\sqrt{3-3b^2}}{2}$ & $-b$ \\[0.7ex]
	37 & $-\frac{\sqrt{1-a^2}}{2}$ & $-\frac{\sqrt{3-3a^2}}{2}$ & $-a$ \\[0.7ex]
	38 & $\sqrt{1-a^2}$ & $0$ & $-a$ \\[0.7ex]
	39 & $-\frac{\sqrt{1-a^2}}{2}$ & $\frac{\sqrt{3-3a^2}}{2}$ & $-a$
\end{longtable}

\noindent
\textit{\textbf{Parameter and Minimal Energy values -}}

The Jacobian matrix, Newton method search has found 50,014 digits of representation for the 10 algebraic parameters $a-j$, but unfortunately the algebraic degree is too high, $>360$ to be able to recover the polynomials for the spherical code.

\begin{longtable}[c]{c|c|c|c}
	\caption{Parameter values for 39 points} \\
	Parameter & log & 1/r & $1/r^2$ \\
	\hline\vspace*{-2.2ex}
	\endfirsthead
	\multicolumn{4}{c}%
	{\tablename\ \thetable\ -- 39 points parameters -- \textit{continued}} \\
	Parameter & log & 1/r & $1/r^2$ \\
	\hline\vspace*{-2.2ex}
	\endhead
	$a$ & 0.9305304495340633877 & 0.9303111905320777064 & 0.9300692717308784119 \\
	$b$ & 0.7781165429733693961 & 0.7777780355447797229 & 0.7773081306470802039 \\
	$c$ & 0.5412896333641035100 & 0.5405636016354160933 & 0.5398022581706319306 \\
	$d$ & 0.2810489018389595961 & 0.2832926783746533173 & 0.2852806007277812605 \\
	$e$ & 0.8268281633462391072 & 0.8276872914815004255 & 0.8285130113618658314 \\
	$f$ & 0.5457792615072795111 & 0.5443946131068471082 & 0.5432324106340845709 \\
	$g$ & 0.3244932998521861406 & 0.3231919362048333817 & 0.3218062684348973968 \\
	$h$ & 0.8481058446013102035 & 0.8475580098711930200 & 0.8470526400321745301 \\
	$i$ & 0.03479057745151933382 & 0.03582306743504348819 & 0.03677398278730272960 \\
	$j$ & 0.8828964220528295373 & 0.8833810773062365082 & 0.8838266228194772597 \\
	\hline\Tstrut
	$energy$ & -183.5092257118552531 & 626.3890090168230113 & 634.4153333806804997
\end{longtable}

\noindent
\textit{\textbf{Symmetries -}}

The symmetry groups for 39 points are identical for all 3 potentials.

\begin{center}
	\begin{tabular}{l|l}
		\multicolumn{2}{c}{Symmetries - 39 points} \\
		\hline\Tstrut
		planes & [[4, 384], [24, 3], [35, 3], [130, 1]] \\[0.2ex]
		\hline\Tstrut
		Gram groups & [[6, 6], [12, 25], [18, 1], [24, 47], [39, 1]] \\
		\hline\Tstrut
		Polygons & [[3, 6], [4, 402], [6, 2], [7, 3], [9, 1]]
	\end{tabular}
\end{center}

\subsection{40 points}

The optimal solution for 40 points contains embedded polygons, in fact 6 octagons and 12 hexagons, which led to finding an optimal axis for orientation. The arrangement is a balanced [2:2:4:2:4:4:2:2:4:4:2:4:2:2] where 2 denotes a dipole and 4 a square.

\begin{figure}[ht]
	\begin{center}
		\includegraphics[type=pdf,ext=pdf,read=pdf,height=1in,width=1in,angle=0]{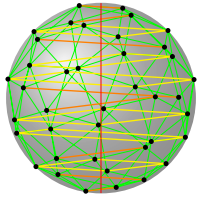}
		\caption{40 points.}
		\label{fig:40pts}
	\end{center}
\end{figure}

\noindent
\textit{\textbf{Parametric Structure -}}

The arrangement was parameterized using 9 algebraic parameters, $a$, $b$, $c$, $d$, $e$, $f$, $g$, $h$, and $i$. These parameters deal with the spacing between the parallel polygons.

\begin{longtable}[c]{r|ccc}
	\caption{Parameterization for 40 points} \\
	pt & $x$ & $y$ & $z$ \\
	\hline\vspace*{-2.2ex}
	\endfirsthead
	\multicolumn{4}{c}%
	{\tablename\ \thetable\ -- 40 points parameters -- \textit{continued\ldots}} \\
	pt & $x$ & $y$ & $z$ \\
	\hline\vspace*{-2.2ex}
	\endhead
	1 & $\sqrt{1-a^2}$ & $0$ & $a$ \\[0.5ex]
	2 & $-\sqrt{1-a^2}$ & $0$ & $a$ \\[0.5ex]
	3 & $0$ & $\sqrt{1-b^2}$ & $b$ \\[0.5ex]
	4 & $0$ & $-\sqrt{1-b^2}$ & $b$ \\[0.5ex]
	5 & $d$ & $\sqrt{1-c^2-d^2}$ & $c$ \\[0.5ex]
	6 & $-d$ & $\sqrt{1-c^2-d^2}$ & $c$ \\[0.5ex]
	7 & $-d$ & $-\sqrt{1-c^2-d^2}$ & $c$ \\[0.5ex]
	8 & $d$ & $-\sqrt{1-c^2-d^2}$ & $c$ \\[0.5ex]
	9 & $\sqrt{\frac{2}{3}}$ & $0$ & $\sqrt{\frac{1}{3}}$ \\[0.5ex]
	10 & $-\sqrt{\frac{2}{3}}$ & $0$ & $\sqrt{\frac{1}{3}}$ \\[0.5ex]
	11 & $f$ & $\sqrt{1-e^2-f^2}$ & $e$ \\[0.5ex]
	12 & $-f$ & $\sqrt{1-e^2-f^2}$ & $e$ \\[0.5ex]
	13 & $-f$ & $-\sqrt{1-e^2-f^2}$ & $e$ \\[0.5ex]
	14 & $f$ & $-\sqrt{1-e^2-f^2}$ & $e$ \\[0.5ex]
	15 & $h$ & $\sqrt{1-g^2-h^2}$ & $g$ \\[0.5ex]
	16 & $-h$ & $\sqrt{1-g^2-h^2}$ & $g$ \\[0.5ex]
	17 & $-h$ & $-\sqrt{1-g^2-h^2}$ & $g$ \\[0.5ex]
	18 & $h$ & $-\sqrt{1-g^2-h^2}$ & $g$ \\[0.5ex]
	19 & $\sqrt{1-i^2}$ & $0$ & $i$ \\[0.5ex]
	20 & $-\sqrt{1-i^2}$ & $0$ & $i$ \\[0.5ex]
	21 & $0$ & $-\sqrt{1-i^2}$ & $-i$ \\[0.5ex]
	22 & $0$ & $\sqrt{1-i^2}$ & $-i$ \\[0.5ex]
	23 & $\sqrt{1-g^2-h^2}$ & $h$ & $-g$ \\[0.5ex]
	24 & $-\sqrt{1-g^2-h^2}$ & $h$ & $-g$ \\[0.5ex]
	25 & $-\sqrt{1-g^2-h^2}$ & $-h$ & $-g$ \\[0.5ex]
	26 & $\sqrt{1-g^2-h^2}$ & $-h$ & $-g$ \\[0.5ex]
	27 & $\sqrt{1-e^2-f^2}$ & $f$ & $-e$ \\[0.5ex]
	28 & $-\sqrt{1-e^2-f^2}$ & $f$ & $-e$ \\[0.5ex]
	29 & $-\sqrt{1-e^2-f^2}$ & $-f$ & $-e$ \\[0.5ex]
	30 & $\sqrt{1-e^2-f^2}$ & $-f$ & $-e$ \\[0.5ex]
	31 & $0$ & $-\sqrt{\frac{2}{3}}$ & $-\sqrt{\frac{1}{3}}$ \\[0.5ex]
	32 & $0$ & $\sqrt{\frac{2}{3}}$ & $-\sqrt{\frac{1}{3}}$ \\[0.5ex]
	33 & $\sqrt{1-c^2-d^2}$ & $d$ & $-c$ \\[0.5ex]
	34 & $-\sqrt{1-c^2-d^2}$ & $d$ & $-c$ \\[0.5ex]
	35 & $-\sqrt{1-c^2-d^2}$ & $-d$ & $-c$ \\[0.5ex]
	36 & $\sqrt{1-c^2-d^2}$ & $-d$ & $-c$ \\[0.5ex]
	37 & $\sqrt{1-b^2}$ & $0$ & $-b$ \\[0.5ex]
	38 & $-\sqrt{1-b^2}$ & $0$ & $-b$ \\[0.5ex]
	39 & $0$ & $-\sqrt{1-a^2}$ & $-a$ \\[0.5ex]
	40 & $0$ & $\sqrt{1-a^2}$ & $-a$
\end{longtable}

\noindent
\textit{\textbf{logarithmic potential -}}

The decimal values for the 9 parameters have been determined to 50,014 digits. However the algebraic degree is $>360$ for this potential.

\begin{longtable}[l]{l}
	$ \qquad a = 0.9608452414555139834\ldots$ \\
	$ \qquad b = 0.8325617222673873581\ldots$ \\
	$ \qquad c = 0.7059050981400227835\ldots$ \\
	$ \qquad d = 0.5403569853489025648\ldots$ \\
	$ \qquad e = 0.3916892764776448136\ldots$ \\
	$ \qquad f = 0.3117438960562255168\ldots$ \\
	$ \qquad g = 0.1959290968338993704\ldots$ \\
	$ \qquad h = 0.8179629789070255120\ldots$ \\
	$ \qquad i = 0.05827507906343504122\ldots$ \\
	$\text{and for the energy:}$ \\
	$energy = -192.3376899173482332\ldots$
\end{longtable}

\noindent
\textit{\textbf{Coulomb $\mathbf{1/r}$ potential -}}

Similarly the 9 parameters have been located for the Coulomb potential to 50,014 digits.

\begin{tabular}{l}
	$ \qquad a = 0.9601769655764723164\ldots$ \\
	$ \qquad b = 0.8329158480857483780\ldots$ \\
	$ \qquad c = 0.7058998933400944309\ldots$ \\
	$ \qquad d = 0.5404423709662100308\ldots$ \\
	$ \qquad e = 0.3913126563296915203\ldots$ \\
	$ \qquad f = 0.3122606114843310841\ldots$ \\
	$ \qquad g = 0.1975603639098160504\ldots$ \\
	$ \qquad h = 0.8186439165125587834\ldots$ \\
	$ \qquad i = 0.05840103736139121117\ldots$ \\
	$ \text{and for the energy:}$ \\
	$ \qquad energy = 660.6752788346224138\ldots$
\end{tabular}

and the algebraic degree for the spherical code is $>360$.

\noindent
\textit{\textbf{Inverse square law $\mathbf{1/r^2}$ potential -}}

The 9 parameters also have been determined for the Inverse square law potential:

\begin{tabular}{l}
	$ \qquad a = 0.9596045207546853771\ldots$ \\
	$ \qquad b = 0.8332509745839302446\ldots$ \\
	$ \qquad c = 0.7058970333118211098\ldots$ \\
	$ \qquad d = 0.5404892082486678872\ldots$ \\
	$ \qquad e = 0.3909557605119346706\ldots$ \\
	$ \qquad f = 0.3127499451566637698\ldots$ \\
	$ \qquad g = 0.1989461783337023559\ldots$ \\
	$ \qquad h = 0.8192190556738142491\ldots$ \\
	$ \qquad i = 0.05847013530974113222\ldots$ \\
	$\text{and for the energy:}$ \\
	$ \qquad energy = 672.3093535034931164\ldots$
\end{tabular}

Once again, the parameters have been determined to 50,014 digits. Searching shows that the polynomials must have degree $>360$, if they are found.

\noindent
\textit{\textbf{Symmetries -}}

The symmetry groups for 40 points are identical for all 3 potentials.

\begin{center}
	\begin{tabular}{l|l}
		\multicolumn{2}{c}{Symmetries - 40 points} \\
		\hline\Tstrut
		planes & [[4, 540], [24, 3], [56, 6], [67, 4]] \\[0.2ex]
		\hline\Tstrut
		Gram groups & [[12, 4], [24, 21], [40, 1], [48, 21]] \\
		\hline\Tstrut
		Polygons & [[3, 28], [4, 558], [6, 12], [8, 6]]
	\end{tabular}
\end{center}

\subsection{41 points}
The optimal solution for 41 points contains embedded polygons, as shown by the symmetry groups. An suitable alignment was found which consisted of 2 poles, and starting at one pole, moving to the opposite, the following polygons: [3:3:6:3:9:3:6:3:3] which is a balanced arrangement.

\begin{figure}[ht]
	\begin{center}
		\includegraphics[type=pdf,ext=pdf,read=pdf,height=1in,width=1in,angle=0]{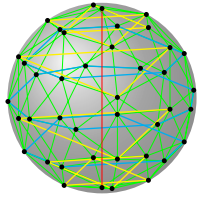}
		\caption{41 points.}
		\label{fig:41pts}
	\end{center}
\end{figure}

\noindent
\textit{\textbf{Parameterized Structure -}}

Again, 10 algebraic parameters are needed to constrain the 41 point minimal solution which is given below:

\begin{longtable}[c]{r|ccc}
	\caption{Parameterization for 41 points} \\
	pt & $x$ & $y$ & $z$ \\
	\hline\vspace*{-2.2ex}
	\endfirsthead
	\multicolumn{4}{c}%
	{\tablename\ \thetable\ -- 41 points parameters -- \textit{continued\ldots}} \\
	pt & $x$ & $y$ & $z$ \\
	\hline\vspace*{-2.2ex}
	\endhead
	1 & $0$ & $0$ & $1$ \\[0.7ex]
	2 & $\frac{\sqrt{3-3a^2}}{2}$ & $\frac{\sqrt{1-a^2}}{2}$ & $a$ \\[0.7ex]
	3 & $0$ & $-\sqrt{1-a^2}$ & $a$ \\[0.7ex]
	4 & $-\frac{\sqrt{3-3a^2}}{2}$ & $\frac{\sqrt{1-a^2}}{2}$ & $a$ \\[0.7ex]
	5 & $\frac{\sqrt{3-3b^2}}{2}$ & $-\frac{\sqrt{1-b^2}}{2}$ & $b$ \\[0.7ex]
	6 & $0$ & $\sqrt{1-b^2}$ & $b$ \\[0.7ex]
	7 & $-\frac{\sqrt{3-3b^2}}{2}$ & $-\frac{\sqrt{1-b^2}}{2}$ & $b$ \\[0.7ex]
	8 & $\sqrt{1-c^2-d^2}$ & $d$ & $c$ \\[0.7ex]
	9 & $\sqrt{1-c^2-f^2}$ & $f$ & $c$ \\[0.7ex]
	10 & $-\sqrt{1-c^2-g^2}$ & $-g$ & $c$ \\[0.7ex]
	11 & $\sqrt{1-c^2-g^2}$ & $-g$ & $c$ \\[0.7ex]
	12 & $-\sqrt{1-c^2-f^2}$ & $f$ & $c$ \\[0.7ex]
	13 & $-\sqrt{1-c^2-d^2}$ & $d$ & $c$ \\[0.7ex]
	14 & $\frac{\sqrt{3-3e^2}}{2}$ & $-\frac{\sqrt{1-e^2}}{2}$ & $e$ \\[0.7ex]
	15 & $0$ & $\sqrt{1-e^2}$ & $e$ \\[0.7ex]
	16 & $-\frac{\sqrt{3-3e^2}}{2}$ & $-\frac{\sqrt{1-e^2}}{2}$ & $e$ \\[0.7ex]
	17 & $-\sqrt{1-h^2}$ & $h$ & $0$ \\[0.7ex]
	18 & $\sqrt{1-h^2}$ & $h$ & $0$ \\[0.7ex]
	19 & $-\frac{\sqrt{3}}{2}$ & $\frac{1}{2}$ & $0$ \\[0.7ex]
	20 & $\frac{\sqrt{3}}{2}$ & $\frac{1}{2}$ & $0$ \\[0.7ex]
	21 & $-i$ & $-\sqrt{1-i^2}$ & $0$ \\[0.7ex]
	22 & $i$ & $-\sqrt{1-i^2}$ & $0$ \\[0.7ex]
	23 & $-\sqrt{1-j^2}$ & $-j$ & $0$ \\[0.7ex]
	24 & $\sqrt{1-j^2}$ & $-j$ & $0$ \\[0.7ex]
	25 & $0$ & $-1$ & $0$ \\[0.7ex]
	26 & $\frac{\sqrt{3-3e^2}}{2}$ & $-\frac{\sqrt{1-e^2}}{2}$ & $-e$ \\[0.7ex]
	27 & $0$ & $\sqrt{1-e^2}$ & $-e$ \\[0.7ex]
	28 & $-\frac{\sqrt{3-3e^2}}{2}$ & $-\frac{\sqrt{1-e^2}}{2}$ & $-e$ \\[0.7ex]
	29 & $\sqrt{1-c^2-d^2}$ & $d$ & $-c$ \\[0.7ex]
	30 & $-\sqrt{1-c^2-g^2}$ & $-g$ & $-c$ \\[0.7ex]
	31 & $\sqrt{1-c^2-g^2}$ & $-g$ & $-c$ \\[0.7ex]
	32 & $\sqrt{1-c^2-f^2}$ & $f$ & $-c$ \\[0.7ex]
	33 & $-\sqrt{1-c^2-f^2}$ & $f$ & $-c$ \\[0.7ex]
	34 & $-\sqrt{1-c^2-d^2}$ & $d$ & $-c$ \\[0.7ex]
	35 & $\frac{\sqrt{3-3b^2}}{2}$ & $-\frac{\sqrt{1-b^2}}{2}$ & $-b$ \\[0.7ex]
	36 & $0$ & $\sqrt{1-b^2}$ & $-b$ \\[0.7ex]
	37 & $-\frac{\sqrt{3-3b^2}}{2}$ & $-\frac{\sqrt{1-b^2}}{2}$ & $-b$ \\[0.7ex]
	38 & $\frac{\sqrt{3-3a^2}}{2}$ & $\frac{\sqrt{1-a^2}}{2}$ & $-a$ \\[0.7ex]
	39 & $0$ & $-\sqrt{1-a^2}$ & $-a$ \\[0.7ex]
	40 & $-\frac{\sqrt{3-3a^2}}{2}$ & $\frac{\sqrt{1-a^2}}{2}$ & $-a$ \\[0.7ex]
	41 & $0$ & $0$ & $-1$
\end{longtable}

The values to 19 digits of the 10 parameters $a-j$ optimized for the minimal solutions of 41 points are:

\begin{longtable}[c]{c|c|c|c}
	\caption{Parameter values for 41 points} \\
	Parameter & log & 1/r & $1/r^2$ \\
	\hline\vspace*{-2.2ex}
	\endfirsthead
	\multicolumn{4}{c}%
	{\tablename\ \thetable\ -- 41 points parameters -- \textit{continued}} \\
	Parameter & log & 1/r & $1/r^2$ \\
	\hline\vspace*{-2.2ex}
	\endhead
	$a$ & 0.8520984319699034101 & 0.8523862771111461391 & 0.8526095822555335960 \\
	$b$ & 0.7875550123025952342 & 0.7889333123109460737 & 0.7904592194498969596 \\
	$c$ & 0.4909732570000036406 & 0.4921990699246269454 & 0.4932559558672481719 \\
	$d$ & 0.1424411984296669891 & 0.1411189147409479565 & 0.1400379041853693472 \\
	$e$ & 0.3472293650970637960 & 0.3453085908033527925 & 0.3443319566507469829 \\
	$f$ & 0.6730856713657296369 & 0.6733284008452331171 & 0.6734969929580199912 \\
	$g$ & 0.8155268697953966259 & 0.8144473155861810736 & 0.8135348971433893384 \\
	$h$ & 0.9002103658277210421 & 0.8994108636176339581 & 0.8988402363895471274 \\
	$i$ & 0.9973326834765621955 & 0.9974647659120622585 & 0.9975566938586787781 \\
	$j$ & 0.8272205140026939919 & 0.8282487897018546242 & 0.8289786408041167464 \\
	\hline\Tstrut
	$energy$ & -201.3592066486634101 & 695.9167443418870270 & 711.5261514819427490
\end{longtable}

All 10 parameters have 50,014 digits precision, for all 3 potentials, but the degree of the algebraic polynomials are $>360$.

\noindent
\textit{\textbf{Minimal Energy values -}}

The decimal values of the 10 parameters, $a-j$, and $energy$ are known to 50,014 digits, nevertheless, the algebraic polynomials were not obtained, they are $>360$ in degree.

The values are given to 19 digits below:
\begin{center}
	\begin{tabular}{l|l}
		\multicolumn{2}{c}{Minimal Energy - 41 points} \\
		\hline\Tstrut
		logarithmic & -201.3592066486634101\ldots \\[0.2ex]
		\hline\Tstrut
		Coulomb & 695.9167443418870270\ldots \\[0.2ex]
		\hline\Tstrut
		Inverse square law & 711.5261514819427490\ldots
	\end{tabular}
\end{center}

\noindent
\textit{\textbf{Symmetries -}}

The symmetry groups for 41 points are identical under all 3 potentials.

\begin{center}
	\begin{tabular}{l|l}
		\multicolumn{2}{c}{Symmetries - 41 points} \\
		\hline\Tstrut
		planes & [[4, 390], [24, 3], [84, 3], [130, 1]] \\[0.2ex]
		\hline\Tstrut
		Gram groups & [[2, 1], [6, 6], [12, 31], [18, 1], [24, 49], [36, 1], [41, 1]] \\
		\hline\Tstrut
		Polygons & [[3, 6], [4, 408], [6, 2], [9, 4]]
	\end{tabular}
\end{center}

\subsection{42 points}
The optimal solution for 42 points is very interesting, it contains 2 poles and 6 pentagons and an decagon on the equator. The arrangement is a balanced [1:5:5:5:10:5:5:5:1] where 1 is a pole, 5 a pentagon and 10 a decagon.

\begin{figure}[ht]
	\begin{center}
		\includegraphics[type=pdf,ext=pdf,read=pdf,height=1in,width=1in,angle=0]{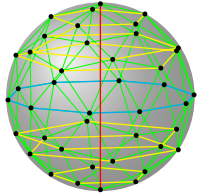}
		\caption{42 points.}
		\label{fig:42pts}
	\end{center}
\end{figure}

\noindent
\textit{\textbf{Parametric Structure -}}

\noindent
Let $c_1$, $c_2$, $s_1$, and $s_2$ be parameters for a pentagon.
\begin{align*}
c_1 = & \frac{\sqrt{5}-1}{4} \\
c_2 = & \frac{\sqrt{5}+1}{4} \\
s_1 = & \frac{\sqrt{10+2\sqrt{5}}}{4} \\
s_2 = & \frac{\sqrt{10-2\sqrt{5}}}{4}
\end{align*}
\noindent
Then using these pentagonal parameters, we find the structure for 42 points:

\begin{longtable}[c]{r|ccc}
	\caption{Parameterization for 42 points} \\
	pt & $x$ & $y$ & $z$ \\
	\hline\vspace*{-2.2ex}
	\endfirsthead
	\multicolumn{4}{c}%
	{\tablename\ \thetable\ -- 42 points parameters -- \textit{continued\ldots}} \\
	pt & $x$ & $y$ & $z$ \\
	\hline\vspace*{-2.2ex}
	\endhead
	1 & $0$ & $0$ & $1$ \\[0.7ex]
	2 & $\sqrt{1-a^2}$ & $0$ & $a$ \\[0.7ex]
	3 & $c_1\sqrt{1-a^2}$ & $s_1\sqrt{1-a^2}$ & $a$ \\[0.7ex]
	4 & $-c_2\sqrt{1-a^2}$ & $s_2\sqrt{1-a^2}$ & $a$ \\[0.7ex]
	5 & $-c_2\sqrt{1-a^2}$ & $-s_2\sqrt{1-a^2}$ & $a$ \\[0.7ex]
	6 & $c_1\sqrt{1-a^2}$ & $-s_1\sqrt{1-a^2}$ & $a$ \\[0.7ex]
	7 & $-\sqrt{1-b^2}$ & $0$ & $b$ \\[0.7ex]
	8 & $-c_1\sqrt{1-b^2}$ & $s_1\sqrt{1-b^2}$ & $b$ \\[0.7ex]
	9 & $c_2\sqrt{1-b^2}$ & $s_2\sqrt{1-b^2}$ & $b$ \\[0.7ex]
	10 & $c_2\sqrt{1-b^2}$ & $-s_2\sqrt{1-b^2}$ & $b$ \\[0.7ex]
	11 & $-c_1\sqrt{1-b^2}$ & $-s_1\sqrt{1-b^2}$ & $b$ \\[0.7ex]
	12 & $\sqrt{1-c^2}$ & $0$ & $c$ \\[0.7ex]
	13 & $c_1\sqrt{1-c^2}$ & $s_1\sqrt{1-c^2}$ & $c$ \\[0.7ex]
	14 & $-c_2\sqrt{1-c^2}$ & $s_2\sqrt{1-c^2}$ & $c$ \\[0.7ex]
	15 & $-c_2\sqrt{1-c^2}$ & $-s_2\sqrt{1-c^2}$ & $c$ \\[0.7ex]
	16 & $c_1\sqrt{1-c^2}$ & $-s_1\sqrt{1-c^2}$ & $c$ \\[0.7ex]
	17 & $d$ & $-\sqrt{1-d^2}$ & $0$ \\[0.7ex]
	18 & $e$ & $-\sqrt{1-e^2}$ & $0$ \\[0.7ex]
	19 & $\frac{\sqrt{2\sqrt{5}+10}\sqrt{1-d^2}}{4}+\frac{d\left(\sqrt{5}-1\right)}{4}$ & $\frac{d\sqrt{2\sqrt{5}+10}}{4}-\frac{\left(\sqrt{5}-1\right)\sqrt{1-d^2}}{4}$ & $0$ \\[0.7ex]
	20 & $\frac{\sqrt{2\sqrt{5}+10}\sqrt{1-e^2}}{4}+\frac{e\left(\sqrt{5}-1\right)}{4}$ & $\frac{e\sqrt{2\sqrt{5}+10}}{4}-\frac{\left(\sqrt{5}-1\right)\sqrt{1-e^2}}{4}$ & $0$ \\[0.7ex]
	21 & $\frac{\sqrt{10-2\sqrt{5}}\sqrt{1-d^2}}{4}-\frac{d\left(\sqrt{5}+1\right)}{4}$ & $\frac{\left(\sqrt{5}+1\right)\sqrt{1-d^2}}{4}+\frac{d\sqrt{10-2\sqrt{5}}}{4}$ & $0$ \\[0.7ex]
	22 & $\frac{\sqrt{10-2\sqrt{5}}\sqrt{1-e^2}}{4}-\frac{e\left(\sqrt{5}+1\right)}{4}$ & $\frac{\left(\sqrt{5}+1\right)\sqrt{1-e^2}}{4}+\frac{e\sqrt{10-2\sqrt{5}}}{4}$ & $0$ \\[0.7ex]
	23 & $\left(-\frac{\sqrt{10-2\sqrt{5}}\sqrt{1-d^2}}{4}\right)-\frac{d\left(\sqrt{5}+1\right)}{4}$ & $\frac{\left(\sqrt{5}+1\right)\sqrt{1-d^2}}{4}-\frac{d\sqrt{10-2\sqrt{5}}}{4}$ & $0$ \\[0.7ex]
	24 & $\left(-\frac{\sqrt{10-2\sqrt{5}}\sqrt{1-e^2}}{4}\right)-\frac{e\left(\sqrt{5}+1\right)}{4}$ & $\frac{\left(\sqrt{5}+1\right)\sqrt{1-e^2}}{4}-\frac{e\sqrt{10-2\sqrt{5}}}{4}$ & $0$ \\[0.7ex]
	25 & $\frac{d\left(\sqrt{5}-1\right)}{4}-\frac{\sqrt{2\sqrt{5}+10}\sqrt{1-d^2}}{4}$ & $\left(-\frac{\left(\sqrt{5}-1\right)\sqrt{1-d^2}}{4}\right)-\frac{d\sqrt{2\sqrt{5}+10}}{4}$ & $0$ \\[0.7ex]
	26 & $\frac{e\left(\sqrt{5}-1\right)}{4}-\frac{\sqrt{2\sqrt{5}+10}\sqrt{1-e^2}}{4}$ & $\left(-\frac{\left(\sqrt{5}-1\right)\sqrt{1-e^2}}{4}\right)-\frac{e\sqrt{2\sqrt{5}+10}}{4}$ & $0$ \\[0.7ex]
	27 & $\sqrt{1-c^2}$ & $0$ & $-c$ \\[0.7ex]
	28 & $c_1\sqrt{1-c^2}$ & $s_1\sqrt{1-c^2}$ & $-c$ \\[0.7ex]
	29 & $-c_2\sqrt{1-c^2}$ & $s_2\sqrt{1-c^2}$ & $-c$ \\[0.7ex]
	30 & $-c_2\sqrt{1-c^2}$ & $-s_2\sqrt{1-c^2}$ & $-c$ \\[0.7ex]
	31 & $c_1\sqrt{1-c^2}$ & $-s_1\sqrt{1-c^2}$ & $-c$ \\[0.7ex]
	32 & $-\sqrt{1-b^2}$ & $0$ & $-b$ \\[0.7ex]
	33 & $-c_1\sqrt{1-b^2}$ & $s_1\sqrt{1-b^2}$ & $-b$ \\[0.7ex]
	34 & $c_2\sqrt{1-b^2}$ & $s_2\sqrt{1-b^2}$ & $-b$ \\[0.7ex]
	35 & $c_2\sqrt{1-b^2}$ & $-s_2\sqrt{1-b^2}$ & $-b$ \\[0.7ex]
	36 & $-c_1\sqrt{1-b^2}$ & $-s_1\sqrt{1-b^2}$ & $-b$ \\[0.7ex]
	37 & $\sqrt{1-a^2}$ & $0$ & $-a$ \\[0.7ex]
	38 & $c_1\sqrt{1-a^2}$ & $s_1\sqrt{1-a^2}$ & $-a$ \\[0.7ex]
	39 & $-c_2\sqrt{1-a^2}$ & $s_2\sqrt{1-a^2}$ & $-a$ \\[0.7ex]
	40 & $-c_2\sqrt{1-a^2}$ & $-s_2\sqrt{1-a^2}$ & $-a$ \\[0.7ex]
	41 & $c_1\sqrt{1-a^2}$ & $-s_1\sqrt{1-a^2}$ & $-a$ \\[0.7ex]
	42 & $0$ & $0$ & $-1$
\end{longtable}

\noindent
\textit{\textbf{logarithmic potential -}}

The decimal values for the parameters have been determined to 50,014 digits. However the algebraic degree is $>360$ for this potential.

\begin{longtable}[l]{l}
	$ \qquad a = 0.8480021192047263213\ldots$ \\
	$ \qquad b = 0.5500987140683140002\ldots$ \\
	$ \qquad c = 0.4199535968733765594\ldots$ \\
	$ \qquad d = 0.9367127737987912556\ldots$ \\
	$ \qquad e = 0.6224239241337443283\ldots$ \\
	$\text{and for the energy:}$ \\
	$ \qquad energy = -210.5845115576339070\ldots$
\end{longtable}

\noindent
\textit{\textbf{Coulomb $\mathbf{1/r}$ potential -}}

Similarly the parameters have been located for the Coulomb potential to 50,014 digits accuracy.

\begin{tabular}{l}
	$ \qquad a = 0.8479934351306249820\ldots$ \\
	$ \qquad b = 0.5475481419910211780\ldots$ \\
	$ \qquad c = 0.4217051127686561307\ldots$ \\
	$ \qquad d = 0.9374180220888377937\ldots$ \\
	$ \qquad e = 0.6208417306215908095\ldots$ \\
	$\text{and for the energy:}$ \\
	$ \qquad energy = 732.0781075436733758\ldots$
\end{tabular}

and the algebraic degree for the spherical code is $>360$.

\noindent
\textit{\textbf{Inverse square law $\mathbf{1/r^2}$ potential -}}

The parameters also have been determined for the Inverse square law potential:

\begin{tabular}{l}
	$ \qquad a = 0.8478870718441163396\ldots$ \\
	$ \qquad b = 0.5453554276036246459\ldots$ \\
	$ \qquad c = 0.4231848715206929912\ldots$ \\
	$ \qquad d = 0.9379732112146647173\ldots$ \\
	$ \qquad e = 0.6195883152605443173\ldots$ \\
	$\text{and for the energy:}$ \\
	$ \qquad energy = 751.8751968168232611\ldots$
\end{tabular}

While the decimal precision of these parameter values is accurately known to 50,014 digits, the algebraic polynomials were not recovered. The polynomials must have degree $>360$, if they are found.

\noindent
\textit{\textbf{Symmetries -}}

The symmetry groups for 42 points are identical under all 3 potentials.

\begin{center}
	\begin{tabular}{l|l}
		\multicolumn{2}{c}{Symmetries - 42 points} \\
		\hline\Tstrut
		planes & [[4, 630], [24, 5], [56, 5], [180, 1]] \\[0.2ex]
		\hline\Tstrut
		Gram groups & [[2, 1], [10, 8], [20, 26], [40, 28], [42, 1]] \\
		\hline\Tstrut
		Polygons & [[4, 660], [5, 6], [8, 5], [10, 1]]
	\end{tabular}
\end{center}

\subsection{43 points}
The optimal arrangement of 43 points has an unusual arrangement of 6 single points and 14 opposing points, all opposing points are parallel to a coordinate axis, and a nonagon, in a [1:2:2:2:1:2:1:2:2:2:9:2:2:2:1:2:1:2:2:2:1] arrangement with a nonagon in the center and 3 single points and 7 opposing points on either side of the nonagon.

While there are embedded squares and one heptagon in the polyhedra, the arrangement with the normal axis through the nonagon was chosen for the preferred alignment. The nonagon is highlighted below in figure \ref{fig:43pts} in yellow. The 14 opposing points are coupled in orange.

\begin{figure}[ht]
	\begin{center}
		\includegraphics[type=pdf,ext=pdf,read=pdf,height=1in,width=1in,angle=0]{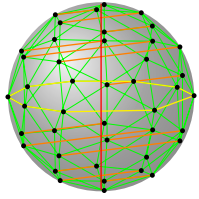}
		\caption{43 points.}
		\label{fig:43pts}
	\end{center}
\end{figure}

\noindent
\textit{\textbf{Algebraic parameterization -}}

It requires 21 parameters $a-u$ to adequately constrain the 43 point set, due to the 14 opposing points and 6 single points in the configuration.

The algebraic parameterized structure is given below:

\begin{longtable}[c]{r|ccc}
	\caption{Parameterization for 43 points} \\
	pt & $x$ & $y$ & $z$ \\
	\hline\vspace*{-2.2ex}
	\endfirsthead
	\multicolumn{4}{c}%
	{\tablename\ \thetable\ -- 43 points parameters -- \textit{continued}} \\
	pt & $x$ & $y$ & $z$ \\
	\hline\vspace*{-2.2ex}
	\endhead
	1 & $\sqrt{1-a^2}$ & $0$ & $a$ \\[0.5ex]
	2 & $-\sqrt{1-b^2-c^2}$ & $c$ & $b$ \\[0.5ex]
	3 & $-\sqrt{1-b^2-c^2}$ & $-c$ & $b$ \\[0.5ex]
	4 & $-\sqrt{1-d^2-e^2}$ & $e$ & $d$ \\[0.5ex]
	5 & $-\sqrt{1-d^2-e^2}$ & $-e$ & $d$ \\[0.5ex]
	6 & $\sqrt{1-f^2-g^2}$ & $g$ & $f$ \\[0.5ex]
	7 & $\sqrt{1-f^2-g^2}$ & $-g$ & $f$ \\[0.5ex]
	8 & $\sqrt{1-h^2}$ & $0$ & $h$ \\[0.5ex]
	9 & $\sqrt{1-i^2-j^2}$ & $j$ & $i$ \\[0.5ex]
	10 & $\sqrt{1-i^2-j^2}$ & $-j$ & $i$ \\[0.5ex]
	11 & $-\sqrt{1-k^2}$ & $0$ & $k$ \\[0.5ex]
	12 & $-\sqrt{1-l^2-m^2}$ & $m$ & $l$ \\[0.5ex]
	13 & $-\sqrt{1-l^2-m^2}$ & $-m$ & $l$ \\[0.5ex]
	14 & $-\sqrt{1-n^2-o^2}$ & $o$ & $n$ \\[0.5ex]
	15 & $-\sqrt{1-n^2-o^2}$ & $-o$ & $n$ \\[0.5ex]
	16 & $\sqrt{1-p^2-q^2}$ & $q$ & $p$ \\[0.5ex]
	17 & $\sqrt{1-p^2-q^2}$ & $-q$ & $p$ \\[0.5ex]
	18 & $1$ & $0$ & $0$ \\[0.5ex]
	19 & $r$ & $\sqrt{1-r^2}$ & $0$ \\[0.5ex]
	20 & $s$ & $\sqrt{1-s^2}$ & $0$ \\[0.5ex]
	21 & $t$ & $\sqrt{1-t^2}$ & $0$ \\[0.5ex]
	22 & $u$ & $\sqrt{1-u^2}$ & $0$ \\[0.5ex]
	23 & $u$ & $-\sqrt{1-u^2}$ & $0$ \\[0.5ex]
	24 & $t$ & $-\sqrt{1-t^2}$ & $0$ \\[0.5ex]
	25 & $s$ & $-\sqrt{1-s^2}$ & $0$ \\[0.5ex]
	26 & $r$ & $-\sqrt{1-r^2}$ & $0$ \\[0.5ex]
	27 & $\sqrt{1-p^2-q^2}$ & $q$ & $-p$ \\[0.5ex]
	28 & $\sqrt{1-p^2-q^2}$ & $-q$ & $-p$ \\[0.5ex]
	29 & $-\sqrt{1-n^2-o^2}$ & $o$ & $-n$ \\[0.5ex]
	30 & $-\sqrt{1-n^2-o^2}$ & $-o$ & $-n$ \\[0.5ex]
	31 & $-\sqrt{1-l^2-m^2}$ & $m$ & $-l$ \\[0.5ex]
	32 & $-\sqrt{1-l^2-m^2}$ & $-m$ & $-l$ \\[0.5ex]
	33 & $-\sqrt{1-k^2}$ & $0$ & $-k$ \\[0.5ex]
	34 & $\sqrt{1-i^2-j^2}$ & $j$ & $-i$ \\[0.5ex]
	35 & $\sqrt{1-i^2-j^2}$ & $-j$ & $-i$ \\[0.5ex]
	36 & $\sqrt{1-h^2}$ & $0$ & $-h$ \\[0.5ex]
	37 & $\sqrt{1-f^2-g^2}$ & $g$ & $-f$ \\[0.5ex]
	38 & $\sqrt{1-f^2-g^2}$ & $-g$ & $-f$ \\[0.5ex]
	39 & $-\sqrt{1-d^2-e^2}$ & $e$ & $-d$ \\[0.5ex]
	40 & $-\sqrt{1-d^2-e^2}$ & $-e$ & $-d$ \\[0.5ex]
	41 & $-\sqrt{1-b^2-c^2}$ & $c$ & $-b$ \\[0.5ex]
	42 & $-\sqrt{1-b^2-c^2}$ & $-c$ & $-b$ \\[0.5ex]
	43 & $\sqrt{1-a^2}$ & $0$ & $-a$
\end{longtable}

The values to 19 digits of the 21 parameters $a-u$ optimized for the minimal solutions of 43 points are:

\begin{longtable}[c]{c|c|c|c}
	\caption{Parameter values for 43 points} \\
	Parameter & log & 1/r & $1/r^2$ \\
	\hline\vspace*{-2.2ex}
	\endfirsthead
	\multicolumn{4}{c}%
	{\tablename\ \thetable\ -- 43 points parameters -- \textit{continued}} \\
	Parameter & log & 1/r & $1/r^2$ \\
	\hline\vspace*{-2.2ex}
	\endhead
$a$ & 0.9982866473863599175 & 0.9982866473863599175 & 0.9984630953851052699 \\
$b$ & 0.8472635797533589187 & 0.8472635797533589187 & 0.8475325598739136651 \\
$c$ & 0.5311293486382919021 & 0.5311293486382919021 & 0.5307272765752298213 \\
$d$ & 0.8400740066768052323 & 0.8400740066768052323 & 0.8399419725010893398 \\
$e$ & 0.2700031255728598295 & 0.2700031255728598295 & 0.2711612073230125316 \\
$f$ & 0.8011237624533662706 & 0.8011237624533662706 & 0.8013390030747537570 \\
$g$ & 0.3149686382777712023 & 0.3149686382777712023 & 0.3143338979088199771 \\
$h$ & 0.5257933545107503321 & 0.5257933545107503321 & 0.5249258057238283796 \\
$i$ & 0.5045367851340684154 & 0.5045367851340684154 & 0.5039541318183206447 \\
$j$ & 0.7922732606467253208 & 0.7922732606467253208 & 0.7934086000053364266 \\
$k$ & 0.4897424218019944371 & 0.4897424218019944371 & 0.4927453299320202459 \\
$l$ & 0.4763943861747329806 & 0.4763943861747329806 & 0.4778957726446077287 \\
$m$ & 0.8451040731358101017 & 0.8451040731358101017 & 0.8437455161174711917 \\
$n$ & 0.4429356904485061573 & 0.4429356904485061573 & 0.4412278165688409847 \\
$o$ & 0.5426636322452562795 & 0.5426636322452562795 & 0.5422887105232160723 \\
$p$ & 0.2966040392756743250 & 0.2966040392756743250 & 0.2975541321855012033 \\
$q$ & 0.4942343699679201443 & 0.4942343699679201443 & 0.4945697569720479639 \\
$r$ & 0.5177710856618404329 & 0.5177710856618404329 & 0.5195541951343300817 \\
$s$ & -0.003442294413174403063 & -0.003442294413174403063 & -0.002581058629101833997 \\
$t$ & -0.5675846094236517838 & -0.5675846094236517838 & -0.5651170132805532067 \\
$u$ & -0.9555580439704250977 & -0.9555580439704250977 & -0.9558789980236047501 \\
\hline\Tstrut
$energy$ & -220.0034770518506974 & 769.1908464591584376 & 793.5218863349545800
\end{longtable}

All 21 parameters have 50,014 digits precision, for all 3 potentials, but the degree of the algebraic polynomials are $>360$.

\noindent
\textit{\textbf{Symmetries -}}

The symmetries are identical for 43 points under all 3 potentials.

\begin{center}
	\begin{tabular}{l|l}
		\multicolumn{2}{c}{Symmetries - 43 points} \\
		\hline\Tstrut
		planes & [[4, 266], [28, 1], [35, 1], [84, 1]] \\[0.2ex]
		\hline\Tstrut
		Gram groups & [[2, 7], [4, 46], [8, 201], [43, 1]] \\
		\hline\Tstrut
		Polygons & [[4, 273], [7, 1], [9, 1]]
	\end{tabular}
\end{center}

\subsection{44 points}
The optimal configuration for 44 points contains extensive symmetries, indicating the high degree of symmetry for the polyhedra.

The figure contains 2 poles, and proceeding from the north pole, the embedded parallel polygons are [1:4:2:4:4:2:10:2:4:4:2:4:1] where 1 is a pole and 2 is a dipole. Clearly the arrangement of polygons is balanced, which is typical with an even number of points.

\begin{figure}[ht]
	\begin{center}
		\includegraphics[type=pdf,ext=pdf,read=pdf,height=1in,width=1in,angle=0]{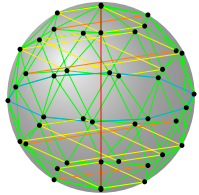}
		\caption{44 points.}
		\label{fig:44pts}
	\end{center}
\end{figure}

\noindent
\textit{\textbf{Algebraic parameterization -}}

After searching and double rotating into an optimal alignment, it became clear that the 44 point set could be parameterized by 3 algebraic parameters $a$, $b$, and $c$.

The algebraic parameterized structure is given below:

\begin{longtable}[c]{r|ccc}
	\caption{Parameterization for 44 points} \\
	pt & $x$ & $y$ & $z$ \\
	\hline\vspace*{-2.2ex}
	\endfirsthead
	\multicolumn{4}{c}%
	{\tablename\ \thetable\ -- 44 points parameters -- \textit{cont.}} \\
	pt & $x$ & $y$ & $z$ \\
	\hline\vspace*{-2.2ex}
	\endhead
	1 & $0$ & $0$ & $1$ \\[0.7ex]
	2 & $b$ & $\sqrt{1-a^2-b^2}$ & $a$ \\[0.7ex]
	3 & $-b$ & $\sqrt{1-a^2-b^2}$ & $a$ \\[0.7ex]
	4 & $-b$ & $-\sqrt{1-a^2-b^2}$ & $a$ \\[0.7ex]
	5 & $b$ & $-\sqrt{1-a^2-b^2}$ & $a$ \\[0.7ex]
	6 & $0$ & $\sqrt{\frac{1}{3}}$ & $\sqrt{\frac{2}{3}}$ \\[0.7ex]
	7 & $0$ & $-\sqrt{\frac{1}{3}}$ & $\sqrt{\frac{2}{3}}$ \\[0.7ex]
	8 & $\frac{1}{2}$ & $\sqrt{\frac{1}{2}}$ & $\frac{1}{2}$ \\[0.7ex]
	9 & $-\frac{1}{2}$ & $\sqrt{\frac{1}{2}}$ & $\frac{1}{2}$ \\[0.7ex]
	10 & $-\frac{1}{2}$ & $-\sqrt{\frac{1}{2}}$ & $\frac{1}{2}$ \\[0.7ex]
	11 & $\frac{1}{2}$ & $-\sqrt{\frac{1}{2}}$ & $\frac{1}{2}$ \\[0.7ex]
	12 & $a$ & $\sqrt{1-a^2-b^2}$ & $b$ \\[0.7ex]
	13 & $-a$ & $\sqrt{1-a^2-b^2}$ & $b$ \\[0.7ex]
	14 & $-a$ & $-\sqrt{1-a^2-b^2}$ & $b$ \\[0.7ex]
	15 & $a$ & $-\sqrt{1-a^2-b^2}$ & $b$ \\[0.7ex]
	16 & $0$ & $\sqrt{1-c^2}$ & $c$ \\[0.7ex]
	17 & $0$ & $-\sqrt{1-c^2}$ & $c$ \\[0.7ex]
	18 & $1$ & $0$ & $0$ \\[0.7ex]
	19 & $\sqrt{\frac{2}{3}}$ & $\sqrt{\frac{1}{3}}$ & $0$ \\[0.7ex]
	20 & $c$ & $\sqrt{1-c^2}$ & $0$ \\[0.7ex]
	21 & $-c$ & $\sqrt{1-c^2}$ & $0$ \\[0.7ex]
	22 & $-\sqrt{\frac{2}{3}}$ & $\sqrt{\frac{1}{3}}$ & $0$ \\[0.7ex]
	23 & $-1$ & $0$ & $0$ \\[0.7ex]
	24 & $-\sqrt{\frac{2}{3}}$ & $-\sqrt{\frac{1}{3}}$ & $0$ \\[0.7ex]
	25 & $-c$ & $-\sqrt{1-c^2}$ & $0$ \\[0.7ex]
	26 & $c$ & $-\sqrt{1-c^2}$ & $0$ \\[0.7ex]
	27 & $\sqrt{\frac{2}{3}}$ & $-\sqrt{\frac{1}{3}}$ & $0$ \\[0.7ex]
	28 & $0$ & $\sqrt{1-c^2}$ & $-c$ \\[0.7ex]
	29 & $0$ & $-\sqrt{1-c^2}$ & $-c$ \\[0.7ex]
	30 & $a$ & $\sqrt{1-a^2-b^2}$ & $-b$ \\[0.7ex]
	31 & $-a$ & $\sqrt{1-a^2-b^2}$ & $-b$ \\[0.7ex]
	32 & $-a$ & $-\sqrt{1-a^2-b^2}$ & $-b$ \\[0.7ex]
	33 & $a$ & $-\sqrt{1-a^2-b^2}$ & $-b$ \\[0.7ex]
	34 & $\frac{1}{2}$ & $\sqrt{\frac{1}{2}}$ & $-\frac{1}{2}$ \\[0.7ex]
	35 & $-\frac{1}{2}$ & $\sqrt{\frac{1}{2}}$ & $-\frac{1}{2}$ \\[0.7ex]
	36 & $-\frac{1}{2}$ & $-\sqrt{\frac{1}{2}}$ & $-\frac{1}{2}$ \\[0.7ex]
	37 & $\frac{1}{2}$ & $-\sqrt{\frac{1}{2}}$ & $-\frac{1}{2}$ \\[0.7ex]
	38 & $0$ & $\sqrt{\frac{1}{3}}$ & $-\sqrt{\frac{2}{3}}$ \\[0.7ex]
	39 & $0$ & $-\sqrt{\frac{1}{3}}$ & $-\sqrt{\frac{2}{3}}$ \\[0.7ex]
	40 & $b$ & $\sqrt{1-a^2-b^2}$ & $-a$ \\[0.7ex]
	41 & $-b$ & $\sqrt{1-a^2-b^2}$ & $-a$ \\[0.7ex]
	42 & $-b$ & $-\sqrt{1-a^2-b^2}$ & $-a$ \\[0.7ex]
	43 & $b$ & $-\sqrt{1-a^2-b^2}$ & $-a$ \\[0.7ex]
	44 & $0$ & $0$ & $-1$
\end{longtable}

\noindent
\textit{\textbf{Algebraic Solutions -}}

The \textit{logarithmic} and \textit{Inverse square $1/r^2$} potentials do have algebraic polynomial solutions.

\noindent
For the \textit{logarithmic} potential:
\begin{align*}
a = & \; 0.8442833840592990183\ldots = \text{ a root of } 691200x^{24} - 3478080x^{22} + 9152467x^{20} - 15187122x^{18} \\ & \quad + 17165426x^{16} - 14250220x^{14} + 8875691x^{12} - 3897914x^{10} + 1120864x^{8} - 220664x^{6} \\ & \quad + 36336x^{4} - 4288x^{2} + 192 \\
b = & \; 0.4632349389083088107\ldots = \text{ a root of } 691200x^{24} - 3478080x^{22} + 9152467x^{20} - 15187122x^{18} \\ & \quad + 17165426x^{16} - 14250220x^{14} + 8875691x^{12} - 3897914x^{10} + 1120864x^{8} - 220664x^{6} \\ & \quad + 36336x^{4} - 4288x^{2} + 192 \\
c = & \; 0.3810484451509902075\ldots = \text{ a root of } 1080x^{12} - 2091x^{10} - 1142x^{8} + 3896x^{6} - 2704x^{4} \\ & \quad + 976x^{2} - 96
\end{align*}

\noindent
For the \textit{Coulomb $1/r$} potential:

The search for the algebraic parameters has failed to discover any polynomial of degree $360$ or less.

\noindent
For the \textit{Inverse square $1/r^2$} potential:
\begin{align*}
a = & \; 0.8443034251069456299\ldots = \text{ a root of } 12043468800x^{52} - 94325299200x^{50} + 429266873280x^{48} \\ & \quad - 1703901418000x^{46} + 5455380810043x^{44} - 12392963212202x^{42} + 19545990279972x^{40} \\ & \quad - 22627456700480x^{38} + 23684144499490x^{36} - 31217165599468x^{34} + 47978366955940x^{32} \\ & \quad - 63666108478696x^{30} + 65311166784555x^{28} - 51305686880890x^{26} + 31120897634080x^{24} \\ & \quad - 14613645697144x^{22} + 5241057684624x^{20} - 1365785078464x^{18} + 217064666944x^{16} \\ & \quad - 686102656x^{14} - 9735799808x^{12} + 2370701824x^{10} - 179171328x^{8} - 28946432x^{6} \\ & \quad + 9895936x^{4} - 1220608x^{2} + 49152 \\
b = & \; 0.4631588728416063474\ldots = \text{ a root of } 12043468800x^{52} - 94325299200x^{50} + 429266873280x^{48} \\ & \quad - 1703901418000x^{46} + 5455380810043x^{44} - 12392963212202x^{42} + 19545990279972x^{40} \\ & \quad - 22627456700480x^{38} + 23684144499490x^{36} - 31217165599468x^{34} + 47978366955940x^{32} \\ & \quad - 63666108478696x^{30} + 65311166784555x^{28} - 51305686880890x^{26} + 31120897634080x^{24} \\ & \quad - 14613645697144x^{22} + 5241057684624x^{20} - 1365785078464x^{18} + 217064666944x^{16} \\ & \quad - 686102656x^{14} - 9735799808x^{12} + 2370701824x^{10} - 179171328x^{8} - 28946432x^{6} \\ & \quad + 9895936x^{4} - 1220608x^{2} + 49152 \\
c = & \; 0.3811445522653392825\ldots = \text{ a root of } 106920x^{26} - 1105110x^{24} + 4390647x^{22} - 11158118x^{20} \\ & \quad + 20541668x^{18} - 24353800x^{16} + 13311840x^{14} + 4501056x^{12} - 11685760x^{10} + 7102720x^{8} \\ & \quad - 1760512x^{6} + 25088x^{4} + 87040x^{2} - 10240
\end{align*}

\noindent
\textit{\textbf{Minimal Energy values -}}

The coordinates for the \textit{Coulomb $1/r$} potential are known to 50,014 digits. The minimal energies have been determined for all 3 potentials as well, but the algebraic polynomial is unknown for all 3.
\begin{center}
	\begin{tabular}{l|l}
		\multicolumn{2}{c}{Minimal Energy - 44 points} \\
		\hline\Tstrut
		logarithmic & -229.6418014875968238\ldots \\[0.2ex]
		\hline\Tstrut
		Coulomb & 807.1742630846280758\ldots \\[0.2ex]
		\hline\Tstrut
		Inverse square law & 836.0418318113427941\ldots
	\end{tabular}
\end{center}

\noindent
\textit{\textbf{Symmetries -}}

The symmetry groups for 44 points are identical under all 3 potentials. From the large count of the groups in both parallel planes and polygons, it is easy to see that this configuration has intensive symmetries embedded in the polyhedra.

\begin{center}
	\begin{tabular}{l|l}
		\multicolumn{2}{c}{Symmetries - 44 points} \\
		\hline\Tstrut
		planes & [[2, 3732], [4, 84], [8, 396], [10, 24], [20, 24], [68, 4], [140, 3], [144, 6]] \\[0.2ex]
		\hline\Tstrut
		Gram groups & [[24, 4], [44, 2], [48, 16], [96, 8], [216, 1]] \\
		\hline\Tstrut
		Polygons & [[3, 7544], [4, 981], [6, 36], [8, 6], [10, 6]]
	\end{tabular}
\end{center}

\subsection{45 points}
The symmetry groups for this configuration has only one entry for polygons [3,15]. Upon examination, this is actually the parameterization indicator, a set of 15 parallel triangles embedded in this polyhedron. The arrangement is [3:3:3:3:3:3:3:3:3:3:3:3:3:3:3] which is balanced but a bit unusual.

\begin{figure}[ht]
	\begin{center}
		\includegraphics[type=pdf,ext=pdf,read=pdf,height=1in,width=1in,angle=0]{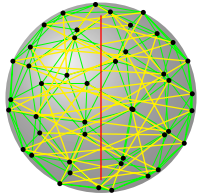}
		\caption{45 points.}
		\label{fig:45pts}
	\end{center}
\end{figure}

\noindent
\textit{\textbf{Algebraic parameterization -}}

It requires 14 parameters to adequately constrain the 45 point set, due to the 15 embedded triangles in the figure.

The algebraic parameterized structure is given below:

\begin{longtable}[c]{r|ccc}
	\caption{Parameterization for 45 points} \\
	pt & $x$ & $y$ & $z$ \\
	\hline\vspace*{-2.2ex}
	\endfirsthead
	\multicolumn{4}{c}%
	{\tablename\ \thetable\ -- 45 points parameters -- \textit{continued}} \\
	pt & $x$ & $y$ & $z$ \\
	\hline\vspace*{-2.2ex}
	\endhead
	1 & $b$ & $\sqrt{1-a^2-b^2}$ & $a$ \\[0.9ex	]
	2 & $\frac{-b-\sqrt{3}\sqrt{1-a^2-b^2}}{2}$ & $\frac{b\sqrt{3}-\sqrt{1-a^2-b^2}}{2}$ & $a$ \\[0.9ex	]
	3 & $\frac{-b+\sqrt{3}\sqrt{1-a^2-b^2}}{2}$ & $\frac{-b\sqrt{3}-\sqrt{1-a^2-b^2}}{2}$ & $a$ \\[0.9ex	]
	4 & $d$ & $-\sqrt{1-c^2-d^2}$ & $c$ \\[0.9ex	]
	5 & $\frac{-d+\sqrt{3}\sqrt{1-c^2-d^2}}{2}$ & $\frac{d\sqrt{3}+\sqrt{1-c^2-d^2}}{2}$ & $c$ \\[0.9ex	]
	6 & $\frac{-d-\sqrt{3}\sqrt{1-c^2-d^2}}{2}$ & $\frac{-d\sqrt{3}+\sqrt{1-c^2-d^2}}{2}$ & $c$ \\[0.9ex	]
	7 & $f$ & $\sqrt{1-e^2-f^2}$ & $e$ \\[0.9ex	]
	8 & $\frac{-f-\sqrt{3}\sqrt{1-e^2-f^2}}{2}$ & $\frac{f\sqrt{3}-\sqrt{1-e^2-f^2}}{2}$ & $e$ \\[0.9ex	]
	9 & $\frac{-f+\sqrt{3}\sqrt{1-e^2-f^2}}{2}$ & $\frac{-f\sqrt{3}-\sqrt{1-e^2-f^2}}{2}$ & $e$ \\[0.9ex	]
	10 & $h$ & $\sqrt{1-g^2-h^2}$ & $g$ \\[0.9ex	]
	11 & $\frac{-h-\sqrt{3}\sqrt{1-g^2-h^2}}{2}$ & $\frac{h\sqrt{3}-\sqrt{1-g^2-h^2}}{2}$ & $g$ \\[0.9ex	]
	12 & $\frac{-h+\sqrt{3}\sqrt{1-g^2-h^2}}{2}$ & $\frac{-h\sqrt{3}-\sqrt{1-g^2-h^2}}{2}$ & $g$ \\[0.9ex	]
	13 & $j$ & $-\sqrt{1-i^2-j^2}$ & $i$ \\[0.9ex	]
	14 & $\frac{-j+\sqrt{3}\sqrt{1-i^2-j^2}}{2}$ & $\frac{j\sqrt{3}+\sqrt{1-i^2-j^2}}{2}$ & $i$ \\[0.9ex	]
	15 & $\frac{-j-\sqrt{3}\sqrt{1-i^2-j^2}}{2}$ & $\frac{-j\sqrt{3}+\sqrt{1-i^2-j^2}}{2}$ & $i$ \\[0.9ex	]
	16 & $l$ & $-\sqrt{1-k^2-l^2}$ & $k$ \\[0.9ex	]
	17 & $\frac{-l+\sqrt{3}\sqrt{1-k^2-l^2}}{2}$ & $\frac{l\sqrt{3}+\sqrt{1-k^2-l^2}}{2}$ & $k$ \\[0.9ex	]
	18 & $\frac{-l-\sqrt{3}\sqrt{1-k^2-l^2}}{2}$ & $\frac{-l\sqrt{3}+\sqrt{1-k^2-l^2}}{2}$ & $k$ \\[0.9ex	]
	19 & $n$ & $\sqrt{1-m^2-n^2}$ & $m$ \\[0.9ex	]
	20 & $\frac{-n-\sqrt{3}\sqrt{1-m^2-n^2}}{2}$ & $\frac{n\sqrt{3}-\sqrt{1-m^2-n^2}}{2}$ & $m$ \\[0.9ex	]
	21 & $\frac{-n+\sqrt{3}\sqrt{1-m^2-n^2}}{2}$ & $\frac{-n\sqrt{3}-\sqrt{1-m^2-n^2}}{2}$ & $m$ \\[0.9ex	]
	22 & $1$ & $0$ & $0$ \\[0.9ex	]
	23 & $-\frac{1}{2}$ & $\sqrt{\frac{3}{2}}$ & $0$ \\[0.9ex	]
	24 & $-\frac{1}{2}$ & $-\sqrt{\frac{3}{2}}$ & $0$ \\[0.9ex	]
	25 & $n$ & $-\sqrt{1-m^2-n^2}$ & $-m$ \\[0.9ex	]
	26 & $\frac{-n+\sqrt{3}\sqrt{1-m^2-n^2}}{2}$ & $\frac{n\sqrt{3}+\sqrt{1-m^2-n^2}}{2}$ & $-m$ \\[0.9ex	]
	27 & $\frac{-n-\sqrt{3}\sqrt{1-m^2-n^2}}{2}$ & $\frac{-n\sqrt{3}+\sqrt{1-m^2-n^2}}{2}$ & $-m$ \\[0.9ex	]
	28 & $l$ & $\sqrt{1-k^2-l^2}$ & $-k$ \\[0.9ex	]
	29 & $\frac{-l-\sqrt{3}\sqrt{1-k^2-l^2}}{2}$ & $\frac{l\sqrt{3}-\sqrt{1-k^2-l^2}}{2}$ & $-k$ \\[0.9ex	]
	30 & $\frac{-l+\sqrt{3}\sqrt{1-k^2-l^2}}{2}$ & $\frac{-l\sqrt{3}-\sqrt{1-k^2-l^2}}{2}$ & $-k$ \\[0.9ex	]
	31 & $j$ & $\sqrt{1-i^2-j^2}$ & $-i$ \\[0.9ex	]
	32 & $\frac{-j-\sqrt{3}\sqrt{1-i^2-j^2}}{2}$ & $\frac{j\sqrt{3}-\sqrt{1-i^2-j^2}}{2}$ & $-i$ \\[0.9ex	]
	33 & $\frac{-j+\sqrt{3}\sqrt{1-i^2-j^2}}{2}$ & $\frac{-j\sqrt{3}-\sqrt{1-i^2-j^2}}{2}$ & $-i$ \\[0.9ex	]
	34 & $h$ & $-\sqrt{1-g^2-h^2}$ & $-g$ \\[0.9ex	]
	35 & $\frac{-h+\sqrt{3}\sqrt{1-g^2-h^2}}{2}$ & $\frac{h\sqrt{3}+\sqrt{1-g^2-h^2}}{2}$ & $-g$ \\[0.9ex	]
	36 & $\frac{-h-\sqrt{3}\sqrt{1-g^2-h^2}}{2}$ & $\frac{-h\sqrt{3}+\sqrt{1-g^2-h^2}}{2}$ & $-g$ \\[0.9ex	]
	37 & $f$ & $-\sqrt{1-e^2-f^2}$ & $-e$ \\[0.9ex	]
	38 & $\frac{-f+\sqrt{3}\sqrt{1-e^2-f^2}}{2}$ & $\frac{f\sqrt{3}+\sqrt{1-e^2-f^2}}{2}$ & $-e$ \\[0.9ex	]
	39 & $\frac{-f-\sqrt{3}\sqrt{1-e^2-f^2}}{2}$ & $\frac{-f\sqrt{3}+\sqrt{1-e^2-f^2}}{2}$ & $-e$ \\[0.9ex	]
	40 & $d$ & $\sqrt{1-c^2-d^2}$ & $-c$ \\[0.9ex	]
	41 & $\frac{-d-\sqrt{3}\sqrt{1-c^2-d^2}}{2}$ & $\frac{d\sqrt{3}-\sqrt{1-c^2-d^2}}{2}$ & $-c$ \\[0.9ex	]
	42 & $\frac{-d+\sqrt{3}\sqrt{1-c^2-d^2}}{2}$ & $\frac{-d\sqrt{3}-\sqrt{1-c^2-d^2}}{2}$ & $-c$ \\[0.9ex	]
	43 & $b$ & $-\sqrt{1-a^2-b^2}$ & $-a$ \\[0.9ex	]
	44 & $\frac{-b+\sqrt{3}\sqrt{1-a^2-b^2}}{2}$ & $\frac{b\sqrt{3}+\sqrt{1-a^2-b^2}}{2}$ & $-a$ \\[0.9ex	]
	45 & $\frac{-b-\sqrt{3}\sqrt{1-a^2-b^2}}{2}$ & $\frac{-b\sqrt{3}+\sqrt{1-a^2-b^2}}{2}$ & $-a$
\end{longtable}

The values to 19 digits of the 14 parameters optimized for the minimal solutions of 45 points are:

\begin{longtable}[c]{c|c|c|c}
	\caption{Parameter values for 45 points} \\
	Parameter & log & 1/r & $1/r^2$ \\
	\hline\vspace*{-2.2ex}
	\endfirsthead
	\multicolumn{4}{c}%
	{\tablename\ \thetable\ -- 45 points parameters -- \textit{continued}} \\
	Parameter & log & 1/r & $1/r^2$ \\
	\hline\vspace*{-2.2ex}
	\endhead
	$a$ & 0.9448012977797529915 & 0.9446090189539352018 & 0.9444249538496973933 \\
	$b$ & 0.2975827807804567773 & 0.2974120108274068113 & 0.2972440703946039624 \\
	$c$ & 0.7866280406501610825 & 0.7863225435219266930 & 0.7859835891091541787 \\
	$d$ & 0.5311607571861045193 & 0.5307863665179054499 & 0.5307061927310007725 \\
	$e$ & 0.6580076535400015539 & 0.6564620371815712186 & 0.6550995497484144041 \\
	$f$ & 0.5150326203808932422 & 0.5165854646291564992 & 0.5179915653250577083 \\
	$g$ & 0.5704617934153301843 & 0.5710451211608503535 & 0.5711378658265501452 \\
	$h$ & 0.8128947994873928471 & 0.8126461199831946751 & 0.8127083999873582932 \\
	$i$ & 0.3593360006397525969 & 0.3591975150162738452 & 0.3587851285984275580 \\
	$j$ & 0.8522973673080205546 & 0.8526016635459988540 & 0.8529259647707977230 \\
	$k$ & 0.2553162973486828234 & 0.2548416787990827249 & 0.2542717550538903899 \\
	$l$ & 0.5522703199981502771 & 0.5531992353715904161 & 0.5538189987679110910 \\
	$m$ & 0.1825666013633008320 & 0.1821998348980555125 & 0.1818035064680277061 \\
	$n$ & 0.8407843910871033489 & 0.8401562364350225808 & 0.8395437670842417963 \\
	\hline\Tstrut
	$energy$ & -239.4536982534525604 & 846.1884010610781359 & 880.3579693212025387
\end{longtable}

All 14 parameters are all known to 50,014 digits precision, for all 3 potentials, but again, attempts to find the algebraic numbers for these 42 parameters failed, the algebraic degree $> 360$.

The symmetry groups for the \textit{Logarithmic}, \textit{Coulomb} and \textit{Inverse square law} potentials are all identical. The polygon group hinted at the 15 parallel triangles, which was the correct alignment to use for parameterization.

\begin{center}
	\begin{tabular}{l|l}
		\multicolumn{2}{c}{Symmetries - 45 points - \textit{Coulomb} or \textit{Inverse sq.}} \\
		\hline\Tstrut
		planes & [[15, 1]] \\[0.2ex]
		\hline\Tstrut
		Gram groups & [[6, 22], [12, 154], [45, 1]] \\
		\hline\Tstrut
		Polygons & [[3, 15]]
	\end{tabular}
\end{center}

\subsection{46 points}
As might be expected, from the optimal 45 point configuration containing 15 co-planar triangles for the optimal solution, adding 1 new point pushes it to a pole, while the other 45 points remain as 15 co-planar triangles with a line through their centroids going through the isolated pole point, but they all slightly descend beneath the isolated point.

\begin{figure}[ht]
	\begin{center}
		\includegraphics[type=pdf,ext=pdf,read=pdf,height=1in,width=1in,angle=0]{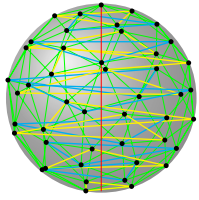}
		\caption{46 points.}
		\label{fig:46pts}
	\end{center}
\end{figure}
The 15 triangles are shown in figure \ref{fig:46pts} in alternating yellow and cyan outlines.

\noindent
\textit{\textbf{Algebraic parameterization -}}

It requires 24 parameters $a$ - $x$ to adequately constrain the 46 point set, due to the 15 embedded triangles in the figure. It should be noted that there are some known points in the configuration, points 8-10 are known, as well as points 32-37, in all 3 potentials.

The algebraic parameterized structure is given below:

\begin{longtable}[c]{r|ccc}
	\caption{Parameterization for 46 points} \\
	pt & $x$ & $y$ & $z$ \\
	\hline\vspace*{-2.2ex}
	\endfirsthead
	\multicolumn{4}{c}%
	{\tablename\ \thetable\ -- 46 points parameters -- \textit{continued}} \\
	pt & $x$ & $y$ & $z$ \\
	\hline\vspace*{-2.2ex}
	\endhead
	1 & $0$ & $0$ & $1$	\\[0.7em]
	2 & $b$ & $-\sqrt{1-a^2-b^2}$ & $a$	\\[0.7em]
	3 & $\frac{+\sqrt{3}\sqrt{1-a^2-b^2}-b}{2}$ & 	$\frac{+b\sqrt{3}+\sqrt{1-a^2-b^2}}{2}$ & $a$	\\[0.7em]
	4 & $\frac{-\sqrt{3}\sqrt{1-a^2-b^2}-b}{2}$ & $\frac{-b\sqrt{3}+\sqrt{1-a^2-b^2}}{2}$ & $a$	\\[0.7em]
	5 & $d$ & $\sqrt{1-c^2-d^2}$ & $c$	\\[0.7em]
	6 & $\frac{+\sqrt{3}\sqrt{1-c^2-d^2}-d}{2}$ & $\frac{-d\sqrt{3}-\sqrt{1-c^2-d^2}}{2}$ & $c$	\\[0.7em]
	7 & $\frac{-\sqrt{3}\sqrt{1-c^2-d^2}-d}{2}$ & $\frac{+d\sqrt{3}-\sqrt{1-c^2-d^2}}{2}$ & $c$	\\[0.7em]
	8 & $\sqrt{\frac{2}{3}}$ & $0$ & $\sqrt{\frac{1}{3}}$	\\[0.7em]
	9 & $-\sqrt{\frac{1}{6}}$ & $+\sqrt{\frac{1}{2}}$ & $\sqrt{\frac{1}{3}}$	\\[0.7em]
	10 & $-\sqrt{\frac{1}{6}}$ & $-\sqrt{\frac{1}{2}}$ & $\sqrt{\frac{1}{3}}$	\\[0.7em]
	11 & $f$ & $-\sqrt{1-e^2-f^2}$ & $e$	\\[0.7em]
	12 & $\frac{+\sqrt{3}\sqrt{1-e^2-f^2}-f}{2}$ & $\frac{+f\sqrt{3}+\sqrt{1-e^2-f^2}}{2}$ & $e$	\\[0.7em]
	13 & $\frac{-\sqrt{3}\sqrt{1-e^2-f^2}-f}{2}$ & $\frac{-f\sqrt{3}+\sqrt{1-e^2-f^2}}{2}$ & $e$	\\[0.7em]
	14 & $h$ & $\sqrt{1-g^2-h^2}$ & $g$	\\[0.7em]
	15 & $\frac{+\sqrt{3}\sqrt{1-g^2-h^2}-h}{2}$ & $\frac{-h\sqrt{3}-\sqrt{1-g^2-h^2}}{2}$ & $g$	\\[0.7em]
	16 & $\frac{-\sqrt{3}\sqrt{1-g^2-h^2}-h}{2}$ & $\frac{+h\sqrt{3}-\sqrt{1-g^2-h^2}}{2}$ & $g$	\\[0.7em]
	17 & $j$ & $-\sqrt{1-i^2-j^2}$ & $i$	\\[0.7em]
	18 & $\frac{+\sqrt{3}\sqrt{1-i^2-j^2}-j}{2}$ & $\frac{+j\sqrt{3}+\sqrt{1-i^2-j^2}}{2}$ & $i$	\\[0.7em]
	19 & $\frac{-\sqrt{3}\sqrt{1-i^2-j^2}-j}{2}$ & $\frac{-j\sqrt{3}+\sqrt{1-i^2-j^2}}{2}$ & $i$	\\[0.7em]
	20 & $l$ & $\sqrt{1-k^2-l^2}$ & $k$	\\[0.7em]
	21 & $\frac{+\sqrt{3}\sqrt{1-k^2-l^2}-l}{2}$ & $\frac{-l\sqrt{3}-\sqrt{1-k^2-l^2}}{2}$ & $k$	\\[0.7em]
	22 & $\frac{-\sqrt{3}\sqrt{1-k^2-l^2}-l}{2}$ & $\frac{+l\sqrt{3}-\sqrt{1-k^2-l^2}}{2}$ & $k$	\\[0.7em]
	23 & $n$ & $-\sqrt{1-m^2-n^2}$ & $m$	\\[0.7em]
	24 & $\frac{+\sqrt{3}\sqrt{1-m^2-n^2}-n}{2}$ & $\frac{+n\sqrt{3}+\sqrt{1-m^2-n^2}}{2}$ & $m$	\\[0.7em]
	25 & $\frac{-\sqrt{3}\sqrt{1-m^2-n^2}-n}{2}$ & $\frac{-n\sqrt{3}+\sqrt{1-m^2-n^2}}{2}$ & $m$	\\[0.7em]
	26 & $p$ & $\sqrt{1-o^2-p^2}$ & $o$	\\[0.7em]
	27 & $\frac{+\sqrt{3}\sqrt{1-o^2-p^2}-p}{2}$ & $\frac{-p\sqrt{3}-\sqrt{1-o^2-p^2}}{2}$ & $o$	\\[0.7em]
	28 & $\frac{-\sqrt{3}\sqrt{1-o^2-p^2}-p}{2}$ & $\frac{+p\sqrt{3}-\sqrt{1-o^2-p^2}}{2}$ & $o$	\\[0.7em]
	29 & $r$ & $-\sqrt{1-q^2-r^2}$ & $q$	\\[0.7em]
	30 & $\frac{+\sqrt{3}\sqrt{1-q^2-r^2}-r}{2}$ & $\frac{+r\sqrt{3}+\sqrt{1-q^2-r^2}}{2}$ & $q$	\\[0.7em]
	31 & $\frac{-\sqrt{3}\sqrt{1-q^2-r^2}-r}{2}$ & $\frac{-r\sqrt{3}+\sqrt{1-q^2-r^2}}{2}$ & $q$	\\[0.7em]
	32 & $\frac{2}{3}\sqrt{2}$ & $0$ & $-\frac{1}{3}$	\\[0.7em]
	33 & $\frac{-\sqrt{2}}{3}$ & $\sqrt{\frac{2}{3}}$ & $-\frac{1}{3}$	\\[0.7em]
	34 & $\frac{-\sqrt{2}}{3}$ & $-\sqrt{\frac{2}{3}}$ & $-\frac{1}{3}$	\\[0.7em]
	35 & $-\sqrt{\frac{2}{3}}$ & $0$ & $-\sqrt{\frac{1}{3}}$	\\[0.7em]
	36 & $\sqrt{\frac{1}{6}}$ & $+\sqrt{\frac{1}{2}}$ & $-\sqrt{\frac{1}{3}}$	\\[0.7em]
	37 & $\sqrt{\frac{1}{6}}$ & $-\sqrt{\frac{1}{2}}$ & $-\sqrt{\frac{1}{3}}$	\\[0.7em]
	38 & $t$ & $\sqrt{1-s^2-t^2}$ & $s$	\\[0.7em]
	39 & $\frac{+\sqrt{3}\sqrt{1-s^2-t^2}-t}{2}$ & $\frac{-t\sqrt{3}-\sqrt{1-s^2-t^2}}{2}$ & $s$	\\[0.7em]
	40 & $\frac{-\sqrt{3}\sqrt{1-s^2-t^2}-t}{2}$ & $\frac{+t\sqrt{3}-\sqrt{1-s^2-t^2}}{2}$ & $s$	\\[0.7em]
	41 & $v$ & $-\sqrt{1-u^2-v^2}$ & $u$	\\[0.7em]
	42 & $\frac{+\sqrt{3}\sqrt{1-u^2-v^2}-v}{2}$ & 		$\frac{+v\sqrt{3}+\sqrt{1-u^2-v^2}}{2}$ & $u$	\\[0.7em]
	43 & $\frac{-\sqrt{3}\sqrt{1-u^2-v^2}-v}{2}$ & 	$\frac{-v\sqrt{3}+\sqrt{1-u^2-v^2}}{2}$ & $u$	\\[0.7em]
	44 & $x$ & $\sqrt{1-w^2-x^2}$ & $w$	\\[0.7em]
	45 & $\frac{+\sqrt{3}\sqrt{1-w^2-x^2}-x}{2}$ & 				$\frac{-x\sqrt{3}-\sqrt{1-w^2-x^2}}{2}$ & $w$	\\[0.7em]
	46 & $\frac{-\sqrt{3}\sqrt{1-w^2-x^2}-x}{2}$ & $\frac{+x\sqrt{3}-\sqrt{1-w^2-x^2}}{2}$ & $w$
\end{longtable}

The values to 19 digits of the 24 parameters optimized for the minimal solutions of 46 points are:

\begin{longtable}[c]{c|c|c|c}
	\caption{Parameter values for 46 points} \\
	Parameter & log & 1/r & $1/r^2$ \\
	\hline\vspace*{-2.2ex}
	\endfirsthead
	\multicolumn{4}{c}%
	{\tablename\ \thetable\ -- 46 points parameters -- \textit{continued}} \\
	Parameter & log & 1/r & $1/r^2$ \\
	\hline\vspace*{-2.2ex}
	\endhead
	$a$ & 0.8675862714330420908 & 0.8677199245471095582 & 0.8678889908992212018 \\
	$b$ & 0.4489102568112173187 & 0.4493240191878207738 & 0.4495077494171742467 \\
	$c$ & 0.8073542358509505115 & 0.8084933810854576747 & 0.8096058113331843383 \\
	$d$ & 0.4709722071870518258 & 0.4696095622181022988 & 0.4684144190532517365 \\
	$e$ & 0.5302074417575879488 & 0.5313632533910031675 & 0.5322654743028723002 \\
	$f$ & 0.6023089241535305098 & 0.6017575319737985536 & 0.6013723144187563307 \\
	$g$ & 0.3911264855983178135 & 0.3902213575213371435 & 0.3895574306235108220 \\
	$h$ & 0.7006540113876593466 & 0.7015599236527649627 & 0.7022821398342988540 \\
	$i$ & 0.1749187766528298038 & 0.1732543475774688254 & 0.1717567457130871582 \\
	$j$ & 0.9181716090488710699 & 0.9187913904859897010 & 0.9194418187267255866 \\
	$k$ & 0.1340412251695277031 & 0.1343867730412334342 & 0.1345036401119759534 \\
	$l$ & 0.9676049333299480676 & 0.9678688634865276264 & 0.9680895041817102717 \\
	$m$ & 0.02657290099507830443 & 0.02606220098186304889 & 0.02560322995709392398 \\
	$n$ & 0.5717648554982521366 & 0.5728102191426471549 & 0.5736048981680224227 \\
	$o$ & -0.2008902683654080937 & -0.2012703286551560381 & -0.2018602701938234770 \\
	$p$ & 0.7853030469670528034 & 0.7863769213430970882 & 0.7873482559694962344 \\
	$q$ & -0.3261261357837803286 & -0.3275212992569916420 & -0.3285473794342779728 \\
	$r$ & 0.8049112026245471598 & 0.8045596983830871199 & 0.8043762462034599139 \\
	$s$ & -0.6755013608187894653 & -0.6745853983313513504 & -0.6738452515769191825 \\
	$t$ & 0.6813884072241316567 & 0.6818540094021473591 & 0.6822950127427863829 \\
	$u$ & -0.7813827441383722216 & -0.7804774000077704620 & -0.7795022868524480195 \\
	$v$ & 0.5800679639436375106 & 0.5815962974107800149 & 0.5831209624307032022 \\
	$w$ & -0.9479068283509840667 & -0.9476468118942033598 & -0.9474261348834770462 \\
	$x$ & 0.2272342431235360469 & 0.2285520961890842169 & 0.2295870670818189546 \\
	\hline\Tstrut
	$energy$ & -249.4558479008570933 & 886.1671136391913397 & 926.0623438450871291
\end{longtable}

The degree of the algebraic codes $>420$. The spherical codes are known to 50,014 digits accuracy for the \textit{logarithmic}, \textit{Coulomb}, and \textit{inverse square} potentials.

\noindent
\textit{\textbf{Symmetries -}}

\begin{center}
	\begin{tabular}{l|l}
		\multicolumn{2}{c}{Symmetries - 46 points} \\
		\hline\Tstrut
		planes & [[4, 75], [15, 4]] \\[0.2ex]
		\hline\Tstrut
		Gram groups & [[6, 1], [12, 10], [24, 81], [46, 1]] \\
		\hline\Tstrut
		Polygons & [[3, 60], [4, 75]]
	\end{tabular}
\end{center}

\subsection{47 points}

No attempts to parameterize 47 points from their optimal solutions were made, although the single embedded heptagon might be an indicator of a possible arrangement. A careful second look showed that 40 to 47 parameters were necessary to constrain this polyhedron.

\begin{figure}[ht]
	\begin{center}
		\includegraphics[type=pdf,ext=pdf,read=pdf,height=1in,width=1in,angle=0]{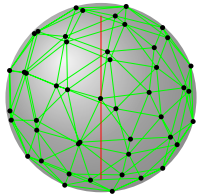}
		\caption{47 points.}
		\label{fig:47pts}
	\end{center}
\end{figure}

\noindent
\textit{\textbf{Minimal Energy values -}}

The coordinates for 47 points are known to 77 digits for the \textit{log} potential and 38 digits for the other two. The minimal energies have been determined for all 3 potentials as well.
\begin{center}
	\begin{tabular}{l|l}
		\multicolumn{2}{c}{Minimal Energy - 47 points} \\
		\hline\Tstrut
		logarithmic & -259.6617598532650787\ldots \\[0.2ex]
		\hline\Tstrut
		Coulomb & 927.0592706797097552\ldots \\[0.2ex]
		\hline\Tstrut
		Inverse square law & 972.8237449081603419\ldots
	\end{tabular}
\end{center}

\noindent
\textit{\textbf{Symmetries -}}

The symmetry groups for 47 points are identical under all 3 potentials.

\begin{center}
	\begin{tabular}{l|l}
		\multicolumn{2}{c}{Symmetries - 47 points} \\
		\hline\Tstrut
		planes & [[4, 190], [35, 1]] \\[0.2ex]
		\hline\Tstrut
		Gram groups & [[2, 41], [4, 520], [47, 1]] \\
		\hline\Tstrut
		Polygons & [[4, 190], [7, 1]]
	\end{tabular}
\end{center}

\subsection{48 points}

The symmetry group polygons hints strongly as squares embedded in the configuration, after determining an optimal axis and rotating the polyhedron twice, it became obvious that 48 points were on 12 different parallel plane squares. This is an 4:4:4:4:4:4:4:4:4:4:4:4 arrangement.

\begin{figure}[ht]
	\begin{center}
		\includegraphics[type=pdf,ext=pdf,read=pdf,height=1in,width=1in,angle=0]{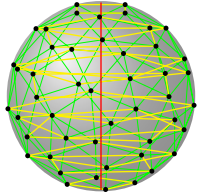}
		\caption{48 points.}
		\label{fig:48pts}
	\end{center}
\end{figure}

\noindent
\textit{\textbf{Algebraic parameterization -}}

Only 18 parameters were required to adequately constrain the 48 point set, because the last 6 squares had their $z$ value the negative of the first 6 squares.

The algebraic parameterized structure is given below:

\begin{longtable}[c]{r|ccc}
	\caption{Parameterization for 48 points} \\
	pt & $x$ & $y$ & $z$ \\
	\hline\vspace*{-2.2ex}
	\endfirsthead
	\multicolumn{4}{c}%
	{\tablename\ \thetable\ -- 48 points parameters -- \textit{continued}} \\
	pt & $x$ & $y$ & $z$ \\
	\hline\vspace*{-2.2ex}
	\endhead
	1 & $b$ & $0$ & $a$ \\[0.5ex]
	2 & $0$ & $b$ & $a$ \\[0.5ex]
	3 & $-b$ & $0$ & $a$ \\[0.5ex]
	4 & $0$ & $-b$ & $a$ \\[0.5ex]
	5 & $d$ & $\sqrt{1-c^2-d^2}$ & $c$ \\[0.5ex]
	6 & $-\sqrt{1-c^2-d^2}$ & $d$ & $c$ \\[0.5ex]
	7 & $-d$ & $-\sqrt{1-c^2-d^2}$ & $c$ \\[0.5ex]
	8 & $\sqrt{1-c^2-d^2}$ & $-d$ & $c$ \\[0.5ex]
	9 & $f$ & $\sqrt{1-e^2-f^2}$ & $e$ \\[0.5ex]
	10 & $-\sqrt{1-e^2-f^2}$ & $f$ & $e$ \\[0.5ex]
	11 & $-f$ & $-\sqrt{1-e^2-f^2}$ & $e$ \\[0.5ex]
	12 & $\sqrt{1-e^2-f^2}$ & $-f$ & $e$ \\[0.5ex]
	13 & $h$ & $\sqrt{1-g^2-h^2}$ & $g$ \\[0.5ex]
	14 & $-\sqrt{1-g^2-h^2}$ & $h$ & $g$ \\[0.5ex]
	15 & $-h$ & $-\sqrt{1-g^2-h^2}$ & $g$ \\[0.5ex]
	16 & $\sqrt{1-g^2-h^2}$ & $-h$ & $g$ \\[0.5ex]
	17 & $j$ & $\sqrt{1-i^2-j^2}$ & $i$ \\[0.5ex]
	18 & $-\sqrt{1-i^2-j^2}$ & $j$ & $i$ \\[0.5ex]
	19 & $-j$ & $-\sqrt{1-i^2-j^2}$ & $i$ \\[0.5ex]
	20 & $\sqrt{1-i^2-j^2}$ & $-j$ & $i$ \\[0.5ex]
	21 & $l$ & $\sqrt{1-k^2-l^2}$ & $k$ \\[0.5ex]
	22 & $-\sqrt{1-k^2-l^2}$ & $l$ & $k$ \\[0.5ex]
	23 & $-l$ & $-\sqrt{1-k^2-l^2}$ & $k$ \\[0.5ex]
	24 & $\sqrt{1-k^2-l^2}$ & $-l$ & $k$ \\[0.5ex]
	25 & $m$ & $\sqrt{1-k^2-m^2}$ & $-k$ \\[0.5ex]
	26 & $-\sqrt{1-k^2-m^2}$ & $m$ & $-k$ \\[0.5ex]
	27 & $-m$ & $-\sqrt{1-k^2-m^2}$ & $-k$ \\[0.5ex]
	28 & $\sqrt{1-k^2-m^2}$ & $-m$ & $-k$ \\[0.5ex]
	29 & $n$ & $\sqrt{1-i^2-n^2}$ & $-i$ \\[0.5ex]
	30 & $-\sqrt{1-i^2-n^2}$ & $n$ & $-i$ \\[0.5ex]
	31 & $-n$ & $-\sqrt{1-i^2-n^2}$ & $-i$ \\[0.5ex]
	32 & $\sqrt{1-i^2-n^2}$ & $-n$ & $-i$ \\[0.5ex]
	33 & $o$ & $-\sqrt{1-g^2-o^2}$ & $-g$ \\[0.5ex]
	34 & $-\sqrt{1-g^2-o^2}$ & $-o$ & $-g$ \\[0.5ex]
	35 & $-o$ & $\sqrt{1-g^2-o^2}$ & $-g$ \\[0.5ex]
	36 & $\sqrt{1-g^2-o^2}$ & $o$ & $-g$ \\[0.5ex]
	37 & $p$ & $\sqrt{1-e^2-p^2}$ & $-e$ \\[0.5ex]
	38 & $-\sqrt{1-e^2-p^2}$ & $p$ & $-e$ \\[0.5ex]
	39 & $-p$ & $-\sqrt{1-e^2-p^2}$ & $-e$ \\[0.5ex]
	40 & $\sqrt{1-e^2-p^2}$ & $-p$ & $-e$ \\[0.5ex]
	41 & $q$ & $\sqrt{1-c^2-q^2}$ & $-c$ \\[0.5ex]
	42 & $-\sqrt{1-c^2-q^2}$ & $q$ & $-c$ \\[0.5ex]
	43 & $-q$ & $-\sqrt{1-c^2-q^2}$ & $-c$ \\[0.5ex]
	44 & $\sqrt{1-c^2-q^2}$ & $-q$ & $-c$ \\[0.5ex]
	45 & $r$ & $\sqrt{1-a^2-r^2}$ & $-a$ \\[0.5ex]
	46 & $-\sqrt{1-a^2-r^2}$ & $r$ & $-a$ \\[0.5ex]
	47 & $-r$ & $-\sqrt{1-a^2-r^2}$ & $-a$ \\[0.5ex]
	48 & $\sqrt{1-a^2-r^2}$ & $-r$ & $-a$
\end{longtable}

The values to 19 digits of the 18 parameters optimized for the minimal solutions of 48 points are:

\begin{longtable}[c]{c|c|c|c}
	\caption{Parameter values for 48 points} \\
	Parameter & log & 1/r & $1/r^2$ \\
	\hline\vspace*{-2.2ex}
	\endfirsthead
	\multicolumn{4}{c}%
	{\tablename\ \thetable\ -- 48 points parameters -- \textit{continued}} \\
	Parameter & log & 1/r & $1/r^2$ \\
	\hline\vspace*{-2.2ex}
	\endhead
	$a$ & 0.9320721754443081266 & 0.9320520077178053482 & 0.9319666247545104778 \\
	$b$ & 0.3622726318154807716 & 0.3623245160201116168 & 0.3625440805525384080 \\
	$c$ & 0.7226497842867026243 & 0.7236545639499285646 & 0.7245542159811955933 \\
	$d$ & 0.4755085525730592076 & 0.4768192066001819708 & 0.4777538453225200198 \\
	$e$ & 0.6261193307867101643 & 0.6251449308091768179 & 0.6242071450014939539 \\
	$f$ & 0.7781389947087470409 & 0.7791559562515492960 & 0.7800655960443368334 \\
	$g$ & 0.3437492491883026209 & 0.3433387122079829769 & 0.3432084341129914291 \\
	$h$ & 0.4027389135091593085 & 0.4074354092793292925 & 0.4108795784146587298 \\
	$i$ & 0.2928342071643686616 & 0.2924344158234941378 & 0.2922099044059918069 \\
	$j$ & 0.8222242513321697287 & 0.8235655084678133153 & 0.8243778575539762615 \\
	$k$ & 0.1143587051561581860 & 0.1157479313360371066 & 0.1168168699178071584 \\
	$l$ & 0.9929257634815719041 & 0.9928952246399484305 & 0.9928484844696680854 \\
	$m$ & 0.6204004765336294144 & 0.6231004337532618481 & 0.6251969494924084789 \\
	$n$ & 0.8833472048742723832 & 0.8854434199494314755 & 0.8870398387587905421 \\
	$o$ & 0.9205139849292685308 & 0.9201892361894765255 & 0.9199018958375438548 \\
	$p$ & 0.5059803484000719089 & 0.5082890237076752701 & 0.5100871017899912362 \\
	$q$ & 0.6865442453396316850 & 0.6858082276624378365 & 0.6850703583720513364 \\
	$r$ & 0.2170228473433140389 & 0.2193654799966509662 & 0.2211732981069128897 \\
	\hline\Tstrut
	$energy$ & -270.1179499592826915 & 968.7134553437876326 & 1019.829580588535623
\end{longtable}

\noindent
\textbf{NOTE:} Due to a problem with the Jacobian, the convergence software was used instead to find 1,021 digits for the parameters which was turned out to be laborious.

\noindent
\textit{\textbf{Symmetries -}}

The symmetry groups for 48 points are identical under all 3 potentials.

\begin{center}
	\begin{tabular}{l|l}
		\multicolumn{2}{c}{Symmetries - 48 points} \\
		\hline\Tstrut
		planes & [[16, 4], [48, 3]] \\[0.2ex]
		\hline\Tstrut
		Gram groups & [[24, 18], [48, 39]] \\
		\hline\Tstrut
		Polygons & [[3, 64], [4, 36]]
	\end{tabular}
\end{center}

\subsection{49 points}

The symmetry group for polygons hints that the structure might be one solitary point and 16 triangles. Indeed, that was the optimal configuration, as rotating a minimal energy point set into alignment with the triangles showed.

\begin{figure}[ht]
	\begin{center}
		\includegraphics[type=pdf,ext=pdf,read=pdf,height=1in,width=1in,angle=0]{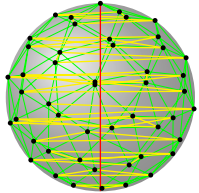}
		\caption{49 points.}
		\label{fig:49pts}
	\end{center}
\end{figure}

\noindent
\textit{\textbf{Algebraic parameterization -}}

It requires 32 parameters to adequately constrain the 49 point set, for the even number of 16 triangles embedded in the polyhedron.

The algebraic parameterized structure is given below:

\begin{longtable}[c]{r|ccc}
	\caption{Parameterization for 49 points} \\
	pt & $x$ & $y$ & $z$ \\
	\hline\vspace*{-2.2ex}
	\endfirsthead
	\multicolumn{4}{c}%
	{\tablename\ \thetable\ -- 49 points parameters -- \textit{continued}} \\
	pt & $x$ & $y$ & $z$ \\
	\hline\vspace*{-2.2ex}
	\endhead
	1 & $0$ & $0$ & $1$ \\[0.5ex]
	2 & $b$ & $0$ & $a$ \\[0.5ex]
	3 & $-\frac{b}{2}$ & $b\frac{\sqrt{3}}{2}$ & $a$ \\[0.5ex]
	4 & $-\frac{b}{2}$ & $-b\frac{\sqrt{3}}{2}$ & $a$ \\[0.5ex]
	5 & $d$ & $\sqrt{1-c^2-d^2}$ & $c$ \\[0.5ex]
	6 & $\frac{-d-\sqrt{3}\sqrt{1-c^2-d^2}}{2}$ & $\frac{d\sqrt{3}-\sqrt{1-c^2-d^2}}{2}$ & $c$ \\[0.5ex]
	7 & $\frac{-d+\sqrt{3}\sqrt{1-c^2-d^2}}{2}$ & $\frac{-d\sqrt{3}-\sqrt{1-c^2-d^2}}{2}$ & $c$ \\[0.5ex]
	8 & $f$ & $-\sqrt{1-e^2-f^2}$ & $e$ \\[0.5ex]
	9 & $\frac{-f-\sqrt{3}\sqrt{1-e^2-f^2}}{2}$ & $\frac{-f\sqrt{3}+\sqrt{1-e^2-f^2}}{2}$ & $e$ \\[0.5ex]
	10 & $\frac{-f+\sqrt{3}\sqrt{1-e^2-f^2}}{2}$ & $\frac{+f\sqrt{3}+\sqrt{1-e^2-f^2}}{2}$ & $e$ \\[0.5ex]
	11 & $h$ & $\sqrt{1-g^2-h^2}$ & $g$ \\[0.5ex]
	12 & $\frac{-h-\sqrt{3}\sqrt{1-g^2-h^2}}{2}$ & $\frac{h\sqrt{3}-\sqrt{1-g^2-h^2}}{2}$ & $g$ \\[0.5ex]
	13 & $\frac{-h+\sqrt{3}\sqrt{1-g^2-h^2}}{2}$ & $\frac{-h\sqrt{3}-\sqrt{1-g^2-h^2}}{2}$ & $g$ \\[0.5ex]
	14 & $j$ & $\sqrt{1-i^2-j^2}$ & $i$ \\[0.5ex]
	15 & $\frac{-j-\sqrt{3}\sqrt{1-i^2-j^2}}{2}$ & $\frac{j\sqrt{3}-\sqrt{1-i^2-j^2}}{2}$ & $i$ \\[0.5ex]
	16 & $\frac{-j+\sqrt{3}\sqrt{1-i^2-j^2}}{2}$ & $\frac{-j\sqrt{3}-\sqrt{1-i^2-j^2}}{2}$ & $i$ \\[0.5ex]
	17 & $l$ & $-\sqrt{1-k^2-l^2}$ & $k$ \\[0.5ex]
	18 & $\frac{-l-\sqrt{3}\sqrt{1-k^2-l^2}}{2}$ & $\frac{-l\sqrt{3}+\sqrt{1-k^2-l^2}}{2}$ & $k$ \\[0.5ex]
	19 & $\frac{-l+\sqrt{3}\sqrt{1-k^2-l^2}}{2}$ & $\frac{+l\sqrt{3}+\sqrt{1-k^2-l^2}}{2}$ & $k$ \\[0.5ex]
	20 & $n$ & $-\sqrt{1-m^2-n^2}$ & $m$ \\[0.5ex]
	21 & $\frac{-n-\sqrt{3}\sqrt{1-m^2-n^2}}{2}$ & $\frac{-n\sqrt{3}+\sqrt{1-m^2-n^2}}{2}$ & $m$ \\[0.5ex]
	22 & $\frac{-n+\sqrt{3}\sqrt{1-m^2-n^2}}{2}$ & $\frac{+n\sqrt{3}+\sqrt{1-m^2-n^2}}{2}$ & $m$ \\[0.5ex]
	23 & $p$ & $\sqrt{1-o^2-p^2}$ & $o$ \\[0.5ex]
	24 & $\frac{-p-\sqrt{3}\sqrt{1-o^2-p^2}}{2}$ & $\frac{p\sqrt{3}-\sqrt{1-o^2-p^2}}{2}$ & $o$ \\[0.5ex]
	25 & $\frac{-p+\sqrt{3}\sqrt{1-o^2-p^2}}{2}$ & $\frac{-p\sqrt{3}-\sqrt{1-o^2-p^2}}{2}$ & $o$ \\[0.5ex]
	26 & $r$ & $\sqrt{1-q^2-r^2}$ & $q$ \\[0.5ex]
	27 & $\frac{-r-\sqrt{3}\sqrt{1-q^2-r^2}}{2}$ & $\frac{r\sqrt{3}-\sqrt{1-q^2-r^2}}{2}$ & $q$ \\[0.5ex]
	28 & $\frac{-r+\sqrt{3}\sqrt{1-q^2-r^2}}{2}$ & $\frac{-r\sqrt{3}-\sqrt{1-q^2-r^2}}{2}$ & $q$ \\[0.5ex]
	29 & $t$ & $-\sqrt{1-s^2-t^2}$ & $s$ \\[0.5ex]
	30 & $\frac{-t-\sqrt{3}\sqrt{1-s^2-t^2}}{2}$ & $\frac{-t\sqrt{3}+\sqrt{1-s^2-t^2}}{2}$ & $s$ \\[0.5ex]
	31 & $\frac{-t+\sqrt{3}\sqrt{1-s^2-t^2}}{2}$ & $\frac{+t\sqrt{3}+\sqrt{1-s^2-t^2}}{2}$ & $s$ \\[0.5ex]
	32 & $v$ & $\sqrt{1-u^2-v^2}$ & $u$ \\[0.5ex]
	33 & $\frac{-v-\sqrt{3}\sqrt{1-u^2-v^2}}{2}$ & $\frac{v\sqrt{3}-\sqrt{1-u^2-v^2}}{2}$ & $u$ \\[0.5ex]
	34 & $\frac{-v+\sqrt{3}\sqrt{1-u^2-v^2}}{2}$ & $\frac{-v\sqrt{3}-\sqrt{1-u^2-v^2}}{2}$ & $u$ \\[0.5ex]
	35 & $x$ & $\sqrt{1-w^2-x^2}$ & $w$ \\[0.5ex]
	36 & $\frac{-x-\sqrt{3}\sqrt{1-w^2-x^2}}{2}$ & $\frac{x\sqrt{3}-\sqrt{1-w^2-x^2}}{2}$ & $w$ \\[0.5ex]
	37 & $\frac{-x+\sqrt{3}\sqrt{1-w^2-x^2}}{2}$ & $\frac{-x\sqrt{3}-\sqrt{1-w^2-x^2}}{2}$ & $w$ \\[0.5ex]
	38 & $z$ & $-\sqrt{1-y^2-z^2}$ & $y$ \\[0.5ex]
	39 & $\frac{-z-\sqrt{3}\sqrt{1-y^2-z^2}}{2}$ & $\frac{-z\sqrt{3}+\sqrt{1-y^2-z^2}}{2}$ & $y$ \\[0.5ex]
	40 & $\frac{-z+\sqrt{3}\sqrt{1-y^2-z^2}}{2}$ & $\frac{+z\sqrt{3}+\sqrt{1-y^2-z^2}}{2}$ & $y$ \\[0.5ex]
	41 & $B$ & $-\sqrt{1-A^2-B^2}$ & $A$ \\[0.5ex]
	42 & $\frac{-B-\sqrt{3}\sqrt{1-A^2-B^2}}{2}$ & $\frac{-B\sqrt{3}+\sqrt{1-A^2-B^2}}{2}$ & $A$ \\[0.5ex]
	43 & $\frac{-B+\sqrt{3}\sqrt{1-A^2-B^2}}{2}$ & $\frac{+B\sqrt{3}+\sqrt{1-A^2-B^2}}{2}$ & $A$ \\[0.5ex]
	44 & $D$ & $\sqrt{1-C^2-D^2}$ & $C$ \\[0.5ex]
	45 & $\frac{-D-\sqrt{3}\sqrt{1-C^2-D^2}}{2}$ & $\frac{D\sqrt{3}-\sqrt{1-C^2-D^2}}{2}$ & $C$ \\[0.5ex]
	46 & $\frac{-D+\sqrt{3}\sqrt{1-C^2-D^2}}{2}$ & $\frac{-D\sqrt{3}-\sqrt{1-C^2-D^2}}{2}$ & $C$ \\[0.5ex]
	47 & $F$ & $-\sqrt{1-E^2-F^2}$ & $E$ \\[0.5ex]
	48 & $\frac{-F-\sqrt{3}\sqrt{1-E^2-F^2}}{2}$ & $\frac{-F\sqrt{3}+\sqrt{1-E^2-F^2}}{2}$ & $E$ \\[0.5ex]
	49 & $\frac{-F+\sqrt{3}\sqrt{1-E^2-F^2}}{2}$ & $\frac{+F\sqrt{3}+\sqrt{1-E^2-F^2}}{2}$ & $E$
\end{longtable}

\noindent
\textit{\textbf{Parameterization values --}}

The values to 19 digits of the 32 parameters optimized for the minimal solutions of 49 points are:

\begin{longtable}[c]{c|c|c|c}
	\caption{Parameter values for 49 points} \\
	Parameter & log & 1/r & $1/r^2$ \\
	\hline\vspace*{-2.2ex}
	\endfirsthead
	\multicolumn{4}{c}%
	{\tablename\ \thetable\ -- 49 points parameters -- \textit{continued}} \\
	Parameter & log & 1/r & $1/r^2$ \\
	\hline\vspace*{-2.2ex}
	\endhead
	$a$ & 0.8793932629824731328 & 0.8792018226377346325 & 0.8791541480195582369 \\
	$b$ & 0.4760960922135770364 & 0.4764495304546803874 & 0.4765374948731785636 \\
	$c$ & 0.8137109170847478634 & 0.8157982924240811064 & 0.8177783481788855993 \\
	$d$ & 0.3161238981617116681 & 0.3128495487032312756 & 0.3099579960467028585 \\
	$e$ & 0.6085309852662297142 & 0.6085508844273354414 & 0.6087174241933589607 \\
	$f$ & 0.7161362722912167638 & 0.7172613587779483886 & 0.7179268099009106912 \\
	$g$ & 0.5600046405690230748 & 0.5613204297264755954 & 0.5625202003388636163 \\
	$h$ & 0.7881249661236182509 & 0.7859324832379948208 & 0.7841722610166736543 \\
	$i$ & 0.4221399993847407942 & 0.4223715856104007130 & 0.4227797911749838521 \\
	$j$ & 0.5475711031770011130 & 0.5430554560191582830 & 0.5396644365686091853 \\
	$k$ & 0.2130873989325435554 & 0.2135152649932475287 & 0.2138041648890469543 \\
	$l$ & 0.9724177353468279008 & 0.9726483617063639642 & 0.9727739053030582458 \\
	$m$ & 0.2021242646392268656 & 0.2012490117326880320 & 0.2003472030848887870 \\
	$n$ & 0.7643353698583752208 & 0.7682282819745653787 & 0.7711428562778039848 \\
	$o$ & 0.09483964870312775635 & 0.09393760913899789848 & 0.09304301855654624734 \\
	$p$ & 0.8784813445859752458 & 0.8767237270049579922 & 0.8756480162916912762 \\
	$q$ & -0.1288980921737543076 & -0.1280584235000077446 & -0.1274073537988331313 \\
	$r$ & 0.5244586197678084657 & 0.5195882780974520738 & 0.5160975032076731476 \\
	$s$ & -0.2586022719455941792 & -0.2600997162632330234 & -0.2618794589369973239 \\
	$t$ & 0.8521992187687367899 & 0.8543171605483494817 & 0.8555517689328800789 \\
	$u$ & -0.2634563076088891862 & -0.2639879441599654627 & -0.2647392774408586690 \\
	$v$ & 0.9624420737103483806 & 0.9619138338644125246 & 0.9614474630667256610 \\
	$w$ & -0.4284095129470309682 & -0.4288251319572060690 & -0.4292433382030535948 \\
	$x$ & 0.7549074785399459162 & 0.7519069282021759227 & 0.7497178690791479332 \\
	$y$ & -0.6134412031364999058 & -0.6145869048352859905 & -0.6156778166220145042 \\
	$z$ & 0.4962353823945488833 & 0.4971013992289302885 & 0.4971382898379245231 \\
	$A$ & -0.6748173185102826475 & -0.6742603371333260692 & -0.6741530469146048714 \\
	$B$ & 0.7246300642092983773 & 0.7257344598675112783 & 0.7261547182502060555 \\
	$C$ & -0.8118942938089824596 & -0.8123279285719047540 & -0.8128547470469438606 \\
	$D$ & 0.5017520060243855969 & 0.5007301850031944294 & 0.4998592320357899631 \\
	$E$ & -0.9476454507644124361 & -0.9476416280064657552 & -0.9477101806997068536 \\
	$F$ & 0.2696875985541636261 & 0.2708687504081964421 & 0.2716010625781302213 \\
	\hline\Tstrut
	$energy$ & -280.7019031182561337 & 1011.557182653571971 & 1069.559739809892139
\end{longtable}

\noindent
\textbf{NOTE:} A bug was found during the run for this configuration with the partial differentiation operator, disabling the ability to find these 32 parameters to more than 38 digits or so. By removing the fixed point and creating and then running the Jacobian as a matrix of all variables, higher precision was available for the parameters (digits: 48,299 - \textit{log}, 28,205 - \textit{Coulomb}, 28,378 - \textit{Inverse Square}).

\noindent
\textit{\textbf{Symmetries -}}

The symmetry groups for 49 points are identical under all 3 potentials.

\begin{center}
	\begin{tabular}{l|l}
		\multicolumn{2}{c}{Symmetries - 49 points} \\
		\hline\Tstrut
		planes & [[16, 1]] \\[0.2ex]
		\hline\Tstrut
		Gram groups & [[6, 392], [49, 1]] \\
		\hline\Tstrut
		Polygons & [[3, 16]]
	\end{tabular}
\end{center}

\subsection{50 points}
The optimal configuration for 50 points is interesting, the embedded polygons actually contain a set of 8 parallel regular hexagons, and 2 poles, thus making a balanced [1:6:6:6:6:6:6:6:6:1] arrangement. See figure \ref{fig:50pts} below, where the regular hexagons are highlighted in yellow.

\begin{figure}[ht]
	\begin{center}
		\includegraphics[type=pdf,ext=pdf,read=pdf,height=1in,width=1in,angle=0]{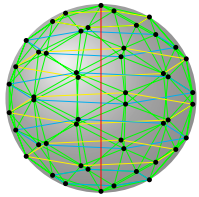}
		\caption{50 points.}
		\label{fig:50pts}
	\end{center}
\end{figure}

The parallel regular hexagons are exploited to create the algebraic structure.

\noindent
\textit{\textbf{Algebraic parameterization -}}

It was discovered, after aligning the 8 regular hexagons, that the 50 point set could be parameterized by 4 algebraic parameters $a$, $b$, $c$ and $d$, which are the distances along the polar axis.

The algebraic parameterized structure is given below:

\begin{longtable}[c]{r|ccc}
	\caption{Parameterization for 50 points} \\
	pt & $x$ & $y$ & $z$ \\
	\hline\vspace*{-2.2ex}
	\endfirsthead
	\multicolumn{4}{c}%
	{\tablename\ \thetable\ -- 50 points parameters -- \textit{continued\ldots}} \\
	pt & $x$ & $y$ & $z$ \\
	\hline\vspace*{-2.2ex}
	\endhead
	1 & $0$ & $0$ & $1$ \\[0.7ex]
	2 & $\sqrt{1-a^2}$ & $0$ & $a$ \\[0.7ex]
	3 & $\frac{\sqrt{1-a^2}}{2}$ & $\frac{\sqrt{3}}{2}\sqrt{1-a^2}$ & $a$ \\[0.7ex]
	4 & $-\frac{\sqrt{1-a^2}}{2}$ & $\frac{\sqrt{3}}{2}\sqrt{1-a^2}$ & $a$ \\[0.7ex]
	5 & $-\sqrt{1-a^2}$ & $0$ & $a$ \\[0.7ex]
	6 & $-\frac{\sqrt{1-a^2}}{2}$ & $-\frac{\sqrt{3}}{2}\sqrt{1-a^2}$ & $a$ \\[0.7ex]
	7 & $\frac{\sqrt{1-a^2}}{2}$ & $-\frac{\sqrt{3}}{2}\sqrt{1-a^2}$ & $a$ \\[0.7ex]
	8 & $0$ & $\sqrt{1-b^2}$ & $b$ \\[0.7ex]
	9 & $\frac{\sqrt{3}}{2}\sqrt{1-b^2}$ & $\frac{\sqrt{1-b^2}}{2}$ & $b$ \\[0.7ex]
	10 & $\frac{\sqrt{3}}{2}\sqrt{1-b^2}$ & $-\frac{\sqrt{1-b^2}}{2}$ & $b$ \\[0.7ex]
	11 & $0$ & $-\sqrt{1-b^2}$ & $b$ \\[0.7ex]
	12 & $-\frac{\sqrt{3}}{2}\sqrt{1-b^2}$ & $-\frac{\sqrt{1-b^2}}{2}$ & $b$ \\[0.7ex]
	13 & $-\frac{\sqrt{3}}{2}\sqrt{1-b^2}$ & $\frac{\sqrt{1-b^2}}{2}$ & $b$ \\[0.7ex]
	14 & $\sqrt{1-c^2}$ & $0$ & $c$ \\[0.7ex]
	15 & $\frac{\sqrt{1-c^2}}{2}$ & $\frac{\sqrt{3}}{2}\sqrt{1-c^2}$ & $c$ \\[0.7ex]
	16 & $-\frac{\sqrt{1-c^2}}{2}$ & $\frac{\sqrt{3}}{2}\sqrt{1-c^2}$ & $c$ \\[0.7ex]
	17 & $-\sqrt{1-c^2}$ & $0$ & $c$ \\[0.7ex]
	18 & $-\frac{\sqrt{1-c^2}}{2}$ & $-\frac{\sqrt{3}}{2}\sqrt{1-c^2}$ & $c$ \\[0.7ex]
	19 & $\frac{\sqrt{1-c^2}}{2}$ & $-\frac{\sqrt{3}}{2}\sqrt{1-c^2}$ & $c$ \\[0.7ex]
	20 & $0$ & $\sqrt{1-d^2}$ & $d$ \\[0.7ex]
	21 & $\frac{\sqrt{3}}{2}\sqrt{1-d^2}$ & $\frac{\sqrt{1-d^2}}{2}$ & $d$ \\[0.7ex]
	22 & $\frac{\sqrt{3}}{2}\sqrt{1-d^2}$ & $-\frac{\sqrt{1-d^2}}{2}$ & $d$ \\[0.7ex]
	23 & $0$ & $-\sqrt{1-d^2}$ & $d$ \\[0.7ex]
	24 & $-\frac{\sqrt{3}}{2}\sqrt{1-d^2}$ & $-\frac{\sqrt{1-d^2}}{2}$ & $d$ \\[0.7ex]
	25 & $-\frac{\sqrt{3}}{2}\sqrt{1-d^2}$ & $\frac{\sqrt{1-d^2}}{2}$ & $d$ \\[0.7ex]
	26 & $\sqrt{1-d^2}$ & $0$ & $-d$ \\[0.7ex]
	27 & $\frac{\sqrt{1-d^2}}{2}$ & $\frac{\sqrt{3}}{2}\sqrt{1-d^2}$ & $-d$ \\[0.7ex]
	28 & $-\frac{\sqrt{1-d^2}}{2}$ & $\frac{\sqrt{3}}{2}\sqrt{1-d^2}$ & $-d$ \\[0.7ex]
	29 & $-\sqrt{1-d^2}$ & $0$ & $-d$ \\[0.7ex]
	30 & $-\frac{\sqrt{1-d^2}}{2}$ & $-\frac{\sqrt{3}}{2}\sqrt{1-d^2}$ & $-d$ \\[0.7ex]
	31 & $\frac{\sqrt{1-d^2}}{2}$ & $-\frac{\sqrt{3}}{2}\sqrt{1-d^2}$ & $-d$ \\[0.7ex]
	32 & $0$ & $\sqrt{1-c^2}$ & $-c$ \\[0.7ex]
	33 & $\frac{\sqrt{3}}{2}\sqrt{1-c^2}$ & $\frac{\sqrt{1-c^2}}{2}$ & $-c$ \\[0.7ex]
	34 & $\frac{\sqrt{3}}{2}\sqrt{1-c^2}$ & $-\frac{\sqrt{1-c^2}}{2}$ & $-c$ \\[0.7ex]
	35 & $0$ & $-\sqrt{1-c^2}$ & $-c$ \\[0.7ex]
	36 & $-\frac{\sqrt{3}}{2}\sqrt{1-c^2}$ & $-\frac{\sqrt{1-c^2}}{2}$ & $-c$ \\[0.7ex]
	37 & $-\frac{\sqrt{3}}{2}\sqrt{1-c^2}$ & $\frac{\sqrt{1-c^2}}{2}$ & $-c$ \\[0.7ex]
	38 & $\sqrt{1-b^2}$ & $0$ & $-b$ \\[0.7ex]
	39 & $\frac{\sqrt{1-b^2}}{2}$ & $\frac{\sqrt{3}}{2}\sqrt{1-b^2}$ & $-b$ \\[0.7ex]
	40 & $-\frac{\sqrt{1-b^2}}{2}$ & $\frac{\sqrt{3}}{2}\sqrt{1-b^2}$ & $-b$ \\[0.7ex]
	41 & $-\sqrt{1-b^2}$ & $0$ & $-b$ \\[0.7ex]
	42 & $-\frac{\sqrt{1-b^2}}{2}$ & $-\frac{\sqrt{3}}{2}\sqrt{1-b^2}$ & $-b$ \\[0.7ex]
	43 & $\frac{\sqrt{1-b^2}}{2}$ & $-\frac{\sqrt{3}}{2}\sqrt{1-b^2}$ & $-b$ \\[0.7ex]
	44 & $0$ & $\sqrt{1-a^2}$ & $-a$ \\[0.7ex]
	45 & $\frac{\sqrt{3}}{2}\sqrt{1-a^2}$ & $\frac{\sqrt{1-a^2}}{2}$ & $-a$ \\[0.7ex]
	46 & $\frac{\sqrt{3}}{2}\sqrt{1-a^2}$ & $-\frac{\sqrt{1-a^2}}{2}$ & $-a$ \\[0.7ex]
	47 & $0$ & $-\sqrt{1-a^2}$ & $-a$ \\[0.7ex]
	48 & $-\frac{\sqrt{3}}{2}\sqrt{1-a^2}$ & $-\frac{\sqrt{1-a^2}}{2}$ & $-a$ \\[0.7ex]
	49 & $-\frac{\sqrt{3}}{2}\sqrt{1-a^2}$ & $\frac{\sqrt{1-a^2}}{2}$ & $-a$ \\[0.7ex]
	50 & $0$ & $0$ & $-1$ \\[0.7ex]
\end{longtable}

Despite obtaining the spherical codes to 50,014 digits, the degree of the algebraic numbers is $>420$.

\noindent
\textit{\textbf{Parameters and Minimal Energy values -}}

\begin{longtable}[c]{c|c|c|c}
	\caption{Parameter values for 50 points} \\
	Parameter & log & 1/r & $1/r^2$ \\
	\hline\vspace*{-2.2ex}
	\endfirsthead
	\multicolumn{4}{c}%
	{\tablename\ \thetable\ -- 50 points parameters -- \textit{continued}} \\
	Parameter & log & 1/r & $1/r^2$ \\
	\hline\vspace*{-2.2ex}
	\endhead
	$a$ & 0.8515838011853908757 & 0.8513832027240754402 & 0.8514411311335073360 \\
	$b$ & 0.5845080765467688786 & 0.5859466784603221965 & 0.5875147416660090262 \\
	$c$ & 0.3823580555306074125 & 0.3828907922740789633 & 0.3834970349379620340 \\
	$d$ & 0.1056903533827164585 & 0.1079740016327089417 & 0.1100866382336972714 \\
	\hline\Tstrut
	$energy$ & -291.5286006577343690 & 1055.182314726296099 & 1119.599506531607568
\end{longtable}

\noindent
\textit{\textbf{Symmetries -}}

The symmetry groups for 50 points are identical under all 3 potentials.

\begin{center}
	\begin{tabular}{l|l}
		\multicolumn{2}{c}{Symmetries - 50 points} \\
		\hline\Tstrut
		planes & [[4, 1008], [120, 6], [160, 1]] \\[0.2ex]
		\hline\Tstrut
		Gram groups & [[2, 1], [12, 4], [24, 40], [48, 30], [50, 1]] \\
		\hline\Tstrut
		Polygons & [[4, 1008], [6, 8], [10, 6]]
	\end{tabular}
\end{center}

\subsection{51 points}
The symmetry groups for polygons shows that there are 17 parallel triangles embedded in the configuration of 51 points, and that is the preferred alignment for parameterization. This is similar to the configuration with 45 points, which required 15 triangles. The arrangement is [3:3:3:3:3:3:3:3:3:3:3:3:3:3:3:3:3].

\begin{figure}[ht]
	\begin{center}
		\includegraphics[type=pdf,ext=pdf,read=pdf,height=1in,width=1in,angle=0]{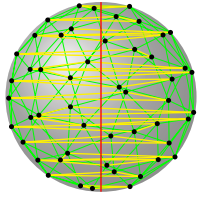}
		\caption{51 points.}
		\label{fig:51pts}
	\end{center}
\end{figure}

\noindent
\textit{\textbf{Algebraic parameterization -}}

It took 16 parameters to adequately constrain the polyhedron. The parameters adjust the height of the triangle along the z-axis, and the rotation of the triangle from the standard [1,0,0] vertex position along the x-axis.

The algebraic parameterized structure is given below:

\begin{longtable}[c]{r|ccc}
	\caption{Parameterization for 51 points} \\
	pt & $x$ & $y$ & $z$ \\
	\hline\vspace*{-2.2ex}
	\endfirsthead
	\multicolumn{4}{c}%
	{\tablename\ \thetable\ -- 51 points parameters -- \textit{continued}} \\
	pt & $x$ & $y$ & $z$ \\
	\hline\vspace*{-2.2ex}
	\endhead
	1 & $b$ & $-\sqrt{1-a^2-b^2}$ & $a$ \\[0.7ex]
	2 & $\frac{-b-\sqrt{3}\sqrt{1-a^2-b^2}}{2}$ & $\frac{-b\sqrt{3}+\sqrt{1-a^2-b^2}}{2}$ & $a$ \\[0.7ex]
	3 & $\frac{-b+\sqrt{3}\sqrt{1-a^2-b^2}}{2}$ & $\frac{+b\sqrt{3}+\sqrt{1-a^2-b^2}}{2}$ & $a$ \\[0.7ex]
	4 & $d$ & $-\sqrt{1-c^2-d^2}$ & $c$ \\[0.7ex]
	5 & $\frac{-d-\sqrt{3}\sqrt{1-c^2-d^2}}{2}$ & $\frac{-d\sqrt{3}+\sqrt{1-c^2-d^2}}{2}$ & $c$ \\[0.7ex]
	6 & $\frac{-d+\sqrt{3}\sqrt{1-c^2-d^2}}{2}$ & $\frac{+d\sqrt{3}+\sqrt{1-c^2-d^2}}{2}$ & $c$ \\[0.7ex]
	7 & $f$ & $\sqrt{1-e^2-f^2}$ & $e$ \\[0.7ex]
	8 & $\frac{-f-\sqrt{3}\sqrt{1-e^2-f^2}}{2}$ & $\frac{+f\sqrt{3}-\sqrt{1-e^2-f^2}}{2}$ & $e$ \\[0.7ex]
	9 & $\frac{-f+\sqrt{3}\sqrt{1-e^2-f^2}}{2}$ & $\frac{-f\sqrt{3}-\sqrt{1-e^2-f^2}}{2}$ & $e$ \\[0.7ex]
	10 & $h$ & $-\sqrt{1-g^2-h^2}$ & $g$ \\[0.7ex]
	11 & $\frac{-h-\sqrt{3}\sqrt{1-g^2-h^2}}{2}$ & $\frac{-h\sqrt{3}+\sqrt{1-g^2-h^2}}{2}$ & $g$ \\[0.7ex]
	12 & $\frac{-h+\sqrt{3}\sqrt{1-g^2-h^2}}{2}$ & $\frac{+h\sqrt{3}+\sqrt{1-g^2-h^2}}{2}$ & $g$ \\[0.7ex]
	13 & $j$ & $-\sqrt{1-i^2-j^2}$ & $i$ \\[0.7ex]
	14 & $\frac{-j-\sqrt{3}\sqrt{1-i^2-j^2}}{2}$ & $\frac{-j\sqrt{3}+\sqrt{1-i^2-j^2}}{2}$ & $i$ \\[0.7ex]
	15 & $\frac{-j+\sqrt{3}\sqrt{1-i^2-j^2}}{2}$ & $\frac{+j\sqrt{3}+\sqrt{1-i^2-j^2}}{2}$ & $i$ \\[0.7ex]
	16 & $l$ & $\sqrt{1-k^2-l^2}$ & $k$ \\[0.7ex]
	17 & $\frac{-l-\sqrt{3}\sqrt{1-k^2-l^2}}{2}$ & $\frac{+l\sqrt{3}-\sqrt{1-k^2-l^2}}{2}$ & $k$ \\[0.7ex]
	18 & $\frac{-l+\sqrt{3}\sqrt{1-k^2-l^2}}{2}$ & $\frac{-l\sqrt{3}-\sqrt{1-k^2-l^2}}{2}$ & $k$ \\[0.7ex]
	19 & $n$ & $\sqrt{1-m^2-n^2}$ & $m$ \\[0.7ex]
	20 & $\frac{-n-\sqrt{3}\sqrt{1-m^2-n^2}}{2}$ & $\frac{+n\sqrt{3}-\sqrt{1-m^2-n^2}}{2}$ & $m$ \\[0.7ex]
	21 & $\frac{-n+\sqrt{3}\sqrt{1-m^2-n^2}}{2}$ & $\frac{-n\sqrt{3}-\sqrt{1-m^2-n^2}}{2}$ & $m$ \\[0.7ex]
	22 & $p$ & $-\sqrt{1-o^2-p^2}$ & $o$ \\[0.7ex]
	23 & $\frac{-p-\sqrt{3}\sqrt{1-o^2-p^2}}{2}$ & $\frac{-p\sqrt{3}+\sqrt{1-o^2-p^2}}{2}$ & $o$ \\[0.7ex]
	24 & $\frac{-p+\sqrt{3}\sqrt{1-o^2-p^2}}{2}$ & $\frac{+p\sqrt{3}+\sqrt{1-o^2-p^2}}{2}$ & $o$ \\[0.7ex]
	25 & $1$ & $0$ & $0$ \\[0.7ex]
	26 & $-\frac{1}{2}$ & $\frac{\sqrt{3}}{2}$ & $0$ \\[0.7ex]
	27 & $-\frac{1}{2}$ & $-\frac{\sqrt{3}}{2}$ & $0$ \\[0.7ex]
	28 & $p$ & $\sqrt{1-o^2-p^2}$ & $-o$ \\[0.7ex]
	29 & $\frac{-p-\sqrt{3}\sqrt{1-o^2-p^2}}{2}$ & $\frac{+p\sqrt{3}-\sqrt{1-o^2-p^2}}{2}$ & $-o$ \\[0.7ex]
	30 & $\frac{-p+\sqrt{3}\sqrt{1-o^2-p^2}}{2}$ & $\frac{-p\sqrt{3}-\sqrt{1-o^2-p^2}}{2}$ & $-o$ \\[0.7ex]
	31 & $n$ & $-\sqrt{1-m^2-n^2}$ & $-m$ \\[0.7ex]
	32 & $\frac{-n-\sqrt{3}\sqrt{1-m^2-n^2}}{2}$ & $\frac{-n\sqrt{3}+\sqrt{1-m^2-n^2}}{2}$ & $-m$ \\[0.7ex]
	33 & $\frac{-n+\sqrt{3}\sqrt{1-m^2-n^2}}{2}$ & $\frac{+n\sqrt{3}+\sqrt{1-m^2-n^2}}{2}$ & $-m$ \\[0.7ex]
	34 & $l$ & $-\sqrt{1-k^2-l^2}$ & $-k$ \\[0.7ex]
	35 & $\frac{-l-\sqrt{3}\sqrt{1-k^2-l^2}}{2}$ & $\frac{-l\sqrt{3}+\sqrt{1-k^2-l^2}}{2}$ & $-k$ \\[0.7ex]
	36 & $\frac{-l+\sqrt{3}\sqrt{1-k^2-l^2}}{2}$ & $\frac{+l\sqrt{3}+\sqrt{1-k^2-l^2}}{2}$ & $-k$ \\[0.7ex]
	37 & $j$ & $\sqrt{1-i^2-j^2}$ & $-i$ \\[0.7ex]
	38 & $\frac{-j-\sqrt{3}\sqrt{1-i^2-j^2}}{2}$ & $\frac{+j\sqrt{3}-\sqrt{1-i^2-j^2}}{2}$ & $-i$ \\[0.7ex]
	39 & $\frac{-j+\sqrt{3}\sqrt{1-i^2-j^2}}{2}$ & $\frac{-j\sqrt{3}-\sqrt{1-i^2-j^2}}{2}$ & $-i$ \\[0.7ex]
	40 & $h$ & $\sqrt{1-g^2-h^2}$ & $-g$ \\[0.7ex]
	41 & $\frac{-h-\sqrt{3}\sqrt{1-g^2-h^2}}{2}$ & $\frac{+h\sqrt{3}-\sqrt{1-g^2-h^2}}{2}$ & $-g$ \\[0.7ex]
	42 & $\frac{-h+\sqrt{3}\sqrt{1-g^2-h^2}}{2}$ & $\frac{-h\sqrt{3}-\sqrt{1-g^2-h^2}}{2}$ & $-g$ \\[0.7ex]
	43 & $f$ & $-\sqrt{1-e^2-f^2}$ & $-e$ \\[0.7ex]
	44 & $\frac{-f-\sqrt{3}\sqrt{1-e^2-f^2}}{2}$ & $\frac{-f\sqrt{3}+\sqrt{1-e^2-f^2}}{2}$ & $-e$ \\[0.7ex]
	45 & $\frac{-f+\sqrt{3}\sqrt{1-e^2-f^2}}{2}$ & $\frac{+f\sqrt{3}+\sqrt{1-e^2-f^2}}{2}$ & $-e$ \\[0.7ex]
	46 & $d$ & $\sqrt{1-c^2-d^2}$ & $-c$ \\[0.7ex]
	47 & $\frac{-d-\sqrt{3}\sqrt{1-c^2-d^2}}{2}$ & $\frac{+d\sqrt{3}-\sqrt{1-c^2-d^2}}{2}$ & $-c$ \\[0.7ex]
	48 & $\frac{-d+\sqrt{3}\sqrt{1-c^2-d^2}}{2}$ & $\frac{-d\sqrt{3}-\sqrt{1-c^2-d^2}}{2}$ & $-c$ \\[0.7ex]
	49 & $b$ & $\sqrt{1-a^2-b^2}$ & $-a$ \\[0.7ex]
	50 & $\frac{-b-\sqrt{3}\sqrt{1-a^2-b^2}}{2}$ & $\frac{+b\sqrt{3}-\sqrt{1-a^2-b^2}}{2}$ & $-a$ \\[0.7ex]
	51 & $\frac{-b+\sqrt{3}\sqrt{1-a^2-b^2}}{2}$ & $\frac{-b\sqrt{3}-\sqrt{1-a^2-b^2}}{2}$ & $-a$
\end{longtable}

The values to 19 digits of the 16 parameters $a-p$ optimized for the minimal solutions of 51 points are:

\begin{longtable}[c]{c|c|c|c}
	\caption{Parameter values for 51 points} \\
	Parameter & log & 1/r & $1/r^2$ \\
	\hline\vspace*{-2.2ex}
	\endfirsthead
	\multicolumn{4}{c}%
	{\tablename\ \thetable\ -- 51 points parameters -- \textit{continued}} \\
	Parameter & log & 1/r & $1/r^2$ \\
	\hline\vspace*{-2.2ex}
	\endhead
	$a$ & 0.9500652207774683361 & 0.9500282116483405921 & 0.9500374260030636337 \\
	$b$ & 0.1646528481400221695 & 0.1639569774313236922 & 0.1630127817167423243 \\
	$c$ & 0.8186390936870304855 & 0.8188298973341696733 & 0.8191941790044680608 \\
	$d$ & 0.5741910519933984701 & 0.5739410414981816943 & 0.5734303120729397311 \\
	$e$ & 0.6832630240617076390 & 0.6834505906070375497 & 0.6837609525896807308 \\
	$f$ & 0.5588323177474157581 & 0.5577214009664225186 & 0.5568250061658443812 \\
	$g$ & 0.6304195432017679592 & 0.6310543252073115408 & 0.6317123645050443373 \\
	$h$ & 0.5770709334741273205 & 0.5784429813501436479 & 0.5794797231422024974 \\
	$i$ & 0.4589774254146343573 & 0.4587811485513831032 & 0.4587428395797030510 \\
	$j$ & 0.8789825263781554468 & 0.8791127250788025677 & 0.8791756144406620436 \\
	$k$ & 0.3201548901453516992 & 0.3212708195372814348 & 0.3219603157388288894 \\
	$l$ & 0.8683386493234643667 & 0.8688030870042235108 & 0.8692352775194469122 \\
	$m$ & 0.2484977046437503842 & 0.2498059292228536024 & 0.2510107591820098522 \\
	$n$ & 0.5452669446733529639 & 0.5461552655574545929 & 0.5469419915411525667 \\
	$o$ & 0.1541160520726872650 & 0.1536691010218526061 & 0.1536525290746748091 \\
	$p$ & 0.8201413317927946452 & 0.8204747398807752826 & 0.8207839753618648537 \\
	\hline\Tstrut
	$energy$ & -302.5336734554331932 & 1099.819290318898223 & 1171.328381382107580
\end{longtable}

All 16 parameters have 50,014 digits precision, for all 3 potentials, but the degree of the algebraic polynomials are $>360$.
The minimal energies also have been found to this precision, but no algebraic polynomials have been found.

\noindent
\textit{\textbf{Symmetries -}}

The symmetry groups for 51 points are identical under all 3 potentials.

\begin{center}
	\begin{tabular}{l|l}
		\multicolumn{2}{c}{Symmetries - 51 points} \\
		\hline\Tstrut
		planes & [[17, 1]] \\[0.2ex]
		\hline\Tstrut
		Gram groups & [[6, 25], [12, 200], [51, 1]] \\
		\hline\Tstrut
		Polygons & [[3, 17]]
	\end{tabular}
\end{center}

\subsection{52 points}
The configuration for 52 points contains 17 parallel triangles also, just as in the configuration for 51 points. However a lone point has been added, and as expected, it migrated to a pole and pushed all 17 parallel triangles a bit closer to the opposite pole to compensate for the minimal energy configuration. The arrangement of polygons is [1:3:3:3:3:3:3:3:3:3:3:3:3:3:3:3:3:3].

The configuration of 52 points shows what happens when 1 solitary point is added to a balanced configuration with even parities. The 51 points still want to stay in the 17 triangles, pushing the lone point away (so to speak).

\begin{figure}[!ht]
	\begin{center}
		\includegraphics[type=pdf,ext=pdf,read=pdf,height=1in,width=1in,angle=0]{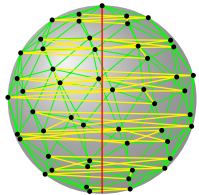}
		\caption{52 points.}
		\label{fig:52pts}
	\end{center}
\end{figure}

\noindent
\textit{\textbf{Algebraic parameterization -}}

It took 34 parameters to adequately constrain the polyhedron, known to 50,014 digits for all 3 power laws.

The algebraic parameterized structure is given below:

\begin{longtable}[c]{r|ccc}
	\caption{Parameterization for 52 points} \\
	pt & $x$ & $y$ & $z$ \\
	\hline\vspace*{-2.2ex}
	\endfirsthead
	\multicolumn{4}{c}%
	{\tablename\ \thetable\ -- 52 points parameters -- \textit{continued}} \\
	pt & $x$ & $y$ & $z$ \\
	\hline\vspace*{-2.2ex}
	\endhead
1 & $0$ & $0$ & $1$ \\[0.7ex]
2 & $b$ & $\sqrt{1-a^2-b^2}$ & $a$ \\[0.7ex]
3 & $\frac{-b-\sqrt{3}\sqrt{1-a^2-b^2}}{2}$ & $\frac{b\sqrt{3}-\sqrt{1-a^2-b^2}}{2}$ & $a$ \\[0.7ex]
4 & $\frac{-b+\sqrt{3}\sqrt{1-a^2-b^2}}{2}$ & $\frac{-b\sqrt{3}-\sqrt{1-a^2-b^2}}{2}$ & $a$ \\[0.7ex]
5 & $d$ & $-\sqrt{1-c^2-d^2}$ & $c$ \\[0.7ex]
6 & $\frac{-d-\sqrt{3}\sqrt{1-c^2-d^2}}{2}$ & $\frac{-d\sqrt{3}+\sqrt{1-c^2-d^2}}{2}$ & $c$ \\[0.7ex]
7 & $\frac{-d+\sqrt{3}\sqrt{1-c^2-d^2}}{2}$ & $\frac{d\sqrt{3}+\sqrt{1-c^2-d^2}}{2}$ & $c$ \\[0.7ex]
8 & $f$ & $\sqrt{1-e^2-f^2}$ & $e$ \\[0.7ex]
9 & $\frac{-f-\sqrt{3}\sqrt{1-e^2-f^2}}{2}$ & $\frac{f\sqrt{3}-\sqrt{1-e^2-f^2}}{2}$ & $e$ \\[0.7ex]
10 & $\frac{-f+\sqrt{3}\sqrt{1-e^2-f^2}}{2}$ & $\frac{-f\sqrt{3}-\sqrt{1-e^2-f^2}}{2}$ & $e$ \\[0.7ex]
11 & $h$ & $-\sqrt{1-g^2-h^2}$ & $g$ \\[0.7ex]
12 & $\frac{-h-\sqrt{3}\sqrt{1-g^2-h^2}}{2}$ & $\frac{-h\sqrt{3}+\sqrt{1-g^2-h^2}}{2}$ & $g$ \\[0.7ex]
13 & $\frac{-h+\sqrt{3}\sqrt{1-g^2-h^2}}{2}$ & $\frac{h\sqrt{3}+\sqrt{1-g^2-h^2}}{2}$ & $g$ \\[0.7ex]
14 & $j$ & $-\sqrt{1-i^2-j^2}$ & $i$ \\[0.7ex]
15 & $\frac{-j-\sqrt{3}\sqrt{1-i^2-j^2}}{2}$ & $\frac{-j\sqrt{3}+\sqrt{1-i^2-j^2}}{2}$ & $i$ \\[0.7ex]
16 & $\frac{-j+\sqrt{3}\sqrt{1-i^2-j^2}}{2}$ & $\frac{j\sqrt{3}+\sqrt{1-i^2-j^2}}{2}$ & $i$ \\[0.7ex]
17 & $l$ & $\sqrt{1-k^2-l^2}$ & $k$ \\[0.7ex]
18 & $\frac{-l-\sqrt{3}\sqrt{1-k^2-l^2}}{2}$ & $\frac{l\sqrt{3}-\sqrt{1-k^2-l^2}}{2}$ & $k$ \\[0.7ex]
19 & $\frac{-l+\sqrt{3}\sqrt{1-k^2-l^2}}{2}$ & $\frac{-l\sqrt{3}-\sqrt{1-k^2-l^2}}{2}$ & $k$ \\[0.7ex]
20 & $n$ & $\sqrt{1-m^2-n^2}$ & $m$ \\[0.7ex]
21 & $\frac{-n-\sqrt{3}\sqrt{1-m^2-n^2}}{2}$ & $\frac{n\sqrt{3}-\sqrt{1-m^2-n^2}}{2}$ & $m$ \\[0.7ex]
22 & $\frac{-n+\sqrt{3}\sqrt{1-m^2-n^2}}{2}$ & $\frac{-n\sqrt{3}-\sqrt{1-m^2-n^2}}{2}$ & $m$ \\[0.7ex]
23 & $p$ & $-\sqrt{1-o^2-p^2}$ & $o$ \\[0.7ex]
24 & $\frac{-p-\sqrt{3}\sqrt{1-o^2-p^2}}{2}$ & $\frac{-p\sqrt{3}+\sqrt{1-o^2-p^2}}{2}$ & $o$ \\[0.7ex]
25 & $\frac{-p+\sqrt{3}\sqrt{1-o^2-p^2}}{2}$ & $\frac{p\sqrt{3}+\sqrt{1-o^2-p^2}}{2}$ & $o$ \\[0.7ex]
26 & $r$ & $-\sqrt{1-q^2-r^2}$ & $q$ \\[0.7ex]
27 & $\frac{-r-\sqrt{3}\sqrt{1-q^2-r^2}}{2}$ & $\frac{-r\sqrt{3}+\sqrt{1-q^2-r^2}}{2}$ & $q$ \\[0.7ex]
28 & $\frac{-r+\sqrt{3}\sqrt{1-q^2-r^2}}{2}$ & $\frac{r\sqrt{3}+\sqrt{1-q^2-r^2}}{2}$ & $q$ \\[0.7ex]
29 & $t$ & $\sqrt{1-s^2-t^2}$ & $s$ \\[0.7ex]
30 & $\frac{-t-\sqrt{3}\sqrt{1-s^2-t^2}}{2}$ & $\frac{t\sqrt{3}-\sqrt{1-s^2-t^2}}{2}$ & $s$ \\[0.7ex]
31 & $\frac{-t+\sqrt{3}\sqrt{1-s^2-t^2}}{2}$ & $\frac{-t\sqrt{3}-\sqrt{1-s^2-t^2}}{2}$ & $s$ \\[0.7ex]
32 & $v$ & $-\sqrt{1-u^2-v^2}$ & $u$ \\[0.7ex]
33 & $\frac{-v-\sqrt{3}\sqrt{1-u^2-v^2}}{2}$ & $\frac{-v\sqrt{3}+\sqrt{1-u^2-v^2}}{2}$ & $u$ \\[0.7ex]
34 & $\frac{-v+\sqrt{3}\sqrt{1-u^2-v^2}}{2}$ & $\frac{v\sqrt{3}+\sqrt{1-u^2-v^2}}{2}$ & $u$ \\[0.7ex]
35 & $x$ & $\sqrt{1-w^2-x^2}$ & $w$ \\[0.7ex]
36 & $\frac{-x-\sqrt{3}\sqrt{1-w^2-x^2}}{2}$ & $\frac{x\sqrt{3}-\sqrt{1-w^2-x^2}}{2}$ & $w$ \\[0.7ex]
37 & $\frac{-x+\sqrt{3}\sqrt{1-w^2-x^2}}{2}$ & $\frac{-x\sqrt{3}-\sqrt{1-w^2-x^2}}{2}$ & $w$ \\[0.7ex]
38 & $z$ & $-\sqrt{1-y^2-z^2}$ & $y$ \\[0.7ex]
39 & $\frac{-z-\sqrt{3}\sqrt{1-y^2-z^2}}{2}$ & $\frac{-z\sqrt{3}+\sqrt{1-y^2-z^2}}{2}$ & $y$ \\[0.7ex]
40 & $\frac{-z+\sqrt{3}\sqrt{1-y^2-z^2}}{2}$ & $\frac{z\sqrt{3}+\sqrt{1-y^2-z^2}}{2}$ & $y$ \\[0.7ex]
41 & $B$ & $\sqrt{1-A^2-B^2}$ & $A$ \\[0.7ex]
42 & $\frac{-B-\sqrt{3}\sqrt{1-A^2-B^2}}{2}$ & $\frac{B\sqrt{3}-\sqrt{1-A^2-B^2}}{2}$ & $A$ \\[0.7ex]
43 & $\frac{-B+\sqrt{3}\sqrt{1-A^2-B^2}}{2}$ & $\frac{-B\sqrt{3}-\sqrt{1-A^2-B^2}}{2}$ & $A$ \\[0.7ex]
44 & $D$ & $-\sqrt{1-C^2-D^2}$ & $C$ \\[0.7ex]
45 & $\frac{-D-\sqrt{3}\sqrt{1-C^2-D^2}}{2}$ & $\frac{-D\sqrt{3}+\sqrt{1-C^2-D^2}}{2}$ & $C$ \\[0.7ex]
46 & $\frac{-D+\sqrt{3}\sqrt{1-C^2-D^2}}{2}$ & $\frac{D\sqrt{3}+\sqrt{1-C^2-D^2}}{2}$ & $C$ \\[0.7ex]
47 & $F$ & $\sqrt{1-E^2-F^2}$ & $E$ \\[0.7ex]
48 & $\frac{-F-\sqrt{3}\sqrt{1-E^2-F^2}}{2}$ & $\frac{F\sqrt{3}-\sqrt{1-E^2-F^2}}{2}$ & $E$ \\[0.7ex]
49 & $\frac{-F+\sqrt{3}\sqrt{1-E^2-F^2}}{2}$ & $\frac{-F\sqrt{3}-\sqrt{1-E^2-F^2}}{2}$ & $E$ \\[0.7ex]
50 & $H$ & $\sqrt{1-G^2-H^2}$ & $G$ \\[0.7ex]
51 & $\frac{-H-\sqrt{3}\sqrt{1-G^2-H^2}}{2}$ & $\frac{H\sqrt{3}-\sqrt{1-G^2-H^2}}{2}$ & $G$ \\[0.7ex]
52 & $\frac{-H+\sqrt{3}\sqrt{1-G^2-H^2}}{2}$ & $\frac{-H\sqrt{3}-\sqrt{1-G^2-H^2}}{2}$ & $G$
\end{longtable}

The values to 19 digits of the 34 parameters $a-H$ optimized for the minimal solutions of 52 points are:

\begin{longtable}[c]{c|c|c|c}
	\caption{Parameter values for 52 points} \\
	Parameter & log & 1/r & $1/r^2$ \\
	\hline\vspace*{-2.2ex}
	\endfirsthead
	\multicolumn{4}{c}%
	{\tablename\ \thetable\ -- 52 points parameters -- \textit{continued}} \\
	Parameter & log & 1/r & $1/r^2$ \\
	\hline\vspace*{-2.2ex}
	\endhead
	$a$ & 0.8706630011276560540 & 0.8723519785822495549 & 0.8738459352037142167 \\
	$b$ & 0.3058683137869242888 & 0.4831238878900578207 & 0.3080817730531670942 \\
	$c$ & 0.8484645284123887382 & 0.8484664910776997699 & 0.8488471726513630728 \\
	$d$ & 0.5225702255327457643 & 0.3361375252612937478 & 0.5207457510779873333 \\
	$e$ & 0.6015576073987396212 & 0.6026753986972674138 & 0.6036865799410745010 \\
	$f$ & 0.7140271626487281802 & 0.04781361365914809882 & 0.7150793584878285535 \\
	$g$ & 0.5833701079946160612 & 0.5856095541115671389 & 0.5874101044771674286 \\
	$h$ & 0.5435720099846360477 & 0.7932534931832303115 & 0.5366078366373861738 \\
	$i$ & 0.5040041862815182345 & 0.5018706754352640984 & 0.5003549756052836641 \\
	$j$ & 0.8427352880201589519 & 0.5854592126228281719 & 0.8443174890015326976 \\
	$k$ & 0.2559195170207020609 & 0.2534462735644271870 & 0.2512709103348453451 \\
	$l$ & 0.6390376129557583582 & -0.3084119808673279526 & 0.6438983071533074998 \\
	$m$ & 0.1908228343504637555 & 0.1897410295901644260 & 0.1885874918858391005 \\
	$n$ & 0.9497525527553429692 & 0.6913853557169998241 & 0.9504214540969711631 \\
	$o$ & 0.1310057787869638767 & 0.1334432919276150316 & 0.1348439491938194020 \\
	$p$ & 0.7473074950384233893 & 0.9364302200454692711 & 0.7476328231305885968 \\
	$q$ & 0.03125817649702673287 & 0.02908229297881019834 & 0.02736435025542482835 \\
	$r$ & 0.9746557858317675961 & 0.2990203430729341725 & 0.9751025032519538915 \\
	$s$ & -0.1828684817110674671 & -0.1869330581786416249 & -0.1903444166730965433 \\
	$t$ & 0.8255355977473824818 & 0.8739681452333602791 & 0.8276992180749213761 \\
	$u$ & -0.2715343734588587209 & -0.2700033045533876398 & -0.2694970936484554191 \\
	$v$ & 0.5000839927501612914 & -0.4619931711964975103 & 0.4973892359385221717 \\
	$w$ & -0.3857234863599258731 & -0.3872818124566919064 & -0.3879682133809512694 \\
	$x$ & 0.9198246249202365391 & 0.5213798463964816191 & 0.9192471047859312750 \\
	$y$ & -0.3963519344729299641 & -0.3938072174446410203 & -0.3920825776355011644 \\
	$z$ & 0.8125509955743114253 & 0.7775882207702676006 & 0.8119100018491940300 \\
	$A$ & -0.6467935116218329224 & -0.6487683444659193073 & -0.6501922799629996047 \\
	$B$ & 0.6401710405118081767 & 0.6771060965321145290 & 0.6410606988265590544 \\
	$C$ & -0.7150341921860741312 & -0.7114755064516761238 & -0.7088680497803384494 \\
	$D$ & 0.4535969026111252927 & 0.6912129949953693095 & 0.4536583943272281180 \\
	$E$ & -0.7968037472875334574 & -0.7970731610485781829 & -0.7968700457690127335 \\
	$F$ & 0.6004272512336641681 & 0.3610786938306778951 & 0.5993422515285554206 \\
	$G$ & -0.9552893441051859324 & -0.9548303569438954492 & -0.9544148962268187706 \\
	$H$ & 0.2019979465132786674 & 0.2900942229392094713 & 0.2014481720952456543 \\
	\hline\Tstrut
	$energy$ & -313.7323719353269477 & 1145.418964319278931 & 1224.478456072020493
\end{longtable}

Although the parameters are known to 50,014 digits, the degree of the algebraic polynomials are $>360$.

\noindent
\textit{\textbf{Symmetries -}}

The symmetry groups for 52 points are identical under all 3 potentials.

\begin{center}
	\begin{tabular}{l|l}
		\multicolumn{2}{c}{Symmetries - 52 points} \\
		\hline\Tstrut
		planes & [[17, 1]] \\[0.2ex]
		\hline\Tstrut
		Gram groups & [[6, 442], [52, 1]] \\
		\hline\Tstrut
		Polygons & [[3, 17]]
	\end{tabular}
\end{center}

\subsection{53 points}
The symmetry groups for this configuration show an embedded nonagon as the possible alignment orientation. After rotation, it was discovered that this was correct, as the arrangement then assumed a balanced configuration with 11 dipoles above and 11 dipoles below the nonagon.
This is a 2:2:2:2:2:2:2:2:2:2:2:9:2:2:2:2:2:2:2:2:2:2:2 arrangement.

\begin{figure}[!ht]
	\begin{center}
		\includegraphics[type=pdf,ext=pdf,read=pdf,height=1in,width=1in,angle=0]{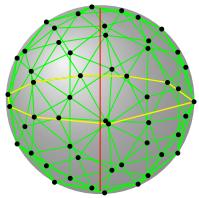}
		\caption{53 points.}
		\label{fig:53pts}
	\end{center}
\end{figure}
\noindent
The nonagon is highlighted in yellow in figure \ref{fig:53pts}.

\noindent
\textit{\textbf{Algebraic parameterization -}}

It requires 27 parameters $a$ - $z$, $A$, to adequately constrain the 53 point set, due to the 22 dipoles and the nonagon in the polyhedron.

The algebraic parameterized structure is given below:

\begin{longtable}[c]{r|ccc}
	\caption{Parameterization for 53 points} \\
	pt & $x$ & $y$ & $z$ \\
	\hline\vspace*{-2.2ex}
	\endfirsthead
	\multicolumn{4}{c}%
	{\tablename\ \thetable\ -- 53 points parameters -- \textit{continued}} \\
	pt & $x$ & $y$ & $z$ \\
	\hline\vspace*{-2.2ex}
	\endhead
	1 & $-\sqrt{1-a^2}$ & $0$ & $a$ \\[0.5em]
	2 & $d$ & $\sqrt{1-c^2-d^2}$ & $c$ \\[0.5em]
	3 & $d$ & $-\sqrt{1-c^2-d^2}$ & $c$ \\[0.5em]
	4 & $-f$ & $\sqrt{1-e^2-f^2}$ & $e$ \\[0.5em]
	5 & $-f$ & $-\sqrt{1-e^2-f^2}$ & $e$ \\[0.5em]
	6 & $-h$ & $\sqrt{1-g^2-h^2}$ & $g$ \\[0.5em]
	7 & $-h$ & $-\sqrt{1-g^2-h^2}$ & $g$ \\[0.5em]
	8 & $\sqrt{1-i^2}$ & $0$ & $i$ \\[0.5em]
	9 & $k$ & $\sqrt{1-j^2-k^2}$ & $j$ \\[0.5em]
	10 & $k$ & $-\sqrt{1-j^2-k^2}$ & $j$ \\[0.5em]
	11 & $m$ & $\sqrt{1-l^2-m^2}$ & $l$ \\[0.5em]
	12 & $m$ & $-\sqrt{1-l^2-m^2}$ & $l$ \\[0.5em]
	13 & $-o$ & $\sqrt{1-n^2-o^2}$ & $n$ \\[0.5em]
	14 & $-o$ & $-\sqrt{1-n^2-o^2}$ & $n$ \\[0.5em]
	15 & $-q$ & $\sqrt{1-p^2-q^2}$ & $p$ \\[0.5em]
	16 & $-q$ & $-\sqrt{1-p^2-q^2}$ & $p$ \\[0.5em]
	17 & $-s$ & $\sqrt{1-r^2-s^2}$ & $r$ \\[0.5em]
	18 & $-s$ & $-\sqrt{1-r^2-s^2}$ & $r$ \\[0.5em]
	19 & $u$ & $\sqrt{1-t^2-u^2}$ & $t$ \\[0.5em]
	20 & $u$ & $-\sqrt{1-t^2-u^2}$ & $t$ \\[0.5em]
	21 & $w$ & $\sqrt{1-v^2-w^2}$ & $v$ \\[0.5em]
	22 & $w$ & $-\sqrt{1-v^2-w^2}$ & $v$ \\[0.5em]
	23 & $x$ & $\sqrt{1-x^2}$ & $0$ \\[0.5em]
	24 & $x$ & $-\sqrt{1-x^2}$ & $0$ \\[0.5em]
	25 & $y$ & $\sqrt{1-y^2}$ & $0$ \\[0.5em]
	26 & $y$ & $-\sqrt{1-y^2}$ & $0$ \\[0.5em]
	27 & $-z$ & $\sqrt{1-z^2}$ & $0$ \\[0.5em]
	28 & $-z$ & $-\sqrt{1-z^2}$ & $0$ \\[0.5em]
	29 & $-A$ & $\sqrt{1-A^2}$ & $0$ \\[0.5em]
	30 & $-A$ & $-\sqrt{1-A^2}$ & $0$ \\[0.5em]
	31 & $-1$ & $0$ & $0$ \\[0.5em]
	32 & $w$ & $\sqrt{1-v^2-w^2}$ & $-v$ \\[0.5em]
	33 & $w$ & $-\sqrt{1-v^2-w^2}$ & $-v$ \\[0.5em]
	34 & $u$ & $\sqrt{1-t^2-u^2}$ & $-t$ \\[0.5em]
	35 & $u$ & $-\sqrt{1-t^2-u^2}$ & $-t$ \\[0.5em]
	36 & $-s$ & $\sqrt{1-r^2-s^2}$ & $-r$ \\[0.5em]
	37 & $-s$ & $-\sqrt{1-r^2-s^2}$ & $-r$ \\[0.5em]
	38 & $-q$ & $\sqrt{1-p^2-q^2}$ & $-p$ \\[0.5em]
	39 & $-q$ & $-\sqrt{1-p^2-q^2}$ & $-p$ \\[0.5em]
	40 & $-o$ & $\sqrt{1-n^2-o^2}$ & $-n$ \\[0.5em]
	41 & $-o$ & $-\sqrt{1-n^2-o^2}$ & $-n$ \\[0.5em]
	42 & $m$ & $\sqrt{1-l^2-m^2}$ & $-l$ \\[0.5em]
	43 & $m$ & $-\sqrt{1-l^2-m^2}$ & $-l$ \\[0.5em]
	44 & $k$ & $\sqrt{1-j^2-k^2}$ & $-j$ \\[0.5em]
	45 & $k$ & $-\sqrt{1-j^2-k^2}$ & $-j$ \\[0.5em]
	46 & $\sqrt{1-i^2}$ & $0$ & $-i$ \\[0.5em]
	47 & $-h$ & $\sqrt{1-g^2-h^2}$ & $-g$ \\[0.5em]
	48 & $-h$ & $-\sqrt{1-g^2-h^2}$ & $-g$ \\[0.5em]
	49 & $-f$ & $\sqrt{1-e^2-f^2}$ & $-e$ \\[0.5em]
	50 & $-f$ & $-\sqrt{1-e^2-f^2}$ & $-e$ \\[0.5em]
	51 & $d$ & $\sqrt{1-c^2-d^2}$ & $-c$ \\[0.5em]
	52 & $d$ & $-\sqrt{1-c^2-d^2}$ & $-c$ \\[0.5em]
	53 & $-\sqrt{1-a^2}$ & $0$ & $-a$
\end{longtable}

The values to 19 digits of the 27 parameters optimized for the minimal solutions of 53 points are:

\begin{longtable}[c]{c|c|c|c}
	\caption{Parameter values for 53 points} \\
	Parameter & log & 1/r & $1/r^2$ \\
	\hline\vspace*{-2.2ex}
	\endfirsthead
	\multicolumn{4}{c}%
	{\tablename\ \thetable\ -- 53 points parameters -- \textit{continued}} \\
	Parameter & log & 1/r & $1/r^2$ \\
	\hline\vspace*{-2.2ex}
	\endhead
	$a$ & 0.9922383091687744256 & 0.9919970752578312165 & 	0.9917562612873874211 \\
	$b$ & 0 & 0 & 0 \\
	$c$ & 0.9141640467130687061 & 0.9137948756226166472 & 0.9134199495737052469 \\
	$d$ & 0.3184397348683263442 & 0.3189184297296864645 & 0.3193806174383299475 \\
	$e$ & 0.8449219588410066192 & 0.8443036628333433686 & 0.8436975246510814288 \\
	$f$ & 0.1810473889784751108 & 0.1799415635599399358 & 0.1790322902624099246 \\
	$g$ & 0.7880792401220870453 & 0.7880317010914089433 & 0.7877851154269705704 \\
	$h$ & 0.5677986568990037053 & 0.5677626888192366040 & 0.5679669190393804607 \\
	$i$ & 0.6966525423673355375 & 0.6946023380292713758 & 0.6926436782490230069 \\
	$j$ & 0.6875712107615540485 & 0.6868567714732249430 & 0.6859067230504037892 \\
	$k$ & 0.2374914528531061836 & 0.2383382588782405128 & 0.2390057734785773905 \\
	$l$ & 0.5532059876067849599 & 0.5524320199170037729 & 0.5516368342462199810 \\
	$m$ & 0.6707627292478623810 & 0.6701920894981987551 & 0.6697095690960534540 \\
	$n$ & 0.4698157495683047133 & 0.4699431853813217994 & 0.4699674631189391416 \\
	$o$ & 0.5813134415847384899 & 0.5827468057332981069 & 0.5840887813626190508 \\
	$p$ & 0.4480791504390199042 & 0.4485072806945607078 & 0.4486284865078888804 \\
	$q$ & 0.1526688121035634178 & 0.1517834562802006705 & 0.1511015231543868104 \\
	$r$ & 0.4133080166336308382 & 0.4136712968057477030 & 0.4135759473295635466 \\
	$s$ & 0.8744719800805626432 & 0.8747384398326747673 & 0.8751431788561740542 \\
	$t$ & 0.2679602460353038208 & 0.2673534741614292354 & 0.2665765364798368368 \\
	$u$ & 0.4154029800247448943 & 0.4141852898838646499 & 0.4133830734016139632 \\
	$v$ & 0.2517580787119970043 & 0.2508448208111062956 & 0.2500433683526314356 \\
	$w$ & 0.9389778880938943014 & 0.9391679047684926740 & 0.9392941178835976341 \\
	$x$ & 0.7671316788838533519 & 0.7673423007824083778 & 0.7675617508389593562 \\
	$y$ & 0.01087226818819071404 & 0.01044232174511778970 & 0.01045987222671940437 \\
	$z$ & 0.4791317312088363128 & 0.4782595721394172755 & 0.4773954386364521234 \\
	$A$ & 0.8394790556632479107 & 0.8404553760294120893 & 0.8412649158721593487 \\
	\hline\Tstrut
	$energy$ & -325.1382346950206751 & 1191.922290416224369 & 1278.652206247138828
\end{longtable}
After working with the parameterization, it was discovered that one parameter, $b$, could be set to 0 as shown. Unfortunately the GP-Pari program encountered a \textit{"matsolve: impossible inverse in Gauss"} error in the Jacobian while enhancing precision for all 3 power laws and the parameters had to be found by direct \textit{descent.3d} search to 1,001 digits accuracy instead.

\noindent
\textit{\textbf{Symmetries -}}

The parallel planes for 53 points are different under each potential, otherwise the Gram matrix and polygon groups are identical.

\begin{center}
	\begin{tabular}{l|l}
		\multicolumn{2}{c}{Symmetries - 53 points - \textit{logarithmic}} \\
		\hline\Tstrut
		planes & [[4, 480], [10, 1], [40, 1], [84, 1]] \\[0.2ex]
		\hline\Tstrut
		Gram groups & [[2, 6], [4, 50], [8, 318], [53, 1]] \\
		\hline\Tstrut
		Polygons & [[4, 490], [5, 1], [9, 1]]
	\end{tabular}
\end{center}

\begin{center}
	\begin{tabular}{l|l}
		\multicolumn{2}{c}{Symmetries - 53 points - \textit{Coulomb}} \\
		\hline\Tstrut
		planes & [[2, 12], [3, 4], [4, 470], [10, 1], [40, 1], [84, 1]] \\[0.2ex]
		\hline\Tstrut
		Gram groups & [[2, 6], [4, 50], [8, 318], [53, 1]] \\
		\hline\Tstrut
		Polygons & [[4, 490], [5, 1], [9, 1]]
	\end{tabular}
\end{center}

\begin{center}
	\begin{tabular}{l|l}
		\multicolumn{2}{c}{Symmetries - 53 points - \textit{Inverse sq.}} \\
		\hline\Tstrut
		planes & [[2, 8], [3, 5], [4, 471], [10, 1], [40, 1], [84, 1]] \\[0.2ex]
		\hline\Tstrut
		Gram groups & [[2, 6], [4, 50], [8, 318], [53, 1]] \\
		\hline\Tstrut
		Polygons & [[4, 490], [5, 1], [9, 1]]
	\end{tabular}
\end{center}

\subsection{54 points}
The optimal configuration for 54 points does not allow for parameterization. This is the fourth such set of points (previous: 26 pts, 35 pts, 36 pts).

\begin{figure}[ht]
	\begin{center}
		\includegraphics[type=pdf,ext=pdf,read=pdf,height=1in,width=1in,angle=0]{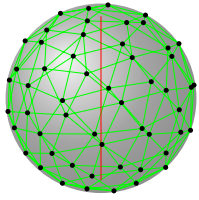}
		\caption{54 points.}
		\label{fig:54pts}
	\end{center}
\end{figure}
\noindent
\textit{\textbf{Minimal Energy values -}}

The coordinates for 54 points are known to 77 digits for the \textit{log} potential and 38 digits for the other two. The minimal energies have been determined for all 3 potentials as well.
\begin{center}
	\begin{tabular}{l|l}
		\multicolumn{2}{c}{Minimal Energy - 54 points} \\
		\hline\Tstrut
		logarithmic & -336.7454643971100850\ldots \\[0.2ex]
		\hline\Tstrut
		Coulomb & 1239.361474729158935\ldots \\[0.2ex]
		\hline\Tstrut
		Inverse square law & 1334.084895358264508\ldots
	\end{tabular}
\end{center}

\noindent
\textit{\textbf{Symmetries -}}

The Gram matrix is slightly different for the \textit{Coulomb $1/r$} potential.

\begin{center}
	\begin{tabular}{l|l}
		\multicolumn{2}{c}{Symmetries - 54 points - \textit{logarithmic}} \\
		\hline\Tstrut
		planes & [] \\[0.2ex]
		\hline\Tstrut
		Gram groups & [[2, 27], [4, 702], [54, 1]] \\
		\hline\Tstrut
		Polygons & []
	\end{tabular}
\end{center}

\begin{center}
	\begin{tabular}{l|l}
		\multicolumn{2}{c}{Symmetries - 54 points - \textit{Coulomb}} \\
		\hline\Tstrut
		planes & [] \\[0.2ex]
		\hline\Tstrut
		Gram groups & [[2, 29], [4, 701], [54, 1]] \\
		\hline\Tstrut
		Polygons & []
	\end{tabular}
\end{center}

\begin{center}
	\begin{tabular}{l|l}
		\multicolumn{2}{c}{Symmetries - 54 points - \textit{Inverse sq}} \\
		\hline\Tstrut
		planes & [] \\[0.2ex]
		\hline\Tstrut
		Gram groups & [[2, 27], [4, 702], [54, 1]] \\
		\hline\Tstrut
		Polygons & []
	\end{tabular}
\end{center}

\subsection{55 points}
Like the case with 54 points, the optimal configurations for 55 points also does not allow for parameterization. This is the fifth such set of points (previous: 26 pts, 35 pts, 36 pts, 54 pts).

\begin{figure}[ht]
	\begin{center}
		\includegraphics[type=pdf,ext=pdf,read=pdf,height=1in,width=1in,angle=0]{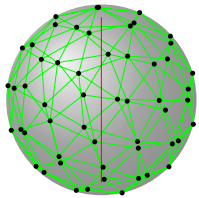}
		\caption{55 points.}
		\label{fig:55pts}
	\end{center}
\end{figure}
\noindent
\textit{\textbf{Minimal Energy values -}}

The coordinates for 55 points are known to 77 digits for the \textit{log} potential and 38 digits for the other two. The minimal energies have been determined for all 3 potentials as well.
\begin{center}
	\begin{tabular}{l|l}
		\multicolumn{2}{c}{Minimal Energy - 55 points} \\
		\hline\Tstrut
		logarithmic & -348.5417962810725430\ldots \\[0.2ex]
		\hline\Tstrut
		Coulomb & 1287.772720782708749\ldots \\[0.2ex]
		\hline\Tstrut
		Inverse square law & 1390.969194237451016\ldots
	\end{tabular}
\end{center}

\noindent
\textit{\textbf{Symmetries -}}

As in the case for 54 points, the Gram matrix is slightly different for the \textit{Coulomb $1/r$} potential.

\begin{center}
	\begin{tabular}{l|l}
		\multicolumn{2}{c}{Symmetries - 55 points - \textit{logarithmic}} \\
		\hline\Tstrut
		planes & [] \\[0.2ex]
		\hline\Tstrut
		Gram groups & [[2, 27], [4, 729], [55, 1]] \\
		\hline\Tstrut
		Polygons & []
	\end{tabular}
\end{center}

\begin{center}
	\begin{tabular}{l|l}
		\multicolumn{2}{c}{Symmetries - 55 points - \textit{Coulomb}} \\
		\hline\Tstrut
		planes & [] \\[0.2ex]
		\hline\Tstrut
		Gram groups & [[2, 29], [4, 728], [55, 1]] \\
		\hline\Tstrut
		Polygons & []
	\end{tabular}
\end{center}

\begin{center}
	\begin{tabular}{l|l}
		\multicolumn{2}{c}{Symmetries - 55 points - \textit{Inverse sq.}} \\
		\hline\Tstrut
		planes & [] \\[0.2ex]
		\hline\Tstrut
		Gram groups & [[2, 27], [4, 729], [55, 1]] \\
		\hline\Tstrut
		Polygons & []
	\end{tabular}
\end{center}

\subsection{56 points}
The optimal configurations for 56 points do not allow for a parameterization. This is the sixth such set (previous: 26, 35, 36, 54, 55 pts). As the number of points in a set increases, more and more optimal configurations occur which cannot be parameterized, thus forcing a direct optimization of the spherical code instead.

\begin{figure}[ht]
	\begin{center}
		\includegraphics[type=pdf,ext=pdf,read=pdf,height=1in,width=1in,angle=0]{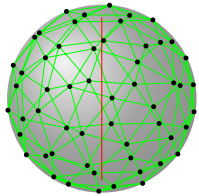}
		\caption{56 points.}
		\label{fig:56pts}
	\end{center}
\end{figure}
\noindent
\textit{\textbf{Minimal Energy values -}}

The coordinates for 56 points are known to 77 digits for the \textit{log} potential and 38 digits for the other two. The minimal energies have been determined for all 3 potentials as well.
\begin{center}
	\begin{tabular}{l|l}
		\multicolumn{2}{c}{Minimal Energy - 56 points} \\
		\hline\Tstrut
		logarithmic & -360.5458992442545471\ldots \\[0.2ex]
		\hline\Tstrut
		Coulomb & 1337.094945275657077\ldots \\[0.2ex]
		\hline\Tstrut
		Inverse square law & 1448.954274110880227\ldots
	\end{tabular}
\end{center}

\noindent
\textit{\textbf{Symmetries -}}

The Gram matrix group for the \textit{logarithmic} potential is slightly different from the other two potentials.

\begin{center}
	\begin{tabular}{l|l}
		\multicolumn{2}{c}{Symmetries - 56 points - \textit{logarithmic}} \\
		\hline\Tstrut
		planes & [] \\[0.2ex]
		\hline\Tstrut
		Gram groups & [[2, 28], [4, 756], [56, 1]] \\
		\hline\Tstrut
		Polygons & []
	\end{tabular}
\end{center}

\begin{center}
	\begin{tabular}{l|l}
		\multicolumn{2}{c}{Symmetries - 56 points - \textit{Coulomb}} \\
		\hline\Tstrut
		planes & [] \\[0.2ex]
		\hline\Tstrut
		Gram groups & [[4, 42], [8, 364], [56, 1]] \\
		\hline\Tstrut
		Polygons & []
	\end{tabular}
\end{center}

\begin{center}
	\begin{tabular}{l|l}
		\multicolumn{2}{c}{Symmetries - 56 points - \textit{Inverse sq.}} \\
		\hline\Tstrut
		planes & [] \\[0.2ex]
		\hline\Tstrut
		Gram groups & [[4, 42], [8, 364], [56, 1]] \\
		\hline\Tstrut
		Polygons & []
	\end{tabular}
\end{center}

\subsection{57 points}
The polygons group of the symmetry groups indicates that the optimal configuration consists of 19 parallel triangles, and rotating a random minimal energy configuration confirms the symmetry group. A constraining parameterization was found.

\begin{figure}[ht]
	\begin{center}
		\includegraphics[type=pdf,ext=pdf,read=pdf,height=1in,width=1in,angle=0]{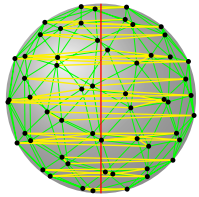}
		\caption{57 points.}
		\label{fig:57pts}
	\end{center}
\end{figure}

\noindent
\textit{\textbf{Algebraic parameterization -}}

It requires 18 parameters to adequately constrain the 57 point set, due to the 19 embedded triangles in the figure. Since the figure is balanced around the equilateral triangle at the equator, only 18 parameters are required.

The algebraic parameterized structure is given below:

\begin{longtable}[c]{r|ccc}
	\caption{Parameterization for 57 points} \\
	pt & $x$ & $y$ & $z$ \\
	\hline\vspace*{-2.2ex}
	\endfirsthead
	\multicolumn{4}{c}%
	{\tablename\ \thetable\ -- 57 points parameters -- \textit{continued}} \\
	pt & $x$ & $y$ & $z$ \\
	\hline\vspace*{-2.2ex}
	\endhead
	1 & $b$ & $-\sqrt{1-a^2-b^2}$ & $a$ \\[0.7ex]
	2 & $\frac{-b-\sqrt{3}\sqrt{1-a^2-b^2}}{2}$ & $\frac{-b\sqrt{3}+\sqrt{1-a^2-b^2}}{2}$ & $a$ \\[0.7ex]
	3 & $\frac{-b+\sqrt{3}\sqrt{1-a^2-b^2}}{2}$ & $\frac{+b\sqrt{3}+\sqrt{1-a^2-b^2}}{2}$ & $a$ \\[0.7ex]
	4 & $d$ & $-\sqrt{1-c^2-d^2}$ & $c$ \\[0.7ex]
	5 & $\frac{-d-\sqrt{3}\sqrt{1-c^2-d^2}}{2}$ & $\frac{-d\sqrt{3}+\sqrt{1-c^2-d^2}}{2}$ & $c$ \\[0.7ex]
	6 & $\frac{-d+\sqrt{3}\sqrt{1-c^2-d^2}}{2}$ & $\frac{+d\sqrt{3}+\sqrt{1-c^2-d^2}}{2}$ & $c$ \\[0.7ex]
	7 & $f$ & $\sqrt{1-e^2-f^2}$ & $e$ \\[0.7ex]
	8 & $\frac{-f-\sqrt{3}\sqrt{1-e^2-f^2}}{2}$ & $\frac{+f\sqrt{3}-\sqrt{1-e^2-f^2}}{2}$ & $e$ \\[0.7ex]
	9 & $\frac{-f+\sqrt{3}\sqrt{1-e^2-f^2}}{2}$ & $\frac{-f\sqrt{3}-\sqrt{1-e^2-f^2}}{2}$ & $e$ \\[0.7ex]
	10 & $h$ & $-\sqrt{1-g^2-h^2}$ & $g$ \\[0.7ex]
	11 & $\frac{-h-\sqrt{3}\sqrt{1-g^2-h^2}}{2}$ & $\frac{-h\sqrt{3}+\sqrt{1-g^2-h^2}}{2}$ & $g$ \\[0.7ex]
	12 & $\frac{-h+\sqrt{3}\sqrt{1-g^2-h^2}}{2}$ & $\frac{+h\sqrt{3}+\sqrt{1-g^2-h^2}}{2}$ & $g$ \\[0.7ex]
	13 & $j$ & $\sqrt{1-i^2-j^2}$ & $i$ \\[0.7ex]
	14 & $\frac{-j-\sqrt{3}\sqrt{1-i^2-j^2}}{2}$ & $\frac{+j\sqrt{3}-\sqrt{1-i^2-j^2}}{2}$ & $i$ \\[0.7ex]
	15 & $\frac{-j+\sqrt{3}\sqrt{1-i^2-j^2}}{2}$ & $\frac{-j\sqrt{3}-\sqrt{1-i^2-j^2}}{2}$ & $i$ \\[0.7ex]
	16 & $l$ & $-\sqrt{1-k^2-l^2}$ & $k$ \\[0.7ex]
	17 & $\frac{-l-\sqrt{3}\sqrt{1-k^2-l^2}}{2}$ & $\frac{-l\sqrt{3}+\sqrt{1-k^2-l^2}}{2}$ & $k$ \\[0.7ex]
	18 & $\frac{-l+\sqrt{3}\sqrt{1-k^2-l^2}}{2}$ & $\frac{+l\sqrt{3}+\sqrt{1-k^2-l^2}}{2}$ & $k$ \\[0.7ex]
	19 & $n$ & $\sqrt{1-m^2-n^2}$ & $m$ \\[0.7ex]
	20 & $\frac{-n-\sqrt{3}\sqrt{1-m^2-n^2}}{2}$ & $\frac{+n\sqrt{3}-\sqrt{1-m^2-n^2}}{2}$ & $m$ \\[0.7ex]
	21 & $\frac{-n+\sqrt{3}\sqrt{1-m^2-n^2}}{2}$ & $\frac{-n\sqrt{3}-\sqrt{1-m^2-n^2}}{2}$ & $m$ \\[0.7ex]
	22 & $p$ & $-\sqrt{1-o^2-p^2}$ & $o$ \\[0.7ex]
	23 & $\frac{-p-\sqrt{3}\sqrt{1-o^2-p^2}}{2}$ & $\frac{-p\sqrt{3}+\sqrt{1-o^2-p^2}}{2}$ & $o$ \\[0.7ex]
	24 & $\frac{-p+\sqrt{3}\sqrt{1-o^2-p^2}}{2}$ & $\frac{+p\sqrt{3}+\sqrt{1-o^2-p^2}}{2}$ & $o$ \\[0.7ex]
	25 & $r$ & $\sqrt{1-q^2-r^2}$ & $q$ \\[0.7ex]
	26 & $\frac{-r-\sqrt{3}\sqrt{1-q^2-r^2}}{2}$ & $\frac{+r\sqrt{3}-\sqrt{1-q^2-r^2}}{2}$ & $q$ \\[0.7ex]
	27 & $\frac{-r+\sqrt{3}\sqrt{1-q^2-r^2}}{2}$ & $\frac{-r\sqrt{3}-\sqrt{1-q^2-r^2}}{2}$ & $q$ \\[0.7ex]
	28 & $1$ & $0$ & $0$ \\[0.7ex]
	29 & $-\frac{1}{2}$ & $\frac{\sqrt{3}}{2}$ & $0$ \\[0.7ex]
	30 & $-\frac{1}{2}$ & $-\frac{\sqrt{3}}{2}$ & $0$ \\[0.7ex]
	31 & $r$ & $-\sqrt{1-q^2-r^2}$ & $-q$ \\[0.7ex]
	32 & $\frac{-r-\sqrt{3}\sqrt{1-q^2-r^2}}{2}$ & $\frac{-r\sqrt{3}+\sqrt{1-q^2-r^2}}{2}$ & $-q$ \\[0.7ex]
	33 & $\frac{-r+\sqrt{3}\sqrt{1-q^2-r^2}}{2}$ & $\frac{+r\sqrt{3}+\sqrt{1-q^2-r^2}}{2}$ & $-q$ \\[0.7ex]
	34 & $p$ & $\sqrt{1-o^2-p^2}$ & $-o$ \\[0.7ex]
	35 & $\frac{-p-\sqrt{3}\sqrt{1-o^2-p^2}}{2}$ & $\frac{+p\sqrt{3}-\sqrt{1-o^2-p^2}}{2}$ & $-o$ \\[0.7ex]
	36 & $\frac{-p+\sqrt{3}\sqrt{1-o^2-p^2}}{2}$ & $\frac{-p\sqrt{3}-\sqrt{1-o^2-p^2}}{2}$ & $-o$ \\[0.7ex]
	37 & $n$ & $-\sqrt{1-m^2-n^2}$ & $-m$ \\[0.7ex]
	38 & $\frac{-n-\sqrt{3}\sqrt{1-m^2-n^2}}{2}$ & $\frac{-n\sqrt{3}+\sqrt{1-m^2-n^2}}{2}$ & $-m$ \\[0.7ex]
	39 & $\frac{-n+\sqrt{3}\sqrt{1-m^2-n^2}}{2}$ & $\frac{+n\sqrt{3}+\sqrt{1-m^2-n^2}}{2}$ & $-m$ \\[0.7ex]
	40 & $l$ & $\sqrt{1-k^2-l^2}$ & $-k$ \\[0.7ex]
	41 & $\frac{-l-\sqrt{3}\sqrt{1-k^2-l^2}}{2}$ & $\frac{+l\sqrt{3}-\sqrt{1-k^2-l^2}}{2}$ & $-k$ \\[0.7ex]
	42 & $\frac{-l+\sqrt{3}\sqrt{1-k^2-l^2}}{2}$ & $\frac{-l\sqrt{3}-\sqrt{1-k^2-l^2}}{2}$ & $-k$ \\[0.7ex]
	43 & $j$ & $-\sqrt{1-i^2-j^2}$ & $-i$ \\[0.7ex]
	44 & $\frac{-j-\sqrt{3}\sqrt{1-i^2-j^2}}{2}$ & $\frac{-j\sqrt{3}+\sqrt{1-i^2-j^2}}{2}$ & $-i$ \\[0.7ex]
	45 & $\frac{-j+\sqrt{3}\sqrt{1-i^2-j^2}}{2}$ & $\frac{+j\sqrt{3}+\sqrt{1-i^2-j^2}}{2}$ & $-i$ \\[0.7ex]
	46 & $h$ & $\sqrt{1-g^2-h^2}$ & $-g$ \\[0.7ex]
	47 & $\frac{-h-\sqrt{3}\sqrt{1-g^2-h^2}}{2}$ & $\frac{+h\sqrt{3}-\sqrt{1-g^2-h^2}}{2}$ & $-g$ \\[0.7ex]
	48 & $\frac{-h+\sqrt{3}\sqrt{1-g^2-h^2}}{2}$ & $\frac{-h\sqrt{3}-\sqrt{1-g^2-h^2}}{2}$ & $-g$ \\[0.7ex]
	49 & $f$ & $-\sqrt{1-e^2-f^2}$ & $-e$ \\[0.7ex]
	50 & $\frac{-f-\sqrt{3}\sqrt{1-e^2-f^2}}{2}$ & $\frac{-f\sqrt{3}+\sqrt{1-e^2-f^2}}{2}$ & $-e$ \\[0.7ex]
	51 & $\frac{-f+\sqrt{3}\sqrt{1-e^2-f^2}}{2}$ & $\frac{+f\sqrt{3}+\sqrt{1-e^2-f^2}}{2}$ & $-e$ \\[0.7ex]
	52 & $d$ & $\sqrt{1-c^2-d^2}$ & $-c$ \\[0.7ex]
	53 & $\frac{-d-\sqrt{3}\sqrt{1-c^2-d^2}}{2}$ & $\frac{+d\sqrt{3}-\sqrt{1-c^2-d^2}}{2}$ & $-c$ \\[0.7ex]
	54 & $\frac{-d+\sqrt{3}\sqrt{1-c^2-d^2}}{2}$ & $\frac{-d\sqrt{3}-\sqrt{1-c^2-d^2}}{2}$ & $-c$ \\[0.7ex]
	55 & $b$ & $\sqrt{1-a^2-b^2}$ & $-a$ \\[0.7ex]
	56 & $\frac{-b-\sqrt{3}\sqrt{1-a^2-b^2}}{2}$ & $\frac{+b\sqrt{3}-\sqrt{1-a^2-b^2}}{2}$ & $-a$ \\[0.7ex]
	57 & $\frac{-b+\sqrt{3}\sqrt{1-a^2-b^2}}{2}$ & $\frac{-b\sqrt{3}-\sqrt{1-a^2-b^2}}{2}$ & $-a$
\end{longtable}

\noindent
\textit{\textbf{Parameterization values --}}

The values to 19 digits of the 18 parameters optimized for the minimal solutions of 57 points are:

\begin{longtable}[c]{c|c|c|c}
	\caption{Parameter values for 57 points} \\
	Parameter & log & 1/r & $1/r^2$ \\
	\hline\vspace*{-2.2ex}
	\endfirsthead
	\multicolumn{4}{c}%
	{\tablename\ \thetable\ -- 57 points parameters -- \textit{continued}} \\
	Parameter & log & 1/r & $1/r^2$ \\
	\hline\vspace*{-2.2ex}
	\endhead
	$a$ & 0.9595125380384837748 & 0.9593982472561573346 & 0.9592430295797528902 \\
	$b$ & 0.1528538042046582041 & 0.1510400201442292428 & 0.1504481057951366924 \\
	$c$ & 0.8064159195044753954 & 0.8061127847218874374 & 0.8051648162909763574 \\
	$d$ & 0.5812803228875366308 & 0.5818017010622614807 & 0.5832136313510178951 \\
	$e$ & 0.7655416542440036844 & 0.7635580248548546735 & 0.7619669959015469652 \\
	$f$ & 0.5343092156110378053 & 0.5362774228532354228 & 0.5377746556112610776 \\
	$g$ & 0.6518605050286107147 & 0.6526101176624277225 & 0.6527511409910811614 \\
	$h$ & 0.5110203688455053209 & 0.5097398748167363052 & 0.5090549759613889785 \\
	$i$ & 0.4658252012187145059 & 0.4653172762165911128 & 0.4644759355926697615 \\
	$j$ & 0.8703353439326157776 & 0.8703913334987289892 & 0.8705680364936731241 \\
	$k$ & 0.4174822430210970745 & 0.4167819456237354241 & 0.4156382026447787826 \\
	$l$ & 0.8512195110578723961 & 0.8514328688417494922 & 0.8519885231386723758 \\
	$m$ & 0.3919407857564659670 & 0.3908468794742435778 & 0.3902284386652641164 \\
	$n$ & 0.6548991798133300362 & 0.6544971461608950337 & 0.6541437514238916804 \\
	$o$ & 0.1839164307976121038 & 0.1833947010632889515 & 0.1828170935972039697 \\
	$p$ & 0.6281900887221845112 & 0.6298855180742855577 & 0.6312350522250580709 \\
	$q$ & 0.04745350350485770277 & 0.04881688478795210130 & 0.05003293613571645555 \\
	$r$ & 0.8919586151422930883 & 0.8926119707663709581 & 0.8930974627550429491 \\
	\hline\Tstrut
	$energy$ & -372.7412006183720172 & 1387.383229252841734 & 1508.368838505948797
\end{longtable}

No algebraic numbers were recovered from the 50,014 digit spherical codes, for all 3 potentials, their degree is $>360$, if found.

\noindent
\textit{\textbf{Symmetries -}}

The symmetry groups for 57 points are identical under all 3 potentials.

\begin{center}
	\begin{tabular}{l|l}
		\multicolumn{2}{c}{Symmetries - 57 points} \\
		\hline\Tstrut
		planes & [[19, 1]] \\[0.2ex]
		\hline\Tstrut
		Gram groups & [[6, 28], [12, 252], [57, 1]] \\
		\hline\Tstrut
		Polygons & [[3, 19]]
	\end{tabular}
\end{center}

\subsection{58 points}

No attempts to parameterize 58 points from their optimal solutions were made. A careful examination showed that 84 parameters are necessary to constrain this polyhedron, well beyond the current computer resources of the correspondence author.

\begin{figure}[ht]
	\begin{center}
		\includegraphics[type=pdf,ext=pdf,read=pdf,height=1in,width=1in,angle=0]{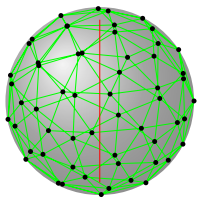}
		\caption{58 points.}
		\label{fig:58pts}
	\end{center}
\end{figure}
\noindent
\textit{\textbf{Minimal Energy values -}}

The coordinates for 58 points are known to 77 digits for the \textit{log} potential and 38 digits for the other two. The minimal energies have been determined for all 3 potentials as well.
\begin{center}
	\begin{tabular}{l|l}
		\multicolumn{2}{c}{Minimal Energy - 58 points} \\
		\hline\Tstrut
		logarithmic & -385.1328297919238166\ldots \\[0.2ex]
		\hline\Tstrut
		Coulomb & 1438.618250640401146\ldots \\[0.2ex]
		\hline\Tstrut
		Inverse square law & 1569.069938532185566\ldots
	\end{tabular}
\end{center}

\noindent
\textit{\textbf{Symmetries -}}

\begin{center}
	\begin{tabular}{l|l}
		\multicolumn{2}{c}{Symmetries - 58 points - \textit{logarithmic}} \\
		\hline\Tstrut
		planes & [[4, 28]] \\[0.2ex]
		\hline\Tstrut
		Gram groups & [[2, 1], [4, 42], [8, 392], [58, 1]] \\
		\hline\Tstrut
		Polygons & [[4, 28]]
	\end{tabular}
\end{center}

\begin{center}
	\begin{tabular}{l|l}
		\multicolumn{2}{c}{Symmetries - 58 points - \textit{Coulomb}} \\
		\hline\Tstrut
		planes & [[4, 28]] \\[0.2ex]
		\hline\Tstrut
		Gram groups & [[2, 1], [4, 42], [8, 392], [58, 1]] \\
		\hline\Tstrut
		Polygons & [[4, 28]]
	\end{tabular}
\end{center}

The Gram group for the \textit{Inverse square law $1/r^2$} potential is slightly different from the previous two.

\begin{center}
	\begin{tabular}{l|l}
		\multicolumn{2}{c}{Symmetries - 58 points - \textit{Inverse sq.}} \\
		\hline\Tstrut
		planes & [[4, 28]] \\[0.2ex]
		\hline\Tstrut
		Gram groups & [[2, 1], [4, 44], [8, 391], [58, 1]] \\
		\hline\Tstrut
		Polygons & [[4, 28]]
	\end{tabular}
\end{center}

\subsection{59 points}
All 3 optimal configurations for 59 points do not allow for a parameterization to be determined. This is the seventh set, in which this occurs (previous: 26, 35, 36, 54, 55, 56).

\begin{figure}[ht]
	\begin{center}
		\includegraphics[type=pdf,ext=pdf,read=pdf,height=1in,width=1in,angle=0]{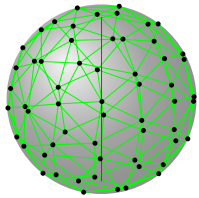}
		\caption{59 points.}
		\label{fig:59pts}
	\end{center}
\end{figure}
\noindent
\textit{\textbf{Minimal Energy values -}}

The coordinates for 59 points are known to 77 digits for the \textit{log} potential and 38 digits for the other two. The minimal energies have been determined for all 3 potentials as well.
\begin{center}
	\begin{tabular}{l|l}
		\multicolumn{2}{c}{Minimal Energy - 59 points} \\
		\hline\Tstrut
		logarithmic & -397.7281496607926605\ldots \\[0.2ex]
		\hline\Tstrut
		Coulomb & 1490.773335278696756\ldots \\[0.2ex]
		\hline\Tstrut
		Inverse square law & 1630.909658338517471\ldots
	\end{tabular}
\end{center}

\noindent
\textit{\textbf{Symmetries -}}

\begin{center}
	\begin{tabular}{l|l}
		\multicolumn{2}{c}{Symmetries - 59 points - \textit{logarithmic}} \\
		\hline\Tstrut
		planes & [] \\[0.2ex]
		\hline\Tstrut
		Gram groups & [[2, 29], [4, 841], [59, 1]] \\
		\hline\Tstrut
		Polygons & []
	\end{tabular}
\end{center}

\begin{center}
	\begin{tabular}{l|l}
		\multicolumn{2}{c}{Symmetries - 59 points - \textit{Coulomb}} \\
		\hline\Tstrut
		planes & [[4, 28]] \\[0.2ex]
		\hline\Tstrut
		Gram groups & [[2, 29], [4, 841], [59, 1]] \\
		\hline\Tstrut
		Polygons & []
	\end{tabular}
\end{center}

\begin{center}
	\begin{tabular}{l|l}
		\multicolumn{2}{c}{Symmetries - 59 points - \textit{Inverse sq.}} \\
		\hline\Tstrut
		planes & [] \\[0.2ex]
		\hline\Tstrut
		Gram groups & [[2, 1711], [59, 1]] \\
		\hline\Tstrut
		Polygons & []
	\end{tabular}
\end{center}

\subsection{60 points}

The polygon group of [3,20] indicates that the configuration consists of 20 parallel triangles on one axis, the only axis of symmetry for the polyhedra. Upon 1 rotation, it was found that this was true, 20 co-planar triangles exists.

\begin{figure}[ht]
	\begin{center}
		\includegraphics[type=pdf,ext=pdf,read=pdf,height=1in,width=1in,angle=0]{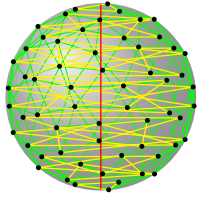}
		\caption{60 points.}
		\label{fig:60pts}
	\end{center}
\end{figure}
\noindent

\noindent
\textit{\textbf{Algebraic parameterization -}}

It requires 30 parameters to adequately constrain the 60 point set, due to the 20 embedded triangles in the figure. Since the figure is balanced around an equator, 10 less parameters than 40 were needed, hence 30 (2 parameters are required per triangle).

The algebraic parameterized structure is given below:

\begin{longtable}[c]{r|ccc}
	\caption{Parameterization for 60 points} \\
	pt & $x$ & $y$ & $z$ \\
	\hline\vspace*{-2.2ex}
	\endfirsthead
	\multicolumn{4}{c}%
	{\tablename\ \thetable\ -- 60 points parameterization -- \textit{continued}} \\
	pt & $x$ & $y$ & $z$ \\
	\hline\vspace*{-2.2ex}
	\endhead
	1 & $b$ & $0$ & $a$ \\[0.7em]
	2 & $\frac{-b}{2}$ & $\frac{+b\sqrt{3}}{2}$ & $a$ \\[0.7em]
	3 & $\frac{-b}{2}$ & $\frac{-b\sqrt{3}}{2}$ & $a$ \\[0.7em]
	4 & $d$ & $-\sqrt{1-c^2-d^2}$ & $c$ \\[0.7em]
	5 & $\frac{-d-\sqrt{3}\sqrt{1-c^2-d^2}}{2}$ & $\frac{-d\sqrt{3}+\sqrt{1-c^2-d^2}}{2}$ & $c$ \\[0.7em]
	6 & $\frac{-d+\sqrt{3}\sqrt{1-c^2-d^2}}{2}$ & $\frac{+d\sqrt{3}+\sqrt{1-c^2-d^2}}{2}$ & $c$ \\[0.7em]
	7 & $f$ & $-\sqrt{1-e^2-f^2}$ & $e$ \\[0.7em]
	8 & $\frac{-f-\sqrt{3}\sqrt{1-e^2-f^2}}{2}$ & $\frac{-f\sqrt{3}+\sqrt{1-e^2-f^2}}{2}$ & $e$ \\[0.7em]
	9 & $\frac{-f+\sqrt{3}\sqrt{1-e^2-f^2}}{2}$ & $\frac{+f\sqrt{3}+\sqrt{1-e^2-f^2}}{2}$ & $e$ \\[0.7em]
	10 & $h$ & $\sqrt{1-g^2-h^2}$ & $g$ \\[0.7em]
	11 & $\frac{-h-\sqrt{3}\sqrt{1-g^2-h^2}}{2}$ & $\frac{+h\sqrt{3}-\sqrt{1-g^2-h^2}}{2}$ & $g$ \\[0.7em]
	12 & $\frac{-h+\sqrt{3}\sqrt{1-g^2-h^2}}{2}$ & $\frac{-h\sqrt{3}-\sqrt{1-g^2-h^2}}{2}$ & $g$ \\[0.7em]
	13 & $j$ & $\sqrt{1-i^2-j^2}$ & $i$ \\[0.7em]
	14 & $\frac{-j-\sqrt{3}\sqrt{1-i^2-j^2}}{2}$ & $\frac{+j\sqrt{3}-\sqrt{1-i^2-j^2}}{2}$ & $i$ \\[0.7em]
	15 & $\frac{-j+\sqrt{3}\sqrt{1-i^2-j^2}}{2}$ & $\frac{-j\sqrt{3}-\sqrt{1-i^2-j^2}}{2}$ & $i$ \\[0.7em]
	16 & $l$ & $-\sqrt{1-k^2-l^2}$ & $k$ \\[0.7em]
	17 & $\frac{-l-\sqrt{3}\sqrt{1-k^2-l^2}}{2}$ & 	$\frac{-l\sqrt{3}+\sqrt{1-k^2-l^2}}{2}$ & $k$ \\[0.7em]
	18 & $\frac{-l+\sqrt{3}\sqrt{1-k^2-l^2}}{2}$ & $\frac{+l\sqrt{3}+\sqrt{1-k^2-l^2}}{2}$ & $k$ \\[0.7em]
	19 & $n$ & $-\sqrt{1-m^2-n^2}$ & $m$ \\[0.7em]
	20 & $\frac{-n-\sqrt{3}\sqrt{1-m^2-n^2}}{2}$ & $\frac{-n\sqrt{3}+\sqrt{1-m^2-n^2}}{2}$ & $m$ \\[0.7em]
	21 & $\frac{-n+\sqrt{3}\sqrt{1-m^2-n^2}}{2}$ & $\frac{+n\sqrt{3}+\sqrt{1-m^2-n^2}}{2}$ & $m$ \\[0.7em]
	22 & $p$ & $\sqrt{1-o^2-p^2}$ & $o$ \\[0.7em]
	23 & $\frac{-p-\sqrt{3}\sqrt{1-o^2-p^2}}{2}$ & $\frac{+p\sqrt{3}-\sqrt{1-o^2-p^2}}{2}$ & $o$ \\[0.7em]
	24 & $\frac{-p+\sqrt{3}\sqrt{1-o^2-p^2}}{2}$ & $\frac{-p\sqrt{3}-\sqrt{1-o^2-p^2}}{2}$ & $o$ \\[0.7em]
	25 & $r$ & $-\sqrt{1-q^2-r^2}$ & $q$ \\[0.7em]
	26 & $\frac{-r-\sqrt{3}\sqrt{1-q^2-r^2}}{2}$ & $\frac{-r\sqrt{3}+\sqrt{1-q^2-r^2}}{2}$ & $q$ \\[0.7em]
	27 & $\frac{-r+\sqrt{3}\sqrt{1-q^2-r^2}}{2}$ & $\frac{+r\sqrt{3}+\sqrt{1-q^2-r^2}}{2}$ & $q$ \\[0.7em]
	28 & $t$ & $-\sqrt{1-s^2-t^2}$ & $s$ \\[0.7em]
	29 & $\frac{-t-\sqrt{3}\sqrt{1-s^2-t^2}}{2}$ & $\frac{-t\sqrt{3}+\sqrt{1-s^2-t^2}}{2}$ & $s$ \\[0.7em]
	30 & $\frac{-t+\sqrt{3}\sqrt{1-s^2-t^2}}{2}$ & $\frac{+t\sqrt{3}+\sqrt{1-s^2-t^2}}{2}$ & $s$ \\[0.7em]
	31 & $u$ & $\sqrt{1-s^2-u^2}$ & $-s$ \\[0.7em]
	32 & $\frac{-u-\sqrt{3}\sqrt{1-s^2-u^2}}{2}$ & 	$\frac{+u\sqrt{3}-\sqrt{1-s^2-u^2}}{2}$ & $-s$ \\[0.7em]
	33 & $\frac{-u+\sqrt{3}\sqrt{1-s^2-u^2}}{2}$ & $\frac{-u\sqrt{3}-\sqrt{1-s^2-u^2}}{2}$ & $-s$ \\[0.7em]
	34 & $v$ & $\sqrt{1-q^2-v^2}$ & $-q$ \\[0.7em]
	35 & $\frac{-v-\sqrt{3}\sqrt{1-q^2-v^2}}{2}$ & $\frac{+v\sqrt{3}-\sqrt{1-q^2-v^2}}{2}$ & $-q$ \\[0.7em]
	36 & $\frac{-v+\sqrt{3}\sqrt{1-q^2-v^2}}{2}$ & $\frac{-v\sqrt{3}-\sqrt{1-q^2-v^2}}{2}$ & $-q$ \\[0.7em]
	37 & $w$ & $-\sqrt{1-o^2-w^2}$ & $-o$ \\[0.7em]
	38 & $\frac{-w-\sqrt{3}\sqrt{1-o^2-w^2}}{2}$ & $\frac{-w\sqrt{3}+\sqrt{1-o^2-w^2}}{2}$ & $-o$ \\[0.7em]
	39 & $\frac{-w+\sqrt{3}\sqrt{1-o^2-w^2}}{2}$ & $\frac{+w\sqrt{3}+\sqrt{1-o^2-w^2}}{2}$ & $-o$ \\[0.7em]
	40 & $x$ & $-\sqrt{1-m^2-x^2}$ & $-m$ \\[0.7em]
	41 & $\frac{-x-\sqrt{3}\sqrt{1-m^2-x^2}}{2}$ & $\frac{-x\sqrt{3}+\sqrt{1-m^2-x^2}}{2}$ & $-m$ \\[0.7em]
	42 & $\frac{-x+\sqrt{3}\sqrt{1-m^2-x^2}}{2}$ & $\frac{+x\sqrt{3}+\sqrt{1-m^2-x^2}}{2}$ & $-m$ \\[0.7em]
	43 & $y$ & $\sqrt{1-k^2-y^2}$ & $-k$ \\[0.7em]
	44 & $\frac{-y-\sqrt{3}\sqrt{1-k^2-y^2}}{2}$ & $\frac{+y\sqrt{3}-\sqrt{1-k^2-y^2}}{2}$ & $-k$ \\[0.7em]
	45 & $\frac{-y+\sqrt{3}\sqrt{1-k^2-y^2}}{2}$ & $\frac{-y\sqrt{3}-\sqrt{1-k^2-y^2}}{2}$ & $-k$ \\[0.7em]
	46 & $z$ & $-\sqrt{1-i^2-z^2}$ & $-i$ \\[0.7em]
	47 & $\frac{-z-\sqrt{3}\sqrt{1-i^2-z^2}}{2}$ & $\frac{-z\sqrt{3}+\sqrt{1-i^2-z^2}}{2}$ & $-i$ \\[0.7em]
	48 & $\frac{-z+\sqrt{3}\sqrt{1-i^2-z^2}}{2}$ & $\frac{+z\sqrt{3}+\sqrt{1-i^2-z^2}}{2}$ & $-i$ \\[0.7em]
	49 & $A$ & $-\sqrt{1-g^2-A^2}$ & $-g$ \\[0.7em]
	50 & $\frac{-A-\sqrt{3}\sqrt{1-g^2-A^2}}{2}$ & $\frac{-A\sqrt{3}+\sqrt{1-g^2-A^2}}{2}$ & $-g$ \\[0.7em]
	51 & $\frac{-A+\sqrt{3}\sqrt{1-g^2-A^2}}{2}$ & $\frac{+A\sqrt{3}+\sqrt{1-g^2-A^2}}{2}$ & $-g$ \\[0.7em]
	52 & $B$ & $-\sqrt{1-e^2-B^2}$ & $-e$ \\[0.7em]
	53 & $\frac{-B-\sqrt{3}\sqrt{1-e^2-B^2}}{2}$ & $\frac{-B\sqrt{3}+\sqrt{1-e^2-B^2}}{2}$ & $-e$ \\[0.7em]
	54 & $\frac{-B+\sqrt{3}\sqrt{1-e^2-B^2}}{2}$ & $\frac{+B\sqrt{3}+\sqrt{1-e^2-B^2}}{2}$ & $-e$ \\[0.7em]
	55 & $C$ & $-\sqrt{1-c^2-C^2}$ & $-c$ \\[0.7em]
	56 & $\frac{-C-\sqrt{3}\sqrt{1-c^2-C^2}}{2}$ & $\frac{-C\sqrt{3}+\sqrt{1-c^2-C^2}}{2}$ & $-c$ \\[0.7em]
	57 & $\frac{-C+\sqrt{3}\sqrt{1-c^2-C^2}}{2}$ & $\frac{+C\sqrt{3}+\sqrt{1-c^2-C^2}}{2}$ & $-c$ \\[0.7em]
	58 & $D$ & $-\sqrt{1-a^2-D^2}$ & $-a$ \\[0.7em]
	59 & $\frac{-D-\sqrt{3}\sqrt{1-a^2-D^2}}{2}$ & 	$\frac{-D\sqrt{3}+\sqrt{1-a^2-D^2}}{2}$ & $-a$ \\[0.7em]
	60 & $\frac{-D+\sqrt{3}\sqrt{1-a^2-D^2}}{2}$ & $\frac{+D\sqrt{3}+\sqrt{1-a^2-D^2}}{2}$ & $-a$
\end{longtable}

\noindent
\textit{\textbf{Parameterization values --}}

The values to 19 digits of the 30 parameters optimized for the minimal solutions of 60 points are:

\begin{longtable}[c]{c|c|c|c}
	\caption{Parameter values for 60 points} \\
	Parameter & log & 1/r & $1/r^2$ \\
	\hline\vspace*{-2.2ex}
	\endfirsthead
	\multicolumn{4}{c}%
	{\tablename\ \thetable\ -- 60 points parameters -- \textit{continued}} \\
	Parameter & log & 1/r & $1/r^2$ \\
	\hline\vspace*{-2.2ex}
	\endhead
	$a$ & 0.9612030957736987686 & 0.9610273821609481894 & 0.9608443350910277295 \\
	$b$ & 0.2758416369496412202 & 0.2764531980948580776 & 0.2770887289722929501 \\
	$c$ & 0.8234204712353256073 & 0.8224141761737171526 & 0.8217230947055163655 \\
	$d$ & 0.2439762150462290542 & 0.2450849419266946147 & 0.2457239541094273887 \\
	$e$ & 0.7455840250177885510 & 0.7460303806901244447 & 0.7463726234745099378 \\
	$f$ & -0.5361015015276512265 & 0.6111264667009431951 & 0.6113411213084060140 \\
	$g$ & 0.7089476929608000620 & 0.7091435077517020016 & 0.7091567382315775358 \\
	$h$ & -0.1722213488026787314 & 0.6780618803338796457 & 0.6780081896535308014 \\
	$i$ & 0.4698382617159265503 & 0.4672482313963320756 & 0.4656700316384325305 \\
	$j$ & 0.6133979280072114038 & 0.6170533981914939001 & 0.6197515650327852649 \\
	$k$ & 0.4441986329257307108 & 0.4460948675469909250 & 0.4475039381692436896 \\
	$l$ & 0.6480765792525520477 & 0.6459914338748474451 & 0.6447102414946823372 \\
	$m$ & 0.3919946602865573912 & 0.3932378910441901092 & 0.3943428730599016595 \\
	$n$ & -0.5949992989412573021 & 0.9040175602878125879 & 0.9031085632585104393 \\
	$o$ & 0.3070616525343694029 & 0.3092036035703416460 & 0.3109695510771343408 \\
	$p$ & -0.1840371215847103444 & 0.9017986899560804310 & 0.9026416358405071867 \\
	$q$ & 0.06448214657458906535 & 0.06334611409114023852 & 0.06267936105544782376 \\
	$r$ & 0.5120179601540782911 & 0.4857562070706075314 & 0.4858788622016544170 \\
	$s$ & 0.05314413404895080137 & 0.05256724407145519520 & 0.05183187901744640248 \\
	$t$ & -0.8531627150792248086 & 0.8748966543883301824 & 0.8742268129934458916 \\
	$u$ & -0.4786918917229729443 & 0.9983900932233055668 & 0.9984699306689336426 \\
	$v$ & -0.8933449497206251001 & 0.8343894807573222920 & 0.8363999986054766405 \\
	$w$ & -0.9241564389899229240 & 0.6592877712136054758 & 0.6596658849521259834 \\
	$x$ & -0.6732220205865382325 & 0.8788231431229936877 & 0.8781347051779220827 \\
	$y$ & 0.8597710007152027251 & 0.8594647332919322702 & 0.8593484806198012411 \\
	$z$ & 0.2534084798238539341 & 0.2534654896181816607 & 0.2535733868746379384 \\
	$A$ & -0.6751607181858807123 & 0.5113567594914695497 & 0.5095418002721745686 \\
	$B$ & 0.6649647202743589516 & 0.6641257929299873136 & 0.6634574151328521656 \\
	$C$ & 0.4519752415520143349 & 0.06814420586753881720 & 0.06605262283446185188 \\
	$D$ & -0.2322132209099418502 & 0.2449912440290526918 & 0.2451035689165447480 \\
	\hline\Tstrut
	$energy$ & -410.5331627932077910 & 1543.830400976378803 & 1693.794611768318406
\end{longtable}

There were some problems encountered with optimizing the parameters, and the highest digit precision obtained for the \textit{logarithmic} potential was 32,829 digits, for the \textit{Coulomb} potential was 45,024 digits and for the \textit{inverse square} potential was 39,118 digits.

No algebraic numbers were recovered from the spherical codes, for all 3 potentials, their degree is $>420$, if found.

\noindent
\textit{\textbf{Symmetries -}}

All symmetry groups for 60 points are identical under all 3 potentials.

\begin{center}
	\begin{tabular}{l|l}
		\multicolumn{2}{c}{Symmetries - 60 points} \\
		\hline\Tstrut
		planes & [[20, 1]] \\[0.2ex]
		\hline\Tstrut
		Gram groups & [[6, 30], [12, 280], [60, 1]] \\
		\hline\Tstrut
		Polygons & [[3, 20]]
	\end{tabular}
\end{center}

\subsection{61 points}
All 3 optimal configurations for 61 points do not allow for a parameterization to be determined. This is the eighth set, the previous are 26, 35, 36, 54, 55, 56, and 59 points.

\begin{figure}[ht]
	\begin{center}
		\includegraphics[type=pdf,ext=pdf,read=pdf,height=1in,width=1in,angle=0]{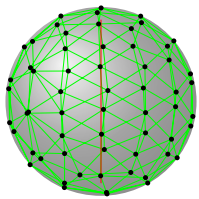}
		\caption{61 points.}
		\label{fig:61pts}
	\end{center}
\end{figure}
\noindent
\textit{\textbf{Minimal Energy values -}}

The coordinates for 61 points are known to 77 digits for the \textit{log} potential and 38 digits for the other two. The minimal energies have been determined for all 3 potentials as well.
\begin{center}
	\begin{tabular}{l|l}
		\multicolumn{2}{c}{Minimal Energy - 61 points} \\
		\hline\Tstrut
		logarithmic & -423.5076359910035272\ldots \\[0.2ex]
		\hline\Tstrut
		Coulomb & 1597.941830198989007\ldots \\[0.2ex]
		\hline\Tstrut
		Inverse square law & 1758.697133957985949\ldots
	\end{tabular}
\end{center}

\noindent
\textit{\textbf{Symmetries -}}

The symmetry groups for 61 points are identical under all 3 potentials.

\begin{center}
	\begin{tabular}{l|l}
		\multicolumn{2}{c}{Symmetries - 61 points} \\
		\hline\Tstrut
		planes & [] \\[0.2ex]
		\hline\Tstrut
		Gram groups & [[2, 1830], [61, 1]] \\
		\hline\Tstrut
		Polygons & []
	\end{tabular}
\end{center}

\subsection{62 points}

The polygon symmetry group indicates that there are 12 pentagons embedded in the configuration. A shrewd guess adds 2 poles and that is exactly what the minimal energy polyhedron is. There is only 1 axis of symmetry in this 62 point set, the dipole axis of the 2 isolated points and the 12 parallel pentagons perpendicular to this axis.

\begin{figure}[ht]
	\begin{center}
		\includegraphics[type=pdf,ext=pdf,read=pdf,height=1in,width=1in,angle=0]{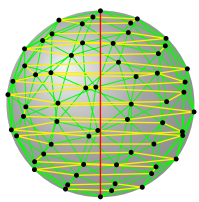}
		\caption{62 points.}
		\label{fig:62pts}
	\end{center}
\end{figure}

\noindent
\textit{\textbf{Algebraic parameterization -}}

Although there are 12 parallel pentagons along a polar axis, only 18 parameters are required to successfully parameterize 62 points.

The algebraic parameterized structure is given below:

\begin{longtable}[c]{r|ccc}
	\caption{Parameterization for 62 points} \\
	pt & $x$ & $y$ & $z$ \\
	\hline\vspace*{-2.2ex}
	\endfirsthead
	\multicolumn{4}{c}%
	{\tablename\ \thetable\ -- 62 points parameters -- \textit{continued}} \\
	pt & $x$ & $y$ & $z$ \\
	\hline\vspace*{-2.2ex}
	\endhead
	1 & $0$ & $0$ & $1$ \\[1.0ex]
	2 & $b$ & $0$ & $a$ \\[1.0ex]
	3 & $\frac{b\left(\sqrt{5}-1\right)}{4}$ & $\frac{b\sqrt{2\sqrt{5}+10}}{4}$ & $a$ \\[1.0ex]
	4 & $-\frac{b\left(\sqrt{5}+1\right)}{4}$ & $\frac{b\sqrt{10-2\sqrt{5}}}{4}$ & $a$ \\[1.0ex]
	5 & $-\frac{b\left(\sqrt{5}+1\right)}{4}$ & $-\frac{b\sqrt{10-2\sqrt{5}}}{4}$ & $a$ \\[1.0ex]
	6 & $\frac{b\left(\sqrt{5}-1\right)}{4}$ & $-\frac{b\sqrt{2\sqrt{5}+10}}{4}$ & $a$ \\[1.0ex]
	7 & $d$ & $-\sqrt{1-c^2-d^2}$ & $c$ \\[1.0ex]
	8 & $\frac{\sqrt{2\sqrt{5}+10}\sqrt{1-c^2}\sqrt{1-c^2-d^2}+d\left(\sqrt{5}-1\right)\sqrt{1-c^2}}{4\sqrt{1-c^2}}$ & $\frac{d\sqrt{2\sqrt{5}+10}\sqrt{1-c^2}-\left(\sqrt{5}-1\right)\sqrt{1-c^2}\sqrt{1-c^2-d^2}}{4\sqrt{1-c^2}}$ & $c$ \\[1.0ex]
	9 & $\frac{\sqrt{10-2\sqrt{5}}\sqrt{1-c^2}\sqrt{1-c^2-d^2}+d\left(-\sqrt{5}-1\right)\sqrt{1-c^2}}{4\sqrt{1-c^2}}$ & $\frac{\left(\sqrt{5}+1\right)\sqrt{1-c^2}\sqrt{1-c^2-d^2}+d\sqrt{10-2\sqrt{5}}\sqrt{1-c^2}}{4\sqrt{1-c^2}}$ & $c$ \\[1.0ex]
	10 & $\frac{\left(-\sqrt{10-2\sqrt{5}}\sqrt{1-c^2}\sqrt{1-c^2-d^2}\right)-d\left(\sqrt{5}+1\right)\sqrt{1-c^2}}{4\sqrt{1-c^2}}$ & $\frac{\left(\sqrt{5}+1\right)\sqrt{1-c^2}\sqrt{1-c^2-d^2}-d\sqrt{10-2\sqrt{5}}\sqrt{1-c^2}}{4\sqrt{1-c^2}}$ & $c$ \\[1.0ex]
	11 & $\frac{\left(-\sqrt{2\sqrt{5}+10}\sqrt{1-c^2}\sqrt{1-c^2-d^2}\right)-d\left(1-\sqrt{5}\right)\sqrt{1-c^2}}{4\sqrt{1-c^2}}$ & $\frac{\left(-\left(\sqrt{5}-1\right)\sqrt{1-c^2}\sqrt{1-c^2-d^2}\right)-d\sqrt{2\sqrt{5}+10}\sqrt{1-c^2}}{4\sqrt{1-c^2}}$ & $c$ \\[1.0ex]
	12 & $f$ & $\sqrt{1-e^2-f^2}$ & $e$ \\[1.0ex]
	13 & $\frac{\sqrt{2\sqrt{5}+10}\sqrt{1-e^2}\sqrt{1-e^2-f^2}+f\left(\sqrt{5}-1\right)\sqrt{1-e^2}}{4\sqrt{1-e^2}}$ & $-\frac{f\sqrt{2\sqrt{5}+10}\sqrt{1-e^2}-\left(\sqrt{5}-1\right)\sqrt{1-e^2}\sqrt{1-e^2-f^2}}{4\sqrt{1-e^2}}$ & $e$ \\[1.0ex]
	14 & $\frac{\sqrt{10-2\sqrt{5}}\sqrt{1-e^2}\sqrt{1-e^2-f^2}+f\left(-\sqrt{5}-1\right)\sqrt{1-e^2}}{4\sqrt{1-e^2}}$ & $-\frac{\left(\sqrt{5}+1\right)\sqrt{1-e^2}\sqrt{1-e^2-f^2}+f\sqrt{10-2\sqrt{5}}\sqrt{1-e^2}}{4\sqrt{1-e^2}}$ & $e$ \\[1.0ex]
	15 & $\frac{\left(-\sqrt{10-2\sqrt{5}}\sqrt{1-e^2}\sqrt{1-e^2-f^2}\right)-f\left(\sqrt{5}+1\right)\sqrt{1-e^2}}{4\sqrt{1-e^2}}$ & $-\frac{\left(\sqrt{5}+1\right)\sqrt{1-e^2}\sqrt{1-e^2-f^2}-f\sqrt{10-2\sqrt{5}}\sqrt{1-e^2}}{4\sqrt{1-e^2}}$ & $e$ \\[1.0ex]
	16 & $\frac{\left(-\sqrt{2\sqrt{5}+10}\sqrt{1-e^2}\sqrt{1-e^2-f^2}\right)-f\left(1-\sqrt{5}\right)\sqrt{1-e^2}}{4\sqrt{1-e^2}}$ & $-\frac{\left(-\left(\sqrt{5}-1\right)\sqrt{1-e^2}\sqrt{1-e^2-f^2}\right)-f\sqrt{2\sqrt{5}+10}\sqrt{1-e^2}}{4\sqrt{1-e^2}}$ & $e$ \\[1.0ex]
	17 & $h$ & $-\sqrt{1-g^2-h^2}$ & $g$ \\[1.0ex]
	18 & $\frac{\sqrt{2\sqrt{5}+10}\sqrt{1-g^2}\sqrt{1-g^2-h^2}+h\left(\sqrt{5}-1\right)\sqrt{1-g^2}}{4\sqrt{1-g^2}}$ & $\frac{h\sqrt{2\sqrt{5}+10}\sqrt{1-g^2}-\left(\sqrt{5}-1\right)\sqrt{1-g^2}\sqrt{1-g^2-h^2}}{4\sqrt{1-g^2}}$ & $g$ \\[1.0ex]
	19 & $\frac{\sqrt{10-2\sqrt{5}}\sqrt{1-g^2}\sqrt{1-g^2-h^2}+h\left(-\sqrt{5}-1\right)\sqrt{1-g^2}}{4\sqrt{1-g^2}}$ & $\frac{\left(\sqrt{5}+1\right)\sqrt{1-g^2}\sqrt{1-g^2-h^2}+h\sqrt{10-2\sqrt{5}}\sqrt{1-g^2}}{4\sqrt{1-g^2}}$ & $g$ \\[1.0ex]
	20 & $\frac{\left(-\sqrt{10-2\sqrt{5}}\sqrt{1-g^2}\sqrt{1-g^2-h^2}\right)-h\left(\sqrt{5}+1\right)\sqrt{1-g^2}}{4\sqrt{1-g^2}}$ & $\frac{\left(\sqrt{5}+1\right)\sqrt{1-g^2}\sqrt{1-g^2-h^2}-h\sqrt{10-2\sqrt{5}}\sqrt{1-g^2}}{4\sqrt{1-g^2}}$ & $g$ \\[1.0ex]
	21 & $\frac{\left(-\sqrt{2\sqrt{5}+10}\sqrt{1-g^2}\sqrt{1-g^2-h^2}\right)-h\left(1-\sqrt{5}\right)\sqrt{1-g^2}}{4\sqrt{1-g^2}}$ & $\frac{\left(-\left(\sqrt{5}-1\right)\sqrt{1-g^2}\sqrt{1-g^2-h^2}\right)-h\sqrt{2\sqrt{5}+10}\sqrt{1-g^2}}{4\sqrt{1-g^2}}$ & $g$ \\[1.0ex]
	22 & $j$ & $\sqrt{1-i^2-j^2}$ & $i$ \\[1.0ex]
	23 & $\frac{\sqrt{2\sqrt{5}+10}\sqrt{1-i^2}\sqrt{1-i^2-j^2}+j\left(\sqrt{5}-1\right)\sqrt{1-i^2}}{4\sqrt{1-i^2}}$ & $-\frac{j\sqrt{2\sqrt{5}+10}\sqrt{1-i^2}-\left(\sqrt{5}-1\right)\sqrt{1-i^2}\sqrt{1-i^2-j^2}}{4\sqrt{1-i^2}}$ & $i$ \\[1.0ex]
	24 & $\frac{\sqrt{10-2\sqrt{5}}\sqrt{1-i^2}\sqrt{1-i^2-j^2}+j\left(-\sqrt{5}-1\right)\sqrt{1-i^2}}{4\sqrt{1-i^2}}$ & $-\frac{\left(\sqrt{5}+1\right)\sqrt{1-i^2}\sqrt{1-i^2-j^2}+j\sqrt{10-2\sqrt{5}}\sqrt{1-i^2}}{4\sqrt{1-i^2}}$ & $i$ \\[1.0ex]
	25 & $\frac{\left(-\sqrt{10-2\sqrt{5}}\sqrt{1-i^2}\sqrt{1-i^2-j^2}\right)-j\left(\sqrt{5}+1\right)\sqrt{1-i^2}}{4\sqrt{1-i^2}}$ & $-\frac{\left(\sqrt{5}+1\right)\sqrt{1-i^2}\sqrt{1-i^2-j^2}-j\sqrt{10-2\sqrt{5}}\sqrt{1-i^2}}{4\sqrt{1-i^2}}$ & $i$ \\[1.0ex]
	26 & $\frac{\left(-\sqrt{2\sqrt{5}+10}\sqrt{1-i^2}\sqrt{1-i^2-j^2}\right)-j\left(1-\sqrt{5}\right)\sqrt{1-i^2}}{4\sqrt{1-i^2}}$ & $-\frac{\left(-\left(\sqrt{5}-1\right)\sqrt{1-i^2}\sqrt{1-i^2-j^2}\right)-j\sqrt{2\sqrt{5}+10}\sqrt{1-i^2}}{4\sqrt{1-i^2}}$ & $i$ \\[1.0ex]
	27 & $l$ & $-\sqrt{1-k^2-l^2}$ & $k$ \\[1.0ex]
	28 & $\frac{\sqrt{2\sqrt{5}+10}\sqrt{1-k^2}\sqrt{1-k^2-l^2}+l\left(\sqrt{5}-1\right)\sqrt{1-k^2}}{4\sqrt{1-k^2}}$ & $\frac{l\sqrt{2\sqrt{5}+10}\sqrt{1-k^2}-\left(\sqrt{5}-1\right)\sqrt{1-k^2}\sqrt{1-k^2-l^2}}{4\sqrt{1-k^2}}$ & $k$ \\[1.0ex]
	29 & $\frac{\sqrt{10-2\sqrt{5}}\sqrt{1-k^2}\sqrt{1-k^2-l^2}+l\left(-\sqrt{5}-1\right)\sqrt{1-k^2}}{4\sqrt{1-k^2}}$ & $\frac{\left(\sqrt{5}+1\right)\sqrt{1-k^2}\sqrt{1-k^2-l^2}+l\sqrt{10-2\sqrt{5}}\sqrt{1-k^2}}{4\sqrt{1-k^2}}$ & $k$ \\[1.0ex]
	30 & $\frac{\left(-\sqrt{10-2\sqrt{5}}\sqrt{1-k^2}\sqrt{1-k^2-l^2}\right)-l\left(\sqrt{5}+1\right)\sqrt{1-k^2}}{4\sqrt{1-k^2}}$ & $\frac{\left(\sqrt{5}+1\right)\sqrt{1-k^2}\sqrt{1-k^2-l^2}-l\sqrt{10-2\sqrt{5}}\sqrt{1-k^2}}{4\sqrt{1-k^2}}$ & $k$ \\[1.0ex]
	31 & $\frac{\left(-\sqrt{2\sqrt{5}+10}\sqrt{1-k^2}\sqrt{1-k^2-l^2}\right)-l\left(1-\sqrt{5}\right)\sqrt{1-k^2}}{4\sqrt{1-k^2}}$ & $\frac{\left(-\left(\sqrt{5}-1\right)\sqrt{1-k^2}\sqrt{1-k^2-l^2}\right)-l\sqrt{2\sqrt{5}+10}\sqrt{1-k^2}}{4\sqrt{1-k^2}}$ & $k$ \\[1.0ex]
	32 & $m$ & $-\sqrt{1-k^2-m^2}$ & $-k$ \\[1.0ex]
	33 & $\frac{\sqrt{2\sqrt{5}+10}\sqrt{1-k^2}\sqrt{1-k^2-m^2}+m\left(\sqrt{5}-1\right)\sqrt{1-k^2}}{4\sqrt{1-k^2}}$ & $\frac{m\sqrt{2\sqrt{5}+10}\sqrt{1-k^2}-\left(\sqrt{5}-1\right)\sqrt{1-k^2}\sqrt{1-k^2-m^2}}{4\sqrt{1-k^2}}$ & $-k$ \\[1.0ex]
	34 & $\frac{\sqrt{10-2\sqrt{5}}\sqrt{1-k^2}\sqrt{1-k^2-m^2}+m\left(-\sqrt{5}-1\right)\sqrt{1-k^2}}{4\sqrt{1-k^2}}$ & $\frac{\left(\sqrt{5}+1\right)\sqrt{1-k^2}\sqrt{1-k^2-m^2}+m\sqrt{10-2\sqrt{5}}\sqrt{1-k^2}}{4\sqrt{1-k^2}}$ & $-k$ \\[1.0ex]
	35 & $\frac{\left(-\sqrt{10-2\sqrt{5}}\sqrt{1-k^2}\sqrt{1-k^2-m^2}\right)-m\left(\sqrt{5}+1\right)\sqrt{1-k^2}}{4\sqrt{1-k^2}}$ & $\frac{\left(\sqrt{5}+1\right)\sqrt{1-k^2}\sqrt{1-k^2-m^2}-m\sqrt{10-2\sqrt{5}}\sqrt{1-k^2}}{4\sqrt{1-k^2}}$ & $-k$ \\[1.0ex]
	36 & $\frac{\left(-\sqrt{2\sqrt{5}+10}\sqrt{1-k^2}\sqrt{1-k^2-m^2}\right)-m\left(1-\sqrt{5}\right)\sqrt{1-k^2}}{4\sqrt{1-k^2}}$ & $\frac{\left(-\left(\sqrt{5}-1\right)\sqrt{1-k^2}\sqrt{1-k^2-m^2}\right)-m\sqrt{2\sqrt{5}+10}\sqrt{1-k^2}}{4\sqrt{1-k^2}}$ & $-k$ \\[1.0ex]
	37 & $n$ & $\sqrt{1-i^2-n^2}$ & $-i$ \\[1.0ex]
	38 & $\frac{\sqrt{2\sqrt{5}+10}\sqrt{1-i^2}\sqrt{1-i^2-n^2}+n\left(\sqrt{5}-1\right)\sqrt{1-i^2}}{4\sqrt{1-i^2}}$ & $-\frac{n\sqrt{2\sqrt{5}+10}\sqrt{1-i^2}-\left(\sqrt{5}-1\right)\sqrt{1-i^2}\sqrt{1-i^2-n^2}}{4\sqrt{1-i^2}}$ & $-i$ \\[1.0ex]
	39 & $\frac{\sqrt{10-2\sqrt{5}}\sqrt{1-i^2}\sqrt{1-i^2-n^2}+n\left(-\sqrt{5}-1\right)\sqrt{1-i^2}}{4\sqrt{1-i^2}}$ & $-\frac{\left(\sqrt{5}+1\right)\sqrt{1-i^2}\sqrt{1-i^2-n^2}+n\sqrt{10-2\sqrt{5}}\sqrt{1-i^2}}{4\sqrt{1-i^2}}$ & $-i$ \\[1.0ex]
	40 & $\frac{\left(-\sqrt{10-2\sqrt{5}}\sqrt{1-i^2}\sqrt{1-i^2-n^2}\right)-n\left(\sqrt{5}+1\right)\sqrt{1-i^2}}{4\sqrt{1-i^2}}$ & $-\frac{\left(\sqrt{5}+1\right)\sqrt{1-i^2}\sqrt{1-i^2-n^2}-n\sqrt{10-2\sqrt{5}}\sqrt{1-i^2}}{4\sqrt{1-i^2}}$ & $-i$ \\[1.0ex]
	41 & $\frac{\left(-\sqrt{2\sqrt{5}+10}\sqrt{1-i^2}\sqrt{1-i^2-n^2}\right)-n\left(1-\sqrt{5}\right)\sqrt{1-i^2}}{4\sqrt{1-i^2}}$ & $-\frac{\left(-\left(\sqrt{5}-1\right)\sqrt{1-i^2}\sqrt{1-i^2-n^2}\right)-n\sqrt{2\sqrt{5}+10}\sqrt{1-i^2}}{4\sqrt{1-i^2}}$ & $-i$ \\[1.0ex]
	42 & $o$ & $-\sqrt{1-g^2-o^2}$ & $-g$ \\[1.0ex]
	43 & $\frac{\sqrt{2\sqrt{5}+10}\sqrt{1-g^2}\sqrt{1-g^2-o^2}+o\left(\sqrt{5}-1\right)\sqrt{1-g^2}}{4\sqrt{1-g^2}}$ & $\frac{o\sqrt{2\sqrt{5}+10}\sqrt{1-g^2}-\left(\sqrt{5}-1\right)\sqrt{1-g^2}\sqrt{1-g^2-o^2}}{4\sqrt{1-g^2}}$ & $-g$ \\[1.0ex]
	44 & $\frac{\sqrt{10-2\sqrt{5}}\sqrt{1-g^2}\sqrt{1-g^2-o^2}+o\left(-\sqrt{5}-1\right)\sqrt{1-g^2}}{4\sqrt{1-g^2}}$ & $\frac{\left(\sqrt{5}+1\right)\sqrt{1-g^2}\sqrt{1-g^2-o^2}+o\sqrt{10-2\sqrt{5}}\sqrt{1-g^2}}{4\sqrt{1-g^2}}$ & $-g$ \\[1.0ex]
	45 & $\frac{\left(-\sqrt{10-2\sqrt{5}}\sqrt{1-g^2}\sqrt{1-g^2-o^2}\right)-o\left(\sqrt{5}+1\right)\sqrt{1-g^2}}{4\sqrt{1-g^2}}$ & $\frac{\left(\sqrt{5}+1\right)\sqrt{1-g^2}\sqrt{1-g^2-o^2}-o\sqrt{10-2\sqrt{5}}\sqrt{1-g^2}}{4\sqrt{1-g^2}}$ & $-g$ \\[1.0ex]
	46 & $\frac{\left(-\sqrt{2\sqrt{5}+10}\sqrt{1-g^2}\sqrt{1-g^2-o^2}\right)-o\left(1-\sqrt{5}\right)\sqrt{1-g^2}}{4\sqrt{1-g^2}}$ & $\frac{\left(-\left(\sqrt{5}-1\right)\sqrt{1-g^2}\sqrt{1-g^2-o^2}\right)-o\sqrt{2\sqrt{5}+10}\sqrt{1-g^2}}{4\sqrt{1-g^2}}$ & $-g$ \\[1.0ex]
	47 & $p$ & $-\sqrt{1-e^2-p^2}$ & $-e$ \\[1.0ex]
	48 & $\frac{\sqrt{2\sqrt{5}+10}\sqrt{1-e^2}\sqrt{1-e^2-p^2}+p\left(\sqrt{5}-1\right)\sqrt{1-e^2}}{4\sqrt{1-e^2}}$ & $\frac{p\sqrt{2\sqrt{5}+10}\sqrt{1-e^2}-\left(\sqrt{5}-1\right)\sqrt{1-e^2}\sqrt{1-e^2-p^2}}{4\sqrt{1-e^2}}$ & $-e$ \\[1.0ex]
	49 & $\frac{\sqrt{10-2\sqrt{5}}\sqrt{1-e^2}\sqrt{1-e^2-p^2}+p\left(-\sqrt{5}-1\right)\sqrt{1-e^2}}{4\sqrt{1-e^2}}$ & $\frac{\left(\sqrt{5}+1\right)\sqrt{1-e^2}\sqrt{1-e^2-p^2}+p\sqrt{10-2\sqrt{5}}\sqrt{1-e^2}}{4\sqrt{1-e^2}}$ & $-e$ \\[1.0ex]
	50 & $\frac{\left(-\sqrt{10-2\sqrt{5}}\sqrt{1-e^2}\sqrt{1-e^2-p^2}\right)-p\left(\sqrt{5}+1\right)\sqrt{1-e^2}}{4\sqrt{1-e^2}}$ & $\frac{\left(\sqrt{5}+1\right)\sqrt{1-e^2}\sqrt{1-e^2-p^2}-p\sqrt{10-2\sqrt{5}}\sqrt{1-e^2}}{4\sqrt{1-e^2}}$ & $-e$ \\[1.0ex]
	51 & $\frac{\left(-\sqrt{2\sqrt{5}+10}\sqrt{1-e^2}\sqrt{1-e^2-p^2}\right)-p\left(1-\sqrt{5}\right)\sqrt{1-e^2}}{4\sqrt{1-e^2}}$ & $\frac{\left(-\left(\sqrt{5}-1\right)\sqrt{1-e^2}\sqrt{1-e^2-p^2}\right)-p\sqrt{2\sqrt{5}+10}\sqrt{1-e^2}}{4\sqrt{1-e^2}}$ & $-e$ \\[1.0ex]
	52 & $q$ & $\sqrt{1-c^2-q^2}$ & $-c$ \\[1.0ex]
	53 & $\frac{\sqrt{2\sqrt{5}+10}\sqrt{1-c^2}\sqrt{1-c^2-q^2}+q\left(\sqrt{5}-1\right)\sqrt{1-c^2}}{4\sqrt{1-c^2}}$ & $-\frac{q\sqrt{2\sqrt{5}+10}\sqrt{1-c^2}-\left(\sqrt{5}-1\right)\sqrt{1-c^2}\sqrt{1-c^2-q^2}}{4\sqrt{1-c^2}}$ & $-c$ \\[1.0ex]
	54 & $\frac{\sqrt{10-2\sqrt{5}}\sqrt{1-c^2}\sqrt{1-c^2-q^2}+q\left(-\sqrt{5}-1\right)\sqrt{1-c^2}}{4\sqrt{1-c^2}}$ & $-\frac{\left(\sqrt{5}+1\right)\sqrt{1-c^2}\sqrt{1-c^2-q^2}+q\sqrt{10-2\sqrt{5}}\sqrt{1-c^2}}{4\sqrt{1-c^2}}$ & $-c$ \\[1.0ex]
	55 & $\frac{\left(-\sqrt{10-2\sqrt{5}}\sqrt{1-c^2}\sqrt{1-c^2-q^2}\right)-q\left(\sqrt{5}+1\right)\sqrt{1-c^2}}{4\sqrt{1-c^2}}$ & $-\frac{\left(\sqrt{5}+1\right)\sqrt{1-c^2}\sqrt{1-c^2-q^2}-q\sqrt{10-2\sqrt{5}}\sqrt{1-c^2}}{4\sqrt{1-c^2}}$ & $-c$ \\[1.0ex]
	56 & $\frac{\left(-\sqrt{2\sqrt{5}+10}\sqrt{1-c^2}\sqrt{1-c^2-q^2}\right)-q\left(1-\sqrt{5}\right)\sqrt{1-c^2}}{4\sqrt{1-c^2}}$ & $-\frac{\left(-\left(\sqrt{5}-1\right)\sqrt{1-c^2}\sqrt{1-c^2-q^2}\right)-q\sqrt{2\sqrt{5}+10}\sqrt{1-c^2}}{4\sqrt{1-c^2}}$ & $-c$ \\[1.0ex]
	57 & $r$ & $-\sqrt{1-a^2-r^2}$ & $-a$ \\[1.0ex]
	58 & $\frac{\sqrt{2\sqrt{5}+10}\sqrt{1-a^2}\sqrt{1-a^2-r^2}+r\left(\sqrt{5}-1\right)\sqrt{1-a^2}}{4\sqrt{1-a^2}}$ & $\frac{r\sqrt{2\sqrt{5}+10}\sqrt{1-a^2}-\left(\sqrt{5}-1\right)\sqrt{1-a^2}\sqrt{1-a^2-r^2}}{4\sqrt{1-a^2}}$ & $-a$ \\[1.0ex]
	59 & $\frac{\sqrt{10-2\sqrt{5}}\sqrt{1-a^2}\sqrt{1-a^2-r^2}+r\left(-\sqrt{5}-1\right)\sqrt{1-a^2}}{4\sqrt{1-a^2}}$ & $\frac{\left(\sqrt{5}+1\right)\sqrt{1-a^2}\sqrt{1-a^2-r^2}+r\sqrt{10-2\sqrt{5}}\sqrt{1-a^2}}{4\sqrt{1-a^2}}$ & $-a$ \\[1.0ex]
	60 & $\frac{\left(-\sqrt{10-2\sqrt{5}}\sqrt{1-a^2}\sqrt{1-a^2-r^2}\right)-r\left(\sqrt{5}+1\right)\sqrt{1-a^2}}{4\sqrt{1-a^2}}$ & $\frac{\left(\sqrt{5}+1\right)\sqrt{1-a^2}\sqrt{1-a^2-r^2}-r\sqrt{10-2\sqrt{5}}\sqrt{1-a^2}}{4\sqrt{1-a^2}}$ & $-a$ \\[1.0ex]
	61 & $\frac{\left(-\sqrt{2\sqrt{5}+10}\sqrt{1-a^2}\sqrt{1-a^2-r^2}\right)-r\left(1-\sqrt{5}\right)\sqrt{1-a^2}}{4\sqrt{1-a^2}}$ & $\frac{\left(-\left(\sqrt{5}-1\right)\sqrt{1-a^2}\sqrt{1-a^2-r^2}\right)-r\sqrt{2\sqrt{5}+10}\sqrt{1-a^2}}{4\sqrt{1-a^2}}$ & $-a$ \\[1.0ex]
	62 & $0$ & $0$ & $-1$
\end{longtable}

\noindent
\textit{\textbf{Algebraic Spherical codes --}}

NOTE: It was not possible to converge the parameters for this set, the problem[PM] has been noted.

\noindent
\textit{\textbf{Parameterization values --}}

The values to 19 digits of the 18 parameters optimized for the minimal solutions of 62 points are:

\begin{longtable}[c]{c|c|c|c}
	\caption{Parameter values for 62 points} \\
	Parameter & log & 1/r & $1/r^2$ \\
	\hline\vspace*{-2.2ex}
	\endfirsthead
	\multicolumn{4}{c}%
	{\tablename\ \thetable\ -- 62 points parameters -- \textit{continued}} \\
	Parameter & log & 1/r & $1/r^2$ \\
	\hline\vspace*{-2.2ex}
	\endhead
	$a$ & 0.8950581099815161237 & 0.8952083481019980206 & 0.8952339116954621106 \\
	$b$ & 0.4459495260187145740 & 0.4456478581666156027 & 0.4455965028480828956 \\
	$c$ & 0.7115887990152248083 & 0.7107706440297346969 & 0.7099490257300468083 \\
	$d$ & 0.5746048585596977335 & 0.5743587404383671248 & 0.5742627517028423773 \\
	$e$ & 0.5864862285863754005 & 0.5862702425810699829 & 0.5857939777097406471 \\
	$f$ & 0.8084360911288725634 & 0.8088034672543563375 & 0.8093093178784396598 \\
	$g$ & 0.3690420338197005303 & 0.3675426452984964825 & 0.3660986012383575400 \\
	$h$ & 0.8610567795529732377 & 0.8608140550571361876 & 0.8606581647024180031 \\
	$i$ & 0.2643284845982289281 & 0.2641801970253393328 & 0.2637019731715880652 \\
	$j$ & 0.8886956868624394914 & 0.8892630517590619250 & 0.8898122357303886774 \\
	$k$ & 0.08506848391343052810 & 0.08497888913923011226 & 0.08475230614790407076 \\
	$l$ & 0.9961679000985086491 & 0.9961195370005955137 & 0.9960832635788161194 \\
	$m$ & 0.8674486547995652933 & 0.8669490839105702982 & 0.8666962829838363170 \\
	$n$ & 0.9145817098405395773 & 0.9152933846401610372 & 0.9159361977219742895 \\
	$o$ & 0.9191367301527322441 & 0.9195351002895553899 & 0.9199689322576672590 \\
	$p$ & 0.6702679237606925788 & 0.6707567191801557171 & 0.6713912319861264496 \\
	$q$ & 0.7004377999404746757 & 0.7013364915873509779 & 0.7022149889616542512 \\
	$r$ & 0.3836906893312154824 & 0.3826056746728261826 & 0.3819068081101356551 \\
	\hline\Tstrut
	$energy$ & -436.7039792383948156 & 1652.909409898300852 & 1824.346926363182267
\end{longtable}

NOTE: The values for the 18 parameters were extended to 1,021 digits precision by the \textit{descent.3d} or \textit{direct} search algorithms for all 3 potentials due.

\noindent
\textit{\textbf{Symmetries -}}

The symmetry groups for 62 points are identical under all 3 potentials.

\begin{center}
	\begin{tabular}{l|l}
		\multicolumn{2}{c}{Symmetries - 62 points} \\
		\hline\Tstrut
		planes & [[120, 1]] \\[0.2ex]
		\hline\Tstrut
		Gram groups & [[2, 1], [10, 30], [20, 174], [62, 1]] \\
		\hline\Tstrut
		Polygons & [[5, 12]]
	\end{tabular}
\end{center}

\subsection{63 points}

The polygon symmetry group of this configuration shows the optimal alignment of 21 parallel triangles. Indeed rotating the randomly oriented minimal set of points into alignment showed the 21 triangles immediately, thus there is only 1 axis of alignment for this set of 63 points.

\begin{figure}[ht]
	\begin{center}
		\includegraphics[type=pdf,ext=pdf,read=pdf,height=1in,width=1in,angle=0]{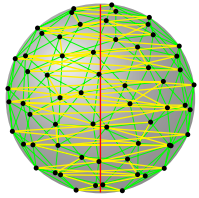}
		\caption{63 points.}
		\label{fig:63pts}
	\end{center}
\end{figure}

\noindent
\textit{\textbf{Algebraic parameterization -}}

It requires 20 parameters to adequately constrain the 63 point set, due to the 21 embedded triangles in the figure. Since the configuration is odd, the parameter count is halved for an even count of parallel polygons.

The algebraic parameterized structure is given below:

\begin{longtable}[c]{r|ccc}
	\caption{Parameterization for 63 points} \\
	pt & $x$ & $y$ & $z$ \\
	\hline\vspace*{-2.2ex}
	\endfirsthead
	\multicolumn{4}{c}%
	{\tablename\ \thetable\ -- 63 points parameters -- \textit{continued}} \\
	pt & $x$ & $y$ & $z$ \\
	\hline\vspace*{-2.2ex}
	\endhead
	1 & $b$ & $-\sqrt{1-a^2-b^2}$ & $a$ \\[0.7ex]
	2 & $\frac{-b+\sqrt{3}\sqrt{1-a^2-b^2}}{2}$ & $\frac{+b\sqrt{3}+\sqrt{1-a^2-b^2}}{2}$ & $a$ \\[0.7ex]
	3 & $\frac{-b-\sqrt{3}\sqrt{1-a^2-b^2}}{2}$ & $\frac{-b\sqrt{3}+\sqrt{1-a^2-b^2}}{2}$ & $a$ \\[0.7ex]
	4 & $d$ & $\sqrt{1-c^2-d^2}$ & $c$ \\[0.7ex]
	5 & $\frac{-d-\sqrt{3}\sqrt{1-c^2-d^2}}{2}$ & $\frac{+d\sqrt{3}-\sqrt{1-c^2-d^2}}{2}$ & $c$ \\[0.7ex]
	6 & $\frac{-d+\sqrt{3}\sqrt{1-c^2-d^2}}{2}$ & $\frac{-d\sqrt{3}-\sqrt{1-c^2-d^2}}{2}$ & $c$ \\[0.7ex]
	7 & $f$ & $\sqrt{1-e^2-f^2}$ & $e$ \\[0.7ex]
	8 & $\frac{-f-\sqrt{3}\sqrt{1-e^2-f^2}}{2}$ & $\frac{+f\sqrt{3}-\sqrt{1-e^2-f^2}}{2}$ & $e$ \\[0.7ex]
	9 & $\frac{-f+\sqrt{3}\sqrt{1-e^2-f^2}}{2}$ & $\frac{-f\sqrt{3}-\sqrt{1-e^2-f^2}}{2}$ & $e$ \\[0.7ex]
	10 & $h$ & $-\sqrt{1-g^2-h^2}$ & $g$ \\[0.7ex]
	11 & $\frac{-h+\sqrt{3}\sqrt{1-g^2-h^2}}{2}$ & $\frac{+h\sqrt{3}+\sqrt{1-g^2-h^2}}{2}$ & $g$ \\[0.7ex]
	12 & $\frac{-h-\sqrt{3}\sqrt{1-g^2-h^2}}{2}$ & $\frac{-h\sqrt{3}+\sqrt{1-g^2-h^2}}{2}$ & $g$ \\[0.7ex]
	13 & $j$ & $\sqrt{1-i^2-j^2}$ & $i$ \\[0.7ex]
	14 & $\frac{-j-\sqrt{3}\sqrt{1-i^2-j^2}}{2}$ & $\frac{+j\sqrt{3}-\sqrt{1-i^2-j^2}}{2}$ & $i$ \\[0.7ex]
	15 & $\frac{-j+\sqrt{3}\sqrt{1-i^2-j^2}}{2}$ & $\frac{-j\sqrt{3}-\sqrt{1-i^2-j^2}}{2}$ & $i$ \\[0.7ex]
	16 & $l$ & $-\sqrt{1-k^2-l^2}$ & $k$ \\[0.7ex]
	17 & $\frac{-l+\sqrt{3}\sqrt{1-k^2-l^2}}{2}$ & $\frac{+l\sqrt{3}+\sqrt{1-k^2-l^2}}{2}$ & $k$ \\[0.7ex]
	18 & $\frac{-l-\sqrt{3}\sqrt{1-k^2-l^2}}{2}$ & $\frac{-l\sqrt{3}+\sqrt{1-k^2-l^2}}{2}$ & $k$ \\[0.7ex]
	19 & $n$ & $-\sqrt{1-m^2-n^2}$ & $m$ \\[0.7ex]
	20 & $\frac{-n+\sqrt{3}\sqrt{1-m^2-n^2}}{2}$ & $\frac{+n\sqrt{3}+\sqrt{1-m^2-n^2}}{2}$ & $m$ \\[0.7ex]
	21 & $\frac{-n-\sqrt{3}\sqrt{1-m^2-n^2}}{2}$ & $\frac{-n\sqrt{3}+\sqrt{1-m^2-n^2}}{2}$ & $m$ \\[0.7ex]
	22 & $p$ & $\sqrt{1-o^2-p^2}$ & $o$ \\[0.7ex]
	23 & $\frac{-p-\sqrt{3}\sqrt{1-o^2-p^2}}{2}$ & $\frac{+p\sqrt{3}-\sqrt{1-o^2-p^2}}{2}$ & $o$ \\[0.7ex]
	24 & $\frac{-p+\sqrt{3}\sqrt{1-o^2-p^2}}{2}$ & $\frac{-p\sqrt{3}-\sqrt{1-o^2-p^2}}{2}$ & $o$ \\[0.7ex]
	25 & $r$ & $\sqrt{1-q^2-r^2}$ & $q$ \\[0.7ex]
	26 & $\frac{-r-\sqrt{3}\sqrt{1-q^2-r^2}}{2}$ & $\frac{+r\sqrt{3}-\sqrt{1-q^2-r^2}}{2}$ & $q$ \\[0.7ex]
	27 & $\frac{-r+\sqrt{3}\sqrt{1-q^2-r^2}}{2}$ & $\frac{-r\sqrt{3}-\sqrt{1-q^2-r^2}}{2}$ & $q$ \\[0.7ex]
	28 & $t$ & $-\sqrt{1-s^2-t^2}$ & $s$ \\[0.7ex]
	29 & $\frac{-t+\sqrt{3}\sqrt{1-s^2-t^2}}{2}$ & $\frac{+t\sqrt{3}+\sqrt{1-s^2-t^2}}{2}$ & $s$ \\[0.7ex]
	30 & $\frac{-t-\sqrt{3}\sqrt{1-s^2-t^2}}{2}$ & $\frac{-t\sqrt{3}+\sqrt{1-s^2-t^2}}{2}$ & $s$ \\[0.7ex]
	31 & $1$ & $0$ & $0$ \\[0.7ex]
	32 & $-\frac{1}{2}$ & $\frac{\sqrt{3}}{2}$ & $0$ \\[0.7ex]
	33 & $-\frac{1}{2}$ & $-\frac{\sqrt{3}}{2}$ & $0$ \\[0.7ex]
	34 & $t$ & $\sqrt{1-s^2-t^2}$ & $-s$ \\[0.7ex]
	35 & $\frac{-t+\sqrt{3}\sqrt{1-s^2-t^2}}{2}$ & $\frac{-t\sqrt{3}-\sqrt{1-s^2-t^2}}{2}$ & $-s$ \\[0.7ex]
	36 & $\frac{-t-\sqrt{3}\sqrt{1-s^2-t^2}}{2}$ & $\frac{+t\sqrt{3}-\sqrt{1-s^2-t^2}}{2}$ & $-s$ \\[0.7ex]
	37 & $r$ & $-\sqrt{1-q^2-r^2}$ & $-q$ \\[0.7ex]
	38 & $\frac{-r-\sqrt{3}\sqrt{1-q^2-r^2}}{2}$ & $\frac{-r\sqrt{3}+\sqrt{1-q^2-r^2}}{2}$ & $-q$ \\[0.7ex]
	39 & $\frac{-r+\sqrt{3}\sqrt{1-q^2-r^2}}{2}$ & $\frac{+r\sqrt{3}+\sqrt{1-q^2-r^2}}{2}$ & $-q$ \\[0.7ex]
	40 & $p$ & $-\sqrt{1-o^2-p^2}$ & $-o$ \\[0.7ex]
	41 & $\frac{-p-\sqrt{3}\sqrt{1-o^2-p^2}}{2}$ & $\frac{-p\sqrt{3}+\sqrt{1-o^2-p^2}}{2}$ & $-o$ \\[0.7ex]
	42 & $\frac{-p+\sqrt{3}\sqrt{1-o^2-p^2}}{2}$ & $\frac{+p\sqrt{3}+\sqrt{1-o^2-p^2}}{2}$ & $-o$ \\[0.7ex]
	43 & $n$ & $\sqrt{1-m^2-n^2}$ & $-m$ \\[0.7ex]
	44 & $\frac{-n+\sqrt{3}\sqrt{1-m^2-n^2}}{2}$ & $\frac{-n\sqrt{3}-\sqrt{1-m^2-n^2}}{2}$ & $-m$ \\[0.7ex]
	45 & $\frac{-n-\sqrt{3}\sqrt{1-m^2-n^2}}{2}$ & $\frac{+n\sqrt{3}-\sqrt{1-m^2-n^2}}{2}$ & $-m$ \\[0.7ex]
	46 & $l$ & $\sqrt{1-k^2-l^2}$ & $-k$ \\[0.7ex]
	47 & $\frac{-l+\sqrt{3}\sqrt{1-k^2-l^2}}{2}$ & $\frac{-l\sqrt{3}-\sqrt{1-k^2-l^2}}{2}$ & $-k$ \\[0.7ex]
	48 & $\frac{-l-\sqrt{3}\sqrt{1-k^2-l^2}}{2}$ & $\frac{+l\sqrt{3}-\sqrt{1-k^2-l^2}}{2}$ & $-k$ \\[0.7ex]
	49 & $j$ & $-\sqrt{1-i^2-j^2}$ & $-i$ \\[0.7ex]
	50 & $\frac{-j-\sqrt{3}\sqrt{1-i^2-j^2}}{2}$ & $\frac{-j\sqrt{3}+\sqrt{1-i^2-j^2}}{2}$ & $-i$ \\[0.7ex]
	51 & $\frac{-j+\sqrt{3}\sqrt{1-i^2-j^2}}{2}$ & $\frac{+j\sqrt{3}+\sqrt{1-i^2-j^2}}{2}$ & $-i$ \\[0.7ex]
	52 & $h$ & $\sqrt{1-g^2-h^2}$ & $-g$ \\[0.7ex]
	53 & $\frac{-h+\sqrt{3}\sqrt{1-g^2-h^2}}{2}$ & $\frac{-h\sqrt{3}-\sqrt{1-g^2-h^2}}{2}$ & $-g$ \\[0.7ex]
	54 & $\frac{-h-\sqrt{3}\sqrt{1-g^2-h^2}}{2}$ & $\frac{+h\sqrt{3}-\sqrt{1-g^2-h^2}}{2}$ & $-g$ \\[0.7ex]
	55 & $f$ & $-\sqrt{1-e^2-f^2}$ & $-e$ \\[0.7ex]
	56 & $\frac{-f-\sqrt{3}\sqrt{1-e^2-f^2}}{2}$ & $\frac{-f\sqrt{3}+\sqrt{1-e^2-f^2}}{2}$ & $-e$ \\[0.7ex]
	57 & $\frac{-f+\sqrt{3}\sqrt{1-e^2-f^2}}{2}$ & $\frac{+f\sqrt{3}+\sqrt{1-e^2-f^2}}{2}$ & $-e$ \\[0.7ex]
	58 & $d$ & $-\sqrt{1-c^2-d^2}$ & $-c$ \\[0.7ex]
	59 & $\frac{-d-\sqrt{3}\sqrt{1-c^2-d^2}}{2}$ & $\frac{-d\sqrt{3}+\sqrt{1-c^2-d^2}}{2}$ & $-c$ \\[0.7ex]
	60 & $\frac{-d+\sqrt{3}\sqrt{1-c^2-d^2}}{2}$ & $\frac{+d\sqrt{3}+\sqrt{1-c^2-d^2}}{2}$ & $-c$ \\[0.7ex]
	61 & $b$ & $\sqrt{1-a^2-b^2}$ & $-a$ \\[0.7ex]
	62 & $\frac{-b+\sqrt{3}\sqrt{1-a^2-b^2}}{2}$ & $\frac{-b\sqrt{3}-\sqrt{1-a^2-b^2}}{2}$ & $-a$ \\[0.7ex]
	63 & $\frac{-b-\sqrt{3}\sqrt{1-a^2-b^2}}{2}$ & $\frac{+b\sqrt{3}-\sqrt{1-a^2-b^2}}{2}$ & $-a$
\end{longtable}

The algebraic degree of the 20 parameters for the 50,014 digit spherical codes is $>420$ for this configuration of 63 points.

\noindent
\textit{\textbf{Parameterization values --}}

The values to 19 digits of the 20 parameters optimized for the minimal solutions of 63 points are:

\begin{longtable}[c]{c|c|c|c}
	\caption{Parameter values for 63 points} \\
	Parameter & log & 1/r & $1/r^2$ \\
	\hline\vspace*{-2.2ex}
	\endfirsthead
	\multicolumn{4}{c}%
	{\tablename\ \thetable\ -- 63 points parameters -- \textit{continued}} \\
	Parameter & log & 1/r & $1/r^2$ \\
	\hline\vspace*{-2.2ex}
	\endhead
	$a$ & 0.9592843536470330307 & 0.9591733374794084792 & 0.9591527382980791461 \\
	$b$ & 0.2784499668446515127 & 0.2786060792466474427 & 0.2785282709608904655 \\
	$c$ & 0.8537430590304877323 & 0.8549097899629047863 & 0.8561151780143505098 \\
	$d$ & 0.3278422333190501673 & 0.3285373986243144746 & 0.3287258106800371257 \\
	$e$ & 0.7322180887100727287 & 0.7322119387367449276 & 0.7323867656935174771 \\
	$f$ & 0.6709648796942250290 & 0.6710701176022733879 & 0.6709284051521776053 \\
	$g$ & 0.7207367785426339164 & 0.7211622602787105148 & 0.7217601318822867451 \\
	$h$ & 0.5920461534016540140 & 0.5906022087084819848 & 0.5894151846510219718 \\
	$i$ & 0.5531407663829349176 & 0.5545741480045833097 & 0.5558831624471157560 \\
	$j$ & 0.5787486675728906156 & 0.5770624199531665560 & 0.5757945992416161897 \\
	$k$ & 0.4314660549086865556 & 0.4320721611641190887 & 0.4325848386049733501 \\
	$l$ & 0.8863782857905197106 & 0.8859588450740977860 & 0.8856730627945439364 \\
	$m$ & 0.3899706322226931493 & 0.3868306768459442495 & 0.3844318610257709651 \\
	$n$ & 0.6401658427302375791 & 0.6410991920123018300 & 0.6418030770225270027 \\
	$o$ & 0.3596758790244882692 & 0.3603258311337392918 & 0.3607371530440930105 \\
	$p$ & 0.8791656688451146780 & 0.8788983656095874008 & 0.8787869935928616277 \\
	$q$ & 0.1325733574371501390 & 0.1329180810841212154 & 0.1332399853868698276 \\
	$r$ & 0.6423333327180779402 & 0.6435121426045185309 & 0.6444888194748329973 \\
	$s$ & 0.07684620483891040093 & 0.07616597530795952632 & 0.07575707619651162835 \\
	$t$ & 0.8889476867017887354 & 0.8902391424612079999 & 0.8913276522889820095 \\
	\hline\Tstrut
	$energy$ & -450.0812391769127035 & 1708.879681503249984 & 1891.632330787479861
\end{longtable}

\noindent
\textit{\textbf{Symmetries -}}

The symmetry groups are identical for 63 points under all 3 potentials.

\begin{center}
	\begin{tabular}{l|l}
		\multicolumn{2}{c}{Symmetries - 63 points} \\
		\hline\Tstrut
		planes & [[21, 1]] \\[0.2ex]
		\hline\Tstrut
		Gram groups & [[6, 31], [12, 310], [63, 1]] \\
		\hline\Tstrut
		Polygons & [[3, 21]]
	\end{tabular}
\end{center}

\subsection{64 points}
All 3 optimal configurations for 64 points do not allow for a parameterization to occur. This is the ninth set, the previous eight are 26, 35, 36, 54, 55, 56, 59 and 61 points.

\begin{figure}[ht]
	\begin{center}
		\includegraphics[type=pdf,ext=pdf,read=pdf,height=1in,width=1in,angle=0]{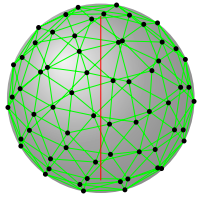}
		\caption{64 points.}
		\label{fig:64pts}
	\end{center}
\end{figure}
\noindent
\textit{\textbf{Minimal Energy values -}}

The coordinates for 64 points are known to 77 digits for the \textit{log} potential and 38 digits for the other two. The minimal energies have been determined for all 3 potentials as well.
\begin{center}
	\begin{tabular}{l|l}
		\multicolumn{2}{c}{Minimal Energy - 64 points} \\
		\hline\Tstrut
		logarithmic & -463.6544329886080460\ldots \\[0.2ex]
		\hline\Tstrut
		Coulomb & 1765.802577927303190\ldots \\[0.2ex]
		\hline\Tstrut
		Inverse square law & 1960.281626974364263\ldots
	\end{tabular}
\end{center}

\noindent
\textit{\textbf{Symmetries -}}

The symmetry groups for 64 points are identical under all 3 potentials.

\begin{center}
	\begin{tabular}{l|l}
		\multicolumn{2}{c}{Symmetries - 64 points} \\
		\hline\Tstrut
		planes & [] \\[0.2ex]
		\hline\Tstrut
		Gram groups & [[4, 48], [8, 480], [64, 1]] \\
		\hline\Tstrut
		Polygons & []
	\end{tabular}
\end{center}

\subsection{65 points}
All 3 optimal configurations for 65 points do not allow for a parameterization to occur. This is the tenth set, the previous nine are 26, 35, 36, 54, 55, 56, 59, 61 and 64 points. In general, the trend points out that trying to parameterize a $n$-point cluster will generally not be possible for $>65$ points, but only for special cases.

\begin{figure}[ht]
	\begin{center}
		\includegraphics[type=pdf,ext=pdf,read=pdf,height=1in,width=1in,angle=0]{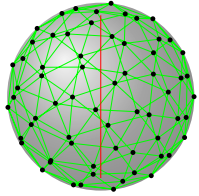}
		\caption{65 points.}
		\label{fig:65pts}
	\end{center}
\end{figure}
\noindent
\textit{\textbf{Minimal Energy values -}}

The coordinates for 65 points are known to 77 digits for the \textit{log} potential and 38 digits for the other two. The minimal energies have been determined for all 3 potentials as well.
\begin{center}
	\begin{tabular}{l|l}
		\multicolumn{2}{c}{Minimal Energy - 65 points} \\
		\hline\Tstrut
		logarithmic & -477.4264260688034746\ldots \\[0.2ex]
		\hline\Tstrut
		Coulomb & 1823.667960263850183\ldots \\[0.2ex]
		\hline\Tstrut
		Inverse square law & 2030.233678508393013\ldots
	\end{tabular}
\end{center}

\noindent
\textit{\textbf{Symmetries -}}

The symmetry groups for 65 points are identical under all 3 potentials.

\begin{center}
	\begin{tabular}{l|l}
		\multicolumn{2}{c}{Symmetries - 65 points} \\
		\hline\Tstrut
		planes & [] \\[0.2ex]
		\hline\Tstrut
		Gram groups & [[2, 32], [4, 1024], [65, 1]] \\
		\hline\Tstrut
		Polygons & []
	\end{tabular}
\end{center}

\subsection{Comments about embedded polygons}

After having found the global minimums and rotating the configurations into desired orientations, it is necessary to make some comments about the polygonal embedding in these solutions.

\begin{itemize}
	\item \textbf{Embedded co-planar equilateral triangles}
\end{itemize}

It is interesting to note that 3 points, 9 points, 12 points, 15 points, 45 points, 51 points, 57 points, 60 points, and 63 points; all consist of parallel triangles sets of 1, 3, 4, 5, 15, 17, 19, 20, and 21 equilateral triangles for their optimal minimal energy arrangement for all 3 potentials. There may be more for points sets $>65$ which contain $n$ co-planar triangles for the optimal arrangement.

Additionally, certain point sets, such as 5, 16, 20 (with a hexagon), 23, 28, 29 (Coulomb and Inverse Square potentials), 31, 46, 49, 52 are composed of co-planar equilateral triangles, with one or two pole points added, indicating that stability is retained by using the co-planar triangle sets, even if a point(s) is deficient or surplus to $n=0\mod 3$. In the case of 2 points added, they become poles.

It is also noted that other point sets contain co-planar triangles + other polygons, for example: 20 points contains 4 triangles, 2 poles and 1 hexagon, 22 points contains 5 triangles, 1 hexagon and one pole, while 39 points is a mix of hexagons and triangles.

41 points is an interesting mix of triangles, hexagons and one nonagon.

\noindent
Other polygons can be involved.

\begin{itemize}
	\item \textbf{Embedded co-planar squares}
\end{itemize}

In the case of squares, 6 points is 1 square plus 2 poles, 8 points is 2 squares, 10 points is two squares plus 2 poles, 18 points has 4 squares plus 2 pole points, 24 points is 6 squares, 40 points is 6 squares plus 8 dipoles, and 48 points is 12 squares for the optimum configuration.

\begin{itemize}
	\item \textbf{Embedded co-planar regular pentagons}
\end{itemize}

Regarding regular pentagons, we have 7 points with 1 pentagon and 2 poles, 17 points with 3 pentagons and 2 poles, 27 points with 5 pentagons and 2 poles, and 37 points with 7 pentagons and 2 poles, and 62 points is 12 pentagons plus 2 poles.

42 points is 6 pentagons plus 1 decagon and 2 poles.

\begin{itemize}
	\item \textbf{Embedded co-planar regular hexagons}
\end{itemize}

Considering regular hexagons, we find 14 points with 2 hexagons and 2 poles, 38 points has 6 hexagons and 2 poles, and 50 points has 8 hexagons and 2 poles.

\begin{itemize}
	\item \textbf{Embedded co-planar regular heptagon or nonagon}
\end{itemize}

The 21 points configuration contains a heptagon, and 53 points contains 1 nonagon and 22 dipoles.

\begin{itemize}
	\item \textbf{Sole dipoles arrangement}
\end{itemize}

It is interesting that 34 points consists of 17 dipoles for the optimum arrangement, although it contains 16 groups of quadrilaterals, but they are not on the preferred orientation.

\textit{It is expected that other larger size ($>65$) point sets might continue this behavior of embedded regular polygons for their global minimum solution.}

\section{Summary of Algebraic Results}

A summary of the searches for global minima is presented below, classified by the potential, then sorted numerically by point count. Both the points and energy are listed.

A key is provided:\hspace*{1em}
\begin{tabular}[c]{rlrl}
	\cellcolor{na} & -- not applicable &
	\cellcolor{solv} & -- algebraic code obtained \\
	\cellcolor{tent} & -- tentative, partially solved code &
	\cellcolor{real} & -- parameterized to 50,014 digits \\
	\cellcolor{newt} & -- Newton's method to 50K digits &
	\cellcolor{impo} & -- impossible to parameterize \\
	\cellcolor{unkn} & -- status unknown &
	\cellcolor{bug} & -- software problem \\
	P & -- coordinates of point &
	E & -- minimal energy \Tstrut
\end{tabular}

\begin{longtable}[c]{r|C{0.15in}C{0.15in}|C{0.15in}C{0.15in}|C{0.15in}C{0.15in}|C{0.15in}C{0.15in}|C{0.15in}C{0.15in}|C{0.15in}C{0.15in}|C{0.15in}C{0.15in}|C{0.15in}C{0.15in}|C{0.15in}C{0.15in}|C{0.15in}C{0.15in}}
	\caption{Algebraic Spherical Codes on S2 - Logarithmic Potential} \\
	& P & E & P & E & P & E & P & E & P & E & P & E & P & E & P & E & P & E & P & E \\
	& \multicolumn{2}{c}{0} & \multicolumn{2}{c}{1} & \multicolumn{2}{c}{2} & \multicolumn{2}{c}{3} & \multicolumn{2}{c}{4} & \multicolumn{2}{c}{5} & \multicolumn{2}{c}{6} & \multicolumn{2}{c}{7} & \multicolumn{2}{c}{8} & \multicolumn{2}{c}{9} \\
	\hline\vspace*{-2.2ex}
	\endfirsthead
	\multicolumn{21}{c}%
	{\tablename\ \thetable\ -- Spherical code -- Logarithmic Potential -- \textit{continued}} \\[0.5ex]
	& P & E & P & E & P & E & P & E & P & E & P & E & P & E & P & E & P & E & P & E \\
	& \multicolumn{2}{c}{0} & \multicolumn{2}{c}{1} & \multicolumn{2}{c}{2} & \multicolumn{2}{c}{3} & \multicolumn{2}{c}{4} & \multicolumn{2}{c}{5} & \multicolumn{2}{c}{6} & \multicolumn{2}{c}{7} & \multicolumn{2}{c}{8} & \multicolumn{2}{c}{9} \\
	\hline\vspace*{-2.2ex}
	\endhead
	0 & \cellcolor{na} & \cellcolor{na} 0 & \cellcolor{solv} 1 & \cellcolor{solv} 1 & \cellcolor{solv} 2 & \cellcolor{solv} 2 & \cellcolor{solv} 3 & \cellcolor{solv} 3 & \cellcolor{solv} 4 & \cellcolor{solv} 4 & \cellcolor{solv} 5 & \cellcolor{solv} 5 & \cellcolor{solv} 6 & \cellcolor{solv} 6 & \cellcolor{solv} 7 & \cellcolor{solv} 7 & \cellcolor{solv} 8 & \cellcolor{solv} 8 & \cellcolor{solv} 9 & \cellcolor{solv} 9 \\
	10 & \cellcolor{solv} 10 & \cellcolor{solv} 10 & \cellcolor{real} 11 & \cellcolor{real} 11 & \cellcolor{solv} 12 & \cellcolor{solv} 12 & \cellcolor{real} 13 & \cellcolor{real} 13 & \cellcolor{solv} 14 & \cellcolor{solv} 14 & \cellcolor{real} 15 & \cellcolor{real} 15 & \cellcolor{solv} 16 & \cellcolor{solv} 16 & \cellcolor{solv} 17 & \cellcolor{solv} 17 & \cellcolor{solv} 18 & \cellcolor{solv} 18 & \cellcolor{real} 19 & \cellcolor{real} 19 \\
	20 & \cellcolor{tent} 20 & \cellcolor{real} 20 & \cellcolor{real} 21 & \cellcolor{real} 21 & \cellcolor{solv} 22 & \cellcolor{solv} 22 & \cellcolor{real} 23 & \cellcolor{real} 23 & \cellcolor{real} 24 & \cellcolor{real} 24 & \cellcolor{real} 25 & \cellcolor{real} 25 & \cellcolor{impo} 26 & \cellcolor{unkn} 26 & \cellcolor{solv} 27 & \cellcolor{real} 27 & \cellcolor{real} 28 & \cellcolor{real} 28 & \cellcolor{impo} 29 & \cellcolor{unkn} 29 \\
	30 & \cellcolor{real} 30 & \cellcolor{real} 30 & \cellcolor{real} 31 & \cellcolor{real} 31 & \cellcolor{solv} 32 & \cellcolor{solv} 32 & \cellcolor{real} 33 & \cellcolor{real} 33 & \cellcolor{real} 34 & \cellcolor{real} 34 & \cellcolor{impo} 35 & \cellcolor{unkn} 35 & \cellcolor{impo} 36 & \cellcolor{unkn} 36 & \cellcolor{bug} 37 & \cellcolor{bug} 37 & \cellcolor{real} 38 & \cellcolor{real} 38 & \cellcolor{real} 39 & \cellcolor{real} 39 \\
	40 & \cellcolor{real} 40 & \cellcolor{real} 40 & \cellcolor{real} 41 & \cellcolor{real} 41 & \cellcolor{real} 42 & \cellcolor{real} 42 & \cellcolor{real} 43 & \cellcolor{real} 43 & \cellcolor{solv} 44 & \cellcolor{real} 44 & \cellcolor{real} 45 & \cellcolor{real} 45 & \cellcolor{bug} 46 & \cellcolor{bug} 46 & \cellcolor{unkn} 47 & \cellcolor{unkn} 47 & \cellcolor{bug} 48 & \cellcolor{bug} 48 & \cellcolor{bug} 49 & \cellcolor{bug} 49 \\
	50 & \cellcolor{real} 50 & \cellcolor{real} 50 & \cellcolor{real} 51 & \cellcolor{real} 51 & \cellcolor{real} 52 & \cellcolor{real} 52 & \cellcolor{bug} 53 & \cellcolor{bug} 53 & \cellcolor{impo} 54 & \cellcolor{unkn} 54 & \cellcolor{impo} 55 & \cellcolor{unkn} 55 & \cellcolor{impo} 56 & \cellcolor{unkn} 56 & \cellcolor{real} 57 & \cellcolor{real} 57 & \cellcolor{unkn} 58 & \cellcolor{unkn} 58 & \cellcolor{impo} 59 & \cellcolor{unkn} 59 \\
	60 & \cellcolor{bug} 60 & \cellcolor{bug} 60 & \cellcolor{impo} 61 & \cellcolor{unkn} 61 & \cellcolor{real} 62 & \cellcolor{real} 62 & \cellcolor{real} 63 & \cellcolor{real} 63 & \cellcolor{impo} 64 & \cellcolor{unkn} 64 & \cellcolor{impo} 65 & \cellcolor{unkn} 65 \\
\end{longtable}

\begin{longtable}[c]{r|C{0.15in}C{0.15in}|C{0.15in}C{0.15in}|C{0.15in}C{0.15in}|C{0.15in}C{0.15in}|C{0.15in}C{0.15in}|C{0.15in}C{0.15in}|C{0.15in}C{0.15in}|C{0.15in}C{0.15in}|C{0.15in}C{0.15in}|C{0.15in}C{0.15in}}
	\caption{Algebraic Spherical Codes on S2 - Coulomb $1/r$ Potential} \\
	& P & E & P & E & P & E & P & E & P & E & P & E & P & E & P & E & P & E & P & E \\
	& \multicolumn{2}{c}{0} & \multicolumn{2}{c}{1} & \multicolumn{2}{c}{2} & \multicolumn{2}{c}{3} & \multicolumn{2}{c}{4} & \multicolumn{2}{c}{5} & \multicolumn{2}{c}{6} & \multicolumn{2}{c}{7} & \multicolumn{2}{c}{8} & \multicolumn{2}{c}{9} \\
	\hline\vspace*{-2.2ex}
	\endfirsthead
	\multicolumn{21}{c}%
	{\tablename\ \thetable\ -- Spherical code -- Coulomb $1/r$ Potential -- \textit{continued}} \\[0.5ex]
	& P & E & P & E & P & E & P & E & P & E & P & E & P & E & P & E & P & E & P & E \\
	& \multicolumn{2}{c}{0} & \multicolumn{2}{c}{1} & \multicolumn{2}{c}{2} & \multicolumn{2}{c}{3} & \multicolumn{2}{c}{4} & \multicolumn{2}{c}{5} & \multicolumn{2}{c}{6} & \multicolumn{2}{c}{7} & \multicolumn{2}{c}{8} & \multicolumn{2}{c}{9} \\
	\hline\vspace*{-2.2ex}
	\endhead
	0 & \cellcolor{na} & \cellcolor{na} 0 & \cellcolor{solv} 1 & \cellcolor{solv} 1 & \cellcolor{solv} 2 & \cellcolor{solv} 2 & \cellcolor{solv} 3 & \cellcolor{solv} 3 & \cellcolor{solv} 4 & \cellcolor{solv} 4 & \cellcolor{solv} 5 & \cellcolor{solv} 5 & \cellcolor{solv} 6 & \cellcolor{solv} 6 & \cellcolor{solv} 7 & \cellcolor{solv} 7 & \cellcolor{solv} 8 & \cellcolor{solv} 8 & \cellcolor{solv} 9 & \cellcolor{solv} 9 \\
	10 & \cellcolor{solv} 10 & \cellcolor{solv} 10 & \cellcolor{real} 11 & \cellcolor{real} 11 & \cellcolor{solv} 12 & \cellcolor{solv} 12 & \cellcolor{real} 13 & \cellcolor{real} 13 & \cellcolor{real} 14 & \cellcolor{real} 14 & \cellcolor{real} 15 & \cellcolor{real} 15 & \cellcolor{real} 16 & \cellcolor{real} 16 & \cellcolor{tent} 17 & \cellcolor{real} 17 & \cellcolor{real} 18 & \cellcolor{real} 18 & \cellcolor{real} 19 & \cellcolor{real} 19 \\
	20 & \cellcolor{real} 20 & \cellcolor{real} 20 & \cellcolor{real} 21 & \cellcolor{real} 21 & \cellcolor{real} 22 & \cellcolor{real} 22 & \cellcolor{real} 23 & \cellcolor{real} 23 & \cellcolor{real} 24 & \cellcolor{real} 24 & \cellcolor{real} 25 & \cellcolor{real} 25 & \cellcolor{impo} 26 & \cellcolor{unkn} 26 & \cellcolor{real} 27 & \cellcolor{real} 27 & \cellcolor{real} 28 & \cellcolor{real} 28 & \cellcolor{real} 29 & \cellcolor{real} 29 \\
	30 & \cellcolor{real} 30 & \cellcolor{real} 30 & \cellcolor{real} 31 & \cellcolor{real} 31 & \cellcolor{solv} 32 & \cellcolor{solv} 32 & \cellcolor{real} 33 & \cellcolor{real} 33 & \cellcolor{real} 34 & \cellcolor{real} 34 & \cellcolor{impo} 35 & \cellcolor{unkn} 35 & \cellcolor{impo} 36 & \cellcolor{unkn} 36 & \cellcolor{bug} 37 & \cellcolor{bug} 37 & \cellcolor{real} 38 & \cellcolor{real} 38 & \cellcolor{real} 39 & \cellcolor{real} 39 \\
	40 & \cellcolor{real} 40 & \cellcolor{real} 40 & \cellcolor{real} 41 & \cellcolor{real} 41 & \cellcolor{real} 42 & \cellcolor{real} 42 & \cellcolor{real} 43 & \cellcolor{real} 43 & \cellcolor{real} 44 & \cellcolor{real} 44 & \cellcolor{real} 45 & \cellcolor{real} 45 & \cellcolor{bug} 46 & \cellcolor{bug} 46 & \cellcolor{unkn} 47 & \cellcolor{unkn} 47 & \cellcolor{bug} 48 & \cellcolor{bug} 48 & \cellcolor{bug} 49 & \cellcolor{bug} 49 \\
	50 & \cellcolor{real} 50 & \cellcolor{real} 50 & \cellcolor{real} 51 & \cellcolor{real} 51 & \cellcolor{real} 52 & \cellcolor{real} 52 & \cellcolor{bug} 53 & \cellcolor{bug} 53 & \cellcolor{impo} 54 & \cellcolor{unkn} 54 & \cellcolor{impo} 55 & \cellcolor{unkn} 55 & \cellcolor{impo} 56 & \cellcolor{unkn} 56 & \cellcolor{real} 57 & \cellcolor{real} 57 & \cellcolor{unkn} 58 & \cellcolor{unkn} 58 & \cellcolor{impo} 59 & \cellcolor{unkn} 59 \\
	60 & \cellcolor{bug} 60 & \cellcolor{bug} 60 & \cellcolor{impo} 61 & \cellcolor{unkn} 61 & \cellcolor{real} 62 & \cellcolor{real} 62 & \cellcolor{real} 63 & \cellcolor{real} 63 & \cellcolor{impo} 64 & \cellcolor{unkn} 64 & \cellcolor{impo} 65 & \cellcolor{unkn} 65 \\
\end{longtable}

\begin{longtable}[c]{r|C{0.15in}C{0.15in}|C{0.15in}C{0.15in}|C{0.15in}C{0.15in}|C{0.15in}C{0.15in}|C{0.15in}C{0.15in}|C{0.15in}C{0.15in}|C{0.15in}C{0.15in}|C{0.15in}C{0.15in}|C{0.15in}C{0.15in}|C{0.15in}C{0.15in}}
	\caption{Algebraic Spherical Codes on S2 - Inverse Square $1/r^2$ Potential} \\
	& P & E & P & E & P & E & P & E & P & E & P & E & P & E & P & E & P & E & P & E \\
	& \multicolumn{2}{c}{0} & \multicolumn{2}{c}{1} & \multicolumn{2}{c}{2} & \multicolumn{2}{c}{3} & \multicolumn{2}{c}{4} & \multicolumn{2}{c}{5} & \multicolumn{2}{c}{6} & \multicolumn{2}{c}{7} & \multicolumn{2}{c}{8} & \multicolumn{2}{c}{9} \\
	\hline\vspace*{-2.2ex}
	\endfirsthead
	\multicolumn{21}{c}%
	{\tablename\ \thetable\ -- Spherical code -- Inverse Square $1/r^2$ Potential -- \textit{continued}} \\[0.5ex]
	& P & E & P & E & P & E & P & E & P & E & P & E & P & E & P & E & P & E & P & E \\
	& \multicolumn{2}{c}{0} & \multicolumn{2}{c}{1} & \multicolumn{2}{c}{2} & \multicolumn{2}{c}{3} & \multicolumn{2}{c}{4} & \multicolumn{2}{c}{5} & \multicolumn{2}{c}{6} & \multicolumn{2}{c}{7} & \multicolumn{2}{c}{8} & \multicolumn{2}{c}{9} \\
	\hline\vspace*{-2.2ex}
	\endhead
	0 & \cellcolor{na} & \cellcolor{na} 0 & \cellcolor{solv} 1 & \cellcolor{solv} 1 & \cellcolor{solv} 2 & \cellcolor{solv} 2 & \cellcolor{solv} 3 & \cellcolor{solv} 3 & \cellcolor{solv} 4 & \cellcolor{solv} 4 & \cellcolor{solv} 5 & \cellcolor{solv} 5 & \cellcolor{solv} 6 & \cellcolor{solv} 6 & \cellcolor{solv} 7 & \cellcolor{solv} 7 & \cellcolor{solv} 8 & \cellcolor{solv} 8 & \cellcolor{solv} 9 & \cellcolor{solv} 9 \\
	10 & \cellcolor{solv} 10 & \cellcolor{solv} 10 & \cellcolor{real} 11 & \cellcolor{real} 11 & \cellcolor{solv} 12 & \cellcolor{solv} 12 & \cellcolor{real} 13 & \cellcolor{real} 13 & \cellcolor{solv} 14 & \cellcolor{solv} 14 & \cellcolor{real} 15 & \cellcolor{real} 15 & \cellcolor{solv} 16 & \cellcolor{solv} 16 & \cellcolor{solv} 17 & \cellcolor{solv} 17 & \cellcolor{solv} 18 & \cellcolor{real} 18 & \cellcolor{newt} 19 & \cellcolor{newt} 19 \\
	20 & \cellcolor{real} 20 & \cellcolor{real} 20 & \cellcolor{real} 21 & \cellcolor{real} 21 & \cellcolor{solv} 22 & \cellcolor{solv} 22 & \cellcolor{real} 23 & \cellcolor{real} 23 & \cellcolor{tent} 24 & \cellcolor{real} 24 & \cellcolor{real} 25 & \cellcolor{real} 25 & \cellcolor{impo} 26 & \cellcolor{unkn} 26 & \cellcolor{real} 27 & \cellcolor{real} 27 & \cellcolor{real} 28 & \cellcolor{real} 28 & \cellcolor{real} 29 & \cellcolor{real} 29 \\
	30 & \cellcolor{real} 30 & \cellcolor{real} 30 & \cellcolor{real} 31 & \cellcolor{real} 31 & \cellcolor{solv} 32 & \cellcolor{solv} 32 & \cellcolor{real} 33 & \cellcolor{real} 33 & \cellcolor{real} 34 & \cellcolor{real} 34 & \cellcolor{impo} 35 & \cellcolor{unkn} 35 & \cellcolor{impo} 36 & \cellcolor{unkn} 36 & \cellcolor{bug} 37 & \cellcolor{bug} 37 & \cellcolor{real} 38 & \cellcolor{real} 38 & \cellcolor{real} 39 & \cellcolor{real} 39 \\
	40 & \cellcolor{real} 40 & \cellcolor{real} 40 & \cellcolor{real} 41 & \cellcolor{real} 41 & \cellcolor{real} 42 & \cellcolor{real} 42 & \cellcolor{real} 43 & \cellcolor{real} 43 & \cellcolor{solv} 44 & \cellcolor{real} 44 & \cellcolor{real} 45 & \cellcolor{real} 45 & \cellcolor{bug} 46 & \cellcolor{bug} 46 & \cellcolor{unkn} 47 & \cellcolor{unkn} 47 & \cellcolor{bug} 48 & \cellcolor{bug} 48 & \cellcolor{bug} 49 & \cellcolor{bug} 49 \\
	50 & \cellcolor{real} 50 & \cellcolor{real} 50 & \cellcolor{real} 51 & \cellcolor{real} 51 & \cellcolor{real} 52 & \cellcolor{real} 52 & \cellcolor{bug} 53 & \cellcolor{bug} 53 & \cellcolor{impo} 54 & \cellcolor{unkn} 54 & \cellcolor{impo} 55 & \cellcolor{unkn} 55 & \cellcolor{impo} 56 & \cellcolor{unkn} 56 & \cellcolor{real} 57 & \cellcolor{real} 57 & \cellcolor{unkn} 58 & \cellcolor{unkn} 58 & \cellcolor{impo} 59 & \cellcolor{unkn} 59 \\
	60 & \cellcolor{bug} 60 & \cellcolor{bug} 60 & \cellcolor{impo} 61 & \cellcolor{unkn} 61 & \cellcolor{real} 62 & \cellcolor{real} 62 & \cellcolor{real} 63 & \cellcolor{real} 63 & \cellcolor{impo} 64 & \cellcolor{unkn} 64 & \cellcolor{impo} 65 & \cellcolor{unkn} 65 \\
\end{longtable}

\vspace{3ex}

\section{Current Status of Unknown Parameters}

There are 1585 unknown parameters in 109 point sets belonging to the 3 potentials. The following table lists what is known about the real decimal values of the parameters. Most (1124) are at least 50,014 digits precision, but 461 exceptions are carefully noted.

The Unknown Status Table provides a challenge for the current mathematical tool sets working with determining the algebraic polynomial given a decimal number. The current tool used is the \textit{GP-Pari algdep(n,d)} command where $n$ is the decimal number in question and $d$ is the tentative degree of the algebraic polynomial we seek.

It has been discovered that these algebraic polynomials have a general form of:
\begin{equation}
	n = a_0 + a_1x^2 + a_2x^4 + _{\dots} + a_n x^{2n}
\end{equation}
where $a_0 , a_1, a_2, _{\dots}, a_n$ are the algebraic number coefficients and $x$ is always an even power.

The "found" column in the table indicates whether or not the algebraic number was discovered for the potential/$n$-point set. A count is given if some of the parameters were resolved.

\begin{longtable}[c]{r|crlrl}
	\caption{\textbf{Status of parameters for Algebraic polynomials}} \\
	potential & pts/parms & Digits & found & Height & Notes \\
	\hline\vspace*{-1.3ex}
	\endfirsthead
	\multicolumn{6}{c}%
	{\tablename\ \thetable\ -- \textbf{Algebraic spherical code parameters} -- \textit{continued}} \\[0.5ex]
	potential & pts/parms & Digits & found & Height & Notes \\
	\hline\vspace*{-1.3ex}
	\endhead
	\multicolumn{6}{c}{\textit{Logarithmic Potential}} \\
	$log$ & [11,5] & 50014 & N & 360 & \\
	$log$ & [13,6] & 50014 & N & 360 & \\
	$log$ & [15,4] & 50014 & N & 360 & \\
	$log$ & [19,9] & 50014 & N & 360 & \\
	$log$ & [20,3] & 50014 & 2 & 500 & 2 of 5 resolved \\
	$log$ & [21,18] & 50014 & N & 360 & \\
	$log$ & [23,6] & 50014 & N & 360 & \\
	$log$ & [24,9] & 22579 & N & 600 & \\
	$log$ & [25,24] & 50014 & N & 420 & \\
	$log$ & [28,16] & 50014 & N & 420 & \\
	$log$ & [30,20] & 50014 & N & 420 & \\
	$log$ & [31,18] & 6049 & N & 420 & \\
	$log$ & [33,32] & 50014 & N & 420 & \\
	$log$ & [34,16] & 50014 & N & 360 & \\
	$log$ & [37,6] & 10018 & N & 420 & \\
	$log$ & [38,3] & 50014 & N & 360 & \\
	$log$ & [39,10] & 50014 & N & 360 & \\
	$log$ & [40,9] & 50014 & N & 360 & \\
	$log$ & [41,10] & 50014 & N & 360 & \\
	$log$ & [42,5] & 50014 & N & 360 & \\
	$log$ & [43,21] & 50014 & N & 360 & \\
	$log$ & [45,14] & 50014 & N & 360 & \\
	$log$ & [46,24] & 50014 & N & 420 & \\
	$log$ & [48,18] & 1021 & N & 60 & \\
	$log$ & [49,32] & 48299 & N & 420 & \\
	$log$ & [50,4] & 50014 & N & 360 & \\
	$log$ & [51,16] & 50014 & N & 360 & \\
	$log$ & [52,34] & 50014 & N & 420 & \\
	$log$ & [53,27] & 1001 & N & 60 & \\
	$log$ & [57,18] & 50014 & N & 360 & \\
	$log$ & [60,30] & 32829 & N & 420 & \\
	$log$ & [62,18] & 1001 & N & 60 & \\
	$log$ & [63,20] & 50014 & N & 360 & \\
	\hline
	\multicolumn{6}{c}{\textit{Coulomb $1/r$ Potential} \mystrut(13,0)} \\[0.5ex]
	$1/r$ & [11,5] & 50014 & N & 360 & \\
	$1/r$ & [13,6] & 50014 & N & 360 & \\
	$1/r$ & [14,3] & 50014 & N & 360 & \\
	$1/r$ & [15,4] & 50014 & N & 360 & \\
	$1/r$ & [16,8] & 50014 & N & 360 & 480 in a,c,e,g \\
	$1/r$ & [17,5] & 50014 & N & 600 & \\
	$1/r$ & [18,2] & 50014 & N & 480 & \\
	$1/r$ & [19,9] & 50014 & N & 360 & \\
	$1/r$ & [20,5] & 50014 & N & 360 & \\
	$1/r$ & [21,18] & 50014 & N & 360 & \\
	$1/r$ & [22,6] & 50014 & N & 360 & \\
	$1/r$ & [23,6] & 50014 & N & 360 & \\
	$1/r$ & [24,9] & 25816 & N & 600 & \\
	$1/r$ & [25,24] & 50014 & N & 420 & \\
	$1/r$ & [27,2] & 50014 & N & 360 & \\
	$1/r$ & [28,16] & 50014 & N & 420 & \\
	$1/r$ & [29,8] & 50014 & N & 420 & \\
	$1/r$ & [30,20] & 50014 & N & 420 & \\
	$1/r$ & [31,18] & 6049 & N & 420 & \\
	$1/r$ & [33,32] & 50014 & N & 420 & \\
	$1/r$ & [34,16] & 50014 & N & 360 & \\
	$1/r$ & [37,6] & 10018 & N & 420 & \\
	$1/r$ & [38,3] & 50014 & N & 360 & \\
	$1/r$ & [39,10] & 50014 & N & 360 & \\
	$1/r$ & [40,9] & 50014 & N & 360 & \\
	$1/r$ & [41,10] & 50014 & N & 360 & \\
	$1/r$ & [42,5] & 50014 & N & 360 & \\
	$1/r$ & [43,21] & 50014 & N & 360 & \\
	$1/r$ & [44,3] & 50014 & N & 360 & \\
	$1/r$ & [45,14] & 50014 & N & 360 & \\
	$1/r$ & [46,24] & 50014 & N & 420 & \\
	$1/r$ & [48,18] & 1021 & N & 60 & \\
	$1/r$ & [49,32] & 28205 & N & 420 & \\
	$1/r$ & [50,4] & 50014 & N & 360 & \\
	$1/r$ & [51,16] & 50014 & N & 360 & \\
	$1/r$ & [52,34] & 50014 & N & 420 & \\
	$1/r$ & [53,27] & 1001 & N & 60 & \\
	$1/r$ & [57,18] & 50014 & N & 360 & \\
	$1/r$ & [60,30] & 45024 & N & 420 & \\
	$1/r$ & [62,18] & 1001 & N & 60 & \\
	$1/r$ & [63,20] & 50014 & N & 360 & \\
	\midrule
	\multicolumn{6}{c}{\textit{Inverse Square Law $1/r^2$ Potential} \mystrut(13,0)} \\
	$1/{r^2}$ & [11,5] & 50033 & N & 360 & \\
	$1/{r^2}$ & [13,6] & 50033 & N & 360 & \\
	$1/{r^2}$ & [15,4] & 50014 & N & 360 & \\
	$1/{r^2}$ & [19,38] & 60013 & N & 360 & \\
	$1/{r^2}$ & [20,5] & 50014 & N & 360 & \\
	$1/{r^2}$ & [21,18] & 50014 & N & 360 & \\
	$1/{r^2}$ & [23,6] & 50014 & N & 360 & \\
	$1/{r^2}$ & [24,18] & 63539 & N & 600 & \\
	$1/{r^2}$ & [25,24] & 50014 & N & 420 & \\
	$1/{r^2}$ & [27,2] & 50014 & N & 360 & \\
	$1/{r^2}$ & [28,16] & 50014 & N & 420 & \\
	$1/{r^2}$ & [29,8] & 50014 & N & 420 & \\
	$1/{r^2}$ & [30,20] & 50014 & N & 420 & \\
	$1/{r^2}$ & [31,18] & 6049 & N & 420 & \\
	$1/{r^2}$ & [33,32] & 50014 & N & 420 & \\
	$1/{r^2}$ & [34,16] & 50014 & N & 360 & \\
	$1/{r^2}$ & [37,6] & 10018 & N & 420 & \\
	$1/{r^2}$ & [38,3] & 50014 & N & 360 & \\
	$1/{r^2}$ & [39,10] & 50014 & N & 360 & \\
	$1/{r^2}$ & [40,9] & 50014 & N & 360 & \\
	$1/{r^2}$ & [41,10] & 50014 & N & 360 & \\
	$1/{r^2}$ & [42,5] & 50014 & N & 360 & \\
	$1/{r^2}$ & [43,21] & 50014 & N & 360 & \\
	$1/{r^2}$ & [45,14] & 50014 & N & 360 & \\
	$1/{r^2}$ & [46,24] & 50014 & N & 420 & \\
	$1/{r^2}$ & [48,18] & 1021 & N & 60 & \\
	$1/{r^2}$ & [49,32] & 28378 & N & 420 & \\
	$1/{r^2}$ & [50,4] & 50014 & N & 360 & \\
	$1/{r^2}$ & [51,16] & 50014 & N & 360 & \\
	$1/{r^2}$ & [52,34] & 50014 & N & 420 & \\
	$1/{r^2}$ & [53,27] & 1001 & N & 60 & \\
	$1/{r^2}$ & [57,18] & 50014 & N & 360 & \\
	$1/{r^2}$ & [60,30] & 39118 & N & 420 & \\
	$1/{r^2}$ & [62,18] & 1001 & N & 60 & \\
	$1/{r^2}$ & [63,20] & 50014 & N & 360 & \\
\end{longtable}

There are 33 logarithmic point sets, 41 $1/r$ point sets, and
35 $1/r^2$ point sets archived and available.

{\color{blue}{The 109 point sets can be obtained from Zenodo\cite{45} for further study. They are located online at \textit{ https://doi.org/10.5281/zenodo.5595337} and can be downloaded individually or as a point set zipped file.}}

80 point sets are 50,014 digits in size for 1092 parameters. However some computational problems were encountered and the full desired precision was not reached for the following 23 point sets with 465 parameters:

\noindent
\textit{\textbf{Spherical codes sets with reduced precision due to computation problems --}}

\begin{longtable}[c]{r|crlrl}
	\caption{\textbf{Spherical code sets with lessened precision}} \\
	potential & pts/parms & Digits & found & Height & Notes \\
	\hline\vspace*{-1.3ex}
	\endfirsthead
	\multicolumn{6}{c}%
	{\tablename\ \thetable\ -- \textbf{Spherical code sets with lessened precision} -- \textit{continued}} \\[0.5ex]
	potential & pts/parms & Digits & found & Height & Notes \\
	\hline\vspace*{-1.3ex}
	\endhead
	\multicolumn{6}{c}{\textit{Logarithmic Potential}} \\
	$log$ & [24,9] & 22579 & N & 420 & \\
	$log$ & [31,18] & 6049 & N & 420 & \\
	$log$ & [37,6] & 10018 & N & 420 & \\
	$log$ & [48,18] & 1021 & N & 60 & \\
	$log$ & [49,32] & 48299 & N & 420 & \\
	$log$ & [53,27] & 1001 & N & 60 & \\
	$log$ & [60,30] & 32829 & N & 420 & \\
	$log$ & [62,18] & 1001 & N & 60 & \\
	\hline
	\multicolumn{6}{c}{\textit{Coulomb $1/r$ Potential} \mystrut(13,0)} \\[0.5ex]
	$1/r$ & [24,9] & 25816 & N & 360 & \\
	$1/r$ & [31,18] & 6049 & N & 420 & \\
	$1/r$ & [37,6] & 10018 & N & 420 & \\
	$1/r$ & [48,18] & 1021 & N & 60 & \\
	$1/r$ & [49,32] & 28205 & N & 420 & \\
	$1/r$ & [53,27] & 1001 & N & 60 & \\
	$1/r$ & [60,30] & 45024 & N & 420 & \\
	$1/r$ & [62,18] & 1001 & N & 60 & \\
	\midrule
	\multicolumn{6}{c}{\textit{Inverse Square Law $1/r^2$ Potential} 	\mystrut(13,0)} \\
	$1/{r^2}$ & [31,18] & 6049 & N & 420 & \\
	$1/{r^2}$ & [37,6] & 10018 & N & 420 & \\
	$1/{r^2}$ & [48,18] & 1021 & N & 60 & \\
	$1/{r^2}$ & [49,32] & 28378 & N & 420 & \\
	$1/{r^2}$ & [53,27] & 1001 & N & 60 & \\
	$1/{r^2}$ & [60,30] & 39118 & N & 420 & \\
	$1/{r^2}$ & [62,18] & 1001 & N & 60 & \\
\end{longtable}
These spherical codes are nevertheless still included in the Zenodo database.

\section{Future Work}

There is more work which could be done in finding the algebraic spherical codes. The following 2 points are parameterizable:

\begin{longtable}[c]{r|cccccccccc}
	\caption{Algebraic Spherical Codes -- Unknown} \\[-2.5ex]
	\hline\vspace*{-2.2ex}
	\endfirsthead
	\multicolumn{11}{c}%
	{\tablename\ \thetable\ -- Spherical codes -- Unknown -- \textit{continued}} \\[0.5ex]
	\hline\vspace*{-2.2ex}
	\endhead
	points & 47 & 58
\end{longtable}
\vspace{-2.0em}
There are 6 unknown spherical codes in these two points. 47 points requires about 47 parameters, and 58 points requires 84 parameters. This was too much for the limited computer resources of the correspondence author, so the parameterization was not done for these two point sets.

Additionally there are 11 point sets which cannot be parameterized by a polygon constraint, they have to be optimized by a direct Jacobian method utilizing 2 parameters per point.

\begin{longtable}[c]{r|cccccccccc}
	\caption{Algebraic Spherical Codes -- Direct Search} \\[-2.5ex]
	\hline\vspace*{-2.2ex}
	\endfirsthead
	\multicolumn{11}{c}%
	{\tablename\ \thetable\ -- Spherical codes -- Direct Search -- \textit{continued}} \\[0.5ex]
	\hline\vspace*{-2.2ex}
	\endhead
	points & 26 & \cellcolor{impo} $29\,\dagger$ & 35 & 36 & 54 & 55 & 56 & 59 & 61 & 64 \\
	& 65
\end{longtable}
\vspace{-2.0em}
There are 31 unknown spherical codes in these 11 points. Each spherical code has to be found by using a Jacobian of $2n$ parameters where $n$ is the point size. $\dagger$The \textit{logarithmic} potential for 29 points cannot be parameterized, and has to use all 58 variables to successfully constrain the configuration.

It is hoped that someone might be motivated to possibly use a direct Jacobian on these point sets in hopes of obtaining the spherical codes and possibly algebraic polynomials for the points directly.

\section{Final Comments}

Having spent 9 months arduously working with the spherical codes to high precision of 50,014 digits, and watching the behavior of the \textit{percolating anneal} and \textit{descent.3d} algorithms, and looking at the solutions, the following comments seem necessary by the correspondence author.

\noindent\textbf{Some things observed and noted --}

\begin{itemize}
	\setlength\itemsep{0.7em}
	\item The pair nature of point distances implies $\mathcal{O}^{\frac{n(n-1)}{2}\;}$ behavior for determining the Jacobian matrix and using Newton's method, rapidly straining computer resources.
	\item The energy polynomials inherit the same $\mathcal{O}^{\frac{n(n-1)}{2}\;}$ behavior, rapidly increasing the size of the polynomial degree to be found. The exception seemed to have occurred with the \textit{logarithmic} potential, but the huge coefficients of those polynomials makes it difficult to determine a small number of digits for the actual value.
	\item Apparently obtaining the exact spherical codes for point sets $>50$ will soon be intractable due to the very large memory needed to process the enormous multi-precision decimal numbers involved. \textbf{This is the main lesson here.}
	\item It is very easy to make manual errors when typesetting polynomials with large coefficients and 1,000+ digit numbers to text. Some type of application working with TeX and high precision numbers is needed for the enormous polynomials involved.
	\item Creating the Hessian matrix for larger point sets will soon become problematic, due to the strain on symbolic algebra computer memory resources.
	\item There is a very remote possibility that one configuration might only be a local minimum, not global, but particular attention was paid to avoid this problem, and 2 algorithms were used as a double check on this. The Hessian matrix was used as the third and final check. \textbf{The energy values obtained matched those of Ridgway \& Cheviakov\cite{47} to the 9 digits listed precision for all 3 potentials.}
	\item Some problems with both Maxima and GP-Pari were encountered when working with matrices of megabyte size. The inverse Gauss routine to invert the matrix for finding the linear algebra solution ran into trouble on several point sets which has known points and the GP-Pari program faulted.
	\item \textbf{Due to the triple checks, it is putative that the \textit{global minimum spherical codes have been located for all 195 point sets}, with 49 algebraic number sets resolved, and 106 other algebraic number sets to be discovered from their 50,014 digit spherical codes.}
\end{itemize}

\large{
\noindent
\textbf{\textit{In summary --}}
}

The $\mathcal{O}^{\frac{n(n-1)}{2}\;}$ behavior of point distances spells trouble for determining the exact algebraic characteristic polynomial for point sets $>50$ and will soon become intractable.

This paper is partially presented as a challenge to improve the current algebraic polynomial mathematics toolset, although GP-Pari does a very admirable job. The 1622 parameters need their algebraic polynomials resolved.

Finally, the discovered 195 spherical codes and 49 algebraic numbers are presented as a starting answer to Problem \#7 of Steve Smale\cite{55,56}, for the first 65 $n-$points configurations on the S2 sphere under 3 potentials or force laws.

\section{Global Minimum Spherical Codes}

{\color{blue}{The spherical codes for all three potentials for 0 to 65 points are archived on Zenodo\cite{46}. They can be accessed at \textit{https://doi.org/10.5281/zenodo.5595366} and downloaded individually or as a zipped file set for a given potential.}}

These high precision spherical codes are accurate to either 77 digits for the \textit{logarithmic} potential or 38 digits for the \textit{Coulomb $1/r$} and \textit{Inverse Square $1/r^2$} potentials. 2 or 3 guard digits were provided for the spherical codes for accuracy.

\noindent
\textbf{NOTE:} It must be mentioned again that the \textit{Coulomb $1/r$} values are initially taken from Neil J.A. Sloane's \textit{"Spherical codes for minimal energy"} database\cite{53}, and improved from 12 digits accuracy to 38 digits accuracy, thus extending the precision of his points. Special care has been taken to keep the order of the points, exactly as in his tables. It is now shown that his putative values are indeed the approximations of the global minimums for the \textit{Coulomb} $1/r$ potential.

The \textit{Logarithmic} and \textit{Inverse Square $1/r^2$} potentials used Sloane's 3d points as an initiating value, so the isometries are kept identical for all 3 potentials.

\noindent
\textit{\textbf{Index of High Precision Spherical Codes --}}

\begin{longtable}[c]{r|ccc}
	\caption{High Precision coordinate point sets} \\
	pts & logarithmic & Coulomb & Inverse Square \\
	\hline\vspace*{-2.2ex}
	\endfirsthead
	\multicolumn{4}{c}%
	{\tablename\ \thetable\ -- High Precision coordinates -- \textit{continued}} \\
	pts & logarithmic & Coulomb & Inverse Square \\
	\hline\vspace*{-2.2ex}
	\endhead
	0 & r1.3.0.normal.80 & r1.3.0.40 & r2.3.0.40 \\
	1 & r1.3.1.normal.80 & r1.3.1.40 & r2.3.1.40 \\
	2 & r1.3.2.normal.80 & r1.3.2.40 & r2.3.2.40 \\
	3 & r1.3.3.normal.80 & r1.3.3.40 & r2.3.3.40 \\
	4 & r1.3.4.normal.80 & r1.3.4.40 & r2.3.4.40 \\
	5 & r1.3.5.normal.80 & r1.3.5.40 & r2.3.5.40 \\
	6 & r1.3.6.normal.80 & r1.3.6.40 & r2.3.6.40 \\
	7 & r1.3.7.normal.80 & r1.3.7.40 & r2.3.7.40 \\
	8 & r1.3.8.normal.80 & r1.3.8.40 & r2.3.8.40 \\
	9 & r1.3.9.normal.80 & r1.3.9.40 & r2.3.9.40 \\
	10 & r1.3.10.normal.80 & r1.3.10.40 & r2.3.10.40 \\
	11 & r1.3.11.normal.80 & r1.3.11.40 & r2.3.11.40 \\
	12 & r1.3.12.normal.80 & r1.3.12.40 & r2.3.12.40 \\
	13 & r1.3.13.normal.80 & r1.3.13.40 & r2.3.13.40 \\
	14 & r1.3.14.normal.80 & r1.3.14.40 & r2.3.14.40 \\
	15 & r1.3.15.normal.80 & r1.3.15.40 & r2.3.15.40 \\
	16 & r1.3.16.normal.80 & r1.3.16.40 & r2.3.16.40 \\
	17 & r1.3.17.normal.80 & r1.3.17.40 & r2.3.17.40 \\
	18 & r1.3.18.normal.80 & r1.3.18.40 & r2.3.18.40 \\
	19 & r1.3.19.normal.80 & r1.3.19.40 & r2.3.19.40 \\
	20 & r1.3.20.normal.80 & r1.3.20.40 & r2.3.20.40 \\
	21 & r1.3.21.normal.80 & r1.3.21.40 & r2.3.21.40 \\
	22 & r1.3.22.normal.80 & r1.3.22.40 & r2.3.22.40 \\
	23 & r1.3.23.normal.80 & r1.3.23.40 & r2.3.23.40 \\
	24 & r1.3.24.normal.80 & r1.3.24.40 & r2.3.24.40 \\
	25 & r1.3.25.normal.80 & r1.3.25.40 & r2.3.25.40 \\
	26 & r1.3.26.normal.80 & r1.3.26.40 & r2.3.26.40 \\
	27 & r1.3.27.normal.80 & r1.3.27.40 & r2.3.27.40 \\
	28 & r1.3.28.normal.80 & r1.3.28.40 & r2.3.28.40 \\
	29 & r1.3.29.normal.80 & r1.3.29.40 & r2.3.29.40 \\
	30 & r1.3.30.normal.80 & r1.3.30.40 & r2.3.30.40 \\
	31 & r1.3.31.normal.80 & r1.3.31.40 & r2.3.31.40 \\
	32 & r1.3.32.normal.80 & r1.3.32.40 & r2.3.32.40 \\
	33 & r1.3.33.normal.80 & r1.3.33.40 & r2.3.33.40 \\
	34 & r1.3.34.normal.80 & r1.3.34.40 & r2.3.34.40 \\
	35 & r1.3.35.normal.80 & r1.3.35.40 & r2.3.35.40 \\
	36 & r1.3.36.normal.80 & r1.3.36.40 & r2.3.36.40 \\
	37 & r1.3.37.normal.80 & r1.3.37.40 & r2.3.37.40 \\
	38 & r1.3.38.normal.80 & r1.3.38.40 & r2.3.38.40 \\
	39 & r1.3.39.normal.80 & r1.3.39.40 & r2.3.39.40 \\
	40 & r1.3.40.normal.80 & r1.3.40.40 & r2.3.40.40 \\
	41 & r1.3.41.normal.80 & r1.3.41.40 & r2.3.41.40 \\
	42 & r1.3.42.normal.80 & r1.3.42.40 & r2.3.42.40 \\
	43 & r1.3.43.normal.80 & r1.3.43.40 & r2.3.43.40 \\
	44 & r1.3.44.normal.80 & r1.3.44.40 & r2.3.44.40 \\
	45 & r1.3.45.normal.80 & r1.3.45.40 & r2.3.45.40 \\
	46 & r1.3.46.normal.80 & r1.3.46.40 & r2.3.46.40 \\
	47 & r1.3.47.normal.80 & r1.3.47.40 & r2.3.47.40 \\
	48 & r1.3.48.normal.80 & r1.3.48.40 & r2.3.48.40 \\
	49 & r1.3.49.normal.80 & r1.3.49.40 & r2.3.49.40 \\
	50 & r1.3.50.normal.80 & r1.3.50.40 & r2.3.50.40 \\
	51 & r1.3.51.normal.80 & r1.3.51.40 & r2.3.51.40 \\
	52 & r1.3.52.normal.80 & r1.3.52.40 & r2.3.52.40 \\
	53 & r1.3.53.normal.80 & r1.3.53.40 & r2.3.53.40 \\
	54 & r1.3.54.normal.80 & r1.3.54.40 & r2.3.54.40 \\
	55 & r1.3.55.normal.80 & r1.3.55.40 & r2.3.55.40 \\
	56 & r1.3.56.normal.80 & r1.3.56.40 & r2.3.56.40 \\
	57 & r1.3.57.normal.80 & r1.3.57.40 & r2.3.57.40 \\
	58 & r1.3.58.normal.80 & r1.3.58.40 & r2.3.58.40 \\
	59 & r1.3.59.normal.80 & r1.3.59.40 & r2.3.59.40 \\
	60 & r1.3.60.normal.80 & r1.3.60.40 & r2.3.60.40 \\
	61 & r1.3.61.normal.80 & r1.3.61.40 & r2.3.61.40 \\
	62 & r1.3.62.normal.80 & r1.3.62.40 & r2.3.62.40 \\
	63 & r1.3.63.normal.80 & r1.3.63.40 & r2.3.63.40 \\
	64 & r1.3.64.normal.80 & r1.3.64.40 & r2.3.64.40 \\
	65 & r1.3.65.normal.80 & r1.3.65.40 & r2.3.65.40
\end{longtable}

\section{Acknowledgments}

I wish to thank the following people and organizations for their help in creating this paper and the mathematical software used to obtain the results.

A special thanks to Wesley Ridgway for carefully explaining the mathematics behind their search procedures, and running the Hessian matrix calculations for that important third and final check and to Professor Alexei F. Cheviakov who graciously consented to allowing Wesley to assist with my questions and to check my work.

A special thanks to Bill Alombert who assisted with some GP-Pari questions, such as the decimal precision representation stored in the software. And I wish to thank the GP-Pari team who made it possible to work with 50,014 digit numbers with repeatable accuracy.

I also want to thank the linux Maxima people for their program to assist in the preparation of the very large size Jacobian matrices involved and doing the differentiation work in multi-variable environments correctly.

Finall, special thanks to Zenodo for hosting the spherical codes and algebraic floating point numbers to very high precision.

With them, none of this would have been possible.

\noindent
-- \textit{Randall L. Rathbun} \\
\hspace*{1.5ex} correspondence author

\section{Declaration of Interest}
There are no conflicts of interest.

\section{Appendix}

\noindent
\textit{\textbf{Algebraic Polynomial or Parameterized Spherical Codes --}}

The 51 algebraic algorithms in GP-Pari code listed below are available also, upon request to the correspondence author, but nearly all require its corresponding parameter file, for full accuracy. The parameter files are large, since they hold variables of size 50,014 digits. The largest holds 34 parameters, so these files can be over 1.7 Megabytes in size.

\begin{longtable}[c]{r|ccc}
	\caption{Algebraic Solution sets} \\
	pts & Logarithmic & Coulomb & Inverse Square \\
	\hline\vspace*{-2.2ex}
	\endfirsthead
	\multicolumn{4}{c}%
	{\tablename\ \thetable\ -- Algebraic solution sets -- \textit{continued}} \\
	pts & Logarithmic & Coulomb & Inverse Square \\
	\hline\vspace*{-2.2ex}
	\endhead
	0 & r1.3.0.normal.alg & r1.3.0.alg & r2.3.0.alg \\
	1 & r1.3.1.normal.alg & r1.3.1.alg & r2.3.1.alg \\
	2 & r1.3.2.normal.alg & r1.3.2.alg & r2.3.2.alg \\
	3 & r1.3.3.normal.alg & r1.3.3.alg & r2.3.3.alg \\
	4 & r1.3.4.normal.alg & r1.3.4.alg & r2.3.4.alg \\
	5 & r1.3.5.normal.alg & r1.3.5.alg & r2.3.5.alg \\
	6 & r1.3.6.normal.alg & r1.3.6.alg & r2.3.6.alg \\
	7 & r1.3.7.normal.alg & r1.3.7.alg & r2.3.7.alg \\
	8 & r1.3.8.normal.alg & r1.3.8.alg & r2.3.8.alg \\
	9 & r1.3.9.normal.alg & r1.3.9.alg & r2.3.9.alg \\
	10 & r1.3.10.normal.alg & r1.3.10.alg & r2.3.10.alg \\
	11 & r1.3.11.normal.alg & r1.3.11.alg & r2.3.11.alg \\
	12 & r1.3.12.normal.alg & r1.3.12.alg & r2.3.12.alg \\
	13 & r1.3.13.normal.alg & r1.3.13.alg & r2.3.13.alg \\
	14 & r1.3.14.normal.alg & r1.3.14.alg & r2.3.14.alg \\
	15 & r1.3.15.normal.alg & r1.3.15.alg & r2.3.15.alg \\
	16 & r1.3.16.normal.alg & r1.3.16.alg & r2.3.16.alg \\
	17 & r1.3.17.normal.alg & r1.3.17.alg & r2.3.17.alg \\
	18 & r1.3.18.normal.alg & r1.3.18.alg & r2.3.18.alg \\
	19 & r1.3.19.normal.alg & r1.3.19.alg & (not applicable) \\
	20 & r1.3.20.normal.alg & r1.3.20.alg & r2.3.20.alg \\
	21 & r1.3.21.normal.alg & r1.3.21.alg & r2.3.21.alg \\
	22 & r1.3.22.normal.alg & r1.3.22.alg & r2.3.22.alg \\
	23 & r1.3.23.normal.alg & r1.3.23.alg & r2.3.23.alg \\
	24 & r1.3.24.normal.alg & r1.3.24.alg & r2.3.24.alg \\
	25 & r1.3.25.normla.alg & r1.3.25.alg & r2.3.25.alg \\
	27 & r1.3.27.normal.alg & r1.3.27.alg & r2.3.27.alg \\
	28 & r1.3.28.normal.alg & r1.3.28.alg & r2.3.28.alg \\
	30 & r1.3.30.normal.alg & r1.3.30.alg & r2.3.30.alg \\
	31 & r1.3.31.normal.alg & r1.3.31.alg & r2.3.31.alg \\
	32 & r1.3.32.normal.alg & r1.3.32.alg & r2.3.32.alg \\
	33 & r1.3.33.normal.alg & r1.3.33.alg & r2.3.33.alg \\
	34 & r1.3.34.normal.alg & r1.3.34.alg & r2.3.34.alg \\
	37 & r1.3.37.normal.alg & r1.3.37.alg & r2.3.37.alg \\
	38 & r1.3.38.normal.alg & r1.3.38.alg & r2.3.38.alg \\
	39 & r1.3.39.normal.alg & r1.3.39.alg & r2.3.39.alg \\
	40 & r1.3.40.normal.alg & r1.3.40.alg & r2.3.40.alg \\
	41 & r1.3.41.normal.alg & r1.3.41.alg & r2.3.41.alg \\
	42 & r1.3.42.normal.alg & r1.3.42.alg & r2.3.42.alg \\
	43 & r1.3.43.normal.alg & r1.3.43.alg & r2.3.43.alg \\
	44 & r1.3.44.normal.alg & r1.3.44.alg & r2.3.44.alg \\
	45 & r1.3.45.normal.alg & r1.3.45.alg & r2.3.45.alg \\
	46 & r1.3.46.normal.alg & r1.3.46.alg & r2.3.46.alg \\
	48 & r1.3.48.normal.alg & r1.3.48.alg & r2.3.48.alg \\
	49 & r1.3.49.normal.alg & r1.3.49.alg & r2.3.49.alg \\
	50 & r1.3.50.normal.alg & r1.3.50.alg & r2.3.50.alg \\
	51 & r1.3.51.normal.alg & r1.3.51.alg & r2.3.51.alg \\
	52 & r1.3.52.normal.alg & r1.3.52.alg & r2.3.52.alg \\
	53 & r1.3.53.normal.alg & r1.3.53.alg & r2.3.53.alg \\
	57 & r1.3.57.normal.alg & r1.3.57.alg & r2.3.57.alg \\
	63 & r1.3.63.normal.alg & r1.3.63.alg & r2.3.63.alg
\end{longtable}

Additionally, serious inquiries into the GP-Pari and Maxima code for finding the algebraic solutions, or discussions about the algorithms is welcomed. Please contact the correspondence author {\it Randall L. Rathbun} at the email address listed at the end of this paper.

\noindent
-- \textit{Randall L. Rathbun} \\
\hspace*{1.5ex} correspondence author

\rule[1ex]{3.2in}{0.4pt}\hspace*{\fill} \\[-0.5em]\hspace*{\fill}


\section{Bibliography}

\begin{thebibliography}{62}
	\bibitem{1} Alishahi, K.; Zamani, M.; The Spherical Ensemble and Uniform Distribution of Points on the Sphere. arXiv:1407.5832 [math-ph] 2014. Available at https://arxiv.org/abs/1407.5832
	\bibitem{2} Altschuler, E. L.; Perez, A.; Stong, R.; A Novel Symmetric Four Dimensional Polytope Found Using Optimization Strategies Inspired by Thomson’s Problem of Charges on a Sphere. Available at https://repositorio.upct.es/bitstream/handle/10317/593/nsf.pdf;
	\bibitem{3} Altschuler, E. L.; Pérez–Garrido, A.; Defect-Free Global Minima in Thomson’s Problem of Charges on a Sphere. Phys. Rev. E 2006, 73 (3), 036108. DOI: https://doi.org/10.1103/PhysRevE.73.036108.
	\bibitem{4} Altschuler, E. L.; Perez-Garrido, A.; New Global Minima for Thomson’s Problem of Charges on a Sphere. arXiv:cond-mat/0408355 2004. Available at https://arxiv.org/abs/cond-mat/0408355
	\bibitem{5} Altschuler, E. L.; Williams, T. J.; Ratner, E. R.; Tipton, R.; Stong, R.; Dowla, F.; Wooten, F.; Possible Global Minimum Lattice Configurations for Thomson’s Problem of Charges on a Sphere. Phys. Rev. Lett. 1997, 78 (14), 2681–2685. DOI: https://doi.org/10.1103/PhysRevLett.78.2681.
	\bibitem{6} Atiyah, M.; Sutcliffe, P.; The Geometry of Point Particles. Proc. R. Soc. Lond. A 2002, 458 (2021), pp1089-1115. DOI: https://doi.org/10.1098/rspa.2001.0913.
	\bibitem{7} Ballinger, B.; Blekherman, G.; Cohn, H.; Giansiracusa, N.; Kelly, E.; Schürmann, A.; Experimental Study of Energy-Minimizing Point Configurations on Spheres. Experimental Mathematics 2009, 18 (3), 257–283. DOI: https://doi.org/10.1080/10586458.2009.10129052.
	\bibitem{8} Batle, J.; Generalized Thomson Problem in Arbitrary Dimensions and Non-Euclidean Geometries. arXiv:1312.1854 [physics, physics:quant-ph] 2013. Available at https://arxiv.org/abs/1312.1854
	\bibitem{9} Băutu, A.; Băutu, E.; Energy Minimization of Point Charges on a Sphere with Particle Swarms. Romanian Reports of Physics 2009, 54, 8. Available at https://rjp.nipne.ro/2009\_54\_1--2/0029\_0037.pdf
	\bibitem{10} Beltrán, C.; Marzo, J.; Ortega-Cerdà, J.; Energy and discrepancy of rotationally invariant determinantal point processes in high dimensional spheres. Journal of Complexity Vol 37, Dec 2016, pp76-109 DOI: https://doi.org/10.1016/j.jco.2016.08.001
	\bibitem{11} Beltrán, C.; The State of the Art in Smale’s 7th Problem. In Foundations of Computational Mathematics, Budapest 2011; Cucker, F., Krick, T., Pinkus, A., Szanto, A., Editors.; Cambridge University Press: Cambridge, 2012; pp 1–15. DOI: https://doi.org/10.1017/CBO9781139095402.002.
	\bibitem{12} Berman, J.; Hanes, K.; Optimizing the Arrangement of Points on the Unit Sphere. Math of Comp, v31, \#140, Oct. 1977, pp1006-1008 Available at https://www.ams.org/journals/mcom/1977-31-140/S0025-5718-1977-0478006-1/S0025-5718-1977-0478006-1.pdf
	\bibitem{13} Bowick, M.; Cacciuto, A.; Nelson, D. R.; Travesset, A.; Crystalline Order on a Sphere and the Generalized Thomson Problem. Phys. Rev. Lett. 2002, 89 (18), 185502. DOI: https://doi.org/10.1103/PhysRevLett.89.185502.
	\bibitem{14} Brauchart, J. S.; Grabner, P. J.; Distributing Many Points on Spheres: Minimal Energy and Designs. Journal of Complexity 2015, 31 (3), 293–326. DOI: https://doi.org/10.1016/j.jco.2015.02.003.
	\bibitem{15} Brown, K. S.; Min-Energy Configurations of Electrons On A Sphere https://www.mathpages.com/home/kmath005/kmath005.htm (accessed Sep 22, 2019).
	\bibitem{16} Cheviakov, A. F.; Ridgway, W. J. M.; Local and Global Optimization of Particle Locations on the Sphere: Models, Applications, Mathematical Aspects, and Computations. (contact author: cheviakov@math.usask.ca)
	\bibitem{17} Cohn, H.; Elkies, N. D.; Kumar, A.; Schürmann, A.; Point Configurations That Are Asymmetric yet Balanced. Proc. Amer. Math. Soc. 2010, 138 (08), pp2863-2863. DOI: https://doi.org/10.1090/S0002-9939-10-10284-6.
	\bibitem{18} Cohn, H.; Kumar, A.; Universally Optimal Distribution of Points on Spheres. Journal AMS, Vol 20, \#1, Jan 2007, pp99--148, S 0894-0347(06)00546-7, pub Sept 5, 2006. Available at https://www.ams.org/journals/jams/2007-20-01/S0894-0347-06-00546-7/S0894-0347-06-00546-7.pdf
	\bibitem{19} Cohn, H.; Stability Configurations of Electrons on a Sphere. Math. Comp. 1956, 10 (55), pp117–120. DOI: https://doi.org/10.1090/S0025-5718-1956-0081133-0.
	\bibitem{20} Costa, S. I. R.; Oggier, F.; Campello, A.; Belfiore, J.-C.; Viterbo, E.; Lattices and Spherical Codes. In Lattices Applied to Coding for Reliable and Secure Communications; Springer International Publishing: Cham, 2017; pp 73–92. DOI: https://doi.org/10.1007/978-3-319-67882-5\_5.
	\bibitem{21} Dragnev, P. D.; Legg, D. A.; Townsend, D. W.; Discrete Logarithmic Energy on the Sphere. Pacific Journal of Mathematics 2002, 207 (2), 345–358.
	\bibitem{22} Dragnev, P. D.; Log-Optimal Configurations on the Sphere. arXiv:1504.02544 [math-ph] 2015. Available at https://arxiv.org/abs/1504.02544
	\bibitem{23} Erber, T.; Hockney, G. M.; Complex Systems: Equilibrium Configurations of N Equal Charges on a Sphere (2 $\le$ N $\le$ 112). In Advances in Chemical Physics; John Wiley \& Sons, Ltd, 2007; pp 495–594. DOI: https://doi.org/10.1002/9780470141571.ch5.
	\bibitem{24} Erber, T.; Hockney, G. M.; Equilibrium Configurations of N Equal Charges on a Sphere. J. Phys. A: Math. Gen. 1991, 24 (23), L1369 – L1377. DOI https://doi.org/10.1088/0305-4470/24/23/008.
	\bibitem{25} Gautam, S.; Vaintrob, D.; A Novel Approach to the Spherical Codes Problem. Available at https://math.mit.edu/research/highschool/rsi/documents/2012Gautam.pdf
	\bibitem{26} Glasser, L.; Every, A. G.; Energies and Spacings of Point Charges on a Sphere. J. Phys. A: Math. Gen. 1992, 25 (9), 2473–2482. DOI: https://doi.org/10.1088/0305-4470/25/9/020.
	\bibitem{27} Global minima for the Thomson problem. Online: http://www-wales.ch.cam.ac.uk/$\sim$wales/CCD/Thomson/table.html (accessed Sep 22, 2019).
	\bibitem{28} Goldberg, M.; Stability Configurations of Electrons on a Sphere. Math. Comp. 1969, 23 (108), pp785–786. DOI: https://doi.org/10.1090/S0025-5718-69-99642-2.
	\bibitem{29} Grabner, P. J.;Point Sets of Minimal Energy. In Applied Algebra and Number Theory; Larcher, G., Pillichshammer, F., Winterhof, A., Xing, C., Editors.; Cambridge University Press: Cambridge, 2014; pp 109–125. DOI: https://doi.org/10.1017/CBO9781139696456.008.
	\bibitem{30} Hardin, D. P.; Michaels, T.; Saff, E. B.; A Comparison of Popular Point Configurations on $S^2$. Available at http://www.math.vanderbilt.edu/saffeb/texts/262.pdf
	\bibitem{31} Heckbert, P. S.; Fast Surface Particle Repulsion. 20. SIGGRAPH'97, New Frontiers in Modeling. Available at https://www.researchgate.net/profile/Paul\-Heckbert/publication/2431573\_Fast\_Surface\_Particle\_Repulsion/links/ 551f58a40cf2f9c1304de221/Fast-Surface-Particle-Repulsion.pdf
	\bibitem{32} Kanimozhi, G.; Rajathy, R.; Kumar, H.; Minimizing Energy of Point Charges on a Sphere Using Symbiotic Organisms Search Algorithm. International Journal on Electrical Engineering and Informatics, Vol 8, \#1, March 2016 DOI: https://doi.org/10.15676/ijeei.2016.8.1.3. Available at https://www.researchgate.net/profile/Kanimozhi-G-2/publication/301535082\_Minimizing\_Energy\_of\_P\_oint\_Charges\_on\_a\_Sphere \_using\_Symbiotic\_Organisms\_Search\_Algorithm/links/5b90d080299bf114b7fd8f2a/ Minimizing-Energy-of-P-oint-Charges-on-a-Sphere-using-Symbiotic-Organisms-Search-Algorithm.pdf
	\bibitem{33} Katanforoush, A.; Shahshahani, M.; Distributing Points on the Sphere, I. Experimental Mathematics 2003, 12 (2), 199–209. DOI: https://doi.org/10.1080/10586458.2003.10504492.
	\bibitem{34} Kuijlaars, A. B. J.; Saff, E. B.; Sun, X.; On Separation of Minimal Riesz Energy Points on Spheres in Euclidean Spaces. arXiv:math-ph/0503063 2005. Available at https://arxiv.org/abs/math-ph/0503063
	\bibitem{35} Kuijlaars, A. B. J.; Saff, E. B.; Sun, X.; On Separation of Minimal Riesz Energy Points on Spheres in Euclidean Spaces. Journal of Computational and Applied Mathematics 2007, 199 (1), 172–180. DOI: https://doi.org/10.1016/j.cam.2005.04.074.
	\bibitem{36} Lakhbab, H.; Bernoussi, S. E.; Harif, A. E.; Energy Minimization of Point Charges on a Sphere with a Spectral Projected Gradient Method. 2012, 3 (5), Available at https://www.ijser.org/onlineResearchPaperViewer.aspx?Energy-Minimization-of-Point-Charges-on-a-Sphere-with-a-Spectral-Projected-Gradient-Method.pdf
	\bibitem{37} Lakhbab, H.; Energy Minimization of Point Charges on a Sphere with a Hybrid Approach. Applied Mathematical Sciences, Vol. 6, 2012, no. 30, 1487 - 1495. Available at http://www.m-hikari.com/ams/ams-2012/ams-29-32-2012/lakhbabAMS29-32-2012
	\bibitem{38} Leech, J.; Equilibrium of Sets of Particles on a Sphere. The Mathematical Gazette 1957, 41 (336), 81–90. DOI: https://doi.org/10.2307/3610579.
	\bibitem{39} Levin, Y.; Arenzon, J. J.; Why Charges Go to the Surface: A Generalized Thomson Problem. Europhys. Lett. 2003, 63 (3), 415–418. DOI: https://doi.org/10.1209/epl/i2003-00546-1.
	\bibitem{40} Melnyk, T. W.; Knop, O.; Smith, W. R.; Extremal Arrangements of Points and Unit Charges on a Sphere: Equilibrium Configurations Revisited. Can. J. Chem. 1977, 55 (10), 1745–1761. DOI: https://doi.org/10.1139/v77-246.
	\bibitem{41} Morris, J. R.; Deaven, D. M.; Ho, K. M.; Genetic-Algorithm Energy Minimization for Point Charges on a Sphere. Phys. Rev. B 1996, 53 (4), R1740–R1743. DOI: https://doi.org/10.1103/PhysRevB.53.R1740.
	\bibitem{42} Paquay S.; Kusumaatmaja H.; Wales, D. J.; Zandi, R; Van der Schoot, P; Energetically favoured defects in dense packings of particles on spherical surfaces. Soft Matter 2016 Jun 29; 12(26), pp5708-17. doi: 10.1039/c6sm00489j. PMID: 27263532 DOI: https://doi.org/10.1039/c6sm00489j
	\bibitem{43} Rakhmanov, E. A.; Saff, E. B.; Zhou, Y.; Electrons on the Sphere. Series in Approximations and Decompositions, Computational Methods and Function Theory 1994, pp. 293-309 (1995) DOI: https://doi.org/10.1142/9789814533232\_0022
	\bibitem{44} Raman, P.; Yang, J.; Optimization on the Surface of the (Hyper)--Sphere. Available at https://www.researchgate.net/profile/Parameswaran--Raman--2/publication/335855051\_Optimization\_on\_the\_Surface\_of\_the\_Hyper--Sphere/links/5d99100ca6fdccfd0e7b2fe7/Optimization--on--the--Surface--of--the--Hyper--Sphere.pdf
	\bibitem{45} Rathbun, Randall L.; Ridgway, Wesley J.M. (2021). 1585 Algebraic numbers to be found for Smale's Problem \#7 [Data set]. Zenodo. https://doi.org/10.5281/zenodo.5595337
	\bibitem{46} Rathbun, Randall L.; Ridgway, Wesley J.M. (2021). 198 Spherical Codes in S2 for Smale's Problem 7 - 0 to 65 points - global minimums [Data set]. Zenodo. https://doi.org/10.5281/zenodo.5595366
	\bibitem{47} Ridgway, W. J. M.; Cheviakov. A. F.; An iterative procedure for finding locally and globally optimal arrangements of particles on the unit sphere. Computer Physics Communications, Volume 233, 2018, pp 84--109, ISSN 0010--4655, DOI: https://doi.org/10.1016/j.cpc.2018.03.029. Available at https://www.sciencedirect.com/science/article/pii/S0010465518301292
	\bibitem{48} Ridgway, W. J. M.; Cheviakov, A. F.; Energy-Minimizing Arrangements of Repelling Particles on the Sphere: Coulombic and Narrow Escape Potentials. Available at http://math.usask.ca/~shevyakov/research/students/wr\_poster\_2016.pdf
	\bibitem{49} Saff, E. B.; Discrete Minimal Energy Problems - Lecture II. Available at http://fourier.math.uoc.gr/saff/lecture-2.pdf
	\bibitem{50} Saff, E. B.; Kuijlaars, A. B. J.; Distributing Many Points on a Sphere. The Mathematical Intelligencer 1997, 19 (1), 5–11. DOI: https://doi.org/10.1007/BF03024331.
	\bibitem{51} Schwartz, R. E.; The 5-Electron Case of Thomson’s Problem. Experimental Mathematics, Volume 22, 2013, Issue 2, pp157-186, Pub online 25 Apr 2013 Available at https://arxiv.org/abs/1001.3702
	\bibitem{52} Schwartz, R. E.; The Triangular Bi-Pyramid Minimizes a Range of Power Law Potentials. Available at https://arxiv.org/abs/1512.04628
	\bibitem{53} Sloane, Neil J. A., "Spherical codes, minimal energy", Database available at http://neilsloane.com/electrons/dim3/
	\bibitem{54} Slosar, A.; Podgornik, R.; On the Connected-Charges Thomson Problem. Euro Physics Letters, 2006, 75 (4), 631. DOI: https://doi.org/10.1209/epl/i2006-10146-1. Available at https://arxiv.org/pdf/cond-mat/0606765
	\bibitem{55} Smale, S.; (1999). "Mathematical problems for the next century". In Arnold, V. I.; Atiyah, M.; Lax, P.; Mazur, B. (editors). Mathematics: frontiers and perspectives. American Mathematical Society. pp. 271-294. ISBN 978-0821820704. Available at http://citeseerx.ist.psu.edu/viewdoc/download?doi=10.1.1.96.9176\&rep=rep1\&type=pdf
	\bibitem{56} Smale, S.; Mathematical Problems for the Next Century. Mathematical Intelligencer 1998, 20, 7–15.
	\bibitem{57} Thomson, J. J.; (1904) XXIV. On the structure of the atom: an investigation of the stability and periods of oscillation of a number of corpuscles arranged at equal intervals around the circumference of a circle; with application of the results to the theory of atomic structure , Philosophical Magazine Series 6, 7:39, 237-265, DOI: https://doi.org/10.1080/14786440409463107
	\bibitem{58} Thomson Problem. Wikipedia; Available at https://en.wikipedia.org/wiki/Thomson\_problem
	\bibitem{59} Wales, D. J.; Ulker, S. Structure and Dynamics of Spherical Crystals Characterized for the Thomson Problem. Phys. Rev. B 2006, 74 (21), 212101. DOI: https://doi.org/10.1103/PhysRevB.74.212101.
	\bibitem{60} Wang, J.; Finding and Investigating Exact Spherical Codes. arXiv:0805.0776 [math] 2008. Available at https://arxiv.org/abs/0805.0776
	\bibitem{61} Wayback Machine, Online: https://web.archive.org/web/20131213172104/http://www.cond-mat.physik.uni-mainz.de/~oettel/ws10/thomson\_PhilMag\_7\_237\_1904.pdf
	\bibitem{62} Yudin, V. A.; The Minimum of Potential Energy of a System of Point Charges. Discrete Mathematics and Applications 2009, 3 (1), 75–82. DOI: https://doi.org/10.1515/dma.1993.3.1.75.
\end{thebibliography}
\end{document}